\pdfoutput=1
\documentclass{amsart}

\usepackage{amssymb,amsmath,stmaryrd,mathrsfs}
\def\definetac{\newif\iftac}    
\ifx\tactrue\undefined
  \definetac
  \ifx\state\undefined\tacfalse\else\tactrue\fi
\fi
\iftac\else\usepackage{amsthm}\fi
\usepackage[all,2cell,pdf]{xy}
\UseAllTwocells
\usepackage{enumitem}
\usepackage{xcolor}
\definecolor{darkgreen}{rgb}{0,0.45,0} 
\usepackage[pagebackref,colorlinks,citecolor=darkgreen,linkcolor=darkgreen]{hyperref}
\usepackage{mathtools}          
\usepackage{tikz}
\usepackage{braket}             

\usepackage{url}                
\usepackage{xspace}             

\makeatletter
\let\ea\expandafter

\def\mdef#1#2{\ea\ea\ea\gdef\ea\ea\noexpand#1\ea{\ea\ensuremath\ea{#2}\xspace}}
\def\alwaysmath#1{\ea\ea\ea\global\ea\ea\ea\let\ea\ea\csname your@#1\endcsname\csname #1\endcsname
  \ea\def\csname #1\endcsname{\ensuremath{\csname your@#1\endcsname}\xspace}}

\DeclareRobustCommand\widecheck[1]{{\mathpalette\@widecheck{#1}}}
\def\@widecheck#1#2{%
    \setbox\z@\hbox{\m@th$#1#2$}%
    \setbox\tw@\hbox{\m@th$#1%
       \widehat{%
          \vrule\@width\z@\@height\ht\z@
          \vrule\@height\z@\@width\wd\z@}$}%
    \dp\tw@-\ht\z@
    \@tempdima\ht\z@ \advance\@tempdima2\ht\tw@ \divide\@tempdima\thr@@
    \setbox\tw@\hbox{%
       \raise\@tempdima\hbox{\scalebox{1}[-1]{\lower\@tempdima\box
\tw@}}}%
    {\ooalign{\box\tw@ \cr \box\z@}}}

\newcount\foreachcount

\def\foreachletter#1#2#3{\foreachcount=#1
  \ea\loop\ea\ea\ea#3\@alph\foreachcount
  \advance\foreachcount by 1
  \ifnum\foreachcount<#2\repeat}

\def\foreachLetter#1#2#3{\foreachcount=#1
  \ea\loop\ea\ea\ea#3\@Alph\foreachcount
  \advance\foreachcount by 1
  \ifnum\foreachcount<#2\repeat}

\def\definescr#1{\ea\gdef\csname s#1\endcsname{\ensuremath{\mathscr{#1}}\xspace}}
\foreachLetter{1}{27}{\definescr}
\def\definecal#1{\ea\gdef\csname c#1\endcsname{\ensuremath{\mathcal{#1}}\xspace}}
\foreachLetter{1}{27}{\definecal}
\def\definebold#1{\ea\gdef\csname b#1\endcsname{\ensuremath{\mathbf{#1}}\xspace}}
\foreachLetter{1}{27}{\definebold}
\def\definebb#1{\ea\gdef\csname l#1\endcsname{\ensuremath{\mathbb{#1}}\xspace}}
\foreachLetter{1}{27}{\definebb}
\def\definefrak#1{\ea\gdef\csname f#1\endcsname{\ensuremath{\mathfrak{#1}}\xspace}}
\foreachletter{1}{9}{\definefrak} 
\foreachletter{10}{27}{\definefrak}
\def\definebar#1{\ea\gdef\csname #1bar\endcsname{\ensuremath{\overline{#1}}\xspace}}
\foreachLetter{1}{27}{\definebar}
\foreachletter{1}{8}{\definebar} 
\foreachletter{9}{15}{\definebar} 
\foreachletter{16}{27}{\definebar}
\def\definetil#1{\ea\gdef\csname #1til\endcsname{\ensuremath{\widetilde{#1}}\xspace}}
\foreachLetter{1}{27}{\definetil}
\foreachletter{1}{27}{\definetil}
\def\definehat#1{\ea\gdef\csname #1hat\endcsname{\ensuremath{\widehat{#1}}\xspace}}
\foreachLetter{1}{27}{\definehat}
\foreachletter{1}{27}{\definehat}
\def\definechk#1{\ea\gdef\csname #1chk\endcsname{\ensuremath{\widecheck{#1}}\xspace}}
\foreachLetter{1}{27}{\definechk}
\foreachletter{1}{27}{\definechk}
\def\defineul#1{\ea\gdef\csname u#1\endcsname{\ensuremath{\underline{#1}}\xspace}}
\foreachLetter{1}{27}{\defineul}
\foreachletter{1}{27}{\defineul}

\def\autofmt@n#1\autofmt@end{\mathrm{#1}}
\def\autofmt@b#1\autofmt@end{\mathbf{#1}}
\def\autofmt@l#1#2\autofmt@end{\mathbb{#1}\mathsf{#2}}
\def\autofmt@c#1#2\autofmt@end{\mathcal{#1}\mathit{#2}}
\def\autofmt@s#1#2\autofmt@end{\mathscr{#1}\mathit{#2}}
\def\autofmt@f#1\autofmt@end{\mathsf{#1}}
\def\autofmt@u#1\autofmt@end{\underline{\smash{\mathsf{#1}}}}
\def\autofmt@U#1\autofmt@end{\underline{\underline{\smash{\mathsf{#1}}}}}
\def\autofmt@h#1\autofmt@end{\widehat{#1}}
\def\autofmt@r#1\autofmt@end{\overline{#1}}
\def\autofmt@t#1\autofmt@end{\widetilde{#1}}
\def\autofmt@k#1\autofmt@end{\check{#1}}

\def\auto@drop#1{}
\def\autodef#1{\ea\ea\ea\@autodef\ea\ea\ea#1\ea\auto@drop\string#1\autodef@end}
\def\@autodef#1#2#3\autodef@end{%
  \ea\def\ea#1\ea{\ea\ensuremath\ea{\csname autofmt@#2\endcsname#3\autofmt@end}\xspace}}
\def\autodefs@end{blarg!}
\def\autodefs#1{\@autodefs#1\autodefs@end}
\def\@autodefs#1{\ifx#1\autodefs@end%
  \def\autodefs@next{}%
  \else%
  \def\autodefs@next{\autodef#1\@autodefs}%
  \fi\autodefs@next}

\DeclareSymbolFont{bbold}{U}{bbold}{m}{n}
\DeclareSymbolFontAlphabet{\mathbbb}{bbold}

\newcommand{\bbone}{\ensuremath{\mathbbb{1}}\xspace}

\mdef\delbar{\overline{\partial}}

\mdef\hf{\textstyle\frac12 }
\mdef\thrd{\textstyle\frac13 }
\mdef\qtr{\textstyle\frac14 }

\SelectTips{cm}{}
\newdir{ >}{{}*!/-10pt/@{>}}    

\mdef\Id{\mathrm{Id}}
\mdef\id{\mathrm{id}}
\alwaysmath{ell}
\alwaysmath{infty}
\alwaysmath{odot}
\def\frc#1/#2.{\frac{#1}{#2}}   
\mdef\ten{\mathrel{\otimes}}

\mdef\sqten{\mathrel{\boxtimes}}

\DeclareRobustCommand\widecheck[1]{{\mathpalette\@widecheck{#1}}}
\def\@widecheck#1#2{%
    \setbox\z@\hbox{\m@th$#1#2$}%
    \setbox\tw@\hbox{\m@th$#1%
       \widehat{%
          \vrule\@width\z@\@height\ht\z@
          \vrule\@height\z@\@width\wd\z@}$}%
    \dp\tw@-\ht\z@
    \@tempdima\ht\z@ \advance\@tempdima2\ht\tw@ \divide\@tempdima\thr@@
    \setbox\tw@\hbox{%
       \raise\@tempdima\hbox{\scalebox{1}[-1]{\lower\@tempdima\box
\tw@}}}%
    {\ooalign{\box\tw@ \cr \box\z@}}}

\mdef\we{\overset{\sim}{\longrightarrow}}
\mdef\leftwe{\overset{\sim}{\longleftarrow}}

\let\xto\xrightarrow

\def\rightarrowtailfill@{\arrowfill@{\Yright\joinrel\relbar}\relbar\rightarrow}
\newcommand\xrightarrowtail[2][]{\ext@arrow 0055{\rightarrowtailfill@}{#1}{#2}}

\def\twoheadrightarrowfill@{\arrowfill@{\relbar\joinrel\relbar}\relbar\twoheadrightarrow}
\newcommand\xtwoheadrightarrow[2][]{\ext@arrow 0055{\twoheadrightarrowfill@}{#1}{#2}}

\def\slashedarrowfill@#1#2#3#4#5{%
  $\m@th\thickmuskip0mu\medmuskip\thickmuskip\thinmuskip\thickmuskip
   \relax#5#1\mkern-7mu%
   \cleaders\hbox{$#5\mkern-2mu#2\mkern-2mu$}\hfill
   \mathclap{#3}\mathclap{#2}%
   \cleaders\hbox{$#5\mkern-2mu#2\mkern-2mu$}\hfill
   \mkern-7mu#4$%
}
\def\rightslashedarrowfill@{%
  \slashedarrowfill@\relbar\relbar\mapstochar\rightarrow}
\newcommand\xslashedrightarrow[2][]{%
  \ext@arrow 0055{\rightslashedarrowfill@}{#1}{#2}}
\mdef\hto{\xslashedrightarrow{}}
\mdef\htoo{\xslashedrightarrow{\quad}}

\def\toiso{\xto{\smash{\raisebox{-.5mm}{$\scriptstyle\sim$}}}}

\long\def\my@drawfill#1#2;{%
\@skipfalse
\fill[#1,draw=none] #2;
\@skiptrue
\draw[#1,fill=none] #2;
}
\newif\if@skip
\newcommand{\skipit}[1]{\if@skip\else#1\fi}
\newcommand{\drawfill}[1][]{\my@drawfill{#1}}

\newif\ifhyperref
\@ifpackageloaded{hyperref}{\hyperreftrue}{\hyperreffalse}
\iftac
  \let\your@state\state
  \def\state#1{\gdef\currthmtype{#1}\your@state{#1}}
  \let\your@staterm\staterm
  \def\staterm#1{\gdef\currthmtype{#1}\your@staterm{#1}}
  \let\defthm\newtheorem
  \def\currthmtype{}
  \ifhyperref
    \def\autoref#1{\ref*{label@name@#1}~\ref{#1}}
  \else
    \def\autoref#1{\ref{label@name@#1}~\ref{#1}}
  \fi
  \AtBeginDocument{%
    \let\old@label\label%
    \def\label#1{%
      {\let\your@currentlabel\@currentlabel%
        \edef\@currentlabel{\currthmtype}%
        \old@label{label@name@#1}}%
      \old@label{#1}}
  }
\else
  \ifhyperref
    \def\defthm#1#2{%
      \newtheorem{#1}{#2}[section]%
      \expandafter\def\csname #1autorefname\endcsname{#2}%
      \expandafter\let\csname c@#1\endcsname\c@thm}
  \else
    \def\defthm#1#2{\newtheorem{#1}[thm]{#2}}
    \ifx\SK@label\undefined\let\SK@label\label\fi
    \let\old@label\label
    \let\your@thm\@thm
    \def\@thm#1#2#3{\gdef\currthmtype{#3}\your@thm{#1}{#2}{#3}}
    \def\currthmtype{}
    \def\label#1{{\let\your@currentlabel\@currentlabel\def\@currentlabel%
        {\currthmtype~\your@currentlabel}%
        \SK@label{#1@}}\old@label{#1}}
    \def\autoref#1{\ref{#1@}}
  \fi
\fi

\newtheorem{thm}{Theorem}[section]

\defthm{cor}{Corollary}
\defthm{prop}{Proposition}
\defthm{lem}{Lemma}
\defthm{sch}{Scholium}
\defthm{assume}{Assumption}
\defthm{claim}{Claim}
\defthm{conj}{Conjecture}
\defthm{hyp}{Hypothesis}
\defthm{fact}{Fact}
\defthm{quest}{Question}
\iftac\theoremstyle{plain}\else\theoremstyle{definition}\fi
\defthm{defn}{Definition}
\defthm{notn}{Notation}
\iftac\theoremstyle{plain}\else\theoremstyle{remark}\fi
\defthm{rmk}{Remark}
\defthm{eg}{Example}
\defthm{egs}{Examples}
\defthm{ex}{Exercise}
\defthm{ceg}{Counterexample}
\defthm{con}{Construction}
\defthm{conv}{Convention}
\defthm{warn}{Warning}
\defthm{digr}{Digression}
\defthm{per}{Perspective}
\defthm{disc}{Discussion}
\defthm{term}{Terminology}
\defthm{heur}{Heuristics}

\def\thmqedhere{\expandafter\csname\csname @currenvir\endcsname @qed\endcsname}

\setitemize[1]{leftmargin=2em}
\setenumerate[1]{leftmargin=*}

\iftac
  \let\c@equation\c@subsection
\else
  \let\c@equation\c@thm
\fi
\numberwithin{equation}{section}

\@ifpackageloaded{mathtools}{\mathtoolsset{showonlyrefs,showmanualtags}}{}

\alwaysmath{alpha}
\alwaysmath{beta}
\alwaysmath{gamma}
\alwaysmath{Gamma}
\alwaysmath{delta}
\alwaysmath{Delta}
\alwaysmath{epsilon}
\mdef\ep{\varepsilon}
\alwaysmath{zeta}
\alwaysmath{eta}
\alwaysmath{theta}
\alwaysmath{Theta}
\alwaysmath{iota}
\alwaysmath{kappa}
\alwaysmath{lambda}
\alwaysmath{Lambda}
\alwaysmath{mu}
\alwaysmath{nu}
\alwaysmath{xi}
\alwaysmath{pi}
\alwaysmath{rho}
\alwaysmath{sigma}
\alwaysmath{Sigma}
\alwaysmath{tau}
\alwaysmath{upsilon}
\alwaysmath{Upsilon}
\alwaysmath{phi}
\alwaysmath{Pi}
\alwaysmath{Phi}
\mdef\ph{\varphi}
\alwaysmath{chi}
\alwaysmath{psi}
\alwaysmath{Psi}
\alwaysmath{omega}
\alwaysmath{Omega}

\makeatother

\tikzset{lab/.style={auto,font=\scriptsize}} 
\usetikzlibrary{arrows}

\usepackage[all]{xy}

\usepackage{xcolor}

\definecolor{fxnote}{rgb}{1.0000,0.0000,0.0000}
\colorlet{fxnotebg}{yellow}

\newcommand{\D}{\sD}

\autodefs{\cDER\cCat\cGrp\cCAT\cDir\cMONCAT\bSet\cProf\bCat
\cPro\cMONDER\cBICAT\ncomp\nid\ncoker\niso\nIm\sProf\ncyl\fC\fZ\fT\nhocolim\nholim\cPsNat\ncoll\nCh\bsSet\cAlg\cIm\cI\cP}

\let\oldboxtimes\boxtimes
\def\boxtimes{\mathrel{\oldboxtimes}}

\newcommand{\fib}{\mathsf{fib}}
\newcommand{\cof}{\mathsf{cof}}

\def\ccsub{_{\mathrm{cc}}}
\def\pdh(#1,#2){\llbracket #1,#2\rrbracket}
\def\ldh(#1,#2){\llbracket #1,#2\rrbracket\ccsub}
\def\pend(#1){\pdh(#1,#1)}
\def\lend(#1){\ldh(#1,#1)}

\def\DTl#1#2#3#4#5#6#7{%
  \xymatrix@C=3pc{{#1} \ar[r]^-{#2} &
    {#3} \ar[r]^-{#4} &
    {#5} \ar[r]^-{#6} &
    {#7}
  }}

\newsavebox{\tvabox}
\savebox\tvabox{\hspace{1mm}\begin{tikzpicture}[>=latex',baseline={(0,-.18)}]
  \draw[->] (0,.1) -- +(1,0);
  \node at (.5,0) {$\scriptscriptstyle\bot$};
  \draw[->] (1,-.1) -- +(-1,0);
  \draw[->] (1,-.2) -- +(-1,0);
\end{tikzpicture}\hspace{1mm}}

\newcommand{\tcof}{\mathsf{tcof}}
\newcommand{\tfib}{\mathsf{tfib}}

\newcommand{\cube}[1]{\square^{#1}}

\newcommand{\A}[1]{{\Vec{A}_{#1}}}

\title{The bivariant parasimplicial $\mathsf{S}_{\bullet}$-construction}
\author{Falk Beckert}

\thanks{%
This research was conducted in the framework of the research training group
\emph{GRK 2240: Algebro-Geometric Methods in Algebra, Arithmetic and Topology} and of the grant HO 4729/2-1, which are funded by the DFG
}

\begin{document}

\begin{abstract}
Coherent strings of composable morphisms play an important role in various important constructions in abstract stable homotopy theory (for example algebraic K-theory or higher Toda brackets) and in the representation theory of finite dimensional algebras (as representations of Dynkin quivers of type A). In a first step we will prove a strong comparison result relating composable strings of morphisms and coherent diagrams on cubes with support on a path from the initial to the final object.

We observe that both structures are equivalent (by passing to higher analogues of mesh categories) to distinguished coherent diagrams on special classes of morphism objects in the 2-category of parasimplices. Furthermore, we show that the notion of distinguished coherent diagrams generalizes well to arbitrary morphism objects in this 2-category. The resulting categories of coherent diagrams lead to higher versions of the $\mathsf{S}_{\bullet}$-construction and are closely related to representations of higher Auslander algebras of Dynkin quivers of type A.

Understanding these categories and the functors relating them in general will require a detailed analysis of the 2-category of parasimplices as well as basic results from abstract cubical homotopy theory (since subcubes of distinguished diagrams very often turn out to be bicartesian). Finally, we show that the previous comparison result extends to a duality theorem on general categories of distinguished coherent diagrams, as a special case leading to some new derived equivalences between higher Auslander algebras.
\end{abstract}

\maketitle

\tableofcontents

\section{Introduction}

\label{sec:intro}

The $\mathsf{S}_{\bullet}$-construction is an important tool for defining algebraic $K$-theory spectra in a general context. Since their first construction for Waldhausen categories \cite{waldhausen:k-theory}, several variants have been considered for other models of higher homotopy theories, like, for instance $\infty$-categories \cite{BGT,lurie:rotation} or derivators \cite{garkusha:I,gst:Dynkin-A}. 

Before giving an overview of the $\mathsf{S}_{\bullet}$-construction, let us remark that we will work with derivators as a model for homotopy theories. This theory has the advantage that for a derivator $\D$ the passage to the homotopy category is given by evaluation at the final category $\D(\bbone)$, whereas evaluation at an arbitrary category $A$ is related to $A$-diagrams in the underlying homotopy category via an underlying diagram functor
\begin{equation}\label{eq:u-dia}
\D(A)\rightarrow\D(\bbone)^A.
\end{equation}

Moreover, this work is part of a project on abstract cubical homotopy theory and global Serre dualities \cite{bg:cubical}, \cite{bg:global} which share the same notations and conventions. Several of our results might be relevant for further research in the context of algebraic $K$-theory. It is known that there cannot be a satisfying construction of algebraic $K$-theory on the 2-category of derivators \cite{muro-raptis:note},\cite{muro-raptis:revisited}, but we highly expect that analogues of the results presented here hold true in the setting of $\infty$-categories. In fact, in \cite[Prop.~9.32]{BGT} it was shown that $K$-theory is a stable invariant. Therefore we will restrict to stable derivators and emphasize that this restriction is required for most of the results presented here.

Given a derivator $\D$, the $\mathsf{S}_{\bullet}$-construction $\mathsf{S}_{\bullet}\D$ is a simplicial object, such that the $n$th level $\mathsf{S}_n\D$ is given by a subderivator $\D^{\mathrm{Ar}[n],ex}$ of the shifted derivator $\D^{\mathrm{Ar}[n]}$ (i.e. presheaves on the arrow category of the $n$-simplex) spanned by objects characterized by certain vanishing and cocartesianess conditions. We observe that the category $\mathrm{Ar}[n]$ sits in the following sequence of embeddings of categories
\begin{equation}\label{eq:meshes}
[n-1] \rightarrow \mathrm{Ar}[n] \rightarrow M_{n-1},
\end{equation}
where the left map is defined by $i\mapsto(0,i+1)$, and the right map is the embedding into the mesh category of the $(n-1)$-simplex \cite[Thm.~4.5]{gst:Dynkin-A}. By considering exponentials and restricting to subobjects, we obtain for every stable derivator a sequence of equivalences in the 2-category of derivators, which are induced by restriction morphisms:
\begin{equation}\label{eq:Sdot-covar}
\D^{[n-1]} \xleftarrow{\simeq} \D^{\mathrm{Ar}[n],ex} \xleftarrow{\simeq} \D^{M_{n-1},ex}.
\end{equation}
Since the equivalences are compatible with the simplicial structure, we in fact obtain equivalences of simplicial objects, thus each of them gives an alternative description of the $\mathsf{S}_{\bullet}$ construction. We refer to the left version as the \textbf{slice model} and to the right one as the \textbf{symmetric model}. Both models have advantages in different situations:
\begin{enumerate}
\item The slice model immediately shows that the simplicial structure of the $\mathsf{S}_{\bullet}$-construction can be described by inverse image functors. Moreover, it is often easier to define morphisms from or to the slice model. And since the categories $[n]$ are homotopically finite, it will follow immediately in many cases, that those morphisms admit adjoints.
\item The categories $M_n$ admit non-trivial symmetries, which carry over to the symmetric model. Moreover, the condition $ex$ forces that one of these symmetries is naturally isomorphic to the suspension functor and a certain composite of those symmetries defines a Serre autoequivalence \cite[Thm.~11.12]{gst:Dynkin-A}.
\end{enumerate}
The interplay between those two pictures is described in detail in \cite{gst:Dynkin-A} and has been proven to be useful in abstract representation theory and abstract homotopy theory.
\begin{itemize}
\item In the special case $\D=\D_{\mathrm{k}}$ of the derivator of a field $\mathrm{k}$, we can identify the value $\D_{\mathrm{k}}(A)\cong D(\mathrm{k}A)$ with the derived category of modules over the category algebra $\mathrm{k}A$. In particular, in the case where $A=[n-1]$ is the $(n-1)$-simplex we obtain 
\[
\D_{\mathrm{k}}([n-1])\cong D(\mathrm{k}\A{n})
\]
the derived category of the path algebra of the $\A{n}$-quiver with linear orientation as an value of the derivator $\D_{\mathrm{k}}$. Hence (for $\D$ general) we can regard $\D^{M_{n-1},ex}$ as a derivator of $\D$-representations of $\A{n}$. In fact, $\D$-representations of general $A_n$-quivers can be obtained from  $\D^{M_{n-1},ex}$ via restriction functors \cite[Thm.~4.14]{gst:Dynkin-A}. Furthermore, one of the autoequivalences induced by the symmetries of $M_{n-1}$ can be identified with the Auslander-Reiten-translation.\\
Since we will work for simplicity with derivators parametrized by all small categories, we obtain unbounded versions of derived categories. In many situations it will be more convenient to work with bounded derived categories, which are typically values of derivators with smaller domains. It is straight forward to generalize our results to this setting. 
\item The simplicial structure morphisms determine a canonical strong triangulation on the underlying category $\D(\bbone)$ \cite[Thm.~13.6]{gst:Dynkin-A}. Hence we obtain compatibility relations for iterated fiber and cofiber constructions generalizing the octahedral axiom as described in detail in \S\ref{sec:triangles}. Since fiber and cofiber constructions belong to the most fundamental operations available in stable homotopy theories, we obtain insight to the general structure of stable homotopy theories.
\end{itemize}

A crucial observation is that these two features of the $\mathsf{S}_{\bullet}$-construction are compatible. The simplicial structure and the equivalences related to the Auslander-Reiten-translation assemble to a parasimplicial structure \cite[Rmk.~4.3.6]{lurie:rotation}, which can be described as follows: the mesh categories are canonically (once we have chosen coordinates on parasimplices \autoref{prop:coord-mor}) isomorphic to certain morphism objects in the 2-category of parasimplices
$$M_n\cong\underline{\Lambda}(\Lambda_1,\Lambda_{n+1}),$$
thus they can be regarded as a parasimplicial version of arrow categories. In \S\ref{sec:vertical} we will make precise in which way the parasimplicial structure on the $\mathsf{S}_{\bullet}$-construction is induced by the functor $\underline{\Lambda}(\Lambda_1,-)$.

For reasons to become clear in a moment, we will now give another model for the $\mathsf{S}_{\bullet}$-construction (c.f. \autoref{thm:dual-Sdot-sl} and \autoref{thm:dual-Sdot}), which is a generalization of the derived equivalence of the $\A{n}$-quiver and the $\A{n}$-quiver with zero relations \cite[Prop.~2.1]{HS-piecewise} to the setting of stable derivators.

\begin{thm}\label{thm:max-path}
Let $\D$ be a stable derivator and $n \geq 0$. Then there are equivalences of stable derivators:
$$\D^{\cube{n+1},ex}_{\rightarrow}\toiso\D^{\cube{n}}_{\rightarrow}\toiso\D^{[n]}.$$
\end{thm}

The conditions $\rightarrow$ and $ex$ define certain full subderivators of the derivator of coherent diagrams of cubical shape. More precisely, a coherent diagram satisfies condition $\rightarrow$, if the support of the diagram is concentrated in a chosen path linking the initial and the final vertex of the cube, and condition $ex$ when the diagram is additionally cocartesian. It is clear that the simplicial structure on the derivators $\D^{[n]}$ transfers to the derivators $\D^{\cube{n}}_{\rightarrow}$ via conjugation with the above equivalences, giving rise to a pseudofunctor $\Delta^{op}\rightarrow Der,[n]\mapsto\D^{\cube{n}}_{\rightarrow}$. Moreover, we observe that this cubical version of the $\mathsf{S}_{\bullet}$-construction also admits a symmetric model which is in fact also related to morphism objects in the 2-category of parasimplices. Analogously to \eqref{eq:Sdot-covar} we have the following equivalences of stable derivators, induced by restriction morphisms
\begin{equation}\label{eq:Sdot-contra}
\D^{\cube{n}}_{\rightarrow}\xleftarrow{\simeq}\D^{\cube{n+1},ex}_{\rightarrow}\xleftarrow{\simeq}\D^{\underline{\Lambda}(\Lambda_{n},\Lambda_{n+1}),ex}.
\end{equation}
It is worth observing that on $\D^{\underline{\Lambda}(\Lambda_1,\Lambda_{n+1})}=\D^{M_n}$ and $\D^{\underline{\Lambda}(\Lambda_{n},\Lambda_{n+1})}$ the conditions $ex$ are special cases of a more general requirement.

The constructions \eqref{eq:Sdot-covar} and \eqref{eq:Sdot-contra} suggest that certain subderivators of $\D^{\underline{\Lambda}(\Lambda_k,\Lambda_n)}$ with $k\leq n$ should be closely related to the $\mathsf{S}_{\bullet}$-construction. We observe that $\underline{\Lambda}(\Lambda_k,\Lambda_n)$ is a subposet of $\mathbb{Z}^{k+1}$, and we say that a coherent diagram $X$ in $\D^{\underline{\Lambda}(\Lambda_k,\Lambda_n)}$ satisfies condition $ex$ if the support of $X$ is concentrated in the subposet of injective maps of parasimplices and moreover, if all $(k+1)$-subcubes of $X$ which are compatible with the embedding into $\mathbb{Z}^{k+1}$ are cocartesian. Then we can show that this condition exactly specializes to the conditions $ex$ appearing in \eqref{eq:Sdot-covar} and \eqref{eq:Sdot-contra}. For this reason we associate to $\D$ and a choice of $n\geq 1$ and $k\geq 2$ another stable derivator
$$\D_{n,k}:=\D^{\underline{\Lambda}(\Lambda_{k-1},\Lambda_{n+k-1}),ex}.$$
We remark, that the subcategory of $\underline{\Lambda}(\Lambda_{k-1},\Lambda_{n+k-1})$ spanned by injective maps of parasimplices is canonically isomorphic to $\underline{\Lambda}(\Lambda_{k-1},\Lambda_{n-1})$, revealing the contravariant nature of simplicial structure of \eqref{eq:Sdot-contra}, and thereby justifying the title \textbf{bivariant parasimplicial} $\mathsf{S}_{\bullet}$-\textbf{construction}.

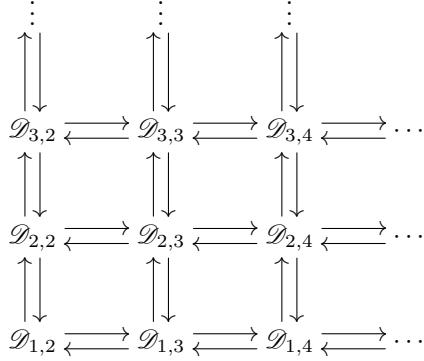
\begin{figure}
\begin{displaymath}
\xymatrix{
\vdots \ar@<1mm>[d] \ar@{<-}@<-1mm>[d]& \vdots \ar@<1mm>[d] \ar@{<-}@<-1mm>[d] & \vdots \ar@<1mm>[d] \ar@{<-}@<-1mm>[d] \\
\D_{3,2} \ar@<1mm>[r] \ar@{<-}@<-1mm>[r] \ar@<1mm>[d] \ar@{<-}@<-1mm>[d]& \D_{3,3} \ar@<1mm>[r] \ar@{<-}@<-1mm>[r] \ar@<1mm>[d] \ar@{<-}@<-1mm>[d] & \D_{3,4} \ar@<1mm>[r] \ar@{<-}@<-1mm>[r] \ar@<1mm>[d] \ar@{<-}@<-1mm>[d] & \cdots\\
\D_{2,2} \ar@<1mm>[r] \ar@{<-}@<-1mm>[r] \ar@<1mm>[d] \ar@{<-}@<-1mm>[d]& \D_{2,3} \ar@<1mm>[r] \ar@{<-}@<-1mm>[r] \ar@<1mm>[d] \ar@{<-}@<-1mm>[d] & \D_{2,4} \ar@<1mm>[r] \ar@{<-}@<-1mm>[r] \ar@<1mm>[d] \ar@{<-}@<-1mm>[d] & \cdots\\
\D_{1,2} \ar@<1mm>[r] \ar@{<-}@<-1mm>[r] & \D_{1,3} \ar@<1mm>[r] \ar@{<-}@<-1mm>[r] & \D_{1,4} \ar@<1mm>[r] \ar@{<-}@<-1mm>[r] & \cdots
}
\end{displaymath}
\caption{The bivariant $\mathsf{S}_{\bullet}$-construction. The arrows indicate infinite chains of adjunctions, the first column is equivalent to \eqref{eq:Sdot-covar}, and the bottom row is equivalent to \eqref{eq:Sdot-contra}.}
\label{fig:intro}
\end{figure}

The main goal of this work, is to analyze in detail the derivators $\D_{n,k}$, as well as the functorialities between them. In a slightly different context, it was already observed in \cite{poguntke-Sdot}, \cite{dyckerhoff-Sdot}, that the columns of \autoref{fig:intro} can be regarded as higher dimensional $\mathsf{S}_{\bullet}$-constructions. Complementary to this we will focus on understanding the interaction of the horizontal and vertical structures in \autoref{fig:intro}, which turns out to be completely symmetric. 

We now outline the internal structure of this work. It consists of four main parts.
\begin{itemize}
\item A short introduction into the most relevant preliminary results and techniques (\S\S\ref{sec:subder}-\ref{sec:cubical}).
\item The definition of the constituents of the bivariant parasimplicial $\mathsf{S}_{\bullet}$-construc-tion and fundamental relations amongs these (\S\S\ref{sec:dnk}-\ref{sec:recollements}).
\item The proof of the symmetry of the bivariant parasimplicial $\mathsf{S}_{\bullet}$-construction (\S\S\ref{sec:contra}-\ref{sec:main}).
\item Applications of the main theorem on higher Toda brackets and higher triangulations (\S\S\ref{sec:Toda}-\ref{sec:triangles}).
\end{itemize}

More precisely, we set up our basic notations and conventions in \S\ref{sec:subder}. Furthermore, we explain in detail in which situations we use the term 'full subderivator'. We consider a few examples, which will become relevant in later sections and point out some useful consequences of this notion.

In \S\ref{sec:para} we review in detail the 2-category $\underline{\Lambda}$ of parasimplices. We show that this 2-category is adjunction complete and describe explicitly how to construct the left (resp. right) adjoints. This leads to the observation that both operations are related by symmetry operations. We introduce \emph{choices of coordinates} on $\underline{\Lambda}$ and use these, to construct embeddings of the simplex category $\Delta$ and duality 2-functors $\underline{\Lambda}^{co}\to\underline{\Lambda}$.

As a preparation for the following chapters, we recollect some elementary facts on the stable homotopy theory of cubes in \S\ref{sec:cubical}. This summarizes parts of \cite{bg:cubical}, although some of the results are stated in a slightly stronger form (due to the earlier use of the stability assumption). 

In \S\ref{sec:dnk} we analyze the properties of the derivators $\D_{n,k}$ for fixed $n$ and $k$. In particular, we show that the symmetries on $\underline{\Lambda}(\Lambda_{k-1},\Lambda_{n+k-1})$ induced by pre- and postcomposition with the paracyclic translation restrict to $\D_{n,k}$ (\autoref{cor:dnk-sym}). We show that they are related to powers of the suspension (\autoref{cor:dnk-sigma}), and that a certain combination of them defines an autoequivalence $\mathsf{s}_3$ which has similar properties as Serre morphisms for $A_n$-quivers (\autoref{cor:fracCY}). Since Serre morphisms and suspension become equal after passing to certain powers, we see that the derivators $\D_{n,k}$ satisfy the fractionally Calabi-Yau property \cite[Thm.~5.19]{gst:Dynkin-A}. Even in this generality we are able to define a slice model (\autoref{thm:slices}), which allows us to use techniques which usually only apply to homotopically finite shifts of derivators. Since the symmetry operations are usually hard to describe on the slice model, we also establish a third model (which to some extend corresponds to the middle terms in \eqref{eq:Sdot-covar} and \eqref{eq:Sdot-contra}). This 'domain model' has the advantage that on the one hand many structure morphisms are still accessible and on the other hand the important symmetry operation $\mathsf{s}_3$ is isomorphic to an inverse image morphism in this situation (\autoref{rmk:slices}). This will turn out to be useful for understanding the interaction of various structure morphisms with symmetry operations (e.g. see \S\ref{sec:horizontal} or \S\ref{sec:contra}). Furthermore, it follows from Iyama's inductive construction of the higher Auslander algebras of $\A{n}$-quivers \cite{iyama-higherA}, that the underlying category of $\D_{n,k}$ is equivalent to the derived category of the (k-1)-Auslander algebra of the $\A{n+1}$-quiver. Hence our results on the structure of the derivators $\D_{n,k}$ can be regarded as a contribution to abstract representation theory. We emphasize that the relation between (para)simplicial combinatorics and higher Auslander algebras of type $A$ is discussed independently in \cite{DJW-Sdot}.

In \S\ref{sec:vertical} we show that the vertical adjoint morphisms in \autoref{fig:intro} are in fact inverse images of the postcomposition functors in $\underline{\Lambda}$. In particular, it follows that many results on the structure of $\underline{\Lambda}$ carry over to these adjunctions. As a consequence, it is immediate that the vertical structure morphisms are compatible with the symmetry operations (\autoref{thm:gsd-Sdot}). As a preparation for the following sections, we show that most of the vertical face morphisms can be described as inverse images on slices. The corresponding statement also holds for vertical degeneracy morphisms if we pass to a slightly larger version of the slice construction (\autoref{cor:slice-restr}).

Unfortunately, the situation in the horizontal direction is significantly more complicated. In \S\ref{sec:horizontal} we construct a fundamental adjoint triple in the horizontal direction of \autoref{fig:intro} by using the slice model. We extend these constructions to the domain model and invoke the techniques established in \S\ref{sec:dnk} to proof the compatibility to the symmetries $\mathsf{s}_3$ (\autoref{thm:Serre-h}). Building on this, we show that the above triples of adjunctions extend to infinite chains of adjunctions and define general horizontal structure morphisms via this infinite chain. However, it seem to be hard to give a direct proof that the horizontal structure morphisms also satisfy the (para)simplicial relations. This statement will be postponed until \S\ref{sec:main} as a corollary of the main theorem.

In \S\ref{sec:recollements} we prove a compatibility result relating the structure morphisms in the horizontal and vertical direction. This will rely on the characterizations of the structure morphisms as inverse images we have established in \S\ref{sec:vertical} and \S\ref{sec:horizontal}. More precisely, if we consider a subsquare
\[
\xymatrix{
\D_{n,k-1} \ar[r] \ar[d] & \D_{n,k} \ar[d]\\
\D_{n-1,k-1} \ar[r] & \D_{n-1,k}
}
\]
of \autoref{fig:intro}, we show that any composition of structure morphisms through the bottom left vertex can be realized as a composition of structure morphisms through the top right vertex. But we observe that there is another composition of structure morphisms through the top right vertex (which is unique up to symmetry), and prove that this composition is isomorphic to the zero-morphism, giving rise to recollements of stable derivators
\[
\xymatrix{
\D_{n,k-1} \ar@{<-}@<1.5mm>[r] \ar@{<-}@<-1.5mm>[r] \ar[r] & \D_{n,k} \ar@{<-}@<1.5mm>[r] \ar@{<-}@<-1.5mm>[r] \ar[r] & \D_{n-1,k}.
}
\]

Paragraph \S\ref{sec:contra} is devoted to the proof of \autoref{thm:max-path}, followed by a closer investigation of the equivalence in a specific example.  Moreover, using the domain model we show that the equivalences from \autoref{thm:max-path} are compatible with the symmetry operation $\mathsf{s}_3$ (\autoref{thm:dual-Serre}). Finally, we show that the equivalences of \autoref{thm:max-path} are in a certain sense self-dual (\autoref{prop:kappa}). This statement is the final important ingredient for the main result in \S\ref{sec:main}, which can be summarized as follows.

There are equivalences of stable derivators
$$\Phi_{n,k}: \D_{n,k}\rightarrow\D_{k-1,n+1}$$ such that:
\begin{enumerate}
\item $\Phi_{n,k}\circ\Phi_{k-1,n+1}\cong id;$
\item the $\Phi_{n,k}$ are compatible with symmetries;
\item they map horizontal structure morphisms to vertical structure morphisms and vice versa;
\item they specialize to the equivalences of \autoref{thm:max-path} in the cases $n=1$ or $k=2$.
\end{enumerate}
In the case $\D=\D_{\mathrm{k}}$ of the derivator of a field $\mathrm{k}$, this result can be regarded as a higher dimensional generalization of the derived equivalence between the $\A{n}$-quiver and the $\A{n}$-quiver with zero relations (\autoref{prop:An-zero}) and, in particular, provides derived equivalences between the $(k-1)$-Auslander algebra of the $\A{n}$-quiver and the $(n-1)$-Auslander algebra of the $\A{k}$-quiver (cf. \autoref{rmk:higher-Auslander}). To the best of the author's knowledge this seems to be a new result even in the case over a field. We emphasize that even in the case $k=n+1$ in which the equivalences $\Phi$ are autoequivalences, they are in general neither identities nor induced by the symmetries considered in \S\ref{sec:dnk}. In fact, for the proof of the main result we consider yet another variant of the slice model which is related more closely to cubical shapes. This allows us to define the morphisms $\Phi_{n,k}$ as certain twisted products of the equivalences of \autoref{thm:max-path}. The most involved part of the proof, that the morphisms $\Phi_{n,k}$ are equivalences, relies on an inductive argument building on the recollements considered in \S\ref{sec:recollements}.

We show in \S\ref{sec:Toda}, that \autoref{thm:max-path} can be used to define a functorial version of higher Toda brackets \cite{toda:brackets} in the context of stable derivators. In particular, we define Toda bracket morphisms on a certain classes of coherent diagrams lifting the classical Toda bracket operations along an underlying (incoherent) diagram functor. This shows that the indeterminacy of Toda brackets corresponds to the failure of the underlying diagram functor of being an equivalence.

As another application, we formulate in \S\ref{sec:triangles} how the derivators $\D_{n,k}$ can be used for an axiomatization of higher triangulations. We start by recalling the definition of a strong triangulation and the construction of these structures in the case of a strong stable derivator $\D$ following \cite{gst:Dynkin-A}. This construction relies heavily on the structure of the derivators $\D_{n,2}$. We explain in detail in which way the axioms of a strong triangulation imply higher versions of the octahedral axioms, which in turn can be regarded as compatibility relations for cofiber sequences. On the other hand, we indicate that the derivators $\D_{1,k}$ for $k\geq3$ encode higher analogues of cofiber sequences and expanding on this we explain how the derivators $\D_{n,k}$ potentially give rise to analogues of higher octahedral axioms in a different way. They encode higher compatibility relations for higher cofiber sequences. However, the underlying diagram functors related to the derivators $\D_{n,k}$ are in important examples neither full nor essentially surjective. A counterexample to the fullness was discussed in \cite{bg:cubical}. Here we show that non-vanishing Toda brackets are obstructions to the essential surjectivity. Therefore, the compatibility relations of higher cofiber sequences cannot descent to a well-behaved axiomatization at the level of underlying categories, and hence rely on having a stable derivator or another model for homotopy theories, which has sufficient information on coherent diagrams. We also indicate some relations to the $n$-angulated categories of \cite{GKO-angulated}.

Finally, in the short appendix we describe a general opertaion on pseudofunctors by conjugation with pointwise equivalences.
\newpage

\begin{center}
{\scshape Acknowledgments \par}
\vspace{0.1cm}
\end{center}

I would like to express my gratitude to my advisor Jens Hornbostel for his support and encouragement, his patience, and his advise. I would like to thank Moritz Groth for countless discussions on stable derivators and related topics. I also wish to thank Tobias Dyckerhoff, Gustavo Jasso and Jan {\v S}{\v t}ov{\'\i}{\v c}ek for answering my questions on representation theory, and Thomas Hudson, Sven Stahn and Sean Tilson for many fruitful discussions. I am very thankful to my family for their continuous support. Finally, I gratefully acknowledge the funding (HO 4729/2-1) received from the German Science Foundation (DFG).

This research was conducted in the framework of the research training group
\emph{GRK 2240: Algebro-Geometric Methods in Algebra, Arithmetic and Topology},
which is funded by the DFG. 

This text is a modified version of the authors doctoral thesis.

\section{Preliminaries on derivators}
\label{sec:subder}

In this section we set up the basic notation and conventions used in this work. Moreover, we make the notion of a full subderivator precise. This will be useful to analyze quivers with zero-relations in the context of abstract representation theory.

\begin{itemize}
\item We denote by $Cat$ the 2-category of small categories and by $CAT$ the 2-category of not necessarily small categories.
\item The terminal category is denoted by $\bbone\in Cat$.
\item Let $A\in Cat$, $a\in A$. Then we denote by $a\colon\bbone\to A$ the unique functor with image $a\in A$.
\item A small category $A$ is called a poset if the cardinality of the set $Hom(a,a')$ is at most 1 for all $a,a'\in A$ and if all isomorphisms in $A$ are identities.
\item For a small category $A$ we denote by $\emptyset$, respectively $\infty$, the initial, respectively final, object in $A$, provided that these objects exist in $A$. Moreover, let $A$ be a category such that $\emptyset,\infty\in A$. If the unique morphism $\emptyset\to\infty$ is an isomorphism, we call $A$ a pointed category and $0:=\emptyset\cong\infty$ a zero object.
\item For a functor $u\colon A\to B$ and $b\in B$ we denote by $A_{b/}$, respectively $A_{/b}$, the slice category of objects $u$-under $b$, respectively of objects $u$-over $b$.  
\item A prederivator is a 2-functor $\D\colon Cat^{op} \to CAT$ and a pseudonatural transformation between prederivators is called a morphism of prederivators.
\item A derivator is a suitable 2-functor $\D\colon Cat^{op} \to CAT$. We refer to \cite{groth:ptstab} for a precise definition of derivators and their pointed and stable variants.
\item Let $\D$ be a derivator and $u\colon A\to B$ a functor between small categories. Then the functor $u^*:=\D(u)\colon\D(B)\to\D(A)$ is called the inverse image along $u$. Since $\D$ is a derivator, the inverse image $u^*$ admits a left and a right adjoint, which we will denote by $u_!$, respectively $u_*$, and call the left, respectively right, Kan extension along $u$.
\item Let $\D,\D'$ be derivators. A morphism of derivators (often called 'a morphism') $\D\to\D'$ is a morphism of prederivators which is compatible with left Kan extensions in the precise sense of \cite[Def.~2.2]{groth:ptstab}. Note, that if $\D,\D'$ are stable derivators, then a morphism of derivators is additionally compatible with homotopically finite right Kan extensions by \cite[Thm.~7.1]{PS-traces}, and is therefore, in particular, exact. We denote by $Der$ the 2-category of derivators and morphisms of derivators, and by $Der^{st}$ the full sub-2-category of stable derivators.
\item Let $\D$ be a derivator and $A\in Cat$ then $\D^A\colon B\mapsto \D(A\times B)$ denotes the derivator shifted by $A$.
\item Let $\D$ be a derivator and $A\in Cat$ then $\mathsf{dia}_A\colon\D(A)\to\D(\bbone)^A$ denotes the underlying diagram functor (c.f. \cite{groth:ptstab}).
\end{itemize}

\begin{con}\label{con:subpre}
Let $\D$ be a prederivator and for $A\in Cat$ full subcategories $\D(A)_0\subseteq\D(A)$ such that for every functor $u\colon A\to B$ between small categories the essential image of $u^*\vert_{\D(B)_0}\colon\D(B)_0\to\D(A)$ is contained in $\D(A)_0$. This yields the existence of a restriction $u^*_0\colon\D(B)_0\to\D(A)_0$ such that the diagram
\begin{equation}\label{eq:restr}
\xymatrix{
\D(B)_0 \ar[r]^{\subseteq}\ar[d]_{u^*_0}&\D(B)\ar[d]^{u^*}\\
\D(A)_0 \ar[r]^{\subseteq}& \D(A)
}
\end{equation}
strictly commutes. Moreover, if $v\colon A\to B$ is a functor and $\alpha\colon u^*\to v^*$ is a natural transformation, then we obtain a natural transformation $\alpha_0\colon u_0^*\to v_0^*$ by restriction, since the fullness of the inclusion $\D(A)_0\subseteq\D(A)$ implies that the map $\alpha_0(X):=\alpha(X)$ for $X\in\D(B)_0$ is a morphism in $\D(A)_0$. This shows that $\D\colon A\mapsto \D(A)_0, u\mapsto u^*_0$ is prederivator and the inclusions $\D_0(A)=\D(A)_0\subseteq\D(A)$ yield a morphism of prederivators $\D_0\xrightarrow{\subseteq}\D$. In this case we call $\D_0$ a full subprederivator of $\D$ and $\D_0\xrightarrow{\subseteq}\D$ the inclusion.
\end{con}

\begin{defn}\label{defn:subder}
Let $\D$ be a derivator and $\D_0\subseteq\D$ a full subprederivator. Then $\D_0$ is called a \textbf{full subderivator} of $\D$ if
\begin{enumerate}
\item $\D_0$ satisfies (Der 1),
\item the essential images of the restricted Kan extensions 
\[
u_!\vert_{\D_0(A)},u_*\vert_{\D_0(A)}\colon\D_0(A)\to\D(B)
\]
are contained in $\D(B)_0$ for all functors $u\colon A\to B$ between small categories $A,B\in Cat$.
\end{enumerate}
\end{defn}

\begin{rmk}
The second condition in the above definition implies that the Kan extensions restrict for all functors $u\colon A \to B$ between small categories to well defined functors $(u_0)_!,(u_0)_*\colon\D_0(A)\to \D_0(B)$. This analogous to the case of inverse images, which was treated in \autoref{con:subpre}.
\end{rmk}

\begin{lem}\label{lem:subder}
Let $\D$ be a derivator and $\D_0\subseteq\D$ be a full subderivator. Then $\D_0$ is a derivator and the inclusion $\D_0\xrightarrow{\subseteq}\D$ is compatible with left and right Kan extensions. Moreover, if $\D$ is pointed or stable then so is $\D_0$.
\end{lem}

\begin{proof} 
The prederivator $\D_0$ satisfies (Der 1) by assumption and (Der 2) follows immediately from the corresponding property of $\D$. Let now $u\colon A\to B$ be a functor between small categories $A,B\in Cat$. For (Der 3) we note that, since the inclusions $\D_0(A)\subseteq\D(A)$ are full, the units and counits of the adjunctions $u_!\dashv u^*\dashv u_*$ give by restriction rise to well defined units and counits of the adjunctions $(u_0)_!\dashv u^*\dashv (u_0)_*$. We note that the triangle equalities are obtained by restricting the corresponding pastings of unrestricted transformations. In the following we will omit the subscript 0 for the functors in the image of $\D_0$. For (Der 4) we have to check that for $b\in B$ the pasting
\[
\xymatrix{
\D_0(\bbone)&&\D_0(A_{/b})\ar[ll]_{p_!} \ar@2{->}[ld]_{\varepsilon}&&\D_0(A) \ar@2{->}[lldd]_{\alpha^*} \ar[ll]_{\pi^*}\\
& & & & & \ar@2{->}[ld]_{\eta}\\
& &\D_0(\bbone)\ar[uu]_{p^*}\ar[lluu]^{\id}&&\D_0(B)\ar[ll]_{b^*}\ar[uu]_{u^*}&&\D_0(A)\ar[ll]_{u_!}\ar[lluu]_{\id}
}
\]
(and its dual) is an isomorphism. Here $\alpha^*$ is obtained by the 2-functoriality of $\D_0$ from the canonical slice square associated to $(u,b)$. For each cell of this pasting we see that it is obtained by restriction from the corresponding cell of the derivator $\D$. For the two triangles this is clear from the proof of (Der 3) and for the middle square since $\D_0\xrightarrow{\subseteq}\D$ is a morphism of prederivators. Therefore, we deduce that the pasting above is an isomorphism from the axiom (Der 4) for $\D$. The same argument applies also to the dual pasting. Let now $\D$ be a pointed derivator. Since $\D_0$ is a derivator we have $\emptyset,\infty\in\D_0(\bbone)$. Since $\D_0(A)\subseteq\D(A)$ is full we conclude $\emptyset\cong\infty\in\D_0(\bbone)$, and hence $\D_0$ is pointed. If $\D$ is stable, then the subcategories $\D(\square)^{cocart}=\D(\square)^{cart}\subseteq\D(\square)$ (cf. \cite[Def.~3.9]{groth:ptstab}) coincide. Since the Kan extensions in $\D_0$ are obtained by restriction from those in $\D$ we have
\[
\D_0(\square)^{cocart}=\D(\square)^{cocart}\cap\D_0(\square)=\D(\square)^{cart}\cap\D_0(\square)=\D_0(\square)^{cart}.
\]
This is exactly the stability for $\D_0$. Next, we show that $\D_0\subseteq\D$ is compatible with left and right Kan extensions. We restrict to the case of left Kan extensions; the dual case is very similar. For $u\colon A\to B$ a functor between small categories $A,B\in Cat$ we consider the pasting
\[
\xymatrix{
\D(B)&&\D(A)\ar[ll]_{u_!} \ar@2{->}[ld]_{\varepsilon}&&\D_0(A) \ar[ll]_{\subseteq}\\
& & & & & \ar@2{->}[ld]_{\eta}\\
& &\D(B)\ar[uu]_{u^*}\ar[lluu]^{\id}&&\D_0(B)\ar[ll]_{\subseteq}\ar[uu]_{u^*}&&\D_0(A)\ar[ll]_{u_!}\ar[lluu]_{\id}.
}
\]
This is obtained by restriction from one of the triangular equalities of $u_!\dashv u^*$ and hence an isomorphism.
\end{proof}

\begin{egs}\label{egs:subders}
\begin{enumerate}
\item Let $\D$ be a pointed derivator, $A$ a small category and $Z\subseteq A$ a set of objects in $A$. Then the full subcategories
$\D^{A,Z}(B)=\lbrace x\in\D^A(B)\vert z^*(x)\cong 0\in\D^A(\bbone)\forall z\in Z\rbrace$ define a full subderivator $\D^{A,Z}$ of $\D^A$.
In fact, in this case both conditions in the definition of a full subderivator follow from the fact that inverse images are always compatible with arbitrary Kan extensions.
\item Let $P\colon\D_0\to\D$ be a fully faithful morphism of derivators. If $P$ additionally is compatible with right Kan extensions, then the inclusion of the essential image of $P$ defines a full subederivator. Typical examples where all these conditions are satisfied are the following.
\begin{itemize}
\item Left Kan extension morphisms (cf. \autoref{rmk:cohoerent}) along fully faithful left adjoint functors, and right Kan extension morphisms along fully faithful left adjoint functors (since these morphisms are isomorphic to inverse image morphisms)
\item Kan extension morphisms along homotopically finite, fully faithful functors for stable derivators \cite[Thm.~2.14]{gs:stable}. This includes in particular the inclusion of the derivator of bicartesian squares into $\D^{\square}$ for a stable derivator $\D$.
\end{itemize}
\end{enumerate}
\end{egs}

\begin{lem}
Let $\D_0\subseteq\D$, $\D_0'\subseteq\D'$ be inclusions of subderivators and $P\colon\D\to\D'$ a morphism of derivators such that  the essential image of $P(A)\vert_{\D_0(A)}$ is contained in $\D_0'(A)$ for all $A\in Cat$. Then $P$ restricts to a morphism of derivators $P\vert_{\D_0}\colon\D_0\to\D_0'$
\end{lem}

\begin{proof}
The pseudonaturality conditions for $P\vert_{\D_0}$ are again (using the fullness of the restrictions for the associativity and unitality constraints) obtained by restriction from those of $P$, hence $P\vert_{\D_0}$ is a morphism of prederivators. To show the compatibility with left Kan extensions, we consider the pasting
\[
\xymatrix{
\D_0'(B)&&\D_0'(A)\ar[ll]_{u_!} \ar@2{->}[ld]_{\varepsilon}&&\D_0(A) \ar@2{->}[lldd]_{\cong}\ar[ll]_{P\vert_{\D_0}(A)}\\
& & & & & \ar@2{->}[ld]_{\eta}\\
& &\D_0'(B)\ar[uu]_{u^*}\ar[lluu]^{\id}&&\D_0(B)\ar[ll]_{P\vert_{\D_0}(B)}\ar[uu]_{u^*}&&\D_0(A)\ar[ll]_{u_!}\ar[lluu]_{\id}.
}
\]
This is an isomorphism, since it is obtained by restriction from an isomorphism. Here we use again the arguments from the proof of \autoref{lem:subder} to show that the triangles are obtained via restriction.
\end{proof}

In the following we use the following additional convention.
\begin{rmk}\label{rmk:cohoerent}
Let $\D$ be a derivator and $u\colon A\to B$ a functor between small categories. Then the inverse image functors along $(u\times C)\colon A\times C\to B\times C$ for $C\in Cat$ define a morphism of derivators $u^*\colon\D^B\to\D^A$ which will be called the inverse image morphisms along $u$. We will work with these morphisms of \textbf{derivators} instead of the functors of \textbf{categories} $u^*\colon\D(B)\to\D(A)$. Note, that we can define similarly a left Kan extension morphism $u_!\colon\D^A\to\D^B$ (\cite[Ex~2.10]{groth:ptstab}). Furthermore, there is morphism of prederivators $u_*\colon\D^A\to\D^B$ defined by right Kan extensions, which is a morphism of derivators if $\D$ is stable and $u$ is homotopically finite (\cite[Thm.~2.14]{gs:stable}).
\end{rmk}

\begin{rmk}\label{rmk:pseudofunctorial}
In \S\S\ref{sec:dnk}-\ref{sec:main} of this work we will usually consider a fixed stable derivator $\D$. Since most of the constructions we will consider are compositions of restrictions of inverse images and Kan extensions between homotopically finite categories (and in the few remaining cases where we also consider more general Kan extensions, we will only consider restrictions of those which are equivalences of derivators) it is immediate that all results in these sections are pseudofunctorial with respect to the 2-category of stable derivators. 
\end{rmk}

\begin{rmk}
Let $\D$ be a derivator and $\D^{c_1}$, $\D^{c_2}$ be full subderivators. If th assignment $A\mapsto\D^{c_1}(A)\cap\D^{c_2}(A)$ defines a full subderivator of $\D$, we will use the notation $\D^{c_1,c_2}$.
\end{rmk}

\section{The 2-category of parasimplices}
\label{sec:para}

\begin{defn}\label{defn:para-s}
\textbf{(The 2-category of parasimplices)}
\begin{enumerate}
\item Let $n \geq 0$. The $n$-\textbf{parasimplex} $\Lambda_n$ is the linearly ordered right $\mathbb{Z}$-set $\mathbb{Z}$ with the $\mathbb{Z}$-operation
\[
(-)+(-)\colon\Lambda_n\times\mathbb{Z}\to\Lambda_n, (\lambda,m)\mapsto \lambda+(n+1)m.
\]
\item The \textbf{2-category of parasimplices} $\underline{\Lambda}$ consists of
\begin{enumerate}
\item objects $\lbrace \Lambda_n \vert n\geq 0 \rbrace$,
\item $\mathbb{Z}$-equivariant maps of linearly ordered sets as 1-morphisms,
\item and natural transformations as 2-morphisms.
\end{enumerate}
\item Let $n,k\geq 0$. We write $\underline{\Lambda}_{n,k}$ for the morphism category $\underline{\Lambda}(\Lambda_k,\Lambda_n)$.
\end{enumerate}
\end{defn}

\begin{rmk}
The 2-category $\underline{\Lambda}$ is a skeletal model of the parasimplex category of \cite[Def.~4.2.1]{lurie:rotation}. Moreover, in \emph{loc.~cit.} the linearly ordered $\mathbb{Z}$-set $\frac{1}{n+1}\mathbb{Z}$ with $\mathbb{Z}$-operation $+1$ is used as a model for $\Lambda_n$.
\end{rmk}

\begin{rmk}
Let $n\geq 0$. The underlying set of $\Lambda_n$ admits a canonical right $\mathbb{Z}$-module structure
\[
(-)\tilde{+}(-)\colon\Lambda_n\times\mathbb{Z}\to\Lambda_n, (\lambda,m)\mapsto\lambda+m.
\]
It will be important to distinguish this module structure from the $\mathbb{Z}$-operation of \autoref{defn:para-s}.
\end{rmk}

\begin{defn}\label{defn:para-t}
\textbf{(Paracyclic translation)} Let $n \geq 0$.
\begin{enumerate}
\item The \textbf{paracyclic operation} $\mathsf{T}\colon\Lambda_n\to\Lambda_n$ is given by
\[
(-)+1=(-)\tilde{+}(n+1)\colon\Lambda_n\to\Lambda_n,\lambda\mapsto\lambda+(n+1).
\]
\item The \textbf{paracyclic translation} $\mathsf{t}\colon\Lambda_n\to\Lambda_n$ is given by
\[
(-)\tilde{+}1\colon\Lambda_n\to\Lambda_n,\lambda\mapsto\lambda+1.
\]
\end{enumerate}
\end{defn}

\begin{rmk}
In the following we use the notation $(-)+(-)$ for the module structure and $\mathsf{T}^{(-)}(-)$ for the $\mathbb{Z}$-operation.
\end{rmk}

\begin{lem}\label{lem:2-transform} 
Let $n\geq 0$.
\begin{enumerate}
\item Paracyclic operations and translations are related by
\[
\mathsf{t}^{n+1}=\mathsf{T}.
\]
\item The 1-morphisms $\mathsf{t}$ and $\mathsf{T}$ are invertible, and conversely every automorphism of $\Lambda_n$ is a power of $\mathsf{t}$. 
\item The paracyclic operations $\mathsf{T}$ assemble into a 2-natural isomorphism
\[
\mathbb{T}:\id_{\underline{\Lambda}}\toiso\id_{\underline{\Lambda}}.
\]
\end{enumerate}
\end{lem}

\begin{proof}
Part (i) is immediate from the definition. The fact that $\mathsf{t}$ and $\mathsf{T}$ come from group operations yields the first statement of (ii). For the second part we consider an automorphism $f\colon\Lambda_n\to\Lambda_n$. Then for $\lambda\in\Lambda_n$ injectivity implies that $f(\lambda+1)\geq f(\lambda)+1$, and by surjectivity $f(\lambda+1)\leq f(\lambda)+1$. By additionally using the dual argument and induction we conclude.
Finally, for (iii) we note that by definition all 1-morphisms in $\underline{\Lambda}$ are $\mathbb{Z}$-equivariant, i.e. they commute with $\mathsf{T}$, for the 2-naturality, and invoke (ii) for the invertibility. 
\end{proof}

\begin{lem}\label{lem:symmetries}
Let $k,n\geq 0$. The automorphisms $\mathsf{t}^*, \mathsf{t}_*: \underline{\Lambda}_{n,k}\to\underline{\Lambda}_{n,k}$ satisfy the relations
\begin{enumerate}
\item $\mathsf{t}^*\circ\mathsf{t}_*=\mathsf{t}_*\circ\mathsf{t}^*$ and
\item $(\mathsf{t}^*)^{k+1}=(\mathsf{t}_*)^{n+1}$.
\end{enumerate}
\end{lem}

\begin{proof}
Let $f\colon\Lambda_k\rightarrow\Lambda_n$ and $\lambda\in\Lambda_k$, then 
\begin{equation}\label{eq:trans-cov-con}
\mathsf{t}^*(f)\colon\lambda\mapsto f(\lambda+1)\qquad \text{ and }\qquad\mathsf{t}_*(f)\colon\lambda\mapsto f(\lambda)+1.
\end{equation}
Hence both sides of (i) describe the assignment
\[
f\mapsto(\lambda\mapsto f(\lambda+1)+1),
\]
which implies (i). For (ii) we invoke \autoref{lem:2-transform} (i) and (iii) to conclude
\[
(\mathsf{t}^*)^{k+1}=\mathsf{T}^*=\mathsf{T}_*=(\mathsf{t}_*)^{n+1}.
\]
\end{proof}

\autoref{lem:symmetries} and the fact, that the notation $(-)^*$ and $(-)_*$ is used later for different purposes, motivate the following definition.

\begin{defn}\label{def:cov-con-sym}
\textbf{(Symmetry operations)} Let $k,n\geq 0$. 
\begin{enumerate}
\item The \textbf{covariant symmetry operation} for $\underline{\Lambda}_{n,k}$ is the automorphism
\[
\mathsf{s}_1:=(\mathsf{t}_*)\colon\underline{\Lambda}_{n,k}\to\underline{\Lambda}_{n,k}.
\]
\item The \textbf{contravariant symmetry operation} for $\underline{\Lambda}_{n,k}$ is the automorphism
\[
\mathsf{s}_2:=(\mathsf{t}^*)\colon\underline{\Lambda}_{n,k}\to\underline{\Lambda}_{n,k}.
\]
\end{enumerate}
\end{defn}

If we replace $\mathsf{T}$ by $\mathsf{t}$ in \autoref{lem:2-transform} (iii), it fails drastically. In the following, we will show that pre- and postcomposition with $\mathsf{t}$, and hence the automorphisms $\mathsf{s}_1$ and $\mathsf{s}_2$, are strongly related to adjunction operations in $\underline{\Lambda}$.

\begin{prop}\label{prop:Lambda-ad-comp}
The 2-category of parasimplices $\underline{\Lambda}$ is adjunction complete (i.e. every 1-morphism $f$ admits both a right adjoint $\mathsf{r}f$ and a left adjoint $\mathsf{l}f$).
\end{prop}

\begin{proof}
Every 1-morphism $f\colon\Lambda_k\to\Lambda_n$ in $\underline{\Lambda}$ gives rise to an endomorphism of the partially ordered set $\mathbb{Z}$ by forgetting the $\mathbb{Z}$-actions. Moreover, $f$ is neither bounded above nor bounded below, since there must be a point in the image, and therefore the entire orbit of this point, which is clearly unbounded, is also contained in the image. For such a map both adjoints clearly exist in the 2-category of posets. More precisely, they are given by the assignments
\begin{equation}\label{eq:left-adjoint}
\mathsf{l}f\colon\mathbb{Z}\to\mathbb{Z},\mu\mapsto\mathrm{min}\lbrace\lambda\in\mathbb{Z}\vert\mu\leq f(\lambda)\rbrace
\end{equation}
and
\begin{equation}\label{eq:right-adjoint}
\mathsf{r}f\colon\mathbb{Z}\to\mathbb{Z},\mu\mapsto\mathrm{max}\lbrace\lambda\in\mathbb{Z}\vert f(\lambda)\leq\mu\rbrace.
\end{equation}
In fact, \eqref{eq:left-adjoint} and \eqref{eq:right-adjoint} are morphisms of parasimplices $\Lambda_n\to\Lambda_k$. We show this only for the left adjoint, since the argument for the right adjoint is very similar. Consider $\mu\in\mathbb{Z}$ and $\lambda=\mathsf{l}f(\mu)$. Since $f\colon\Lambda_k\to\Lambda_n$ is $\mathbb{Z}$-equivariant, the inequalities
\[
\mu+(n+1)\leq f(\lambda+(k+1))
\]
and
\[
f(\lambda+(k+1)-1)=f((\lambda-1)+(k+1))\leq(\mu-1)+(n+1)=\mu+(n+1)-1
\]
hold true, thereby proving that $\mathsf{l}f(\mu+(n+1))=\mathsf{l}f(\mu)+(k+1)$. Finally, we conclude by observing that the inclusion of $\underline{\Lambda}$ into the 2-category of posets is locally fully faithful, such that the units and counits of the adjunctions $\mathsf{l}f\dashv f$ and $f\dashv\mathsf{r}f$ in the 2-category of posets also define 2-morphisms in~$\underline{\Lambda}$.
\end{proof}

\begin{con}
Let $\mathscr{C}$ be an adjunction complete 2-category. Then the essential uniqueness of adjoints and their compatibility with compositions guarantee the existence of pseudofunctors
\[
\mathsf{L},\mathsf{R}\colon\mathscr{C}\to\mathscr{C}^{coop}
\]
which are the identity on objects and such that for every 1-morphism $f\in\mathscr{C}$ there is a triple of adjoint morphisms
\[
\mathsf{L}f\dashv f \dashv\mathsf{R}f.
\]
The category of all such pseudofunctors is a contractible groupoid.
Since the 2-category $\underline{\Lambda}$ is in fact a category enriched over posets, adjoint morphisms are unique. Therefore, in this case the pseudofunctors $\mathsf{L}$ and $\mathsf{R}$ are uniquely determined, and furthermore they are even 2-functors.
\end{con}

\begin{defn}\label{defn:adj-fun}
\textbf{(Adjunction functors)} Let $k,n\geq 0$.
\begin{enumerate}
\item The \textbf{right adjunction functor} $\mathsf{r}$ for $\underline{\Lambda}_{n,k}$ is the structure morphism defined by $\mathsf{R}$
\[
\mathsf{R}_{n,k}:\underline{\Lambda}_{n,k}\to\underline{\Lambda}_{k,n}^{op}.
\]
\item The \textbf{left adjunction functor} $\mathsf{l}$ for $\underline{\Lambda}_{n,k}$ is the structure morphism defined by $\mathsf{L}$
\[
\mathsf{L}_{n,k}:\underline{\Lambda}_{n,k}\to\underline{\Lambda}_{k,n}^{op}.
\]
\end{enumerate}
\end{defn}

The following proposition addresses the compatibility of symmetry operations and adjunction functors.

\begin{prop}\label{prop:adj-comp}
Let $k,n \geq 0$. There are equalities of functors $\underline{\Lambda}_{n,k}\to\underline{\Lambda}_{k,n}^{op}$
\begin{equation}
\mathsf{l}\circ\mathsf{s}_2=(\mathsf{s}_1^{op})^{-1}\circ\mathsf{l}=(\mathsf{s}_2^{op})^{-1}\circ\mathsf{r}=\mathsf{r}\circ\mathsf{s}_1.
\end{equation}
\end{prop}

\begin{proof}
Let $f\colon\Lambda_k\to\Lambda_n$ be a morphism. The compatibility of taking left (resp. right) adjoints with compositions applied to $f\circ\mathsf{t}$ (resp. $\mathsf{t}\circ f$) immediately yields the first (resp. third) equality. For the second equality we consider $\mu\in\Lambda_n$. Then \eqref{eq:trans-cov-con} and \eqref{eq:left-adjoint} show that $(\mathsf{s}_1^{op})^{-1}\circ\mathsf{l}$ is given by the assignment
\begin{equation}\label{eq:standart-ad}
\mu\mapsto\mathrm{min}\lbrace\lambda\in\Lambda_k\vert\mu\leq f(\lambda)\rbrace-1.
\end{equation}
The values of the assignment can be re-expressed as the maxima of all elements, which are smaller than the respective minima
\begin{equation}
\mathrm{min}\lbrace\lambda\in\Lambda_k\vert\mu\leq f(\lambda)\rbrace-1=\mathrm{max}\lbrace\lambda\in\Lambda_k\vert f(\lambda)\leq\mu-1\rbrace.
\end{equation}
By using \eqref{eq:trans-cov-con} and \eqref{eq:right-adjoint} we conclude that \eqref{eq:standart-ad} also describes $(\mathsf{s}_2^{op})^{-1}\circ\mathsf{r}$.
\end{proof}

We are now ready to discuss the 2-naturality of the paracyclic translations $\mathsf{t}$.

\begin{cor}\label{cor:gsd-Lambda}
The paracyclic translations $\mathsf{t}$ assemble into a 2-natural isomorphism
\[
\mathbb{S}\colon\mathsf{R}\to\mathsf{L}\colon\underline{\Lambda}\to\underline{\Lambda}^{coop}.
\]
\end{cor}

\begin{proof}
Let $f\colon\Lambda_k\to\Lambda_n$ be a 1-morphism in $\underline{\Lambda}$. The statement follows from the commutativity of the squares
\[
\xymatrix{
\Lambda_n \ar[r]^{\mathsf{r}f} \ar[d]^{\mathsf{t}} & \Lambda_k \ar[d]^{\mathsf{t}}\\
\Lambda_n \ar[r]^{\mathsf{l}f} & \Lambda_k,
}
\]
which in turn is a reformulation of the second equality in \autoref{prop:adj-comp}.
\end{proof}

\begin{rmk}
The 2-natural isomorphism $\mathbb{S}$ is an example of a global Serre duality. These structures will be investigated in more detail in \cite{bg:global}.
\end{rmk}

\begin{defn}
A family of basepoints $\lambda_\bullet=\lbrace\lambda_n\in\Lambda_n, n \geq 0\rbrace$ is called a \textbf{choice of coordinates} in $\underline{\Lambda}$.
\end{defn}

\begin{prop}\label{prop:coord-mor}
Let $n,k\geq 0$, $\lambda_\bullet$ a choice of coordinates, and $\mu_0,\cdots,\mu_k\in\Lambda_n$ such that
\begin{equation}\label{eq:coord1}
\mu_0\leq\mu_1\leq\cdots\leq\mu_k\leq\mu_0+(n+1).
\end{equation}
Then there is a unique morphism $f\colon\Lambda_k\to\Lambda_n$ such that for $i\in\lbrace 0,\cdots,k\rbrace$
\begin{equation}\label{eq:coord2}
f(\lambda_k+i)=\mu_i.
\end{equation}
\end{prop}

\begin{proof}
Let $\lambda\in\Lambda_k$. Then there are uniquely determined $i\in\lbrace 0,\cdots,k\rbrace,p\in\mathbb{Z}$ such that $\lambda=\lambda_k+i+p(k+1)$. Hence the assignment
\[
\lambda\mapsto\mu_i+p(n+1)
\]
is the only possible definition of $f$ that satisfies equivariance and \eqref{eq:coord2}. This assignment is order preserving by \eqref{eq:coord1}, thus it defines a map of parasimplices.
\end{proof}

\begin{defn}
Let $k,n\geq 0$, $\lambda_\bullet$ a choice of coordinates, $f\in\underline{\Lambda}_{n,k}$. The $\lambda_\bullet$\textbf{-coordinate representation} $f_{\lambda_\bullet}$ of $f$ is the $k+1$-tuple
\[
(f(\lambda_k),f(\lambda_k+1),\cdots,f(\lambda_k+k))\in\mathbb{Z}^{k+1}.
\]
\end{defn}

\begin{cor}\label{cor:coord-rep}
Let $n,k\geq 0$ and $\lambda_\bullet$ a choice of coordinates. Then the $\lambda_\bullet$-coordinate representation induces an embedding of posets
\[
\underline{\Lambda}_{n,k}\to\mathbb{Z}^{k+1},f\mapsto f_{\lambda_\bullet}.
\]
\end{cor}

\begin{proof}
We see immediately, that the $\lambda_\bullet$-coordinate representation induces a morphism of posets
\[
\underline{\Lambda}_{n,k}\to\mathbb{Z}^{k+1}
\]
and that the image is characterized by \eqref{eq:coord1}. The injectivity follows from \autoref{prop:coord-mor}.
\end{proof}

\begin{cor}
Let $\lambda_\bullet$ be a choice of coordinates. Then there is a locally fully faithful 2-functor $\mathsf{i}=\mathsf{i}_{\lambda_{\bullet}}\colon\Delta\to\underline{\Lambda}$
defined by
\begin{enumerate}
\item $\Delta_n\mapsto\Lambda_n$ for $n\geq 0$ and
\item $(g\colon\Delta_k\to\Delta_n)\mapsto(\lambda_n+g(0),\cdots,\lambda_n+g(k))$ for $k,n\geq 0$.
\end{enumerate}
\end{cor}

\begin{proof}
By \autoref{prop:coord-mor} the map $\Delta(\Delta_k,\Delta_n)\to\underline{\Lambda}_{n,k}$ defined by (ii) is order preserving and injective. The compatibility with identities is obvious. Finally, consider $g\colon\Delta_k\to\Delta_n$ and $h\colon\Delta_n\to\Delta_l$. The composition of the images of $g$ and $h$ is determined by its values on $\lambda_k+i$ for $i\in\lbrace 0,\cdots,k\rbrace$ (\autoref{prop:coord-mor})
\[
\lambda_k+i\mapsto\lambda_n+g(i)\mapsto\lambda_l+h(g(i))=\lambda_l+(h\circ g)(i),
\]
which is by definition the image of the composition of $g$ and $h$.
\end{proof}

\begin{cor}
Let $n,k \geq 0$ and $\lambda_{\bullet}$ be choice of coordinates. Then the functor
\begin{equation}
\mathsf{d}_{n,k}\colon\underline{\Lambda}_{n,k}\to\underline{\Lambda}_{n,k}^{op},(f_i)_{0 \leq i \leq k}\mapsto (2\lambda_n+n-f_{k-i})_{0 \leq i \leq k}
\end{equation}
is a self-inverse isomorphism of posets. Moreover, the functors $\mathsf{d}_{n,k}\colon\underline{\Lambda}_{n,k}\to\underline{\Lambda}_{n,k}^{op}$ assemble into a strictly self-inverse 2-functor
\[
\mathsf{D}=\mathsf{D}_{\lambda_{\bullet}}\colon\underline{\Lambda}\to\underline{\Lambda}^{co}.
\]
\end{cor}

\begin{proof}
The assignment $\mathsf{d}_{n,k}$ is clearly order reversing. Moreover, the property $f_0\leq \cdots \leq f_k\leq f_0+n+1$ yields
\[
2\lambda_n+n-f_k\leq\cdots\leq 2\lambda_n+n-f_0\leq 2\lambda_n+n-f_k+n+1,
\]
where the last inequality is obtained from the chain of implications
\begin{align}
 f_k\leq f_0+n+1 \Rightarrow  &-(f_0+n+1)\leq -f_k \\
\Rightarrow &-f_0\leq -f_k+n+1\\
\Rightarrow & 2\lambda_n+n-f_0\leq 2\lambda_n+n-f_k+n+1.
\end{align}
Hence $\mathsf{d}_{n,k}\colon\underline{\Lambda}_{n,k}\to\underline{\Lambda}_{n,k}^{op}$ defines indeed a functor. Furthermore,  $\mathsf{d}_{n,k}$ is self inverse since
\[
\mathsf{d}_{n,k}^2((f_i)_{0 \leq i \leq k})=(2\lambda_n+n-(2\lambda_n+n-f_{k-(k-i)}))_{0 \leq i \leq k}=(f_i)_{0 \leq i \leq k},
\]
and therefore an equivalence.
For the functoriality statement we have to understand how $\mathsf{d}_{n,k}(f)$ operates on $\lambda_k+j$ for $j\in\mathbb{Z}$ arbitrary. More precisely, we claim that for a map of parasimplices $f\colon\Lambda_k\to\Lambda_n,\lambda_k+j\mapsto f_j$ the image under $\mathsf{d}_{n,k}$ is given by the assignment $f'\colon\lambda_k+j\mapsto 2\lambda_n+n-f_{k-j}$. For this, we only have to show that this assignment is equivariant. Then $f'=\mathsf{d}_{n,k}(f)$ will follow from \autoref{prop:coord-mor} by using that both sides agree on $\lbrace \lambda_k,\lambda_k+1,\cdots,\lambda_k+k\rbrace$. The equivariance follows from the corresponding property of $f$
\[
f'(\lambda_k+j+k+1)=2\lambda_n+n-f_{k-(j+k+1)}=2\lambda_n+n-(f_{k-j}-(n+1))=f'(\lambda_k+j)+n+1.
\]
For the unitality, we compute
\[
\mathsf{d}_{k,k}(\id)(\lambda_k+j)=2\lambda_k+k-\id(\lambda_k+k-j)=\lambda_k+j.
\]
For the compatibility with compositions, we consider $f\colon\Lambda_k\to\Lambda_n,\lambda_k+j\mapsto f_j$ and $g\colon\Lambda_n\to\Lambda_m,\lambda_n+j\mapsto g_j$. Then the composition of $f$ and $g$ is given by
\begin{equation}\label{eq:para-composition}
g\circ f\colon\Lambda_k\to\Lambda_n,\lambda_k+j\xmapsto{f}f_j=\lambda_n+(f_j-\lambda_n)\xmapsto{g}g_{(f_j-\lambda_n)}.
\end{equation}
On the other hand we have by definition
\[
\mathsf{D}(f)\colon\lambda_k+j\mapsto 2\lambda_n+n-f_{k-j}\qquad\text{and}\qquad\mathsf{D}(g)\colon\lambda_n+j\mapsto 2\lambda_m+m-g_{n-j}.
\]
In particular, we can compute the composition
\begin{align}
\mathsf{D}(g)\circ\mathsf{D}(f)\colon\lambda_k+j\xmapsto{\mathsf{D}(f)}&2\lambda_n+n-f_{k-j}=\lambda_n+(\lambda_n+n-f_{k-j})\\
\xmapsto{\mathsf{D}(g)}&2\lambda_m+m+g_{(n-(\lambda_n+n-f_{k-j}))}=2\lambda_m+m+g_{(f_{k-j}-\lambda_n)},
\end{align}
and invoke \eqref{eq:para-composition} to identify the last term with $\mathsf{D}(g\circ f)(\lambda_k+j)$.
\end{proof}

In the following, we show that for different choices of coordinates the resulting induced simplicial embeddings and duality 2-functors are related via 2-natural isomorphisms. Hence we can fix one choice of coordinates, and work without loss of generality with the structures $\mathsf{i}$ and $\mathsf{D}$ induced by this fixed choice. In particular, we will use the choice of coordinates $0_{\bullet}=\lbrace 0\in\Lambda_n,n\geq 0\rbrace$. Since the element $0$ in the $\mathbb{Z}$-set $\mathbb{Z}$ is not characterized by a "special property", the choice of the coordinates $0_{\bullet}$ is not "more natural" than any other possible choice. It just leads in many occasions to a simpler notation.

\begin{prop}\label{prop:change-coord}
Let $\lambda_{\bullet}$ be a choice of coordinates.
\begin{enumerate}
\item There is a 2-functor $F^{\lambda_{\bullet}}\colon\underline{\Lambda}\to\underline{\Lambda}$ which is the identity on objects and such that for $n,k\geq 0$
\[
F^{\lambda_{\bullet}}_{n,k}=\mathsf{s}_1^{\lambda_n}\circ\mathsf{s}_2^{-\lambda_k}\colon\underline{\Lambda}_{n,k}\to\underline{\Lambda}_{n,k}.
\]
\item There is a 2-natural isomorphism $\phi^{\lambda_{\bullet}}\colon\mathrm{id}\to F^{\lambda_{\bullet}}$ such that for $n\geq 0$
\[
\phi^{\lambda_{\bullet}}_n=\mathsf{t}^{\lambda_n}\colon\Lambda_n\to\Lambda_n.
\]
\end{enumerate}
\end{prop}

\begin{proof}
We invoke \autoref{prop:conjugation} and define $F^{\lambda_{\bullet}}=\id[\mathsf{S}]$ and $\phi^{\lambda_{\bullet}}=\alpha[\mathsf{S}]$ for $\mathsf{S}_{\Lambda_n}=\mathsf{t}^{\lambda_n}$. Moreover, $F^{\lambda_{\bullet}}$ is a 2-functor and $\phi^{\lambda_{\bullet}}$ is a 2-natural isomorphism since $\id\colon\underline{\Lambda}\to\underline{\Lambda}$ is a 2-functor and $\mathsf{t}\colon\Lambda_n\to\Lambda_n$ is an isomorphism for all $n\geq 0$.
\end{proof}

\begin{cor}
Let $\lambda_{\bullet}$ be a choice of coordinates.
\begin{enumerate}
\item There is a 2-natural isomorphism $\mathsf{i}_{0_{\bullet}}\to\mathsf{i}_{\lambda_{\bullet}}$.
\item There is a 2-natural isomorphism $\mathsf{D}_{0_{\bullet}}\to\mathsf{D}_{\lambda_{\bullet}}$.
\end{enumerate}
\end{cor}

\begin{proof}
For (i) we observe by plugging in the definitions, that $\mathsf{i}_{\lambda_{\bullet}}=F^{\lambda_{\bullet}}\circ\mathsf{i}_{0_{\bullet}}$. Similarly, for (ii) we have $\mathsf{D}_{\lambda_{\bullet}}=F^{\lambda_{\bullet}}\circ\mathsf{D}_{0_{\bullet}}\circ(F^{\lambda_{\bullet}})^{-1}$. Hence, by \autoref{prop:change-coord}, the desired 2-natural isomorphisms are induced by $\phi^{\lambda_{\bullet}}$. 
\end{proof}

In the following we discuss some applications of the structures $\mathsf{i}$ and $\mathsf{D}$ with respect to the choice of coordinates $0_{\bullet}$. We start by describing the symmetry operations with respect to this choice of coordinates. For this we define for $n,k\geq 0$ the following automorphisms of the poset $\mathbb{Z}^{k+1}$.
\begin{enumerate}
\item $\mathsf{s}_1^{\mathbb{Z}}\colon(\lambda_0,\cdots,\lambda_k)\mapsto(\lambda_0+1,\cdots,\lambda_{k}+1)$,
\item $\mathsf{s}_2^{\mathbb{Z}}\colon(\lambda_0,\cdots,\lambda_k)\mapsto(\lambda_1,\cdots,\lambda_{k},\lambda_0+n+1)$.
\end{enumerate}

\begin{prop}\label{prop:coord-sym}
Let $n,k\geq 0$ and $f\in\underline{\Lambda}_{n,k}$. Then
\begin{enumerate}
\item $\mathsf{s}_1^{\mathbb{Z}}(f_{0_{\bullet}})=(\mathsf{s}_1(f))_{0_{\bullet}}$,
\item $\mathsf{s}_2^{\mathbb{Z}}(f_{0_{\bullet}})=(\mathsf{s}_2(f))_{0_{\bullet}}$.
\end{enumerate}
\end{prop}

\begin{proof}
Let $f\in\underline{\Lambda}_{n,k}$, then $f_{0_{\bullet}}=(f(0),\cdots,f(k))$. Then we conclude by observing
\begin{enumerate}
\item $(\mathsf{s}_1(f))_{0_{\bullet}}=(\mathsf{t}\circ f)_{0_{\bullet}}=(f(0)+1,\cdots,f(k)+1)$,
\item $(\mathsf{s}_2(f))_{0_{\bullet}}=(f\circ\mathsf{t})_{0_{\bullet}}=(f(1),\cdots,f(k),f(0)+n+1)$,
\end{enumerate}
where we have used equivariance to identify the last coordinate in (ii).
\end{proof}

From now on we will omit the $\mathbb{Z}$-superscript on $\mathsf{s}_1^{\mathbb{Z}}$ and $\mathsf{s}_2^{\mathbb{Z}}$.

\begin{prop}\label{prop:lambda-shift}
Let $n,k\geq 0$ and $f\in\underline{\Lambda}_{n,k}$ Then there exists a unique $l\in\mathbb{Z}$ and $g\colon\Delta_k\to\Delta_n$ such that 
\[
f=\mathsf{s}_2^l\circ\mathsf{i}(g).
\]
\end{prop}

\begin{proof}
We first show the existence of $l$. Let $f=(f_0,\cdots,f_k)\in\underline{\Lambda}_{n,k}$ not in the image of $\mathsf{i}$. We assume $0\leq f_0$ and observe directly that then $f_k\geq n+1$. The case $f_0\leq-1$ is completely dual to this one. Then there exist $r\in\lbrace 0,\cdots,n\rbrace$, $q\geq 0$ such that $f_0=q(n+1)+r$. Then we consider $f'=(f'_0,\cdots,f'_k):=\mathsf{s}_2^{-q(k+1)}(f)=\mathsf{s}_1^{-q(n+1)}(f)$. From the latter description of $f'$ (which follows from \autoref{lem:symmetries}) we deduce that $0\leq f'_0\leq n$ and hence $f'_k\leq f'_0+n+1\leq 2n+1$. If $f'_k\leq n$, then $f'$ is in the image of $\mathsf{i}$. Otherwise consider $j$ minimal with $n+1\leq f'_j$. Then we have by construction that $\mathsf{s}_2^{-k-j+1}(f')$ is in the image of $\mathsf{i}$. 

\noindent For the converse direction, we note that the image of $\mathsf{i}$ has a minimal element $\underline{0}:=(0,\cdots,0)$ and a maximal element $\underline{n}:=(n,\cdots,n)$. Comparing the last coordinates yields
\[
\mathsf{s}_2(\underline{0}) \nleq \underline{n}.
\]
Since $\mathsf{s}_2$ is order preserving we conclude by using the minimality of $\underline{0}$ and the maximality of $\underline{n}$ that $\mathsf{s}_2(\mathsf{i}(g)) \nleq \mathsf{i}(g')$ for all $g,g'\colon\Delta_k\to\Delta_n$. In particular, the images of $\mathsf{i}$ and $\mathsf{s}_2\circ\mathsf{i}$ are disjoint. This yields inductively the injectivity of 
\[
\Delta(\Delta_k,\Delta_n)^{\mathbb{Z}},(g_i)_{i\in\mathbb{Z}}\mapsto(\mathsf{s}_2^i\circ\mathsf{i}(g))_{i\in\mathbb{Z}}
\]
and hence the uniqueness of $l$.
\end{proof}

\begin{defn}
Let $n\geq 0$.
\begin{enumerate}
\item Let $0\geq i\geq n+1$. We call $\mathsf{i}(\mathrm{d}_i)\colon\Lambda_{n}\to\Lambda_{n+1}$ the $i$\textbf{th parasimplicial face map} and denote it also by $\mathrm{d}_i$.
\item Let $0\geq i\geq n$. We call $\mathsf{i}(\mathrm{s}_i)\colon\Lambda_{n+1}\to\Lambda_{n}$ the $i$\textbf{th parasimplicial degeneracy map} and denote it also by $\mathrm{s}_i$.
\end{enumerate}
\end{defn}

\begin{rmk}\label{rmk:para-fd}
In the 2-category of simplices there is the following chain of adjunctions
\[
\mathrm{d}_{n+1}\dashv\mathrm{s}_{n}\dashv\mathrm{d}_n\dashv\cdots\dashv\mathrm{d}_1\dashv\mathrm{s}_{0}\dashv\mathrm{d}_0
\]
of 1-morphisms relating $\Delta_n$ and $\Delta_{n+1}$. The 2-functoriality of $\mathsf{i}$ yields a chain of adjunctions relating the  parasimplicial face and degeneracy maps. Moreover, by \autoref{cor:gsd-Lambda} this can be reformulated as
\begin{equation}\label{eq:para-fd}
\mathrm{d}_{i}=\mathsf{t}^i\circ\mathrm{d}_0\circ\mathsf{t}^{-i}\qquad\text{and}\qquad\mathrm{s}_{i}=\mathsf{t}^i\circ\mathrm{s}_0\circ\mathsf{t}^{-i}.
\end{equation}
Furthermore, \eqref{eq:para-fd} can be taken as a definition in cases where the left hand sides are not defined.
\end{rmk}

\begin{cor}\label{cor:Lambda-generators}
Let $n,k\geq 0$ and $f\in\underline{\Lambda}_{n,k}$. Then $f$ is a composition of morphisms of the form $\mathsf{t}$, $\mathsf{t}^{-1}$, $\mathrm{d}_0$ and $\mathrm{s}_0$.
\end{cor}

\begin{proof}
By \autoref{prop:lambda-shift} there is $l\in\mathbb{Z}$ and $g\colon\Delta_k\to\Delta_n$ such that there is a decomposition $f=\mathsf{i}(g)\circ\mathsf{t}^l$. Since $\Delta$ is generated by the simplicial face and degeneracy maps, we can decompose in terms of those. The 2-functoriality of $\mathsf{i}$ yields a decomposition of $\mathsf{i}(g)$ into parasimplicial face and degeneracy maps. Finally, we obtain the desired decomposition by applying \eqref{eq:para-fd} to each of these.
\end{proof}

\begin{rmk}\label{adjunction-duality}
Since the duality $\mathsf{D}\colon\underline{\Lambda}\to\underline{\Lambda}^{co}$ is a 2-functor, it preserves adjunctions, but reverses the order. In other words, the diagrams
\[
\xymatrix{
\underline{\Lambda}\ar[r]^{\mathsf{L}}\ar[d]_{\mathsf{D}}&\underline{\Lambda}^{coop}\ar[d]^{\mathsf{D}}\\
\underline{\Lambda}^{co}\ar[r]_{\mathsf{R}}&\underline{\Lambda}^{op}
}
\qquad\text{and}\qquad
\xymatrix{
\underline{\Lambda}\ar[r]^{\mathsf{R}}\ar[d]_{\mathsf{D}}&\underline{\Lambda}^{coop}\ar[d]^{\mathsf{D}}\\
\underline{\Lambda}^{co}\ar[r]_{\mathsf{L}}&\underline{\Lambda}^{op}
}
\]
commute. We emphasize that this construction generalizes to arbitrary adjunction complete 2-categories, and indicates that for such a 2-category $\mathcal{C}$, there is a duality between equivalences $\mathcal{C}\to\mathcal{C}^{op}$ and equivalences $\mathcal{C}\to\mathcal{C}^{co}$. We consider, the right square and $n,k\geq 0$. In particular, we obtain isomorphisms
\[
\mathrm{ad}_{n,k}:=\mathsf{d}_{k,n}\circ\mathsf{r}_{n,k}=\mathsf{l}_{n,k}\circ\mathsf{d}_{n,k}\colon\underline{\Lambda}_{n,k}\toiso\underline{\Lambda}_{k,n}.
\]
The following result describes the compatibility of $\mathrm{ad}_{n,k}$ with choices of coordinates.
\end{rmk}

\begin{prop}
Let $n,k\geq 0$. Then $\mathrm{ad}_{n,k}(0,\cdots,0)=(0,\cdots,0)$.
\end{prop}

\begin{proof}
By definition $\mathsf{d}_{k,n}(0,\cdots,0)=(n,\cdots,n)$ and \eqref{eq:left-adjoint} yields $\mathsf{l}_{k,n}(n,\cdots,n)=(0,\cdots,0)$.
\end{proof}

Another interesting feature of the 2-category of parasimplices is the behavior of injective 1-morphisms. More precisely, we will see that for given $n,k \geq 0$ the poset $\underline{\Lambda}_{n,k}$ exactly describes the injective morphisms for a certain shift in the codomain.

\begin{defn}
Let $n,k \geq 0$. The poset $\underline{\Lambda}_{n,k}^{inj}$ is the subposet of $\underline{\Lambda}_{n,k}$ consisting of those 1-morphisms $f\colon\Lambda_k\to\Lambda_n$ whose underlying maps $f\colon\mathbb{Z}\to\mathbb{Z}$ are injective.
\end{defn}

\begin{prop}\label{prop:para-inj}
Let $n,k \geq 0$. The functor $\mathsf{inj}=\mathsf{inj}_{n,k}\colon\underline{\Lambda}_{n,k}\to\underline{\Lambda}_{n+k+1,k}^{inj}$ which is defined on coordinate representations by the assignment
\[
(f_0,f_1,\cdots,f_k)\mapsto(f_0,f_1+1,\cdots,f_k+k)
\]
is an isomorphism.
\end{prop}

\begin{proof}
Note that a morphism $g\colon\Lambda_k\to\Lambda_{n+k+1}$ with coordinate representation $(g_0,\cdots,g_k)$ is injective if and only if the relations
\begin{equation}\label{eq:coord-inj}
g_0<\cdots <g_k < g_0+n+k+2
\end{equation}
hold true. Since the relations
\[
f_0\leq\cdots\leq f_k\leq f_0+n+1
\]
imply the inequalities (for the last inequality we use $f_k\leq f_0+n+1 \Rightarrow f_k+k\leq f_0+n+1+k < f_0+n+k+2$)
\[
f_0 < f_1+1 < \cdots < f_k+k < f_0+n+k+2,
\]
we conclude that $\mathsf{inj}$ is well defined. Moreover, $\mathsf{inj}$ is a restriction of the translation map $+(0,\cdots,k)\colon\mathbb{Z}\to\mathbb{Z}$. Hence $\mathsf{inj}$ is order preserving and injective. It remains to prove the surjectivity of $\mathsf{inj}$. For this we consider $g\in\underline{\Lambda}_{n+k+1,k}^{inj}$ with coordinate representation $(g_0,\cdots,g_k)$. Then \eqref{eq:coord-inj} yields the inequalities
\[
g_0\leq g_1-1\leq\cdots\leq g_k-k\leq g_0+n+1,
\]
which exhibit $(g_0,g_1-1,\cdots,g_k-k)$ as a preimage of $g$.
\end{proof}

\begin{prop}\label{prop:sym-inj}
Let $n,k\geq 0$. The symmetries $\mathsf{s}_1, \mathsf{s}_2\colon\underline{\Lambda}_{n+k+1,k}\to\underline{\Lambda}_{n+k+1,k}$ restrict to automorphisms of $\underline{\Lambda}_{n+k+1,k}^{inj}$ and furthermore there are equalities of isomorphisms $\underline{\Lambda}_{n,k}\to\underline{\Lambda}_{n+k+1,k}^{inj}$
\begin{enumerate}
\item $\mathsf{inj}\circ\mathsf{s}_1=\mathsf{s}_1\circ\mathsf{inj},$
\item $\mathsf{inj}\circ\mathsf{s}_2=\mathsf{s}_2\circ\mathsf{s}_1^{-1}\circ\mathsf{inj}.$
\end{enumerate}
\end{prop}

\begin{proof}
We invoke \autoref{prop:coord-sym} to describe all compositions explicitly in terms of coordinate representations. In the case of (i), we observe immediately that both sides are given by the assignment
\[
(f_0,f_1,\cdots,f_k)\mapsto(f_0+1,f_1+2,\cdots,f_k+k+1).
\]
For (ii) we compute both sides separately and conclude by comparing the results.
\begin{align}
(f_0,f_1,\cdots,f_k)\xmapsto{\mathsf{s}_2}&(f_1,f_2,\cdots,f_k,f_0+n+1)\\
\xmapsto{\mathsf{inj}}&(f_1,f_2+1,\cdots,f_k+k-1,f_0+n+k+1)\\
&\text{and}\\
(f_0,f_1,\cdots,f_k)\xmapsto{\mathsf{inj}}&(f_0,f_1+1,\cdots,f_k+k)\\
\xmapsto{\mathsf{s}_1^{-1}}&(f_0-1,f_1,\cdots,f_k+k-1)\\
\xmapsto{\mathsf{s}_2}&(f_1,f_2+1,\cdots,f_k+k-1,f_0-1+n+k+2).
\end{align}
\end{proof}

\section{The stable homotopy theory of cubes}
\label{sec:cubical}

In this short section we recollect some important facts about the stable homotopy theory of cubes. Building on these, we will construct stable derivators associated to the categories $\underline{\Lambda}_{n,k}$ in \S\ref{sec:dnk}. More precisely, the proof of \autoref{thm:slices} will rely on \autoref{prop:stability} and \autoref{prop:frankes-lemma} and the proof of \autoref{prop:slnk-sigma} will rely on \autoref{cor:bicart-concat}.

\begin{defn}
Let $n\geq 0$. The $n$-\textbf{cube} $\cube{n}\in Cat$ is the category $[1]^n$.
\end{defn}

\begin{rmk}
Let $n\geq 0$.
\begin{enumerate}
\item The category $\cube{n}$ is a poset.
\item Let $\mathbf{n}=\lbrace 0, \cdots, n-1 \rbrace$. Then the power set $\mathcal{P}(\mathbf{n})$ inherits the structure of a poset via the inclusion of subsets. Then the assignment 
\[
\chi_{(-)}\colon\mathcal{P}(\mathbf{n})\to [1]^n,
\]
defined by the characteristic functions of subsets, is an isomorphism of posets. In the following, we will sometimes use elements of $\mathcal{P}(\mathbf{n})$ to specify objects of $\cube{n}$, without mentioning the implicit use of the isomorphism $\chi_{(-)}$.
\item In particular, there is a functor $\#\colon\cube{n}\to[n]$ defined by the cardinality of subsets of $\mathbf{n}$.
\item In \cite{bg:cubical} the set $\lbrace 1, \cdots, n \rbrace$ was used to parametrize the coordinates of $\cube{n}$. The reason for the different convention here is the fact that objects in $\underline{\Lambda}_{n,k}$ can be described as functions on $\mathbf{k+1}$ via choices of coordinates (\autoref{cor:coord-rep}).
\end{enumerate}
\end{rmk}

\begin{defn}
Let $0\leq k\leq l\leq n$ and $i\in\mathbf{n}$.
\begin{enumerate}
\item We denote by
\[
\mathrm{d}_0^i,\mathrm{d}_1^i\colon\cube{n-1}\to\cube{n}\qquad\text{and}\qquad\mathrm{s}_0^i\colon\cube{n}\to\cube{n-1}
\]
the face maps and the projection with respect to the $i$th coordinate.
\item We denote by 
\[
\iota_{k,l}\colon\cube{n}_{k,l}\to\cube{n}
\]
the inclusion of $\#^{-1}(\lbrace k,k+1,\cdots,l\rbrace)$.
\end{enumerate}
\end{defn}

\begin{rmk}
Let $M\subset\mathbf{n}$ be a subset. Then we use the more general notation
\begin{enumerate}
\item $\mathrm{d}^M_0=\prod_{m\in M}\mathrm{d}^m_0\colon\cube{n-\#(M)}\to\cube{n}$,
\item $\mathrm{d}^M_1=\prod_{m\in M}\mathrm{d}^m_1\colon\cube{n-\#(M)}\to\cube{n}$,
\item $\mathrm{s}^M_0=\prod_{m\in M}\mathrm{s}^m_0\colon\cube{n}\to\cube{n-\#(M)}$.
\end{enumerate}
\end{rmk}

\begin{con}
We recall the construction of cofiber and fiber sequences for a stable derivator $\D$ and some of their properties. We refer to \cite{groth:ptstab} and in the cubical case to \cite{bg:cubical} for a more detailed treatment of these constructions. Let $j_1\colon[1]\to\cube{2}_{0,1}$ be the inclusion induced by $\mathrm{d}_1^0$ and $j_2\colon[2]\to\cube{2},i\mapsto\mathbf{i}$. Then we define for a stable derivator $\D$
\begin{itemize}
\item the \textbf{cofiber sequence morphism} $\mathsf{Cof}=j_2^*\circ(\iota_{0,1})_!\circ(j_1)_*\colon\D^{[1]}\to\D^{[2]}$,
\item the \textbf{cofiber morphism} $\cof=\mathrm{d}_0^*\circ\mathsf{Cof}\colon\D^{[1]}\to\D^{[1]}$,
\item the \textbf{cone morphism} $C=\infty^*\circ\cof\colon\D^{[1]}\to\D$,
\item the \textbf{suspension morphism} $\Sigma=C\circ(\mathrm{d}_1)_*\colon\D\to\D$.
\end{itemize}
Let now $n\geq 1$. Then we define
\begin{itemize}
\item the \textbf{iterated cofiber sequence morphism} $\mathsf{Cof}^{\underline{1}}\colon\D^{\cube{n}}\to\D^{[2]^n}$ by applying the cofiber sequence morphism in every variable separately.
\item the \textbf{iterated cofiber morphism} $\cof^{\underline{1}}=(\mathrm{d}_0\times\cdots\times\mathrm{d}_0)^*\circ\mathsf{Cof}^{\underline{1}}\colon\D^{\cube{n}}\to\D^{\cube{n}}$.
\item the \textbf{total cofiber morphism} $\tcof=\infty^*\circ\cof^{\underline{1}}\colon\D^{\cube{n}}\to\D$.
\end{itemize}
Dually, let $j'_1\colon[1]\to\cube{2}_{1,2}$ be the inclusion induced by $\mathrm{d}_0^1$. Then we define 
\begin{itemize}
\item the \textbf{fiber sequence morphism} $\mathsf{Fib}=j_2^*\circ(\iota_{1,2})_*\circ(j'_1)_!\colon\D^{[1]}\to\D^{[2]}$,
\item the \textbf{fiber morphism} $\fib=\mathrm{d}_2^*\circ\mathsf{Fib}\colon\D^{[1]}\to\D^{[1]}$,
\item the \textbf{cocone morphism} $F=\emptyset^*\circ\cof\colon\D^{[1]}\to\D$.
\item the \textbf{loop morphism} $\Omega=F\circ(\mathrm{d}_0)_!\colon\D\to\D$.
\end{itemize}
Let now $n\geq 1$. Then we define
\begin{itemize}
\item the \textbf{iterated fiber sequence morphism} $\mathsf{Fib}^{\underline{1}}\colon\D^{\cube{n}}\to\D^{[2]^n}$ by applying the fiber sequence morphism in every variable separately.
\item the \textbf{iterated fiber morphism} $\fib^{\underline{1}}=(\mathrm{d}_2\times\cdots\times\mathrm{d}_2)^*\circ\mathsf{Fib}^{\underline{1}}\colon\D^{\cube{n}}\to\D^{\cube{n}}$.
\item the \textbf{total fiber morphism} $\tcof=\emptyset^*\circ\cof^{\underline{1}}\colon\D^{\cube{n}}\to\D$.
\end{itemize}
Then the following properties are satisfied
\begin{itemize}
\item $(\Sigma,\Omega)$ is a pair of mutually inverse equivalences.
\item $(\cof^{\underline{1}},\fib^{\underline{1}})$  is a pair of mutually inverse equivalences.
\end{itemize}
\end{con}

\begin{rmk}\label{rmk:elementary}
We summarize some of the well-known isomorphisms between the fiber and cofiber constructions for a stable derivator $\D$, which are elementary consequences of the definitions.
\begin{enumerate}
\item $\cof^3\cong\Sigma\colon\D^{[1]}\to\D^{[1]}$,
\item $\fib^3\cong\Omega\colon\D^{[1]}\to\D^{[1]}$,
\item $\Omega\circ 2^*\circ\mathsf{Cof}\cong \Omega\circ\mathrm{d}_0^*\circ\cof\cong F \cong\mathrm{d}_1^*\circ\fib \cong 0^*\circ\mathsf{Fib}\colon\D^{[1]}\to\D$,
\item $0^*\circ\mathsf{Cof}\cong F\circ\cof\cong \mathrm{d}_1^* \cong\mathrm{d}_0^*\circ\fib \cong 1^*\circ\mathsf{Fib}\colon\D^{[1]}\to\D$,
\item $1^*\circ\mathsf{Cof}\cong \mathrm{d}_1^*\circ\cof\cong \mathrm{d}_0^* \cong C\circ\fib \cong 2^*\circ\mathsf{Fib}\colon\D^{[1]}\to\D$,
\item $2^*\circ\mathsf{Cof}\cong \mathrm{d}_0^*\circ\cof\cong C \cong \Sigma\circ\mathrm{d}_1^*\circ\fib \cong \Sigma\circ 0^*\circ\mathsf{Fib}\colon\D^{[1]}\to\D$,
\item $\mathsf{Cof}\circ(\mathrm{d}_0)_!\cong(\mathrm{d}_0)_!\circ\mathrm{s}_0^*\cong(\mathrm{d}_2)_!\circ(\mathrm{d}_0)_!\cong\mathsf{Fib}\circ\mathrm{s}_0^*\colon\D\to\D^{[2]}$,
\item $\mathsf{Cof}\circ\mathrm{s}_0^*\cong(\mathrm{d}_0)_*\circ(\mathrm{d}_1)_*\cong(\mathrm{d}_2)_*\circ\mathrm{s}_0^*\cong\mathsf{Fib}\circ(\mathrm{d}_1)_*\colon\D\to\D^{[2]}$,
\item $\mathsf{Cof}\circ(\mathrm{d}_1)_*\cong\Sigma\circ\mathsf{Fib}\circ(\mathrm{d}_0)_!\colon\D\to\D^{[2]}$
\item $\mathrm{d}_0^*\circ\mathsf{Cof}\circ(\mathrm{d}_1)_*\cong\Sigma\circ(\mathrm{d}_0)_!\colon\D\to\D^{[1]}$,
\item $\mathrm{d}_2^*\circ\mathsf{Fib}\circ(\mathrm{d}_0)_!\cong\Omega\circ(\mathrm{d}_1)_!\colon\D\to\D^{[1]}$.
\end{enumerate}
Let now $n\geq 0$ and $M_0\cup M_1\cup M_2=\mathbf{n}$ be a partition. Let $m\colon\mathbf{n}\to\mathbf{3}$ be the map with $m^{-1}(j)=M_j$ for $j\in\mathbf{3}$. Then (i)-(vi) above generalize to the following isomorphisms (c.f. section 8 of \cite{bg:cubical}). In each case we use that morphisms of stable derivators commute with compositions of inverse images and homotopically finite Kan extensions.
\begin{enumerate}
\item $(\cof^{\underline{1}})^3\cong\Sigma^n\colon\D^{\cube{n}}\to\D^{\cube{n}}$,
\item $(\fib^{\underline{1}})^3\cong\Omega^n\colon\D^{\cube{n}}\to\D^{\cube{n}}$,
\item $m^*\circ\mathsf{Cof}^{\underline{1}}\cong\tcof\circ(\mathrm{d}_1^{M_0}\times\mathrm{d}_0^{M_1}\times\id^{M_2})^* \cong\tfib\circ(\id^{M_0}\times\mathrm{d}_1^{M_1}\times\mathrm{d}_0^{M_2})^*\circ\cof^{\underline{1}}\colon$ $\D^{\cube{n}}\to\D$,
\item $m^*\circ\mathsf{Fib}^{\underline{1}}\cong\tcof\circ(\mathrm{d}_1^{M_0}\times\mathrm{d}_0^{M_1}\times\id^{M_2})^* \circ\fib^{\underline{1}}\cong\tfib\circ(\id^{M_0}\times\mathrm{d}_1^{M_1}\times\mathrm{d}_0^{M_2})^*\colon$ $\D^{\cube{n}}\to\D$.
\end{enumerate}
\end{rmk}

\begin{prop}\label{prop:bicart-obstruction}
Let $n\geq 0$ and $\D$ a stable derivator. Then for $X\in\D^{\cube{n}}$ the following properties are equivalent.
\begin{enumerate}
\item $X\in\mathrm{essim}((\iota_{0,n-1})_!\colon\D^{\cube{n}_{0,n-1}}\to\D^{\cube{n}})$
\item $\tcof(X)=0\in\D$
\end{enumerate}
\end{prop}

\begin{proof}
This is \cite[Prop.~9.2]{bg:cubical}.
\end{proof}

\begin{cor}\label{prop:stability}
Let $n\geq 0$ and $\D$ a stable derivator. Then
\[
\mathrm{essim}((\iota_{0,n-1})_!\colon\D^{\cube{n}_{0,n-1}}\to\D^{\cube{n}})=\mathrm{essim}((\iota_{1,n})_*\colon\D^{\cube{n}_{1,n}}\to\D^{\cube{n}}).
\]
\end{cor}

\begin{proof}
This is \cite[Cor.~9.3]{bg:cubical}. The first equivalence below follows from \autoref{prop:bicart-obstruction}, the third equivalence from the dual version of this statement and for the second equivalence we invoke \cite[Rem.~8.27]{bg:cubical} for the relation $\tcof\cong\Sigma^n\circ\tfib$.
\begin{align}
X\in\mathrm{essim}((\iota_{0,n-1})_!)\iff&\tcof(X)=0\\
\iff&\tfib(X)=0\\
\iff& X\in\mathrm{essim}((\iota_{1,n})_*).
\end{align}
\end{proof}

\begin{rmk}
Let $n\geq 0$ and $\D$ a derivator. Then objects in $\mathrm{essim}((\iota_{0,n-1})_!)$ are called \textbf{cocartesian} $n$\textbf{-cubes} while objects in $\mathrm{essim}((\iota_{1,n})_*)$ are called \textbf{cartesian} $n$\textbf{-cubes} in $\D$. Hence in the case $n=2$ \autoref{prop:stability} just states the coincidence of cocartesian and cartesian squares for stable derivators, which is in fact the defining property of stability.
\end{rmk}

\begin{defn}
Let $n\geq 0$ and $\D$ a stable derivator.
\begin{enumerate}
\item The stable derivator $\D^{\cube{n},ex}:=\mathrm{essim}((\iota_{0,n-1})_!)=\mathrm{essim}((\iota_{1,n})_*)$ is called the \textbf{derivator of bicartesian} $n$\textbf{-cubes} in $\D$.
\item Objects in $\D^{\cube{n},ex}$ are called \textbf{bicartesian} $n$\textbf{-cubes} in $\D$.
\end{enumerate}
\end{defn}

\begin{rmk}
Let $n\geq 0$ and $\D$ a stable derivator.
\begin{enumerate}
\item To see that $\D^{\cube{n},ex}$ is in fact a stable derivator we observe that $\iota_{0,n-1}$ and $\iota_{1,n}$ are fully faithful and invoke \cite[Prop.~1.20]{groth:ptstab} to obtain equivalences of prederivators
\[
\D^{\cube{n}_{0,n-1}}\toiso\D^{\cube{n},ex}\qquad\text{and}\qquad\D^{\cube{n}_{1,n}}\toiso\D^{\cube{n},ex},
\]
respectively.
\item In \cite{bg:cubical} the derivator of bicartesian $n$-cubes in $\D$ is denoted by $\D^{\cube{n},n-1}$ and called the derivator of $(n-1)$-determined $n$-cubes.
\end{enumerate}
\end{rmk}

\begin{eg}
Let $n\geq 0$, $\D$ a stable derivator and $X\in\D^{\cube{n}}$ such that $\iota_{1,n-1}^*(X)=0$. Then
\[
X\in\D^{\cube{n},ex}\Rightarrow\infty^*(X)\cong\Sigma^n\circ\emptyset^*(X).
\]
This is \cite[Ex.~6.7]{bg:cubical}.
\end{eg}

\begin{prop}\label{prop:tcof-concat}
Let $n\geq 0$ and $\D$ a stable derivator. Let $X\in\D^{\cube{n-1}\times[2]}$. Then there is a natural cofiber sequence
\begin{equation}\label{eq:tcof-seq}
\tcof\circ\mathrm{d}_2^*(X)\to\tcof\circ\mathrm{d}_1^*(X)\to\tcof\circ\mathrm{d}_0^*(X)
\end{equation}
of total cofibers of $n$-cubes in $\D$.
\end{prop}

\begin{proof}
There is a functor $f\colon\cube{3}\to[2]$ defined by the diagram
\[
\xymatrix{
& 0 \ar[rr]\ar[dd] & & 2 \ar[dd]\\
0 \ar[rr]\ar[dd]\ar[ur] & & 1 \ar[dd]\ar[ur]\\
& 1 \ar[rr] & & 2\\
1 \ar[rr]\ar[ur] & & 1, \ar[ur]
}
\]
where the coordinates associated to $0,1,2\in\mathbf{3}$ are displayed in the horizontal, vertical and diagonal direction, respectively. We observe that
\begin{equation}\label{eq:3-cube}
f\circ\mathrm{d}_1^0=\mathrm{d}_2\circ\mathrm{s}_0^2\qquad\text{and}\qquad f\circ\mathrm{d}_0^0=\mathrm{d}_0\circ\mathrm{s}_0^1.
\end{equation}
We now consider a stable derivator $\D_0$, and define $\D_1=\D_0^{[1]}$. From \eqref{eq:3-cube} and \cite[Thm.~8.11]{gst:basic} we deduce that the essential image of $f^*\colon\D_0^{[2]}\to\D_0^{\cube{3}}$ is contained in $\D_1^{\cube{2},ex}$, where here the coordinates of $\cube{2}$ correspond to the coordinates 1 and 2 of $\cube{3}$. Furthermore, since the cone functor $C:\D_1\to\D_0$ is exact, it commutes with homotopically finite left Kan extensions and hence induces $C\colon\D_1^{\cube{2},ex}\to\D_0^{\cube{2},ex}$. Moreover, for $X_0\in\D_0^{[2]}$ the object $Y_0=C\circ f^*(X_0)\in\D_0^{\cube{2},ex}$ has by construction the underlying diagram
\begin{equation}\label{eq:cof-seq}
\xymatrix{
C\circ\mathrm{d}_2^*(X_0) \ar[r] \ar[d] & C\circ\mathrm{d}_1^*(X_0) \ar[d] \\
0 \ar[r] & C\circ\mathrm{d}_0^*(X_0).
}
\end{equation}
Finally, we set $\D_0=\D^{\cube{n-1}}$ and $X_0=X$ and apply $\tcof\colon\D^{\cube{n-1}}\to\D$ to $Y_0$. Using \eqref{eq:cof-seq} we conclude, that the bicartesian square $\tcof(Y_0)$ exhibits \eqref{eq:tcof-seq} as a cofiber sequence.
\end{proof}

\begin{cor}\label{cor:bicart-concat}
Let $n\geq 0$ and $\D$ a stable derivator. Let $X\in\D^{\cube{n-1}\times[2]}$. Then, if two of the three $n$-cubes $\mathrm{d}_2^*(X), \mathrm{d}_1^*(X)$ and $\mathrm{d}_0^*(X)$ are bicartesian, also the third one is bicartesian.
\end{cor}

\begin{proof}
This is immediate from \autoref{prop:bicart-obstruction}, \autoref{prop:tcof-concat} and the 2-out-of-3 property for zero-objects in cofiber sequences.
\end{proof}

\begin{rmk}
An alternative strategy to prove \autoref{cor:bicart-concat} relies on the corresponding unstable statements \cite[Prop.~7.20]{bg:cubical} and \autoref{prop:stability}.
\end{rmk}

\begin{prop}\label{prop:frankes-lemma}
Let $n\geq 0$, $\D$ be a derivator and $f\colon A\to B$ and $g\colon\cube{n}\to B$ functors between small categories such that there is a full subcategory $B'\subseteq B$ with
\begin{enumerate}
\item $f(A)\subseteq B'$
\item $g(\cube{n}_{0,n-1})\subseteq B'$ and $g(\infty)\notin B'$,
\item the functor $\cube{n}_{0,n-1}\to B'_{/g(\infty)}$ is a right adjoint
\end{enumerate}
Then the essential image of $f_!\colon\D^A\to\D^B$ consists objects $X\in\D^B$ such that $g^*(X)$ is cocartesian.
\end{prop}

\begin{proof}
This is a special case of \cite[Lem.~7.6]{bg:cubical}.
\end{proof}

We conclude this section with a new characterization of the total cofiber morphism.

\begin{prop}\label{prop:tcof-new}
Let $n\geq 0$ and $\D$ a stable derivator. Consider the factorization of the inclusion $\mathrm{d}_1^{0}\colon\cube{n}\to\cube{n+1}$
\[
\cube{n}\xrightarrow{\alpha}\cube{n+1}_{0,n}\xrightarrow{\beta}\cube{n+1}.
\]
Then there is an natural isomorphism $\tcof\cong\infty^*\circ\beta_!\circ\alpha_*\colon\D^{\cube{n}}\to\D$.
\end{prop}

\begin{proof}
The essential image of $\beta_!\circ\alpha_*$ is contained in the full subderivator of bicartesian $n+1$-cubes. As a consequence, \cite[Cor.~9.9]{bg:cubical} yields an isomorphism $\tcof\circ(\mathrm{d}_1^{0})^*\circ\beta_!\circ\alpha_*\toiso\tcof\circ(\mathrm{d}_0^{0})^*\circ\beta_!\circ\alpha_*$. Since $\alpha$ and $\beta$ are fully faithful we invoke \cite[Prop.~1.20]{groth:ptstab} for isomorphisms $(\mathrm{d}_1^{0})^*\circ\beta_!\circ\alpha_*\cong\id_{\cube{n}}$. On the other hand $\alpha_*$ is an extension-by-zero morphism \cite[Prop.~1.23]{groth:ptstab}. We conclude that essential image of $(\mathrm{d}_0^{0})^*\circ\beta_!\circ\alpha_*$ coincides with the essential image of $\infty'_!\colon\D\to\D^{\cube{n}}$ (where $\infty'$ denotes the final object of $\cube{n}$) and using the previously mentioned results of \cite{groth:ptstab} again we obtain  the first isomorphism in
\[
\tcof\circ(\mathrm{d}_0^{0})^*\circ\beta_!\circ\alpha_*\cong\tcof\circ\infty'_!\circ\infty'^*\circ(\mathrm{d}_0^{0})^*\circ\beta_!\circ\alpha_*\cong\infty^*\circ\beta_!\circ\alpha_*.
\]
The second isomorphism is \cite[Lem.~8.19]{bg:cubical}.
\end{proof}

\section{The derivators \texorpdfstring{$\D_{n,k}$}{Dnk}}
\label{sec:dnk}

In this chapter we introduce for a given stable derivators $\D$ a family of stable derivator $\D_{n,k}$ parametrized by pairs of natural numbers. These will be the objects of main interest in this work, and their properties will be analyzed in the forthcoming chapters.

\begin{defn}
Let $P \subset \mathbb{Z}^{k+1}$ be a subposet, $x=(x_0,\cdots,x_k) \in P$ and $\D$ a stable derivator.
\begin{enumerate}
\item If $\cube{}_x := \lbrace (x_0+\delta_0,\cdots,x_k+\delta_k) \vert \delta_i \in \lbrace 0,1 \rbrace \text{ for } i \in \lbrace 0,\cdots,k\rbrace \rbrace \subset P$ we call this the elementary subcube of $P$ starting in $x$.
\item If $\cube{}_x \not\subset P$ we say that $P$ does not contain the subcube starting in $x$.
\item Let $x \in P$ such that $\cube{}_x \subset P$ and $X \in \D^P$, then we call 
\[
\cube{}_x(X) := X \vert_{\cube{}_x} \in \D^{\cube{k+1}}
\]
the elementary subcube of $X$ starting in $x$.
\end{enumerate}
\end{defn}

\begin{eg}
The subposet $\underline{\Lambda}_{n,k}\subset\mathbb{Z}^{k+1}$ contains $\cube{}_f$ if and only if $f$ is injective.
\end{eg}

\begin{defn}
Let $P \subset \mathbb{Z}^{k+1}$ be a subposet. We say an object $X \in \D^P$ satisfies property (P1), resp. property (P2) at a point $x \in P$ if the following condition holds:

(P1) $\cube{}_x \subset P$ and $\cube{}_x(X)$ is bicartesian.

(P2) $X \vert_x = 0$.
\end{defn}

Recall from \autoref{cor:coord-rep} that the choice of coordinates $0_\bullet$ induces embeddings of posets $\underline{\Lambda}_{n,k}\subset\mathbb{Z}^{k+1}$.

\begin{defn}\label{defn:dnk}
Let $n\geq -k+1, k\geq 2$. The derivator $\D_{n,k}$ is the full subderivator of $\D^{\underline{\Lambda}_{n+k-1,k-1}}$ spanned by those objects $X \in \D^{\underline{\Lambda}_{n+k-1,k-1}}$ satisfying
\begin{enumerate}
\item property (P1) for all $f \in \underline{\Lambda}_{n+k-1,k-1}$ that are injective,
\item property (P2) for all $f \in \underline{\Lambda}_{n+k-1,k-1}$ that are not injective.
\end{enumerate}
\end{defn}

We note, that \autoref{lem:subder} implies that $\D_{n,k}$ is stable, as $\D^{\underline{\Lambda}_{n+k-1,k-1}}$ is stable. 

\begin{rmk}
The index shift appearing in \autoref{defn:dnk} is motivated by the following observations.
\begin{enumerate}
\item The index $n$ of $\D_{n,k}$ refers to the dimension of a maximal subsimplex. More precisely, for all $l \in \mathbf{k}$ the natural number $n$ is exactly the maximal number such that there exists an injective coordinate $(x_0,\cdots,x_k) \in \underline{\Lambda}_{k-1,n+k}$ with $(x_0,\cdots,x_{l-1},x_l+m,x_{l+1},\cdots,x_k) \in \underline{\Lambda}_{n+k-1,k-1}$ injective for all $m \in \mathbf{n}$. In particular, by \cite[Thm.~4.5.]{gst:Dynkin-A} we have equivalences
\[
\D_{n,2}\cong\D^{[n]}.
\]
\item The index $k$ of $\D_{n,k}$ refers to the dimension of the subcubes which are forced to be bicartesian.
\end{enumerate}
\end{rmk}

\begin{rmk}
The property (P2) for non-injective objects implies immediately, that there are equivalences $\D_{n,k}\cong 0$ for $n\leq -1$ (since this assumption implies that all objects in $\underline{\Lambda}_{n+k-1,k-1}$ are not injective).
\end{rmk}

\begin{rmk}
From the discussion in \autoref{egs:subders} it is clear that the inclusions $\D_{n,k}\subseteq\D^{\underline{\Lambda}_{n+k-1,k-1}}$ are indeed inclusions of full subderivators in the sense of \autoref{defn:subder}. Nevertheless, we carry out the details for the case of property (P1) for an injective object $x\in\underline{\Lambda}_{n+k-1,k-1}$ and the left Kan extension morphism along a functor $u\colon A\to B$ between small categories $A,B\in Cat$. The diagram
\[
\xymatrix{
\D(\underline{\Lambda}_{n+k-1,k-1}\times A) \ar[rr]^{(\square_x\times\id)^*}\ar[d]_{(\id\times u)_!} && \D(\cube{k}\times A) \ar[d]_{(\id\times u)_!}&& \D(\cube{k}_{0,k-1}\times A)\ar[d]_{(\id\times u)_!}\ar[ll]_{(\iota_{0,k-1}\times\id)_!}\\
\D(\underline{\Lambda}_{n+k-1,k-1}\times B) \ar[rr]_{(\square_x\times\id)^*} && \D(\cube{k}\times B) && \D(\cube{k}_{0,k-1}\times B)\ar[ll]^{(\iota_{0,k-1}\times\id)_!}
}
\]
commutes up to natural isomorphisms, since $u_!$ is a morphism of derivators. We assume that $X\in\D(\underline{\Lambda}_{n+k-1,k-1}\times A)$ satisfies property (P1) in $x$. This means exactly that there exists $Y\in\D(\cube{k}_{0,k-1}\times A)$ such that 
\[
(\square_x\times\id)^*(X)\cong(\iota_{0,k-1}\times\id)_!(Y).
\]
Then the above diagram yields
\begin{align}
(\square_x\times\id)^*\circ(\id\times u)_!(X)\cong&(\id\times u)_!\circ(\square_x\times\id)^*(X)\\
\cong&(\id\times u)_!\circ(\iota_{0,k-1}\times\id)_!(Y)\\
\cong&(\iota_{0,k-1}\times\id)_!\circ(\id\times u)_!(Y).
\end{align}
This is property (P1) for $(\id\times u)_!(X)$ at $x$.
\end{rmk}

In the following, we will discuss the most important automorphisms of the derivators $\D_{n,k}$.

\begin{cor}\label{cor:dnk-sym}
Let $n\ge -k+1, k\geq 2$ and $\D$ a stable derivator. The inverse images of the symmetry automorphisms 
\[
\mathsf{s}_1,\mathsf{s}_2\colon\underline{\Lambda}_{n+k-1,k-1}\to\underline{\Lambda}_{n+k-1,k-1}
\]
restrict to automorphisms
\[
\mathsf{s}_1^*,\mathsf{s}_2^*\colon\D_{n,k}\to\D_{n,k}.
\]
\end{cor}

\begin{proof}
\autoref{prop:sym-inj} implies that all of the four operations $\mathsf{s}_1,\mathsf{s}_2,\mathsf{s}_1^{-1}$ and $\mathsf{s}_2^{-1}$ preserve injective objects. Hence all of the maps also preserve non-injective objects. This yields the compatibility of $\mathsf{s}_1^*$ and $\mathsf{s}_2^*$ with the property (P2) on non-injective objects. Let $e_j\in\mathbb{Z}^k$, $j\in\lbrace 0,\cdots,k-1\rbrace$ denote the $j$th basis vector. Then we invoke \autoref{prop:coord-sym} for $\mathsf{s}_1(f+e_j)=\mathsf{s}_1(f)+e_j$, $\mathsf{s}_2(f+e_j)=\mathsf{s}_2(f)+e_{j-1}$ if $j\geq1$ and $\mathsf{s}_2(f+e_0)=\mathsf{s}_2(f)+e_{k-1}$. This shows that postcomposition with $\mathsf{s}_1$ and $\mathsf{s}_2$ maps elementary subcubes to cyclic permutations of elementary subcubes. Since bicartesian cubes are stable under permutation of coordinates \cite[Cor.~4.10]{bg:cubical}, we conclude that $\mathsf{s}_1^*$ and $\mathsf{s}_2^*$ are compatible with the property (P1) on injective objects.
\end{proof}

\begin{rmk}
\autoref{prop:para-inj} implies that $\underline{\Lambda}_{n+k-1,k-1}^{inj}\cong\underline{\Lambda}_{n-1,k-1}$. It will become clear later, that the derivator $\D_{n,k}$ has many properties one would expect from an object associated to $\underline{\Lambda}_{n-1,k-1}$. Moreover, by \autoref{prop:sym-inj} the automorphism $\mathsf{s}_1^{-1}\circ\mathsf{s}_2\colon\underline{\Lambda}_{n+k-1,k-1}\to\underline{\Lambda}_{n+k-1,k-1}$ restricts to $\mathsf{s}_2$ on the subposet of injective morphisms. This motivates us to use the notation 
\[
\mathsf{s}_3:=\mathsf{s}_1^{-1}\circ\mathsf{s}_2\colon\underline{\Lambda}_{n+k-1,k-1}\to\underline{\Lambda}_{n+k-1,k-1}.
\]
\end{rmk}

\begin{defn}
Let $n\geq -k+1,k\geq 2$ and $\xi=(0,1,\cdots,k-1)\in\underline{\Lambda}_{n+k-1,k-1}$. Consider the following full subcategories of $\underline{\Lambda}_{n+k-1,k-1}$.
\begin{enumerate}
\item The \textbf{fundamental domain} $Do_{n,k} := (\underline{\Lambda}_{n+k-1,k-1})_{\xi / \mathsf{s}_3^k(\xi)}$ with inclusion \[
do_{n,k}\colon Do_{n,k} \longrightarrow \underline{\Lambda}_{n+k-1,k-1}.
\]
\item The \textbf{fundamental slice} $Sl_{n,k} := (\underline{\Lambda}_{n+k-1,k-1})_{\xi / \mathsf{s}_3^{k-1}(\xi)}$ with inclusion \[
sl_{n,k}\colon Sl_{n,k} \longrightarrow \underline{\Lambda}_{n+k-1,k-1}.
\]
\end{enumerate}
\end{defn}

\begin{defn}
Let $n \geq -k+1,k\geq 2$.
\begin{enumerate}
\item The \textbf{derivator of fundamental domains} $do\D_{n,k}$ is the full subderivator of $\D^{Do_{n,k}}$ spanned by those objects $X \in \D^{Do_{n,k}}$ that satisfy
\begin{enumerate}
\item property (P1) for all $x \in Do_{n,k}$ such that $\cube{}_x\subset Do_{n,k}$,
\item property (P2) for all $x \in Do_{n,k}$ such that $do_{n,k}(x)$ is not injective.
\end{enumerate}
\item The \textbf{derivator of fundamental slices} $sl\D_{n,k}$ is the full subderivator of $\D^{Sl_{n,k}}$ spanned by those objects $X \in \D^{Do_{n,k}}$ that satisfy property (P2) for all $x\in Sl_{n,k}$ such that $sl_{n,k}(x)$ is not injective.
\end{enumerate}
\end{defn}

\begin{thm}\label{thm:slices}
Let $n\geq -k+1,k\geq 2$ and $\D$ a stable derivator. The inverse images of $do_{n,k}$ and $sl_{n,k}$ restrict to equivalences
\[
do_{n,k}^*\colon\D_{n,k}\to do\D_{n,k} \qquad\text{and}\qquad sl_{n,k}^*\colon\D_{n,k}\to sl\D_{n,k}.
\]
\end{thm}

\begin{proof}
We show the statement for $sl_{n,k}^*$. The statement for $do_{n,k}^*$ is very similar.
We consider the following subposets of $\underline{\Lambda}_{n+k-1,k-1}$
\begin{itemize}
\item $A_1=Sl_{n,k}\cup\lbrace x\in\underline{\Lambda}_{n+k-1,k-1}\text{ non-injective}\vert \exists y\in sl_{n,k}\colon y\leq x\rbrace$,
\item $A_2=A_1\cup\lbrace x\in\underline{\Lambda}_{n+k-1,k-1}\vert \exists y\in sl_{n,k}\colon y\leq x\rbrace$,
\item $A_3=A_2\cup\lbrace x\in\underline{\Lambda}_{n+k-1,k-1}\text{ non-injective}\vert \exists y\in sl_{n,k}\colon x\leq y\rbrace$
\end{itemize}
with inclusions
\[
Sl_{n,k}\xrightarrow{u_1}A_1\xrightarrow{u_2}A_2\xrightarrow{u_3}A_3\xrightarrow{u_4}\underline{\Lambda}_{n+k-1,k-1}.
\]
We claim that $u:=(u_4)_*\circ(u_3)_!\circ(u_2)_!\circ(u_1)_*$ is inverse to $sl_{n,k}^*$. We observe that the inclusions $u_i$, $1\leq i\leq 4$ are fully faithful. Then \cite[Prop.~1.20]{groth:ptstab} yields that Kan extensions along $u_i$ are also fully faithful and that the units of $(u_i)_!\dashv u_i^*$ and counits of $u_i^*\dashv (u_i)_*$ are invertible. Hence it is sufficient to show the essential image of $u$ is contained in $\D_{n,k}$. Since $u_1$ is a sieve, the corresponding right Kan extension is an extension-by-zero morphism \cite[Prop.~1.23]{groth:ptstab}. Using that the units of $(u_i)_!\dashv u_i^*$ and counits of $u_i^*\dashv (u_i)_*$ for $i\geq 2$ are invertible, we conclude that objects in the essential image of $u$ satisfy property (P2) for non-injective objects in $A_1$. We claim that objects in the essential image of $(u_2)_!\colon\D^{A_1}\to\D^{A_2}$ satisfy property (P1) for all injective objects in $A_2$. Let $x$ be an injective objects in $A_2$. We define $B_x=A_1\cup\lbrace y\in A_2\vert y\leq \mathsf{s}_1(x)\rbrace$ and consider the inclusions $A_1\xrightarrow{v_x}B_x\xrightarrow{w_x}A_2$. We observe that $u_2=w_x\circ v_x$. Hence $(u_2)_!=(w_x)_!\circ (v_x)_!$ by the pseudofunctoriality of right Kan extensions.
In particular we have $\square_x^*\circ u_!=\square_x^*\circ (v_x)_!\circ (u_1)_*$. We show that we are now in a situation where \autoref{prop:frankes-lemma} applies. For this we consider $B_x'=B_x\setminus\lbrace \mathsf{s}_1(x)\rbrace$. We observe that $(B'_x)_{/\mathsf{s}_1(x)}=\lbrace y\in B'_x\vert y\leq \mathsf{s}_1(x)\rbrace$. To show that the restriction $\square_x\vert_{\cube{k}_{0,k-1}}\colon\cube{k}_{0,k-1}\to(B'_x)_{/\mathsf{s}_1(x)}$ is a right adjoint, we use the fact that an embedding of posets $q\colon P\to S$ is a right adjoint if for every element $s$ not in the image of $q$ the set $\lbrace p\in P\vert q(p)\geq s\rbrace$ admits a unique minimal element. To establish this in our situation we use that for $b\in(B'_x)_{/\mathsf{s}_1(x)}$ the implication
\[
b\leq\square_x(M)\land b\leq\square_x(N) \Rightarrow b\leq\square_x(M\cap N) 
\]
holds by construction of $\square_x$ for $M,N\in\cube{k}$. Therefore, \autoref{prop:frankes-lemma} yields the cocartesianess of objects in the image of $\square_x^*\circ(v_x)_!$. Moreover, \autoref{prop:stability} yields the bicartesianess of these squares. By using again that the units of $(u_i)_!\dashv u_i^*$ and counits of $u_i^*\dashv (u_i)_*$ are invertible and the dual arguments for $u_3$ and $u_4$, we conclude that the essential image of $u$ is indeed contained in $\D_{n,k}$.
\end{proof}

\begin{rmk}\label{rmk:slices}
 Let $n\geq 0,k\geq 2$ and $\D$ a stable derivator.
\begin{enumerate}
\item Let $i\in\lbrace 1,2,3\rbrace$. \autoref{thm:slices} allows us to define the symmetry operations on derivators of slices and domains via
\begin{itemize}
\item $\mathsf{s}_i^*=sl_{n,k}^*\circ\mathsf{s}_i^*\circ(sl_{n,k}^*)^{-1}\colon sl\D_{n,k}\to sl\D_{n,k}$,
\item $\mathsf{s}_i^*=do_{n,k}^*\circ\mathsf{s}_i^*\circ(do_{n,k}^*)^{-1}\colon do\D_{n,k}\to do\D_{n,k}$.
\end{itemize}
We emphasize that, although the notation might suggest that these morphisms are restrictions of inverse images, this is in general not the case.
\item Let $sd_{n,k}\colon Sl_{n,k}\to Do_{n,k}$ be the inclusion. Since $sl_{n,k}=do_{n,k}\circ sd_{n,k}$ it follows from \autoref{thm:slices} that $sd_{n,k}^*\colon do\D_{n,k}\to sl\D_{n,k}$ is invertible. Moreover, by construction (recall that the units of $(u_i)_!\dashv u_i^*$ and counits of $u_i^*\dashv (u_i)_*$ are invertible) the inverse can be constructed as $(u'_2)_!\circ(u'_1)_*$, where
\[
Sl_{n,k}\xrightarrow{u'_1}Do_{n,k}\setminus\lbrace f=(f_0,\cdots,f_k)\in Do_{n,k}\vert f \text{ injective}, f_0\geq 1\rbrace\xrightarrow{u'_2}Do_{n,k}
\]
are the inclusions.
\item Consider the map $\tilde{sd}_{n,k}\colon Sl_{n,k}\to Do_{n,k},f\mapsto\mathsf{s}_3(f)$. Then the diagram
\[
\xymatrix{
&&\D_{n,k}\ar[lld]_{sl_{n,k}^*}\ar[d]^{do_{n,k}^*}\ar[rrd]^{\mathsf{s}_3^*\circ sl_{n,k}^*}\\
sl\D_{n,k}&&do\D_{n,k}\ar[ll]^{sd_{n,k}*}\ar[rr]_{\tilde{sd}_{n,k}^*}&&sl\D_{n,k}
}
\
\]
commutes. This yield an isomorphism
\[
\tilde{sd}_{n,k}^*\circ (sd_{n,k}^*)^{-1}\cong \mathsf{s}_3^*\circ sl_{n,k}^*\circ(sl_{n,k}^*)^{-1}\cong\mathsf{s}_3^*.
\]
Thus we regard $sl\D_{n,k}$ as a minimal model of $\D_{n,k}$, and $do\D_{n,k}$ as minimal among those models where $\mathsf{s}_3^*$ is computable as an inverse image.
\end{enumerate}
\end{rmk}

\begin{egs}\label{egs:simplex-slice}
Let $\D$ be a stable derivator.
\begin{enumerate}
\item Let $n\geq 0$. We consider the case $k=2$. Then
\[
\xi=(0,1)\qquad\text{and}\qquad\mathsf{s}_3(\xi)=(0,n+1).
\]
Hence the assignment $c_n\colon[n]\to Sl_{n,2}, i\mapsto (0,i+1)$ defines an isomorphism of categories. As a consequence, \autoref{thm:slices} provides us with an equivalence of derivators
\[
\D_{n,2}\xrightarrow{sl_{n,2}^*} sl\D_{n,2}\xrightarrow{c_n^*}\D^{[n]}.
\]
We emphasize that this special case is exactly \cite[Thm.~4.5]{gst:Dynkin-A}.
\item Let $n=0$ and $k\geq 2$. In this case we have $\mathsf{s}_3=\id_{\underline{\Lambda}(\Lambda_{k-1},\Lambda_{k-1})}$. In particular, there are isomorphisms $\bbone\cong Sl_{0,k}\cong Do_{0,k}$. As a consequence, \autoref{thm:slices} implies that
\[
\xi^*\colon\D_{0,k}\toiso sl\D_{0,k}\cong \D
\]
is an equivalence of derivators.
\item Let now $n=2$ and $k=3$. Then
\[
\xi=(0,1,2),\qquad\mathsf{s}_3(\xi)=(0,1,4)\qquad\text{and}\qquad\mathsf{s}_3^2(\xi)=(0,3,4).
\]
We define the category $X=[2]\times[2]\setminus(2,0)$. Then the assignment 
\[
X\to Sl_{2,3},(i,j)\mapsto(0,i+1,j+2)
\]
defines an isomorphism of categories. Moreover, by \autoref{thm:slices} the derivator $\D_{2,3}$ is equivalent to the full subderivator of $\D^X$ spanned by those objects $x$ with $(1,0)^*x=0$ and $(2,1)^*x=0$, i.e. objects such that the underlying diagram is of the form
\[
\xymatrix{
x_{0,0}\ar[r]\ar[d] & x_{0,1} \ar[r]\ar[d] & x_{0,2} \ar[d]\\
0 \ar[r] & x_{1,1} \ar[r]\ar[d] & x_{1,2} \ar[d]\\
& 0\ar[r] & x_{2,2}.
}
\]
\end{enumerate}
\end{egs}

\begin{rmk}\label{rmk:higher-Auslander}
Let $\mathrm{k}$ be a field. It is immediate from Iyama's inductive construction (\cite[Thm.~1.18]{iyama-higherA}, c.f. also \cite[Thm.~3.4]{OT-cluster}) of the $k$-Auslander algebra $T^{(k)}_n(\mathrm{k})$ of the $\A{n}$-quiver, that $T^{(k)}_n(\mathrm{k})$ can be described by the quiver generated by the injective objects in $Sl_{n+1,k-1}$ and two types of relations. The first type corresponds exactly to the commutativity of all existing elementary subcubes and the second type type of relations deals with the vanishing of compositions, which factor through a non-injective object in $Sl_{n-1,k+1}$. In particular, there are exact equivalences of triangulated categories
\[
(\D_{\mathrm{k}})_{n-1,k+1}(\bbone)\cong sl(\D_{\mathrm{k}})_{n-1,k+1}(\bbone)\cong D(T^{(k)}_n(\mathrm{k})).
\]
Furthermore, the automorphism $\mathsf{s}_1^{-1}\colon(\D_{\mathrm{k}})_{n-1,k+1}\to(\D_{\mathrm{k}})_{n-1,k+1}$ corresponds under these equivalences to the $k$-Auslander-Reiten-translate (\cite[\S1]{iyama-higherA}). This follows from the relation $\mathsf{s}_1^{-1}=\mathsf{s}_3\circ\mathsf{s}_2^{-1}$ by identifying $\mathsf{s}_2\cong\Sigma^k$ (\autoref{cor:dnk-sigma}) and relating $\mathsf{s}_3$ to the Serre-functor (\autoref{rmk:Serre}).
\end{rmk}

\begin{prop}\label{prop:slnk-sigma}
Let $n\geq 1,k\geq 2$ and $\D$ a stable derivator. There is a natural isomorphism
\[
\Sigma^{k-1}\cong sl_{n,k}^*\circ\mathsf{s}_2^*\circ(sl_{n,k}^*)^{-1}\colon sl\D_{n,k}\to sl\D_{n,k}.
\]
\end{prop}

\begin{proof}
Consider the map $J\colon Sl_{n,k}\times\cube{k}\to\underline{\Lambda}_{n+k-1,k-1}$ defined by the assignment $((f_0,\cdots,f_{k-1}),(\delta_0,\cdots,\delta_{k-1}))
\mapsto$
\[
(f_0+\delta_0(f_1-f_0),\cdots,f_{k-2}+\delta_{k-2}(f_{k-1}-f_{k-2}),f_{k-1}+\delta_{k-1}(f_0+n+k-f_{k-1}).
\]
It is easy to see that $J$ is a well defined morphism of posets (and hence a functor). Let $f=(f_0,\cdots f_{k-1})\in Sl_{n,k}$ and $\delta=(\delta_0,\cdots,\delta_{k-1})\in\cube{k}\setminus\lbrace\emptyset,\infty\rbrace$. We claim, that in this case $J(f,\delta)\in \underline{\Lambda}_{n+k-1,k-1}$ is not injective. If $\delta_0=1$, there exists by assumption $i\leq k-1$ minimal with $\delta_i=0$. Hence $J(f,\delta)_{i-1}=J(f,\delta)_i$. If $\delta_0=0$, there exists by  assumption $i\leq k-1$ maximal with $\delta_i=1$. If $i\leq k-2$, we have $J(f,\delta)_{i}=J(f,\delta)_{i+1}$. In the remaining case $i=k-1$ we observe $J(f,\delta)_{k-1}=J(f,\delta)_{0}+n+k$ to conclude the claim.\\
Moreover, we observe 
\begin{equation}\label{eq:J-empty}
J\circ(\id\times\emptyset)=sl_{n,k}
\end{equation}
and 
\begin{equation}\label{eq:J-infty}
J\circ(\id\times\infty)=\mathsf{s}_2\circ sl_{n,k}.
\end{equation}
As a consequence, the inverse image of $J$ restricts to a morphism
\[
J^*\colon\D_{n,k}\to (sl\D_{n,k})^{\cube{k}}.
\]
We claim furthermore, that for all $f\in Sl_{n,k}$ and $X\in\D_{n,k}$ the $k$-cube $(f\times\id)^*(X)$ is bicartesian. In fact, if $f$ is not injective, the cube in question is constantly zero (in particular bicartesian). Otherwise it is a concatination of the bicartesian cubes $\cube{}_g(X)$ for $g=f+(\kappa_0,\cdots,\kappa_{k-1})$ with $0\leq\kappa_i\leq f_{i+1}-f_i-1$ for $0\leq i\leq k-2$ and $0\leq\kappa_{k-1}\leq f_0+n+k-f_{k-1}-1$, and hence by \autoref{cor:bicart-concat} bicartesian. Therefore the essential image of $J^*$ is contained in $(sl\D_{n,k})^{\cube{k},ex}$. As a consequence of this, the first claim and \eqref{eq:J-empty}, we obtain by \cite[Ex.~6.7]{bg:cubical} the identification 
\[
(\id\times\infty)^*\circ J^*\cong\Sigma^{k-1}\circ sl_{n,k}^*.
\]
Finally, \eqref{eq:J-infty} yields $sl_{n,k}^*\circ\mathsf{s}_2^*\cong\Sigma^{k-1}\circ sl_{n,k}^*$ and we conclude by precomposing with $(sl_{n,k}^*)^{-1}$.
\end{proof}

\begin{cor}\label{cor:dnk-sigma}
Let $n\geq 1, k\geq 2$ and $\D$ a stable derivator. There is a natural isomorphism
\[
\Sigma^{k-1}\cong\mathsf{s}_2^*\colon\D_{n,k}\to\D_{n,k}.
\]
\end{cor}

\begin{proof}
\autoref{prop:slnk-sigma} implies the second isomorphism in
\[
\Sigma^{k-1}\cong(sl_{n,k}^*)^{-1}\circ\Sigma^{k-1}\circ sl_{n,k}^*\cong\mathsf{s}_2^*\colon \D_{n,k}\to\D_{n,k},
\]
and the first isomorhism is induced by the exactness of $sl_{n,k}^*$.
\end{proof}

\begin{rmk}
Let $n\geq 1,k\geq 2$ and $\D$ a stable derivator. \autoref{thm:slices} implies immediately, that an object $X\in\D_{n,k}$ is completely determined by $do_{n,k}^*(X)$. Using \autoref{cor:dnk-sigma} we can relate the objects $f^*(X)$ for $f\in\underline{\Lambda}_{n+k-1,k-1}$ arbitrary directly to $do_{n,k}^*(X)$. For this we note that $Do_{n,k}$ contains $\mathsf{i}(g)$ for every injective map $g\colon\Delta_{k-1}\to\Delta_{n+k-1}$. Hence, for an arbitrary injective $f\in\underline{\Lambda}_{n+k-1,k-1}$, we obtain
\[
f^*(X)\cong(\Sigma^{k-1})^l\circ\mathsf{i}(g)^*\circ do_{n,k}^*(X)
\]
for the unique $l\in\mathbb{Z}$ induced by \autoref{prop:lambda-shift}.
\end{rmk}

\begin{cor}\label{cor:fracCY}
Let $n\geq 1,k\geq 2$ and $\D$ a stable derivator. There are natural isomorphisms
\begin{equation}\label{eq:fracCY}
(\mathsf{s}_3^*)^{n+k}\cong(\mathsf{s}_2^*)^n\cong\Sigma^{n(k-1)}\colon\D_{n,k}\to\D_{n,k}.
\end{equation}
\end{cor}

\begin{proof}
The second isomorphism is \autoref{cor:dnk-sigma}. For the first isomorphism we invoke \autoref{lem:symmetries} for the relation $\mathsf{s}_1^{n+k}=\mathsf{s}_2^k$ of automorphisms of $\underline{\Lambda}_{n+k-1,k-1}$. Hence
\[
(\mathsf{s}_3^*)^{n+k}=(\mathsf{s}_2^*\circ(\mathsf{s}_1^*)^{-1})^{n+k}=(\mathsf{s}_2^*)^{n+k}\circ(\mathsf{s}_2^*)^{-k}=(\mathsf{s}_2^*)^{n}.
\]
\end{proof}

\begin{rmk}\label{rmk:Serre}
For a field $\mathrm{k}$ and $\D=\D_{\mathrm{k}}$ the autoequivalences $\mathsf{s}_3$ induce Serre equivalences on underlying categories. Moreover, for $\D$ general it was shown in the special case of $k=2$ that the autoequivalences $\mathsf{s}_3^*$ can be considered as Serre equivalences in a derivator theoretic sense \cite[Thm.~11.12]{gst:Dynkin-A}. Hence \autoref{cor:fracCY} suggests that the derivators $\D_{n,k}$ have a fractionally Calabi--Yau dimension of $\frac{n(k-1)}{n+k}$. We refer to \cite{Ladkani-derived} for some related examples of fractionally Calabi--Yau categories. We also observe that the enumerator and denominator of this fraction both are invariant under the assignment
\begin{equation}\label{eq:assignment}
(n,k)\mapsto(k-1,n+1).
\end{equation}
We will strengthen this observation with \autoref{thm:main} by constructing equivalences 
\[
\Phi_{n,k}\colon\D_{n,k}\to\D_{k-1,n+1},
\]
which commute with $\mathsf{s}_3^*$ and $\Sigma$. In fact, the first step towards this result, the special case $k=2$ (or equivalently $n=1$) will be the main content of \S\ref{sec:contra}.

\noindent Furthermore, it is worth to mention that infinite chains of adjunctions very often turn out to be 2-periodic with respect to Serre equivalence. This observation will be investigated more closely in \cite{bg:global}. We will encounter examples of such infinite chains of adjunctions in \autoref{cor:gsdver} and \autoref{cor:gsd-dual}. In these cases the fractionally Calabi--Yau property will be very useful for the understanding of iterated adjoints.
\end{rmk}

\noindent For technical reasons it will become useful in the proof of \autoref{thm:main} to extend objects of $\D_{n,k}$ to all of $\mathbb{Z}^k$ by 0. The next proposition shows, that this is always possible.

\begin{prop}\label{prop:dnk-hat}
Let $n\geq -k+1, k\geq2$. Let $\hat{\D}_{n,k}$ be the full subderivator of $\D^{\mathbb{Z}^k}$ spanned by those objects $X \in \D^{\mathbb{Z}}$ that satisfy
\begin{enumerate}
\item property (P1) for all $x \in \mathbb{Z}^k$ representing an injective object in $\underline{\Lambda}_{n+k-1,k-1}$,
\item property (P2) for all other $x \in \mathbb{Z}^k$.
\end{enumerate}
Then the Kan extensions $\mathsf{j}_!,\mathsf{j}_*\colon\D^{\underline{\Lambda}_{n+k-1,k-1}}\to\D^{\mathbb{Z}^k}$ restrict to morphisms \[
\mathsf{j}_!,\mathsf{j}_*\colon\D_{n,k}\to\hat{\D}_{n,k}.
\]
Moreover, the restrictions coincide and are equivalences.
\end{prop}

\begin{proof}
We show that the restriction of $\mathsf{j}_*$ is well defined and an equivalence. Since $\mathsf{j}\colon\underline{\Lambda}_{n+k-1,k-1}\to\mathbb{Z}^k$ is fully faithful, the same is true for $\mathsf{j}_*$ \cite[Prop.~1.20]{groth:ptstab}. Hence it is sufficient to identify the essential image of $\mathsf{j}_*$ with $\hat{\D}_{n,k}$. For this it is enough (\cite[Lem.~1.21]{groth:ptstab}) to show that for $x=(x_0,\cdots,x_{k-1}) \in \mathbb{Z}^k\setminus\underline{\Lambda}_{n+k-1,k-1}$ and $X \in\D_{n,k}$ we have $x^*\mathsf{j}_*(X)\cong 0$. Since $\underline{\Lambda}_{n+k-1,k-1}$ and $\mathbb{Z}^k$ are posets, we can identify $\underline{\Lambda}_{n+k-1,k-1,x/}$ with the full subposet of $\underline{\Lambda}_{n+k-1,k-1}$ on those objects which are pointwise greater or equal then $x$. Let $i_x$ denote the inclusion of these posets. Axiom (Der4) implies 
\begin{equation}\label{eq:dhat}
x^*\circ\mathsf{j}_*\cong(\pi_{\underline{\Lambda}_{n+k-1,k-1,x/}})_*\circ i_x^*.
\end{equation}
We claim that $\underline{\Lambda}_{n+k-1,k-1,x/}$ admits a minimal object. To see this we note, that if $x\leq x^1$ and $x\leq x^2$ for $x^1=(x_0^1,\cdots,x_{k-1}^1),x^2=(x_0^2,\cdots,x_{k-1}^2)\in\underline{\Lambda}_{n+k-1,k-1}$, then also $\mathrm{min}(x^1,x^2)=(\mathrm{min}(x^1_0,x^2_0),\cdots,\mathrm{min}(x^1_{k-1},x^2_{k-1}))\in\underline{\Lambda}_{n+k-1,k-1}$ and by construction $x\leq\mathrm{min}(x^1,x^2)$.
Next, we claim that the minimal object $x^0=(x_0^0,\cdots,x_{k-1}^0)\in\underline{\Lambda}_{n+k-1,k-1,x/}$ is not injective in $\underline{\Lambda}_{k-1,n+k-1}$. For this we assume that $x^0$ is injective. Since $x\notin\underline{\Lambda}_{n+k-1,k-1}$ there is either $i\in\lbrace 1,\cdots,k-1\rbrace$ such that $x_i<x_{i-1}$ or $x_0+n+k<x_{k-1}$. Then we define $x^{-1}\in\mathbb{Z}^k$ in the first case by $(x_0^0,\cdots,x^0_i-1,\cdots,x^0_{k-1})$ and in the second case by $(x_0^0-1,x_1^0,\cdots,x_{k-1}^0)$. Then we have by construction
\[
x\leq x^{-1}<x^0\qquad\text{and}\qquad x^{-1}\in\underline{\Lambda}_{n+k-1,k-1}.
\]
This contradicts the minimality of $x^0$ and proves the second claim. Hence property (P2) holds for $x^0$ which leads to the last isomorphism in
\[
x^*\circ\mathsf{j}_*(X)\cong(\pi_{\underline{\Lambda}_{n+k-1,k-1,x/}})_*\circ i_x^*(X)\cong (x^0)^*\circ i_x^*(X)\cong 0.
\]
The first isomorphism above follows from \eqref{eq:dhat}, and the second from the first claim.
The proof of the statement for $\mathsf{j}_!$ is completely dual to the above. Finally, the isomorphism $\mathsf{j}_!\cong\mathsf{j}_*\colon\D_{n,k}\to\hat{\D}_{n,k}$ follows from the observation, that both functors are inverse to the restriction of the inverse image morphism $\mathsf{j}^*$.
\end{proof}

\begin{rmk}\label{rmk:slnk-hat}
We define the cubical slice $Sl_{n,k}^{\square}=\mathbb{Z}^k_{\xi/\mathsf{s}_3^{k-1}{\xi}}$ and the derivator $sl\hat{\D}_{n,k}$ to be the full subderivator of $\D^{Sl_{n,k}^{\square}}$ spanned by those objects $X\in\D^{Sl_{n,k}^{\square}}$ that satisfy property (P2) for all $x\in Sl_{n,k}^{\square}$, which are not the image of an injective object in $Sl_{n,k}$. We can show with the same strategy as in the proof of \autoref{prop:dnk-hat}, that 
\[
(\mathsf{j}\vert_{Sl_{n,k}})_!\colon sl\D_{n,k}\toiso sl\hat{\D}_{n,k}
\]
is an equivalence. Let $sl_{n,k}^{\square}\colon Sl_{n,k}^{\square}\to\mathbb{Z}^k$ be the inclusion. The square
\[
\xymatrix{
\D_{n,k}\ar[r]^{sl_{n,k}^*}\ar[d]_{\mathsf{j}_!}&sl\D_{n,k}\ar[d]^{(\mathsf{j}\vert_{Sl_{n,k}})_!}\\
\hat{\D}_{n,k}\ar[r]_{(sl_{n,k}^{\square})^*}&sl\hat{\D}_{n,k}
}
\]
commutes up to natural isomorphism (since the inverses of the vertical morphisms are given by restrictions of inverse image morphisms). We note that the top morphism is an equivalence by \autoref{thm:slices}. Hence the entire square consists of equivalences.
\end{rmk}

\section{Vertical functoriality}
\label{sec:vertical}

The main objective of this section is the construction of canonical morphisms relating the derivators $\D_{n,k}$ for fixed $k\geq 2$. For this we will show that for a stable derivator $\D$ the inverse images associated to the postcomposition functors
\[
(\mathrm{s}_i)_*\colon\underline{\Lambda}_{n+k,k-1}\rightleftarrows\underline{\Lambda}_{n+k-1,k-1}\colon(\mathrm{d}_i)_*.
\]
restrict to morphisms between $\D_{n,k}$ and $\D_{n+1,k}$. As a consequence we obtain for $k\geq 2$ fixed a $2$-functor
\[
(-)^*\colon\underline{\Lambda}^{op}\to Der,\Lambda_n\mapsto\D_{n-k+1,k}.
\]
We describe these operations also on fundamental slices, which will be useful for later applications. In particular in the case $k=2$ we will obtain an extended version of and hence recover the standard simplicial structure on the derivators $\D^{[n]}, n\geq 0$.

\begin{con}
Let $m,k\geq 0$ and $\D$ a stable derivator. Consider the adjunction \[
\mathrm{s}_0\colon\Lambda_{m+1}\rightleftarrows\Lambda_{m}\colon\mathrm{d}_0
\]
(c.f.\autoref{rmk:para-fd}) in the 2-category $\underline{\Lambda}$ of parasimplices. We now apply the 2-functor $\underline{\Lambda}(\Lambda_k,-)\colon\underline{\Lambda}\to Cat$ to this adjunction, and hence obtain the adjunction 
\[
\underline{\Lambda}(\Lambda_k,\mathrm{s}_0)\colon\underline{\Lambda}_{m+1,k}\rightleftarrows\underline{\Lambda}_{m,k}\colon\underline{\Lambda}(\Lambda_k,\mathrm{d}_0)
\]
in $Cat$. Finally, we apply $\D\colon Cat^{op}\to Der$ to obtain an adjunction
\[
\mathsf{d}\colon\D^{\underline{\Lambda}_{m+1,k}}\rightleftarrows\D^{\underline{\Lambda}_{m,k}}\colon\mathsf{s}.
\]
\end{con}

\begin{prop}\label{prop:d0s0-v}
Let $n\geq 0$, $k\geq 2$, $m=n+k-1$ and $\D$ a stable derivator. Then the adjunction $\mathsf{d}\dashv\mathsf{s}$ restricts to an adjunction
\[
\mathsf{d}\colon\D_{n+1,k}\rightleftarrows\D_{n,k}\colon\mathsf{s}.
\]
\end{prop}

\begin{proof}
Since postcomposition functors automatically preserve non-injective objects, we deduce that objects in the image of  $\mathsf{d}$ and $\mathsf{s}$ satisfy property (P2) on all non-injective objects. It remains to verify property (P1) for all injective objects.

\noindent First we take care of the morphism $\mathsf{d}$. For this let $f=(f_0,\cdots,f_{k-1})$ be an injective object in $\underline{\Lambda}_{n+k-1,k-1}$. Let $M=\lbrace i\in\mathbf{k}\vert\exists j\in\mathbb{Z}\colon f_i+1=j\cdot(n+k)\rbrace\subseteq\mathbf{k}$ and define
\[
\tilde{\square}_f\colon\cube{k}\to\underline{\Lambda}_{n+k,k-1},(\delta_0,\cdots,\delta_{k-1})\mapsto(f_0+\mu(0)\cdot\delta_0,\cdots,f_{k-1}+\mu(k-1)\cdot\delta_{k-1}),
\]
where $\mu\colon\mathbf{k}\to\mathbb{Z}$ is defined by $\mu(i)=1$ if $i\notin M$ and $\mu(i)=2$ if $i\in M$. Then the elementary subcube $\square_f$ starting in $f$ satisfies $\underline{\Lambda}(\Lambda_{k-1},\mathrm{d}_0)\circ\square_f=\tilde{\square}_f$. For $x\in\D_{n+1,k}$ we obtain therefore $\square_f^*(\mathsf{d}^v(x))=\tilde{\square}_f^*(x)$, which is bicartesian by assumption on $x$ and \autoref{cor:bicart-concat}.

\noindent Now we consider the morphism $\mathsf{s}$. For this let $g=(g_0,\cdots,g_{k-1})$ be an injective object in $\underline{\Lambda}_{n+k,k-1}$. Let $N=\lbrace i\in\mathbf{k}\vert\exists j\in\mathbb{Z}\colon g_i=j\cdot(n+k+1)\rbrace\subseteq\mathbf{k}$ and define
\[
\tilde{\square}_g\colon\cube{k}\to\underline{\Lambda}_{n+k-1,k-1},(\delta_0,\cdots,\delta_{k-1})\mapsto(g_0+\nu(0)\cdot\delta_0,\cdots,g_{k-1}+\nu(k-1)\cdot\delta_{k-1}),
\]
where $\nu\colon\mathbf{k}\to\mathbb{Z}$ is defined by $\nu(i)=1$ if $i\notin N$ and $\nu(i)=0$ if $i\in N$. Then the elementary subcube $\square_g$ starting in $g$ satisfies $\underline{\Lambda}(\Lambda_{k-1},\mathrm{s}_0)\circ\square_g=\tilde{\square}_g$. For $x\in\D_{n,k}$ we obtain therefore $\square_g^*(\mathsf{s}^v(x))=\tilde{\square}_g^*(x)$, which is bicartesian by assumption on $x$ if $N=\emptyset$. If $N\neq\emptyset$ the cube $\tilde{\square}_g^*(x)$ is in the essential image of the inverse image associated to the canonical projection $\cube{k}\to\cube{\mathbf{k}\setminus N}$, and hence bicartesian by \autoref{prop:bicart-obstruction}.
\end{proof}

\begin{cor}\label{cor:higher-Sdot-I}
Let $m,m'\geq 0, k\geq 2$, $f\in\underline{\Lambda}(\Lambda_m,\Lambda_{m'})$, and $\D$ a stable derivator. Then the inverse image $\underline{\Lambda}(\Lambda_{k-1},f)^*\colon\D^{\underline{\Lambda}_{m,k-1}}\to\D^{\underline{\Lambda}_{m',k-1}}$ restricts to a morphism of derivators
\[
\D_{m-k+1,k}\to\D_{m'-k+1,k}.
\]
\end{cor}

\begin{proof}
The statement is clear whenever $f$ is of the form $\mathrm{d}_0$ or $\mathrm{s}_0$ by \autoref{prop:d0s0-v} or of the form $\mathsf{t}$ by \autoref{cor:dnk-sym}. For the general case we invoke \autoref{cor:Lambda-generators} to factor $f$ as a composition of morphisms of the above form and use the 2-functoriality of $\underline{\Lambda}(\Lambda_{k-1},-)$ and $\D$.
\end{proof}

\begin{cor}\label{cor:higher-Sdot-II}
Let $k\geq 2$ and $\D$ a stable derivator. Then there is 2-functor
\[
\mathsf{S}_{\bullet}^{(k-1)}(\D)\colon\underline{\Lambda}^{op}\to Der, \Lambda_m\mapsto\D_{m-k+1,k},f\mapsto\underline{\Lambda}(\Lambda_{k-1},f)^*.
\]
Moreover, $\mathsf{S}_{\bullet}^{(k-1)}(-)$ is 2-functorial with respect to morphisms of derivators.
\end{cor}

\begin{proof}
The assignment is well defined on 1-morphisms by \autoref{cor:higher-Sdot-I} and on 2-morphisms since the derivators $\D_{m-k+1,k}$ are full subderivators of $\D^{\underline{\Lambda}_{m,k-1}}$. The 2-functoriality follows from the one of $\underline{\Lambda}(\Lambda_{k-1},-)$ and $\D$. Furthermore, morphisms of prederivators commute by definition with inverse images. This implies the naturality statement.
\end{proof}

\begin{defn}
Let $k\geq 1$. The 2-functor $\mathsf{S}_{\bullet}^{(k)}(-)\colon Der^{st}\to 2-Fun(\underline{\Lambda}^{op},Der)$ is called the $k$\textbf{th higher parasimplicial} $\mathsf{S}_{\bullet}$\textbf{-construction}.
\end{defn}

\begin{rmk}
By passing to derivators of fundamental domains in the case $k=1$ we recover a parasimplicial enhancement of the standard simplicial $\mathsf{S}_{\bullet}$-construction (c.f. \cite{waldhausen:k-theory,garkusha:I}). For $k\geq 2$ a construction very similar to $\mathsf{S}_{\bullet}^{(k)}$, but in a slightly different context, was considered recently in \cite{poguntke-Sdot}. Another variant (in the context of $\infty$-categories) thereof was shown in \cite{dyckerhoff-Sdot} to appear naturally in the categorified Dold--Kan correspondence.
\end{rmk}

\begin{thm}\label{thm:gsd-Sdot}
Let $k\geq 1$ and $\D$ a stable derivator. Then the image of $\mathsf{S}_{\bullet}^{(k)}(\D)$ is contained in $Der^{st,\infty-ad}$. Moreover, there is a pseudonatural equivalence 
\[
\mathcal{S}\colon\mathsf{L}\mathsf{S}_{\bullet}^{(k)}(\D)\toiso\mathsf{R}\mathsf{S}_{\bullet}^{(k)}(\D),
\]
defined by $\mathcal{S}_{\Lambda_m}=\mathsf{s}_3^*\colon\D_{m-k,k+1}\toiso\D_{m-k,k+1}$ for $m\geq 0$.
\end{thm}

\begin{proof}
The 2-category $\Lambda$ is adjunction complete by \autoref{prop:Lambda-ad-comp}. Since 2-functors preserve adjunctions, images of adjunction complete 2-categories under 2-functors are again adjunction complete. Hence the image of $\mathsf{S}_{\bullet}^{(k)}$ is forced to be contained in the largest adjunction complete sub-2-category of $Der^{st}$, which is $Der^{st,\infty-ad}$.
For the pseudonatural equivalence, we remember \autoref{cor:gsd-Lambda}, which states that the parasimplicial translations $\mathsf{t}$ define a 2-natural isomorphism
\[
\mathbb{S}\colon\mathsf{R}\toiso\mathsf{L}\colon\underline{\Lambda}\to\underline{\Lambda}^{coop}.
\]
We again use the compatibility of 2-functors with adjunctions to deduce that the whiskering of $\mathbb{S}$ with $\mathsf{S}_{\bullet}^{(k)}(\D)$ defines a pseudonatural equivalence
\begin{equation}\label{eq:inverse-gsd}
\mathsf{R}\mathsf{S}_{\bullet}^{(k)}(\D)\toiso\mathsf{L}\mathsf{S}_{\bullet}^{(k)}(\D),
\end{equation}
which is locally defined by $\underline{\Lambda}(\Lambda_k,\mathsf{t})^*=\mathsf{s}_1^*$. Finally, we use the natural equivalence $\mathsf{s}_3^*=\Sigma^k\circ(\mathsf{s}_1^*)^{-1}$ in $\D_{m-k+2,k+1}$ (c.f \autoref{cor:dnk-sigma}) to exhibit $\mathcal{S}$ as the inverse of the pasting of \eqref{eq:inverse-gsd} with $\Omega^k\colon\id_{Der^{st,\infty-ad}}\toiso\id_{Der^{st,\infty-ad}}$.
\end{proof}

\begin{rmk}
The proof of \autoref{thm:gsd-Sdot} suggests that the pseudonatural equivalence \eqref{eq:inverse-gsd} might be the more important, or at least more natural construction. But, however, there are various reasons to prefer the equivalence $\mathcal{S}$, which are all incarnations of the fact that $\mathsf{s}_3^*$ admits more useful properties then $\mathsf{s}_1^*$.
\begin{itemize}
\item In the case $n=0$, where $\D_{n,k}\cong\D$ holds true, we can identify $\mathsf{s}_3^*=\id$.
\item In the case $k=2$ it is known that $\mathsf{s}_3\colon\D_{n,2}\to\D_{n,2}$ define Serre equivalences (\cite[Thm.~11.12]{gst:Dynkin-A}).
\item The relation involving $\mathsf{s}_3^*$ and $\Sigma$ in \autoref{cor:fracCY} is invariant under $(n,k)\mapsto(k-1,n+1)$.
\item The duality morphisms $\Psi_n$ are compatible with $\mathsf{s}_3^*$ (\autoref{thm:dual-Serre}).
\end{itemize}
And we will see even more reasons in the following chapters.
\end{rmk}


For later computations it will be useful to know how the generalized face and degeneracy morphisms (i.e. the morphisms appearing in the infinite chain of adjunctions generated by $\mathsf{d}\dashv\mathsf{s}$) interact with the restrictions to fundamental slices, which will be the content of the remainder of this chapter. For this we observe the following.

\begin{prop}\label{prop:slice-comp}
Let $n\geq 0, k\geq 2$ and $\D$ a stable derivator. Then
\[
\underline{\Lambda}(\Lambda_{k-1},\mathrm{d_i})(Sl_{n,k})\subseteq Sl_{n+1,k}\text{ for }1\leq i\leq n+k.
\]
\end{prop}

\begin{proof}
Under the assumption $1\leq i\leq n+k$, we have 
\begin{equation}\label{eq:slice-comp}
\mathrm{d}_i(j)\in\lbrace j,j+1\rbrace \text{ for }0\leq j\leq n+k-1.
\end{equation}
Moreover, for $f=(0,f_1,\cdots,f_{k-1})\in\underline{\Lambda}_{k-1,n+k-1}$ the image is of the form
\[
\underline{\Lambda}(\Lambda_{k-1},\mathrm{d_i})(f)=(\mathrm{d_i}(0),\mathrm{d_i}(f_1),\cdots,\mathrm{d_i}(f_{k-1})).
\]
In particular, for $\xi_{n,k}=(0,1,\cdots,k-1)\in\underline{\Lambda}_{k-1,n+k-1}$ it is now a consequence of \eqref{eq:slice-comp} that
\begin{equation}\label{eq:slice-initial}
\xi_{n+1,k}=(0,1,\cdots,k-1)\leq\underline{\Lambda}(\Lambda_{k-1},\mathrm{d_i})(\xi_{n,k})\in\underline{\Lambda}_{k-1,n+k}.
\end{equation}
By using additionally that $\mathrm{d_i}(0)=0$ for $1\leq i\leq n+k$ we conclude
\begin{equation}\label{eq:slice-final}
\underline{\Lambda}(\Lambda_{k-1},\mathrm{d_i})(\mathsf{s}_3^{k-1}\xi_{n,k})\leq(0,n+2,\cdots,n+k)=\mathsf{s}_3^{k-1}\xi_{n+1,k}\in\underline{\Lambda}_{k-1,n+k}.
\end{equation}
The inequalities \eqref{eq:slice-initial} and \eqref{eq:slice-final} together yield the first statement. For the second statement, we use that
\begin{enumerate}
\item $\mathrm{s}_i(j)\in\lbrace j-1,j\rbrace$ for $0\leq j\leq n+k-1$,
\item $\mathrm{s}_i(0)=0$
\end{enumerate}
holds for $1\leq i\leq n+k-1$ and conclude with a very similar strategy.
\end{proof}

Unfortunately, the analogue of \autoref{prop:slice-comp} fails in many cases if one replaces the face maps $\mathrm{d}_i$ by degeneracy maps $\mathrm{s}_i$, since for $0\leq i\leq k-2$ we have $\underline{\Lambda}(\Lambda_{k-1},\mathrm{s_i})(\xi_{n+1,k})\leq\xi_{n,k}$. Therefore, we have to consider a slightly larger version of the slice.

\begin{defn}
Let $n\geq 0,k\geq 2$ and $\D$ a stable derivator.
\begin{enumerate}
\item The triangular slice $Sl^{\triangle}_{n,k}$ is the full subcategory 
\[
(\underline{\Lambda}_{n+k-1,k-1})_{(0,\cdots,0)/(0,n+k-1,\cdots,n+k-1)}\subseteq\underline{\Lambda}_{n+k-1,k-1}
\]
with inclusion $sl^{\triangle}_{n,k}\colon Sl^{\triangle}_{n,k}\to\underline{\Lambda}_{n+k-1,k-1}$.
\item The derivator of triangular slices $sl^{\triangle}\D_{n,k}$ is the full subderivator of $\D^{Sl^{\triangle}_{n,k}}$ spanned by those objects $X\in\D^{Sl^{\triangle}_{n,k}}$ that satisfy property (P2) for all $x\in Sl^{\triangle}_{n,k}$ such that $sl^{\triangle}_{n,k}(x)$ is non-injective.
\end{enumerate}
\end{defn}

\begin{con}\label{con:triangle}
Let $n\geq 0,k\geq 2$ and $\D$ a stable derivator. The $\triangle^{n,k}\colon Sl_{n,k}\to Sl^{\triangle}_{n,k}$ inclusion admits a factorization
\[
Sl_{n,k}\xrightarrow{a}(Sl^{\triangle}_{n,k})_{\xi/}\xrightarrow{b}Sl^{\triangle}_{n,k}.
\]
We observe $a$ is a sieve and $b$ is a cosieve. Since the complement of $\triangle^{n,k}$ consists of non-injective objects, it follows from \cite[Prop.~1.23]{groth:ptstab} that $b_!\circ a_*\colon sl\D_{n,k}\to sl^{\triangle}\D_{n,k}$ is an equivalence inverse to $(\triangle^{n,k})^*$. Together with \autoref{thm:slices} this implies that also $(sl^{\triangle_{n,k}})^*\colon\D_{n,k}\to sl^{\triangle}\D_{n,k}$ is an equivalence.

\noindent On the other hand, let $\Delta(\Delta_{k-1},\Delta_{n+k-1})_0\subseteq\Delta(\Delta_{k-1},\Delta_{n+k-1})$ be the full subcategory on those objects $g\colon\Delta_{k-1}\to\Delta_{n+k-1}$ with $g(0)$. Then the functor
\begin{equation}\label{eq:triang}
\Delta(\Delta_{k-1},\Delta_{n+k-1})_0\to Sl_{n,k}^{\triangle},g\mapsto(g(0),\cdots,g(k-1))
\end{equation}
is an isomorphism. This yields for every morphism $f\colon\Lambda_{n+k-1}\to\Lambda_{n'+k-1}$ with 
\begin{itemize}
\item $f(0)=0$,
\item $f(n+k-1)\leq n'+k-1$
\end{itemize}
a commutative diagram
\[
\xymatrix{
Sl^{\triangle}_{n,k}\ar[r]_{sl^{\triangle}_{n,k}}\ar[d]_{(\underline{\Lambda}(\Lambda_{k-1},f))\vert_{Sl_{n,k}^{\triangle}}}&\underline{\Lambda}_{k-1,n+k-1}\ar[d]^{\underline{\Lambda}(\Lambda_{k-1},f)}\\
Sl^{\triangle}_{n',k}\ar[r]_{sl^{\triangle}_{n',k}}&\underline{\Lambda}_{k-1,n'+k-1},
}
\]
since the left vertical morphism is (up to composition with \eqref{eq:triang}) $\Delta(\Delta_{k-1},f)_0$, which is well-defined by the assumptions on $f$.
\end{con}

\begin{cor}\label{cor:slice-restr}
Let $n\geq 0$, $k\geq 2$ and $\D$ a stable derivator. Then
\begin{enumerate}
\item for $1\leq i\leq n+k$ there is a strictly commutative diagram
\[
\xymatrix{
\D_{n+1,k}\ar[r]^{sl_{n+1,k}^*}\ar[d]_{\underline{\Lambda}(\Lambda_{k-1},\mathrm{d}_i)^*}& sl\D_{n+1,k}\ar[d]^{(\underline{\Lambda}(\Lambda_{k-1},\mathrm{d}_i)\vert_{Sl_{n,k}})^*}\\
\D_{n,k}\ar[r]_{sl_{n,k}^*}& sl\D_{n,k}
}
\]
\item for $0\leq i\leq n+k-1$ there is a strictly commutative diagram
\[
\xymatrix{
\D_{n,k}\ar[r]^{(sl_{n+1,k}^{\triangle})^*}\ar[d]_{\underline{\Lambda}(\Lambda_{k-1},\mathrm{s}_i)^*}& sl^{\triangle}\D_{n,k}\ar[d]^{(\underline{\Lambda}(\Lambda_{k-1},\mathrm{s}_i)\vert_{Sl^{\triangle}_{n+1,k}})^*}\\
\D_{n+1,k}\ar[r]_{(sl_{n,k}^{\triangle})^*}& sl^{\triangle}\D_{n+1,k}
}
\]
\end{enumerate}
\end{cor}

\begin{proof}
\autoref{prop:slice-comp} and \autoref{con:triangle} yield that
\begin{enumerate}
\item $\underline{\Lambda}(\Lambda_{k-1},\mathrm{d}_i)\vert_{Sl_{n,k}}\colon Sl_{n,k}\to Sl_{n+1,k}$ and
\item $\underline{\Lambda}(\Lambda_{k-1},\mathrm{s}_i)\vert_{Sl^{\triangle}_{n+1,k}}\colon Sl^{\triangle}_{n+1,k}\to Sl^{\triangle}_{n,k}$
\end{enumerate}
are well defined under the respective assumptions on $i$. As a consequence, the 2-functoriality of $\D$ implies the statements immediately.
\end{proof}

From \autoref{cor:slice-restr} and the compatibility of 2-functors with adjunctions we deduce that there is chain of adjunctions
\[
(\underline{\Lambda}(\Lambda_{k-1},\mathrm{d}_1))^*\dashv(\underline{\Lambda}(\Lambda_{k-1},\mathrm{s}_1))^*\dashv\cdots
\]
\[
\cdots\dashv(\underline{\Lambda}(\Lambda_{k-1},\mathrm{s}_{n+k-1}))^*\dashv(\underline{\Lambda}(\Lambda_{k-1},\mathrm{d}_{n+k}))^*
\]
relating $\D_{n,k}$ and $\D_{n+1,k}$. Moreover, \autoref{thm:gsd-Sdot} implies that that his chain extends to an infinite chain of adjunctions. In the following we show that the left adjoint of $(\underline{\Lambda}(\Lambda_{k-1},\mathrm{d}_{1}))^*$ and the right adjoint of $(\underline{\Lambda}(\Lambda_{k-1},\mathrm{d}_{n+k}))^*$ also admit a simple description. For simplicity we use in the following the notation
\[
d^v:=\underline{\Lambda}(\Lambda_{k-1},\mathrm{d}_{1})\vert_{Sl_{n,k}}\colon Sl_{n,k}\to Sl_{n+1,k},(0,f_1,\cdot,f_{k-1})\mapsto(0,f_1+1,\cdots,f_{k-1}+1),
\]
\[
d^{v\vee}:=\underline{\Lambda}(\Lambda_{k-1},\mathrm{d}_{n+k})\vert_{Sl_{n,k}}\colon Sl_{n,k}\to Sl_{n+1,k}(0,f_1,\cdot,f_{k-1})\mapsto(0,f_1,\cdots,f_{k-1}).
\]

\begin{prop}\label{prop:the-sieve}
Let $n\geq 1,k\geq 2$ and $\D$ a stable derivator. Then the adjunctions
\[
d^v_!\colon\D^{Sl_{n,k}}\rightleftarrows\D^{Sl_{n+1,k}}\colon(d^v)^*\qquad\text{and}\qquad(d^{v\vee})^*\colon\D^{Sl_{n+1,k}}\leftrightarrows\D^{Sl_{n,k}}\colon d^{v\vee}_*
\]
restrict to an adjunctions
\[
d^v_!\colon sl\D_{n,k}\rightleftarrows sl\D_{n+1,k}\colon(d^v)^*\qquad\text{and}\qquad(d^{v\vee})^*\colon sl\D_{n+1,k}\leftrightarrows sl\D_{n,k}\colon d^{v\vee}_*.
\]
\end{prop}

\begin{proof}
The statements are completely dual to each other. We show the statement for $d^{v\vee}$.
We have to show that the image of $(d^{v\vee}_*)\vert_{sl\D_{n,k}}$ is contained in $sl\D_{n+1,k}$. Since $d$ is fully faithful, the counit $(d^{v\vee})^*d^{v\vee}_*\toiso\id$ is invertible by \cite[Prop.~1.26.]{groth:ptstab}. This shows property (P2) for non-injective objects of $Sl_{n+1,k}$, which are in the image of $d^{v\vee}$. Moreover, we observe that 
\[
d^{v\vee}\colon(0,f_1,\cdots,f_{k-1})\mapsto(0,f_1,\cdots,f_{k-1})
\]
is a sieve. Hence the corresponding right Kan extension morphism is an extension-by-zero morphism \cite[Prop.~3.6.]{groth:ptstab}, which shows the property (P2) for non-injective objects of $Sl_{n+1,k}$ which are not in the image of $d^{v\vee}$.
\end{proof}

In \S\ref{sec:recollements} it will be important to have a systematic notation for all generalized face and degeneracy morphisms.

\begin{notn}
Let $F\colon X\to Y$ be a morphism in a 2-category, such that all iterated adjoints of $F$ exist. Then we denote the $n$th iterated right adjoint by $F[n]$ and the $n$th iterated left adjoint by $F[-n]$. Occasionally we use the convention $F[0]=F$.
\end{notn}

\begin{defn}
Let $n\geq 0, k\geq 2$ and $\D$ a stable derivator. Then
\[
\mathsf{d}^v:=\mathsf{d}^v_{n,k}:=\underline{\Lambda}(\Lambda_{k-1},\mathrm{d}_{k})^*\colon\D_{n+1,k}\to\D_{n,k}
\]
is called the \textbf{standard vertical face morphism}.
\end{defn}

\begin{eg}\label{eg:ver-fd}
Let $n\geq 0$, $k\geq 2$ and $\D$ a stable derivator. With this notation we have
\begin{enumerate}
\item $\underline{\Lambda}(\Lambda_{k-1},\mathrm{d}_i)^*\cong\mathsf{d}^v[2(i-k)]$ for  $0\leq i \leq n+k$,
\item $\underline{\Lambda}(\Lambda_{k-1},\mathrm{s}_i)^*\cong\mathsf{d}^v[2(i-k)+1]$ for  $0\leq i \leq n+k-1$.
\end{enumerate}
\end{eg}

Even more explicitely, we have the following.

\begin{eg}\label{eg:explicit}
We consider the case $k=2$.
\begin{enumerate}
\item Recall from \autoref{egs:simplex-slice} that $[n]\to Sl_{n,2},i\mapsto(0,i+1)$ is an isomorphism of categories. Since all objects in $Sl_{n,2}$ are injective, these isomorphisms lead to equivalences of derivators $\D^{[n]}\cong sl\D_{n,2}$. By plugging in the definition we obtain therefore the commutativity of the diagram
\begin{equation}\label{eq:fd-comparison}
\xymatrix{
\D^{[n+1]}\ar[r]^{\sim}\ar[d]_{\mathrm{d}_i^*}&\D^{\underline{\Lambda}(\Lambda_1,\Lambda_{n+2})}\ar[r]^{=}\ar[d]^{\underline{\Lambda}(\Lambda_1,\mathrm{d}_{i+1})^*}&sl\D_{n+1,2}\ar[d]^{\mathsf{d}^v[2i-2]}\\
\D_{[n]}\ar[r]_{\sim}&\D^{\underline{\Lambda}(\Lambda_1,\Lambda_{n+1})}\ar[r]_{=}&sl\D_{n,2}.
}
\end{equation}
In particular we observe that, due to the index shift between the left and middle vertical morphism, we have one face and degeneracy morphism (those defined by $\underline{\Lambda}(\Lambda_1,\mathrm{d}_{0})^*$ and $\underline{\Lambda}(\Lambda_1,\mathrm{s}_{0})^*$) more then a priori expected. These extra morphisms also satisfy the simplicial relations.
\item If we specialize even further to the case $n=0$, the commutativity of \eqref{eq:fd-comparison} yields the coincidence of the following sequences of adjoint morphisms, which correspond to the vertical morphisms in \eqref{eq:fd-comparison}
\[
C\dashv(\mathrm{d}_0)_!\dashv \mathrm{d}_0^*\dashv\mathrm{s}_0^* \dashv\mathrm{d}_1^*\dashv(\mathrm{d}_1)_*\dashv F,
\]
\[
\underline{\Lambda}(\Lambda_1,\mathrm{d}_{0})^*\dashv\underline{\Lambda}(\Lambda_1,\mathrm{s}_{0})^*\dashv \underline{\Lambda}(\Lambda_1,\mathrm{d}_{1})^*\dashv\underline{\Lambda}(\Lambda_1,\mathrm{s}_{1})^*
\]
\[\dashv \underline{\Lambda}(\Lambda_1,\mathrm{d}_{2})^*\dashv\underline{\Lambda}(\Lambda_1,\mathrm{s}_{2})^*\dashv \underline{\Lambda}(\Lambda_1,\mathrm{d}_{3})^*,
\]
\[
\mathsf{d}^v[-4]\dashv\mathsf{d}^v[-3]\dashv\mathsf{d}^v[-2]\dashv\mathsf{d}^v[-1]\dashv\mathsf{d}^v\dashv\mathsf{d}^v[1]\dashv\mathsf{d}^v[2].
\]
We note, that in the second chain of adjunctions the parasimplicial maps $\mathrm{s}_{2}\colon\Lambda_2\leftrightarrows\Lambda_1\colon\mathrm{d}_{3}$ are not in the image of the embedding $\mathsf{i}$ of the 2-category of simplices $\Delta$.
\end{enumerate}
\end{eg}

\begin{cor}\label{cor:gsdver}
Let $n\geq 1,k\geq 2$ and $\D$ a stable derivator. Then the standard vertical face morphism $\mathsf{d}^v\colon\D_{n+1,k}\to\D_{n,k}$ generates an infinite chain of adjunctions such that for $p\in\mathbb{Z}$
\begin{enumerate}
\item $\mathsf{d}^v[2p]=(\mathsf{s}_3^*)^p\circ\mathsf{d}^v\circ(\mathsf{s}_3^*)^{-p}$,
\item $\mathsf{d}^v[2p+1]=(\mathsf{s}_3^*)^p\circ\mathsf{d}^v[1]\circ(\mathsf{s}_3^*)^{-p}=(\mathsf{s}_3^*)^{p+1}\circ\mathsf{d}^v[-1]\circ(\mathsf{s}_3^*)^{-p-1}$.
\end{enumerate}
\end{cor}

\begin{proof}
The isomorphisms
\[
\mathsf{d}^v[-1]=(\mathsf{s}_3^*)^{-1}\circ\mathsf{d}^v[1]\circ\mathsf{s}_3^*\qquad\text{and}\qquad\mathsf{d}^v[1]=\mathsf{s}_3^*\circ\mathsf{d}^v[-1]\circ(\mathsf{s}_3^*)^{-1},
\]
are immediate consequences of \autoref{thm:gsd-Sdot}. The general statement follows via induction.
\end{proof}

\section{Horizontal functoriality}
\label{sec:horizontal}

In this section we construct canonical morphisms relating the derivators $\D_{n,k}$ for fixed $n\geq 1$. For this we define, in a first step, morphisms of derivators $sl\D_{n,k}\to sl\D_{n,k'}$ and invoke \autoref{thm:slices} to transfer them to the desired morphisms $\D_{n,k}\to\D_{n,k'}$. This has the advantage, that we do not have to worry about cartesianess conditions for subcubes. More precisely, we first define a morphism $\mathsf{d}^h\colon\D_{n,k+1}\to\D_{n,k}$, which will be exhibited as the horizontal analogue of the standard vertical face morphisms in \S\ref{sec:main}. In the next step we will construct the left and the right adjoint of $\mathsf{d}^h$. We show that the resulting adjoint triples are periodic with respect to the autoequivalences $\mathsf{s}_3^*$, and therefore extend to infinite chains of adjunctions.

\begin{defn}\label{defn:hormor}
Let $n\geq 1,k\geq 2$ and $\D$ a stable derivator.
\begin{enumerate}
\item The \textbf{standard horizontal face map} is the inclusion of poets 
\[
d^h\colon Sl_{n,k}\to Sl_{n,k+1}, (0,f_1,\cdots,f_{k-1})\mapsto(0,1,f_1+1,\cdots,f_{k-1}+1).
\]
\item The \textbf{standard horizontal face morphism} $\mathsf{d}^h:=\mathsf{d}^h_{n,k}$ is the restriction of the inverse image of $d^h$
\[
\mathsf{d}^h=\mathsf{d}^h_{n,k}\colon sl\D_{n,k+1}\to sl\D_{n,k}.
\]
\end{enumerate}
\end{defn}

\begin{rmk}\label{rmk:dh-welldef}
The inclusion $d^h$ obviously preserves non-injective objects. This implies directly that the restriction of the inverse image of $d^h$ above is well defined.
\end{rmk}

\begin{prop}\label{prop:dh-ext-zero}
Let $n\geq 1,k\geq 2$ and $\D$ a stable derivator. Then the adjunction $(d^h)^*\dashv d^h_*\colon\D^{Sl_{n,k+1}}\leftrightarrows\D^{Sl_{n,k}}$ restricts to an adjunction
\[
\mathsf{d}^h\dashv\mathsf{d}^h[1]\colon sl\D_{n,k+1}\leftrightarrows sl\D_{n,k}.
\]
\end{prop}

\begin{proof}
We have to show that the image of $(d^h_*)\vert_{sl\D_{n,k}}$ is contained in $sl\D_{n,k+1}$. Since $d^h$ is a sieve, the corresponding right Kan extension morphism is an extension by zero, which shows the property (P2) for non-injective objects of $Sl_{n,k+1}$ which are not in the image of $d^h$. On the other hand, for non-injective objects in the image of $d^v$ the property (P2) follows from the invertibility of the counit $(d^h)^*d^h_*\toiso\id$.
\end{proof}

Unfortunately, it turns out that the left Kan extension $d^h_!\colon\D^{Sl_{n,k}}\to\D^{Sl_{n,k+1}}$ does not restrict to a morphism of the form $sl\D_{n,k+1}\to sl\D_{n,k}$. Because of this, the description of $\mathsf{d}^h[-1]$ will be slightly more involved. More precisely, we consider the subposet 
\[
B_{n,k}=Sl_{n,k+1}\setminus\lbrace (0,g_1,\cdots,g_k)\in Sl_{n,k+1}\vert (0,g_1,\cdots,g_k) \text{ injective with } g_1\geq 1\rbrace
\]
of $Sl_{n,k+1}$. Since the image of $d^h$ is contained in $B_{n,k}$, we obtain the following factorization of $d^h$
\[
Sl_{n,k}\xrightarrow{i} B_{n,k}\xrightarrow{j} Sl_{n,k+1}.
\]
But now we are in a situation where standard techniques from the theory of pointed derivators apply. More explicitly, we use the following result.

\begin{lem}\label{lem:adjoint-lem}
Let $\D$ be a pointed derivator and $u\colon A\to B$ be a functor such that there is a factorization $u=w\circ v$
\[
A\xrightarrow{v}C\xrightarrow{w}B
\]
with $v$ a sieve and $w$ fully faithful. Then the restriction of the inverse image
\[
u^*\colon\D^{B,w(C\setminus v(A))}\to\D^A
\]
is a right adjoint and the left adjoint is given by $w_!\circ v_*\colon\D^A\to\D^{B,w(C\setminus v(A))}$.
\end{lem}

\begin{proof}
Since $v$ is a sieve, \cite[Prop.~3.6]{groth:ptstab} implies that the adjunction $v^*\dashv v_*$ restricts to an equivalence of derivators
\begin{equation}\label{eq:adjoint-lem1}
v_*\colon\D^A\leftrightarrows\D^{C,C\setminus v(A)}\colon v^*.
\end{equation}
Moreover, since $w$ is fully faithful, by \cite[Prop.~1.20]{groth:ptstab} the same is true for $w_!$ and the unit of the adjunction $w_!\dashv w^*$ is an isomorphism.
As a consequence $w_!\dashv w^*$ restricts to an adjunction
\begin{equation}\label{eq:adjoint-lem2}
w_!\colon\D^{C,C\setminus v(A)}\leftrightarrows\D^{B,w(C\setminus v(A))}\colon w^*.
\end{equation}
We conclude by composing \eqref{eq:adjoint-lem1} and \eqref{eq:adjoint-lem2}.
\end{proof}

\begin{prop}
Let $n\geq 1,k\geq 2$ and $\D$ a stable derivator. Then the composition $j_!\circ i_*\colon\D^{Sl_{n,k}}\to\D^{Sl_{n,k+1}}$ restricts to a morphism $\mathsf{d}^h[-1]\colon sl\D_{n,k}\to sl\D_{n,k+1}$, which is left adjoint to $\mathsf{d}^h\colon sl\D_{n,k+1}\to sl\D_{n,k}$.
\end{prop}

\begin{proof}
The inclusion $i$ and $j$ are both fully faithful. Hence by \cite[Prop.~1.20]{groth:ptstab} the Kan extensions $i_*$ and $j_!$ are also fully faithful and the counit of $i^*\dashv i_*$ and the unit of $j_!\dashv j^*$ are isomorphisms. Hence, for $X\in sl\D_{n,k}$, the condition (P2) for non-injective objects implies the condition (P2) for $j_!\circ i_*(X)$ on all non-injective objects in $d^h(Sl_{n,k})$. For the condition (P2) for remaining non-injective objects we note that $i$ is a sieve and invoke \autoref{lem:adjoint-lem} which also yields the statement about the adjunction. 
\end{proof}

\begin{cor}\label{prop:preparations-h}
Let $n\geq 1, k\geq2$ and $\D$ a stable derivator. The adjoint triple 
\[
\xymatrix{\mathsf{d}^h[-1]\dashv\mathsf{d}^h\dashv\mathsf{d}^h[1]\colon sl\D_{n,k+1}\ar@{<-}@<1.5mm>[r] \ar@{<-}@<-1.5mm>[r] \ar[r] & sl\D_{n,k}}
\]
induces adjoint triples
\[
\xymatrix{\mathsf{d}^h[-1]\dashv\mathsf{d}^h\dashv\mathsf{d}^h[1]\colon do\D_{n,k+1}\ar@{<-}@<1.5mm>[r] \ar@{<-}@<-1.5mm>[r] \ar[r] & do\D_{n,k}}
\]
and
\[
\xymatrix{\mathsf{d}^h[-1]\dashv\mathsf{d}^h\dashv\mathsf{d}^h[1]\colon\D_{n,k+1}\ar@{<-}@<1.5mm>[r] \ar@{<-}@<-1.5mm>[r] \ar[r] & \D_{n,k}}.
\]
Moreover, in all three cases the units of the adjunctions $\mathsf{d}^h[-1]\dashv\mathsf{d}^h$ and the counits of the adjunctions $\mathsf{d}^h\dashv\mathsf{d}^h[1]$ are isomorphisms.
\end{cor}

\begin{proof}
The existence of the adjoint triples follows directly from \autoref{thm:slices}. The statement on the units and counits follows from the fully faithfulness of $d^h$ and \cite[Prop.~1.20]{groth:ptstab}. 
\end{proof}

\begin{rmk}\label{rmk:factorization}
Let $n\geq 0,k\geq 2$ and $\D$ a stable derivator. There are mutually inverse isomorphisms of categories
\[
\phi_{n,k}\colon Do_{n,k}\rightleftarrows Sl_{n,k+1}\colon\psi_{n,k}
\]
defined by $\phi_{n,k}(f_0,\cdots,f_{k-1})=(0,f_0+1,\cdots,f_{k-1}+1)$ and $\psi_{n,k}(0,g_1,\cdots,g_k)=(g_1-1,\cdots,g_k-1)$. Moreover, the inverse image of $\psi_{n,k}$ restricts to an embedding $do\D_{n,k}\to sl\D_{n,k+1}$. We denote by $sl\D_{n,k+1}^{\simeq}\xrightarrow{i_{n,k}}sl\D_{n,k+1}$ the inclusion of the essential image. In particular, the inverse images associated to $\psi_{n,k}$ and $\phi_{n,k}$ restrict to mutualls inverse morphisms
\[
(\psi_{n,k}^*)^{\simeq}\colon do\D_{n,k}\rightleftarrows sl\D_{n,k+1}^{\simeq}\colon(\phi_{n,k}^*)^{\simeq},
\]
where we use use the superscripts $\simeq$ to distinguish the above functors from the unrestricted inverse images.
\end{rmk}

\begin{prop}\label{prop:dh-domain}
Let $n\geq 0, k\geq 2$, $\D$ a stable derivator and $x\in\D_{n,k}$. Then there is an isomorphism $\mathsf{d}^h[-1]\cong\psi_{n,k}^*\circ(sd^*)^{-1}$. In particular, the essential image of $\mathsf{d}^h[-1]\colon sl\D_{n,k}\to sl\D_{n,k+1}$  is contained in $sl\D_{n,k+1}^{\simeq}$.
\end{prop}

\begin{proof}
We consider the diagram
\[
\xymatrix{
\D^{Sl_{n,k}}\ar[r]^{(u_1)_!}\ar@{=}[d]&\D^{A_1}\ar[r]^{(u_2)_*}\ar[d]_{i_1^*}&\D^{A_2}\ar[r]^{(u_3)_*}\ar[d]_{i_2^*}&\D^{A_3}\ar[r]^{(u_4)_!}\ar[d]_{i_3^*}&\D^{\underline{\Lambda}_{n+k-1,k-1}}\ar[d]_{do_{n,k}^*}\\
\D^{Sl_{n,k}}\ar[r]^{i'_!}\ar@{=}[d]&\D^{Do_{n,k}\cap A_1}\ar[r]^{j'_*}\ar[d]_{(\psi_{n,k}\vert_{Do_{n,k}\cap A_1})^*}&\D^{Do_{n,k}}\ar@{=}[r]\ar[d]_{\psi_{n,k}^*}&\D^{Do_{n,k}}\ar@{=}[r]&\D^{Do_{n,k}}\\
\D^{Sl_{n,k}}\ar[r]^{i^*}&\D^{B_{n,k}}\ar[r]^{j^*}&\D^{Sl_{n,k+1}},
}
\]
where $i_1,i_2,i_3,i'$ and $j'$ are the respective obvious inclusions. We claim that the diagram above is commutative.
\begin{itemize}
\item The first square in the top row commutes because $u_1=i_1\circ i'$, the fully faithfulness of $i_1$ and \cite[Prop.~1.20]{groth:ptstab}.
\item The second square in the top row commutes because
\[
\xymatrix{
Do_{n,k}\cap A_1\ar[r]^{j'}\ar[d]_{i_1}&Do_{n,k}\ar[d]^{i_2}\\
A_1\ar[r]_{u_2}&A_2
}
\]
is a strict pull-back, $i_2$ is a sieve and \cite[Prop.~1.24]{groth:ptstab}.
\item The third square in the top row commutes because $i_3=u_3\circ i_2$, the fully faithfullness of $u_3$ and \cite[Prop.~1.20]{groth:ptstab}.
\item The third square in the top row commutes because $do_{n,k}=u_4\circ i_3$, the fully faithfullness of $u_4$ and \cite[Prop.~1.20]{groth:ptstab}.
\item The squares in the bottom row commute because they are mates of inverse image squares where the vertical maps are isomorphisms.
\end{itemize}
Hence $\mathsf{d}^h[-1]\cong\psi_{n,k}^*\circ do_{n,k}^*\circ(sl_{n,k}^*)^{-1}$.
\end{proof}

\begin{rmk}
As a consequence, we note, that an object in $sl\D_{n,k+1}$ is contained in $sl\D_{n,k+1}^{\simeq}$ if it satisfies property (P1) for all $g=(0,g_1,\cdots,g_k)\in Sl_{n,k+1}$ with $g_k\leq n+k-2$.
\end{rmk}

\begin{thm}\label{thm:Serre-h}
Let $n\geq 1,k\geq 2$ and $\D$ a stable derivator. Then there is a natural isomorphism
\[
\mathsf{s}_3^*\circ\mathsf{d}^h[-1]\cong\mathsf{d}^h[1]\circ\mathsf{s}_3^*.
\]
\end{thm}

\begin{proof}
We show the equivalent statement, that there is an isomorphism
\begin{equation}\label{eq:hor-Serre}
(\mathsf{s}_3^*)^{-1}\circ\mathsf{d}^h[1]\cong\mathsf{d}^h[-1]\circ(\mathsf{s}_3^*)^{-1}.
\end{equation}
\autoref{prop:dh-domain} implies that the essential image of $\mathsf{d}^h[-1]\circ(\mathsf{s}_3^*)^{-1}$ is contained in $sl\D_{n,k+1}^{\simeq}$. In the following we show that also the essential image of $(\mathsf{s}_3^*)^{-1}\circ\mathsf{d}^h[1]$ is contained in $sl\D_{n,k+1}^{\simeq}$.
Let $sd_{n,k}\colon Sl_{n,k}\to Do_{n,k}$ be the inclusion and $\tilde{sd}_{n,k}\colon Sl_{n,k}\to Do_{n,k}$ be defined by $f\mapsto\mathsf{s}_3(f)$. Recall, that \autoref{thm:slices} implies that $sd_{n,k}^*\colon do\D_{n,k}\to sl\D_{n,k}$ is invertible.
Let $x\in sl\D_{n,k}$ and consider 
\[
y=(\tilde{sd}_{n,k+1}^*)^{-1}\circ\mathsf{d}^h[1](x)\in do\D_{n,k+1}.
\]
Then for $g=(g_0,\cdots,g_k)\in Do_{n,k+1}$ we have by \autoref{prop:dh-ext-zero} that $g^*y=0$ if $g_k=n+k$ and $g_0\geq 1$. Let now $h=(h_0,\cdots,h_k)\in Do_{n,k+1}$ injective with $h_0\geq 1$ and $h_k\leq n+k-1$. Then we define $\tilde{\square}_h\colon\cube{k+1}\to Do_{n,k+1}$ by 
\scalebox{.9}{\parbox{.95\hsize}{%
\begin{equation}
(\delta_0,\cdots,\delta_k)\mapsto
\begin{cases}
(h_0+\delta_0(h_1-h_0),\cdots,h_{k-1}+\delta_{k-1}(h_k-h_{k-1}),h_k) \text{ if } \delta_k=0\\
(h_0+\delta_0(h_1-h_0),\cdots,h_{k-1}+\delta_{k-1}(h_k-h_{k-1}),n+k) \text{ if } \delta_k=1.
\end{cases}
\end{equation}
}}

By construction, $\tilde{\square}_h$ is well-defined and a concatination of elementary subcubes. Therefore we can conclude by \autoref{cor:bicart-concat} that $\tilde{\square}_h^*(y)$ is bicartesian. We observe that $\iota_{1,k+1}^*\circ\tilde{\square}_h^*(y)\cong 0$, where we use non-injectivity for $\delta_k=0$ and our assumption on $y$ for $\delta_k=1$. Hence, we obtain 
\begin{equation}\label{eq:vanishing}
h^*(y)=\emptyset^*\circ\tilde{\square}_h^*(y)\cong 0.
\end{equation}
In the next step we consider $sd_{n,k+1}^*(y)\in sl\D_{n,k+1}$. Let now $e=(0,e_1,\cdots,e_{k})\in Sl_{n,k+1}$ injective with $e_{k}\leq n+k-1$. Then we consider the elementary subcube $\tilde{\square}_e\colon\cube{k+1}\to Do_{n,k+1}$ by 
\[
(\delta_0,\cdots,\delta_k)\mapsto(\delta_0,e_1+\delta_1,\cdots,e_k+\delta_k).
\]
Hence $\tilde{\square}_e(y)$ is bicartesian. Now \eqref{eq:vanishing} implies $(\mathrm{d}_0^0)^*\circ\tilde{\square}_e^*(y)\cong 0$. On the other hand
\[
(\mathrm{d}_0^0)^*\circ\tilde{\square}_e^*(y)=\square_e^*\circ sd_{n,k+1}^*y,
\]
where $\square_e$ is the elementary subcube of $Sl_{n,k+1}$ starting in $e$. We invoke \cite[Thm.~8.11]{gst:basic} to see that $\square_e^*\circ sd_{n,k+1}^*(y)$ is bicartesian, and therefore $sd_{n,k+1}^*(y)$ is in fact an object in $sl\D_{n,k+1}^{\simeq}$, and \autoref{rmk:slices} to identify $sd_{n,k+1}^*(y)\cong(\mathsf{s}_3^*)^{-1}\circ\mathsf{d}^h[1](x)$.

\noindent We conclude by showing that both sides of \eqref{eq:hor-Serre} are, when considered as morphisms $sl\D_{n,k}\to sl\D_{n,k+1}^{\simeq}$, inverse to the equivalence $\tilde{sd}_{n,k}^*\circ(\phi_{n,k}^*)^{\simeq}\colon sl\D_{n,k+1}^{\simeq}\to sl\D_{n,k}$.
For the left hand side we consider the composition of isomorphisms
\begin{align}
 &\tilde{sd}_{n,k}^*\circ(\phi_{n,k}^*)^{\simeq}\circ(\mathsf{s}_3^*)^{-1}\circ\mathsf{d}^h[1]\\
\toiso & \tilde{sd}_{n,k}^*\circ(\phi_{n,k}^*)^{\simeq}\circ sd_{n,k+1}^*\circ(\tilde{sd}_{n,k+1}^*)^{-1}\circ\mathsf{d}^h[1]\\
= & sd_{n,k}^*\circ(\phi_{n,k}^*)^{\simeq}\circ\tilde{sd}_{n,k+1}^*\circ(\tilde{sd}_{n,k+1}^*)^{-1}\circ\mathsf{d}^h[1] \label{eq:gsd-hor-I}\\
\toiso & sd_{n,k}^*\circ(\phi_{n,k}^*)^{\simeq}\circ\mathsf{d}^h[1]\\
= & \mathsf{d}^h\circ\mathsf{d}^h[1]\\
\toiso &\id_{sl\D_{n,k}},
\end{align}
where
\begin{enumerate}
\item the first step is induced by \autoref{rmk:slices},
\item the second step follows from $sd_{n,k+1}\circ\phi_{n,k}\circ\tilde{sd}_{n,k}=\tilde{sd}_{n,k+1}\circ\phi_{n,k}\circ sd_{n,k}$,
\item the third step is induced by the equivalence $\tilde{sd}_{n,k+1}$ \autoref{rmk:slices},
\item the fourth step follows from $\phi_{n,k}\circ sd_{n,k}=d^h$,
\item the fifth step is the invertibility of the counit of the adjunction $\mathsf{d}^h\dashv\mathsf{d}^h[1]$.
\end{enumerate}
Finally, for the right hand side the composition of isomorphisms
\begin{align}
&\tilde{sd}_{n,k}^*\circ(\phi_{n,k}^*)^{\simeq}\circ\mathsf{d}^h[-1]\circ(\mathsf{s}_3^*)^{-1}\\
\toiso & \tilde{sd}_{n,k}^*\circ(sd_{n,k}^*)^{-1}\circ sd_{n,k}^*\circ(\phi_{n,k}^*)^{\simeq}\circ\mathsf{d}^h[-1]\circ(\mathsf{s}_3^*)^{-1}\\
= & \tilde{sd}_{n,k}^*\circ(sd_{n,k}^*)^{-1}\circ\mathsf{d}^h\circ\mathsf{d}^h[-1]\circ(\mathsf{s}_3^*)^{-1} \label{eq:gsd-hor-II}\\
\toiso & \tilde{sd}_{n,k}^*\circ(sd_{n,k}^*)^{-1}\circ(\mathsf{s}_3^*)^{-1}\\
\toiso & \mathsf{s}_3^*\circ(\mathsf{s}_3^*)^{-1}\\
\toiso & \id_{sl\D_{n,k}},
\end{align}
where the single step are
\begin{enumerate}
\item first, induced by the equivalence $sd_{n,k}^*$ \autoref{rmk:slices},
\item second, the equality $\phi_{n,k}\circ sd_{n,k}=d^h$,
\item third, the invertibility of the unit of the adjunction $\mathsf{d}^h[-1]\dashv\mathsf{d}^h$ \autoref{prop:preparations-h},
\item fourth, induced by \autoref{rmk:slices},
\item fifth, induced by the equivalence $\mathsf{s}_3^*$,
\end{enumerate}
completes the proof.
\end{proof}

As an application, we obtain the horizontal analogue of \autoref{cor:gsdver}.

\begin{cor}\label{cor:gsd-dual}
Let $n\geq 1,k\geq 2$ and $\D$ a stable derivator. Then the standard horizontal face morphism $\mathsf{d}^h\colon\D_{n,k+1}\to\D_{n,k}$ generates an infinite chain of adjunctions such that for $p\in\mathbb{Z}$
\begin{enumerate}
\item $\mathsf{d}^h[2p]=(\mathsf{s}_3^*)^p\circ\mathsf{d}^h\circ(\mathsf{s}_3^*)^{-p}$,
\item $\mathsf{d}^h[2p+1]=(\mathsf{s}_3^*)^p\circ\mathsf{d}^h[1]\circ(\mathsf{s}_3^*)^{-p}=(\mathsf{s}_3^*)^{p+1}\circ\mathsf{d}^h[-1]\circ(\mathsf{s}_3^*)^{-p-1}$.
\end{enumerate}
\end{cor}

\begin{proof}
This follows inductively from the relations
\[
\mathsf{d}^h[-1]=(\mathsf{s}_3^*)^{-1}\circ\mathsf{d}^h[1]\circ\mathsf{s}_3^*\qquad\text{and}\qquad\mathsf{d}^h[1]=\mathsf{s}_3^*\circ\mathsf{d}^h[-1]\circ(\mathsf{s}_3^*)^{-1},
\]
which are consequences of \autoref{thm:Serre-h}.
\end{proof}

\section{The structure of a local square}
\label{sec:recollements}

In \S\ref{sec:vertical} and \S\ref{sec:horizontal} the main result was the existence of a good supply of morphisms relating derivators of the form $\D_{n,k}$ in the situation where $k$, respectively $n$, is fixed. In this section the main goal is to understand how the structure morphisms in the vertical and horizontal direction interact. For this we will consider those cases, where the relevant maps can be described as inverse images on (triangular) slices. More precisely, we will first show that the horizontal face maps $\mathsf{d}^h$ assemble into a map of certain \emph{simplicial} derivators associated to higher $\mathsf{S}_{\bullet}$-constructions. We will build on this result in \S\ref{sec:main}. Second, by additionally analyzing some boundary cases, we will obtain a full understanding of commutative  (up to natural isomorphism) squares of the form
\[
\xymatrix{
\D_{n+1,k}\ar[r]\ar[d]&\D_{n+1,k+1}\ar[d]\\
\D_{n,k}\ar[r]&\D_{n,k+1},
}
\]
 in the 2-category of derivators, where all displayed morphisms are generalized face and degeneracy morphisms.

\noindent First we observe that standard horizontal face morphism can also on triangular slices be defined as an inverse image morphism. More precisely, the functor
\[
d^h\colon Sl^{\triangle}_{n,k}\to Sl^{\triangle}_{n,k+1},(0,f_1,\cdot,f_{k-1})\mapsto(0,1,f_1+1,\cdots,f_{k-1}+1)
\]
satisfies $\triangle^{n,k+1}\circ d^h=d^h\circ\triangle^{n,k}\colon Sl_{n,k}\to Sl^{\triangle}_{n,k+1}$ and therefore \autoref{thm:slices} and \autoref{con:triangle} yield
\[
\xymatrix{
\D_{n,k+1}\ar[r]^{\mathsf{d}^h}\ar[d]_{sl^{\triangle}_{n,k+1}}&\D_{n,k}\ar[d]^{sl^{\triangle}_{n,k}}\\
sl^{\triangle}\D_{n,k+1}\ar[r]_{(d^h)^*}&sl^{\triangle}\D_{n,k}.
}
\]

\begin{prop}\label{prop:relations}
Let $n\geq 0, k\geq 2$. The following squares commute for $1\leq i \leq n+k$.
\begin{enumerate}
\item
\begin{equation}\label{eq:III}
\xymatrix{
Sl^{\triangle}_{n+1,k}\ar[r]^{d^h}\ar[d]_{\underline{\Lambda}(\Lambda_{k-1},\mathrm{s}_{i-1})}&Sl^{\triangle}_{n+1,k+1}\ar[d]^{\underline{\Lambda}(\Lambda_k,\mathrm{s}_{i})}\\
Sl^{\triangle}_{n,k}\ar[r]_{d^h}&Sl^{\triangle}_{n,k+1}.
}
\end{equation}
\item
\begin{equation}\label{eq:IV}
\xymatrix{
Sl_{n+1,k}\ar[r]^{d^h}&Sl_{n+1,k+1}\\
Sl_{n,k}\ar[r]_{d^h}\ar[u]^{\underline{\Lambda}(\Lambda_{k-1},\mathrm{d}_{i})}&Sl_{n,k+1}\ar[u]_{\underline{\Lambda}(\Lambda_k,\mathrm{d}_{i+1})}.
}
\end{equation}
\end{enumerate}
\end{prop}

\begin{proof}
For the first statement we consider $f=(0,f_1,\cdots,f_{k-1})\in Sl^{\triangle}_{n+1,k}$ such that $j$ is maximal in $\mathbf{k}$ with $f_j\leq i-1$. Then we compute the composition through the upper right vertex as
\begin{align}
& (0,f_1,\cdots,f_j,f_{j+1},\cdots,f_{k-1})\\
\mapsto & (0,1,f_1+1,\cdots,f_j+1,f_{j+1}+1,\cdots,f_{k-1}+1)\\
\mapsto & (0,1,f_1+1,\cdots,f_j+1,f_{j+1},\cdots,f_{k-1})
\end{align}
and the composition through the lower left vertex as
\begin{align}
& (0,f_1,\cdots,f_j,f_{j+1},\cdots,f_{k-1})\\
\mapsto & (0,f_1,\cdots,f_j,f_{j+1}-1,\cdots,f_{k-1}-1)\\
\mapsto & (0,1,f_1+1,\cdots,f_j+1,f_{j+1},\cdots,f_{k-1}).
\end{align}
For the second statement we consider $f=(0,f_1,\cdots,f_{k-1})\in Sl_{n,k}$ such that $j$ is maximal in $\mathbf{k}$ with $f_j< i$. Then we compute the composition through the lower right vertex as
\begin{align}
& (0,f_1,\cdots,f_j,f_{j+1},\cdots,f_{k-1})\\
\mapsto & (0,1,f_1+1,\cdots,f_j+1,f_{j+1}+1,\cdots,f_{k-1}+1)\\
\mapsto & (0,1,f_1+1,\cdots,f_j+1,f_{j+1}+2,\cdots,f_{k-1}+2)
\end{align}
and the composition through the lower left vertex as
\begin{align}
& (0,f_1,\cdots,f_j,f_{j+1},\cdots,f_{k-1})\\
\mapsto & (0,f_1,\cdots,f_j,f_{j+1}+1,\cdots,f_{k-1}+1)\\
\mapsto & (0,1,f_1+1,\cdots,f_j+1,f_{j+1}+2,\cdots,f_{k-1}+2).
\end{align}
\end{proof}

Let $k\geq 0$. In the following we denote by 
\[
k^+\colon\Delta\to\Delta,\Delta_n\mapsto\Delta_{n+k},\mathrm{d}_i\mapsto\mathrm{d}_{i+k},\mathrm{s}_i\mapsto\mathrm{s}_{i+k}
\]
and observe the obvious transitivity property $l^+\circ k^+=(k+l)^+$.

\begin{cor}\label{cor:dh-comp}
Let $k\geq 1$ and $\D$ be a stable derivator.
\begin{enumerate}
\item The standard horizontal face morphisms $\mathsf{d}^h\colon\D_{n,k+2}\to\D_{n,k+1}$ assemble into  a pseudonatural transformation 
\[
\mathsf{S}^{(k+1)}_{\bullet}(\D)\circ\mathsf{i}\circ 2^+\to\mathsf{S}^{(k)}_{\bullet}(\D)\circ\mathsf{i}\circ 1^+\colon\Delta^{op}\to Der.
\]
\item Let $0\leq a\leq 2k$. If $a$ is even, then the morphisms $\mathsf{d}^h[a]\colon\D_{n,k+2}\to\D_{n,k+1}$ assemble into  a pseudonatural transformation 
\[
\mathsf{S}^{(k+1)}_{\bullet}(\D)\circ\mathsf{i}\circ (k+2)^+\to\mathsf{S}^{(k)}_{\bullet}(\D)\circ\mathsf{i}\circ (k+1)^+\colon\Delta^{op}\to Der.
\]
If $a$ is odd, then the morphisms $\mathsf{d}^h[a]\colon\D_{n,k+1}\to\D_{n,k+2}$ assemble into  a pseudonatural transformation 
\[
\mathsf{S}^{(k)}_{\bullet}(\D)\circ\mathsf{i}\circ (k+1)^+\to\mathsf{S}^{(k+1)}_{\bullet}(\D)\circ\mathsf{i}\circ (k+2)^+\colon\Delta^{op}\to Der.
\]
\end{enumerate}
\end{cor}

\begin{proof}
We invoke \autoref{prop:conjugation} applied to the equivalences of \autoref{thm:slices} and \autoref{con:triangle} for a pseudonatural equivalence
\[
\mathsf{S}^{(k)}_{\bullet}(\D)\circ\mathsf{i}\circ 1^+\to\mathsf{S}_{\bullet}^{(k),\triangle}(\D),
\]
where $\mathsf{S}_{\bullet}^{(k),\triangle}(\D)\colon\Delta\to Der$ is the 2-functor defined by 
\begin{itemize}
\item $\Delta_n\mapsto sl^{\triangle}\D_{n-k+1,k+1}$,
\item $(f\colon\Delta_n\to\Delta_{n'})\mapsto(f_0\colon\Delta_{n+1}\to\Delta_{n'+1},0\mapsto 0,i+1\mapsto f(i)+1)\mapsto\Delta(\Delta_k,f_0)_0^*.$
\end{itemize}
To show the first statement we claim that the morphism $\mathsf{d}^h\colon\D_{n,k+2}\to\D_{n,k+1}$ assemble into a 2-natural transformation $\mathsf{S}_{\bullet}^{(k+1),\triangle}(\D)\circ 1^+\to\mathsf{S}_{\bullet}^{(k),\triangle}(\D)$. For this it is sufficient to check the naturality condition on the generators of $\Delta$, which is exactly the 2-functoriality of $\D$ applied to \autoref{prop:relations}.

\noindent For the second statement we additionally use the pseudofunctoriality of right adjoints, since in this case the pseudonaturality squares are obtained by passing to the $a$th right adjoints of the inverse image squares associated to \autoref{prop:relations}.
\end{proof}

\begin{prop}\label{prop:square_I}
Let $n\geq 0,k\geq 2,1\leq i\leq n+k+1$ and $\D$ a stable derivator. Then there are natural isomorphisms
\begin{enumerate}
\item $\mathsf{d}^h_{n+1,k}\circ\mathsf{d}^v_{n,k+1}[-2k-1]\cong0\colon\D_{n,k+1}\to\D_{n+1,k}$,
\item $\mathsf{d}^h_{n+1,k}\circ\mathsf{d}^v_{n,k+1}[2(i-k)-1]\cong\mathsf{d}^v_{n,k}[2(i-k)-1]\circ\mathsf{d}^h_{n,k}\colon\D_{n,k+1}\to\D_{n+1,k}$.
\end{enumerate}
\end{prop}

\begin{proof}
We use the notation from \autoref{prop:the-sieve}
\[
d^v=\underline{\Lambda}(\Lambda_k,\mathrm{d}_1)\vert_{Sl_{n,k+1}}\colon(0,f_1,\cdots,f_k)\mapsto(0,f_1+1,\cdots,f_k+1).
\] 
Since the image of $d^h$ is contained in the complement of the image of $d^v$ and the left Kan extension $d^v_!$ is an extension-by-zero morphism, the composition 
\begin{equation}\label{eq:square-I}
\D^{Sl_{n,k+1}}\xrightarrow{d^v_!}\D^{Sl_{n+1,k+1}}\xrightarrow{(d^h)^*}\D^{Sl_{n+1,k}}
\end{equation}
vanishes. We apply \autoref{prop:the-sieve} and \autoref{eg:explicit} to the left arrow and \autoref{defn:hormor} to the right arrow, to conclude that \eqref{eq:square-I} restricts to 
\begin{equation}\label{eq:square-II}
sl\D_{n,k+1}\xrightarrow{\mathsf{d}_{n,k+1}^v[-2k-1]}sl\D_{n+1,k+1}\xrightarrow{\mathsf{d}_{n+1,k}^h}sl\D_{n+1,k}.
\end{equation}
Hence also the composition \eqref{eq:square-II} has to vanish. Now \autoref{thm:slices} yields the first statement.

For the second statement, we consider for $1\leq i\leq n+k$ the following square of slices (and observe that the vertical morphisms are well defined by \autoref{prop:slice-comp})
\begin{equation}\label{eq:square-III}
\xymatrix{
Sl^{\triangle}_{n+1,k}\ar[r]^{d^h}\ar[d]_{\underline{\Lambda}(\Lambda_{k-1},\mathrm{s}_{i-1})}&Sl^{\triangle}_{n+1,k+1}\ar[d]^{\underline{\Lambda}(\Lambda_k,\mathrm{s}_{i})}\\
Sl^{\triangle}_{n,k}\ar[r]_{d^h}&Sl^{\triangle}_{n,k+1}.
}
\end{equation}
This square is commutative by \autoref{prop:relations}. We apply the derivator $\D$ to the square \eqref{eq:square-III} and obtain by \autoref{cor:slice-restr} and \autoref{rmk:dh-welldef} a restricted square
\[
\xymatrix{
sl^{\triangle}\D_{n+1,k} & sl^{\triangle}\D_{n+1,k+1}\ar[l]_{\mathsf{d}^h}\\
sl^{\triangle}\D_{n,k} \ar[u]^{\mathsf{d}^v[2(i-k)-1]} & sl^{\triangle}\D_{n,k+1} \ar[l]^{\mathsf{d}^h} \ar[u]_{\mathsf{d}^v[2(i-k)-1]}.
}
\]
\autoref{thm:slices} yields the second statement for $1\leq i\leq n+k$.

\noindent It remains to show the case $i=n+k+1$. For this we consider the square (which is again well defined by \autoref{prop:slice-comp})
\begin{equation}\label{eq:square-IV}
\xymatrix{
Sl_{n+1,k}\ar[r]^{d^h}&Sl_{n+1,k+1}\\
Sl_{n,k}\ar[r]_{d^h}\ar[u]^{d^{v\vee}=\underline{\Lambda}(\Lambda_{k-1},\mathrm{d}_{n+k})}&Sl_{n,k+1}\ar[u]_{\underline{\Lambda}(\Lambda_k,\mathrm{d}_{n+k+1})=d^{v\vee}}.
}
\end{equation}
We claim that also this square commutes. However, this is clear as the effect of the vertical morphisms on coordinate representations is trivial. Moreover, this square is a pullback of sieve and hence homotopy exact by \cite[Prop.~1.24]{groth:ptstab}. This yields the commutativity of mate (obtained by passing to the right adjoints in the vertical direction) of the associated inverse image square. However, this in turn restricts by \autoref{prop:the-sieve} and \autoref{rmk:dh-welldef} to 
\[
\xymatrix{
sl\D_{n+1,k} & sl\D_{n+1,k+1}\ar[l]_{\mathsf{d}^h}\\
sl\D_{n,k} \ar[u]^{\mathsf{d}^v[2n+1]} & sl\D_{n,k+1} \ar[l]^{\mathsf{d}^h} \ar[u]_{\mathsf{d}^v[2n+1]}.
}
\]
A final application of \autoref{thm:slices} completes the proof.
\end{proof} 

\begin{thm}\label{thm:full-square}
Let $n\geq 0, k\geq 2$, $a\in\mathbb{Z}$ even, $b\in\mathbb{Z}$ odd, and $\D$ a stable derivator.
\begin{enumerate}
\item Let $p\in\mathbb{Z}$ and $i\in \lbrace -2k-1,-2k+1,\cdots,2n+1\rbrace$ the unique elements, such that $b-a=2p(n+k+2)+i$. Then there are natural equivalences
\[
\mathsf{d}_{n+1,k}^h[a]\circ\mathsf{d}^v_{n,k+1}[b]\cong
\begin{cases}
0 \text{ for } i=-2k-1,\\
\mathsf{d}^v_{n,k}[b-2p]\circ\mathsf{d}^h_{n,k}[a+2p] \text{ for } i\geq -2k+1,
\end{cases}
\]
\item Let $q\in\mathbb{Z}$ and $j\in \lbrace -2n-1,-2n+1,\cdots,2k+1\rbrace$ the unique elements, such that $b-a=2q(n+k+2)+j$. Then there are natural equivalences
\[
\mathsf{d}_{n,k+1}^v[a]\circ\mathsf{d}^h_{n+1,k}[b]\cong
\begin{cases}
0 \text{ for } j=2k+1,\\
\mathsf{d}^h_{n,k}[b-2q]\circ\mathsf{d}^h_{n,k}[a+2q] \text{ for } j\leq 2k-1.
\end{cases}
\]
\item There are natural equivalences
\[
\mathsf{d}_{n+1,k}^h[a+2(n+k+1)]\circ\mathsf{d}^v_{n,k+1}[b+2(n+k+1)]\cong\Sigma^{k-n-1}\circ\mathsf{d}_{n+1,k}^h[a]\circ\mathsf{d}^v_{n,k+1}[b]
\]
and
\[
\mathsf{d}_{n,k+1}^v[a+2(n+k+1)]\circ\mathsf{d}^h_{n+1,k}[b+2(n+k+1)]\cong\Omega^{k-n-1}\circ\mathsf{d}_{n,k+1}^v[a]\circ\mathsf{d}^h_{n+1,k}[b].
\]
\end{enumerate}
\end{thm}

\begin{proof}
We use the notation $l=n+k+2$ and show the first statement in the case $a=0$. To achieve this, we consider the following composition of natural isomorphisms

\scalebox{.9}{\parbox{.95\hsize}{%
\begin{align}
& \mathsf{d}_{n+1,k}^h\circ\mathsf{d}_{n,k+1}^v[b]\\
= &\mathsf{d}_{n+1,k}^h\circ\mathsf{d}_{n,k+1}^v[i+2pl]\\
\cong &\mathsf{d}_{n+1,k}^h\circ(\mathsf{s}_3^*)^{pl}\circ\mathsf{d}_{n,k+1}^v[i]\circ(\mathsf{s}_3^*)^{-pl}\\
\cong &\mathsf{d}_{n+1,k}^h\circ(\Sigma^{(n+1)k})^{p}\circ\mathsf{d}_{n,k+1}^v[i]\circ(\mathsf{s}_3^*)^{-p}\circ(\Sigma^{nk})^{-p}\\
\cong &\mathsf{d}_{n+1,k}^h\circ\mathsf{d}_{n,k+1}^v[i]\circ(\mathsf{s}_3^*)^{-p}\circ(\Sigma^{(n+1)k})^{p}\circ(\Sigma^{nk})^{-p}\\
\cong &\mathsf{d}_{n,k}^v[i]\circ\mathsf{d}_{n,k}^h\circ(\mathsf{s}_3^*)^{-p}\circ(\Sigma^{(n+1)k})^{p}\circ(\Sigma^{nk})^{-p}\\
\cong &\mathsf{d}_{n,k}^v[i]\circ(\mathsf{s}_3^*)^{-p}\circ\mathsf{d}_{n,k}^h[2p]\circ(\Sigma^{(n+1)k})^{p}\circ(\Sigma^{nk})^{-p}\\
\cong &(\mathsf{s}_3^*)^{-p(l-1)}\circ\mathsf{d}_{n,k}^v[i+2p(l-1)]\circ(\mathsf{s}_3^*)^{p(l-1)}\circ(\mathsf{s}_3^*)^{-p}\circ\mathsf{d}_{n,k}^h[2p]\circ(\Sigma^{(n+1)k})^{p}\circ(\Sigma^{nk})^{-p}\\
\cong &(\Sigma^{(n+1)(k-1)})^{-p}\circ\mathsf{d}_{n,k}^v[i+2p(l-1)]\circ(\Sigma^{n(k-1)})^{p}\circ\mathsf{d}_{n,k}^h[2p]\circ(\Sigma^{(n+1)k})^{p}\circ(\Sigma^{nk})^{-p}\\
\cong &\mathsf{d}_{n,k}^v[i+2p(l-1)]\circ\mathsf{d}_{n,k}^h[2p]\\
= &\mathsf{d}_{n,k}^v[b-2p]\circ\mathsf{d}_{n,k}^h[2p],
\end{align}
}}

where the single steps are induced by
\begin{enumerate}
\item first and tenth (last), the assumption on $b-a$,
\item second and seventh, \autoref{thm:gsd-Sdot},
\item third and eighth, \autoref{cor:fracCY},
\item fourth and ninth, the exactness of morphisms of derivators, and in the latter case also the equality $(n+1)k+n(k-1)-(n+1)(k-1)-nk=0$ for $n,k\in\mathbb{Z}$.
\item fifth, \autoref{prop:square_I},
\item sixth, \autoref{cor:gsd-dual},
\end{enumerate}
respectively. In the case $i=-2k-1$, we consider only the first four steps, and apply \autoref{prop:square_I} to obtain the desired vanishing.

\noindent For the general case, we consider the natural isomorphism (for $i\geq -2k+1$)
\begin{equation}\label{eq:full-square-I}
\mathsf{d}^h_{n+1,k}\circ\mathsf{d}^v_{n,k+1}[b-a]\cong\mathsf{d}^v_{n,k}[b-a-2p]\circ\mathsf{d}^h_{n,k}[2p],
\end{equation}
which exists by the special case above. The uniqueness of adjoints, and their compatibility with composition yields inductively natural equivalences between the $a$th right adjoints in \eqref{eq:full-square-I}, and hence the statement for $a>0$. For $a<0$ we consider left adjoints of \eqref{eq:full-square-I} instead. The arguments for the case $i=-2k-1$ are very similar (and even simpler).

\noindent For the second statement we consider $a'\in\mathbb{Z}$ even, $b'\in\mathbb{Z}$ odd with $b'-a'=2pl+i$ as in the first statement. 
Hence we have a nautral isomorphism
\begin{equation}\label{eq:full-square-II}
\mathsf{d}^h_{n+1,k}[a']\circ\mathsf{d}^v_{n,k+1}[b']\cong\mathsf{d}^v_{n,k}[b'-2p]\circ\mathsf{d}^h_{n,k}[a'+2p].
\end{equation}
We again invoke the uniqueness of adjoints, and their compatibility with composition to obtain a natural isomorphism relating the left adjoints of \eqref{eq:full-square-II}
\begin{equation}\label{eq:full-square-III}
\mathsf{d}^v_{n,k+1}[b'-1]\circ\mathsf{d}^h_{n+1,k}[a'-1]\cong\mathsf{d}^h_{n,k}[a'-1+2p]\circ\mathsf{d}^v_{n,k}[b'-1-2p].
\end{equation}
By substituting $a=b'-1$ and $b=a'-1$, we have $b-a=-(b'-a')=2ql+j$ with $q=-p$ and $i=-j$. Plugging this into \eqref{eq:full-square-III} leads exactly to the second statement in the case $j\leq 2k-1$. The case $j=2k+1$ is again very similar (and even simpler).
For the first part of statement (iii) we consider the following composition of equivalences
\begin{align}
&\mathsf{d}_{n+1,k}^h[a+2(n+k+1)]\circ\mathsf{d}^v_{n,k+1}[b+2(n+k+1)]\\
\cong&(\mathsf{s}_3^*)^{n+k+1}\circ\mathsf{d}_{n+1,k}^h[a]\circ\mathsf{d}^v_{n,k+1}[b]\circ(\mathsf{s}_3^*)^{-(n+k+1)}\\
\cong&\Sigma^{(n+1)(k-1)}\circ\mathsf{d}_{n+1,k}^h[a]\circ\mathsf{d}^v_{n,k+1}[b]\circ\Omega^{nk}\\
\cong&\Sigma^{k-n-1}\circ\mathsf{d}_{n+1,k}^h[a]\circ\mathsf{d}^v_{n,k+1}[b],
\end{align}
where the first equivalence is \autoref{cor:gsdver} and \autoref{cor:gsd-dual}, the second equivalence in \autoref{cor:fracCY}, and the third equivalence is the exactness of morphisms of derivators.
\end{proof}

\begin{prop}\label{prop:dual-square}
Let $n\geq 0,k\geq 2$, $a,b\in\mathbb{Z}$ such that $b-a$ is even, and $\D$ a stable derivator. If there are $p\in\mathbb{Z}$ and $-k+1\leq i \leq n$, such that $b-a=2p(n+k+1)+2i$, then there are natural isomorphisms
\begin{enumerate}
\item 
\[
\mathsf{d}^h_{n,k}[a]\circ\mathsf{d}^v_{n,k+1}[b]\cong\Sigma^{p(n-k+1)}\circ\mathsf{d}^v_{n,k}[b-2p(n+k+1)]\circ\mathsf{d}^h_{n+1,k}[a+2p(n+k+1)]
\]
and
\[
\mathsf{d}^v_{n,k}[b]\circ\mathsf{d}^h_{n+1,k}[a]\cong\Sigma^{p(n-k+1)}\circ\mathsf{d}^h_{n,k}[a+2p(n+k+1)]\circ\mathsf{d}^v_{n,k+1}[b-2p(n+k+1)]
\]
if $a,b$ are even,
\item
\[
\mathsf{d}^h_{n+1,k}[a]\circ\mathsf{d}^v_{n,k}[b]\cong\Omega^{p(n-k+1)}\circ\mathsf{d}^v_{n,k+1}[b-2p(n+k+1)]\circ\mathsf{d}^h_{n,k}[a+2p(n+k+1)]
\]
and
\[
\mathsf{d}^v_{n,k+1}[b]\circ\mathsf{d}^h_{n,k}[a]\cong\Omega^{p(n-k+1)}\circ\mathsf{d}^h_{n+1,k}[a+2p(n+k+1)]\circ\mathsf{d}^v_{n,k}[b-2p(n+k+1)]
\]
if $a,b$ are odd.
\end{enumerate}
\end{prop}

\begin{proof}
We apply inverse images to \autoref{prop:relations} (i) and reformulate using the notation introduced in \autoref{eg:ver-fd} and \autoref{defn:hormor} to obtain isomorphisms
\[
\mathsf{d}^h_{n,k}\circ\mathsf{d}^v_{n,k+1}[2i]\cong\mathsf{d}^v_{n,k}[2i]\circ\mathsf{d}^h_{n+1,k}
\]
for $-k+1\leq i\leq n$. By passing to adjoint isomorphisms, we conclude the first part of statement (i) in the case $p=0$. For $p\in\mathbb{Z}$ general, we use the notation $l=n+k+1$ and invoke the following composition of isomorphisms
\begin{align}
&\mathsf{d}^h_{n,k}[a]\circ\mathsf{d}^v_{n,k+1}[b]\\
\cong&\mathsf{d}^h_{n,k}[a]\circ\mathsf{s}_3^{pl}\circ\mathsf{d}^v_{n,k+1}[b-2pl]\circ\mathsf{s}_3^{-pl}\\
\cong&\mathsf{d}^h_{n,k}[a]\circ\Sigma^{pnk}\circ\mathsf{d}^v_{n,k+1}[b-2pl]\circ\mathsf{s}_3^{-pl}\\
\cong&\Sigma^{pnk}\circ\mathsf{d}^h_{n,k}[a]\circ\mathsf{d}^v_{n,k+1}[b-2pl]\circ\mathsf{s}_3^{-pl}\\
\cong&\Sigma^{pnk}\circ\mathsf{d}^v_{n,k}[b-2pl]\circ\mathsf{d}^h_{n+1,k}[a]\circ\mathsf{s}_3^{-pl}\\
\cong&\Sigma^{p(n-k+1)}\circ\mathsf{d}^v_{n,k}[b-2pl]\circ\Sigma^{p(n+1)(k-1)}\circ\mathsf{d}^h_{n+1,k}[a]\circ\mathsf{s}_3^{-pl}\\
\cong&\Sigma^{p(n-k+1)}\circ\mathsf{d}^v_{n,k}[b-2pl]\circ\mathsf{s}_3^{pl}\circ\mathsf{d}^h_{n+1,k}[a]\circ\mathsf{s}_3^{-pl}\\
\cong&\Sigma^{p(n-k+1)}\circ\mathsf{d}^v_{n,k}[b-2pl]\circ\mathsf{d}^h_{n+1,k}[a+2pl],
\end{align}
where the first isomorphisms is \autoref{cor:gsdver}, the second and sixth one is \autoref{cor:fracCY}, the third and fifth one is the exactness of morphisms of derivators, the fourth one is statement (i) in the case $p=0$, and the seventh one is \autoref{cor:gsd-dual}.

\noindent For the second part of (i) we use the first part with $a'=a+2pl$ and $b'=b-2pl$. We note that then $b'-a'=-2pl+2i$. 

The first (resp. second) isomorphism of (ii) is obtained from the second (resp. first) isomorphism of (i) by passing to adjoints.
\end{proof}

\begin{rmk}
It is worth to note that under the substitution $(n',k')=(k-1,n+1)$ the suspension exponent tranforms as $n'-k'+1=k-n-1=-(n-k+1)$. In fact one can show that the two equivalences of \autoref{thm:full-square} (iii), \autoref{prop:dual-square} (i) and (ii), respectively, can be transformed into each other by conjugating with the equivalences of \autoref{thm:main} and applying the substitution above.
\end{rmk}

\section{The symmetry theorem - The case \texorpdfstring{$k=2$}{k=2}}
\label{sec:contra}

In this section we establish a comparison result between the derivators $\D_{n,2}$ and $\D_{1,n+1}$ for $n\geq 1$. This can be regarded as a generalization of the following classical result from representation theory. Let $\mathrm{k}$ be a field and $R_n$ be the quotient of the path algebra $\mathrm{k}\A{n}$ by the ideal generated by paths of length two.

\begin{prop}\label{prop:An-zero}
Let $n\geq 1$ and $\mathrm{k}$ be a field. Then there is an exact equivalence of triangulated categories
\[
D^b(\mathrm{k}\A{n})\toiso D^b(R_n).
\]
\end{prop}

\begin{proof}
This is \cite[Prop.~2.1]{HS-piecewise}.
\end{proof}

Let us consider the special case $n=4$. In this case an $R_n$-module is is given by a functor $F\colon[2]\times[1]\to\mathrm{k}-Mod$ such that $F(0,1)\cong F(2,0)\cong 0$. Using the universal property of the zero-object we can extend $F$ to a functor $\cube{3}\to\mathrm{k}-Mod$ as indicated by the diagram
\begin{equation}\label{an-rel}
\xymatrix{
&0\ar@{-->}[rr]\ar@{-->}[dd]&&0\ar[dd]\\
F(0,0)\ar@{-->}[ur]\ar[rr]\ar[dd]&&F(1,0)\ar[ur]\ar[dd]\\
&0\ar@{-->}[rr]&&F(2,1)\\
0\ar@{-->}[ur]\ar[rr]&&F(1,1).\ar[ur]
}
\end{equation}
Hence we obtain an equivalence between the corresponding functor categories with the respective zero-conditions. We will show in \autoref{thm:dual-Sdot-sl} that the description of the relations based on cubical shapes \eqref{an-rel} generalizes well to the setting of stable derivators. This motivates us to introduce the  notion of a cube supported on a maximal path in a stable derivator $\D$. The resulting derivators are of the form $sl\D_{1,n}$ (\autoref{prop:higher-cof}).

\begin{defn}
Let $n,k\geq 0$. A non-degenerate $k$-simplex $s\colon[k]\to\cube{n}$ is called a \textbf{path} in $\cube{n}$. A path is called \textbf{maximal} if $k=n$.
\end{defn}

\begin{eg}
Let $n\geq 0$. The standard maximal path is defined by the functor
\[
\rightarrow\colon [n]\to\cube{n},i\mapsto\lbrace n-i,\cdots,n-1\rbrace.
\]
Moreover, all other maximal paths in $\cube{n}$ are obtained from $\rightarrow$ and a permutation $\sigma\in\mathrm{Aut}(\mathbf{n})$ by $\rightarrow_{\sigma}\colon [n]\to\cube{n},i\mapsto\sigma(\rightarrow(i))$. Furthermore, this construction defines a bijection between $\mathrm{Aut}(\mathbf{n})$ and the set of maximal paths in $\cube{n}$.
\end{eg}

\begin{defn}
Let $n\geq 0$, $\D$ a stable derivator and $\rightarrow_{\sigma}$ a maximal path in $\cube{n}$. The \textbf{derivator of} $n$\textbf{-cubes with} $\rightarrow_{\sigma}$\textbf{-support} $\D^{\cube{n}}_{\rightarrow_{\sigma}}$ is the full subderivator of $\D^{\cube{n}}$ spanned by those objects $X$ such that for all $M\subseteq\mathbf{n}$ with $M$ not in the image of $\rightarrow_{\sigma}$ the property (P2) holds, i.e. $M^*X\cong 0$.
\end{defn}

\begin{prop}\label{prop:higher-cof}
Let $n\geq 2$ and $\D$ a stable derivator. Then there are equivalences of derivators
\begin{equation}\label{eq:max-path}
sl\D_{1,n}\toiso\D^{\cube{n-1}}_{\rightarrow}\qquad\text{and}\qquad do\D_{1,n}\toiso\D^{\cube{n},ex}_{\rightarrow}.
\end{equation}
\end{prop}

\begin{proof}
We note that in this case $\xi=(0,1,\cdots,n-1)$ and by \autoref{lem:symmetries} 
\[
\mathsf{s}_3^n=\mathsf{s}_2^n\circ(\mathsf{s}_1)^{-n}=\mathsf{s}_1.
\]
Therefore, we conclude that $(\mathsf{s}_3^*)^n(\xi)=(1,2,\cdots,n)$ and $(\mathsf{s}_3^*)^{n-1}(\xi)=(0,2,\cdots,n)$ that hence $\square_{\xi}\colon\cube{n-1}\to Sl_{1,n}\subset\mathbb{Z}^{n-1}$ and $\square_{\xi}\colon\cube{n}\to Do_{1,n}\subset\mathbb{Z}^{n}$ are isomorphisms of categories. Moreover, under these isomorphisms the objects of $\cube{n-1}$ and $\cube{n}$ with non-decreasing coordinates, i.e. the objects in the image of the maximal path $\rightarrow$, correspond exactly to the injective objects in $Sl_{1,n}$ and $Do_{1,n}$, respectively. As a consequence, we obtain that the inverse images $\square_{\xi}^*$ restrict to an equivalence $sl\D_{1,n}\toiso\D^{\cube{n-1}}_{\rightarrow}$ and an embedding $do\D_{1,n}\to\D^{\cube{n}}_{\rightarrow}$, respectively. In the latter case, the property (P1) for $\xi$ (which is unique in $Do_{1,n}$ with $\square_{\xi}\subseteq Do_{1,n}$) allows us to identify the essential image of this embedding with $\D^{\cube{n},ex}_{\rightarrow}$.
\end{proof}

\begin{defn}\label{defn:higher-cof}
Let $n\geq 1$ and $\D$ a stable derivator.
\begin{enumerate}
\item An $n$\textbf{-cofiber sequence} in $\D$ is an object in $\D_{1,n+1}$.
\item Let $X$ be an $n$-cofiber sequence in $D$. Then the \textbf{base} of $X$ is the object $sl_{1,n+1}^*(X)\in sl\D_{1,n+1}.$
\item Let $X$ be an $n$-cofiber sequence in $D$. Then the $n$\textbf{-cone} of $X$ is the object $\mathsf{s}_1(\xi)^*(X)\in\D$.
\end{enumerate}
\end{defn}

\begin{rmk}\label{rmk:highercof}
Let $n\geq 1$, $\D$ a stable derivator and $X$ an $n$-cofiber sequence in $\D$. Then by \autoref{prop:higher-cof} and \autoref{cor:dnk-sigma} the underling diagram of $X$ is determined by a sequence
\[
\cdots\to x_0\to x_1 \to \cdots\to x_{n+1}\to\Sigma^nx_0\to\Sigma^nx_1\to\cdots
\]
in $\D(\bbone)$ such that
\begin{itemize}
\item all consecutive compositions vanish,
\item $x_0\to x_1\to\cdots\to x_n$ gives rise to the underlying diagram of the base of $X$,
\item $x_{n+1}$ is the $n$-cone of $X$.
\end{itemize}
Moreover, if $\D$ is a strong derivator the underlying diagrams of 1-cofiber sequences give rise to distinguished triangles. We will come back to this notion in \S\ref{sec:Toda} and \S\ref{sec:triangles}.
\end{rmk}

\begin{thm}\label{thm:dual-Sdot-sl}
Let $n\geq 0$ and $\D$ a stable derivator. Then there is an equivalence of derivators
\[
\Psi_n^{\square}\colon\D^{[n]}\toiso\D^{\cube{n}}_{\rightarrow}.
\]
\end{thm}

\begin{proof}
Let $\tau\in\mathrm{Aut}(\mathbf{n}),\tau(i)=n-1-i$ be the flip permutation. We claim that the desired equivalence of derivators is defined by the composition
\begin{equation}\label{eq:dual-Sdot-sl}
\D^{[n]}\xrightarrow{(\rightarrow_{\tau})_!}\D^{\cube{n}}\xrightarrow{\fib^{\underline{1}}}\D^{\cube{n}}.
\end{equation}
Since $\rightarrow_{\tau}\colon[n]\to\cube{n}$ is fully faithful the same is true for the associated Kan extension morphism. And since $\fib^{\underline{1}}\colon\D^{\cube{n}}\to\D^{\cube{n}}$ is an equivalence (\cite[Prop.~8.9]{bg:cubical}), we deduce that \eqref{eq:dual-Sdot-sl} is an embedding. As a consequence, to conclude it is sufficient to identify the essential image of \eqref{eq:dual-Sdot-sl} with $\D^{\cube{n}}_{\rightarrow}$.
To do this, we proceed in three steps.
\begin{enumerate}
\item For $0\leq i\leq n$ we denote by $L_i\subseteq\cube{n}$ the subcategories spanned by the objects
\[
L_i=\lbrace M \subseteq\mathbf{n}\vert\mathbf{i}\subseteq M, i\notin M\rbrace.
\]
We observe, that the collections $L_i$ are $(n-1-i)$-dimensional subcubes of $\cube{n}$ (here we use the convention that $(-1)$-cubes are singletons) and that the $L_i$ for $0\leq i\leq n$ together define a partition of the set of objects of $\cube{n}$. We denote by
\[
\iota_i\colon L_i\to\cube{n}\qquad\text{and}\qquad \pi_i\colon L_i\to\bbone
\]
the inclusion and the canonical projection, respectively. Moreover, we see that for $0\leq i\leq n$
\[
L_i=\lbrace M\in\cube{n}\vert i\in[n]\text{ is maximal with }\exists (\rightarrow_{\tau}(i)\to M)\in\cube{n}\rbrace.
\]
In particular, the right adjoint $p\colon\cube{n}\to[n]$ of $\rightarrow_{\tau}\colon[n]\to\cube{n}$ exists and can be described by
\[
p(M)=i\iff M\in L_i.
\]
\item We can characterize the essential image of $(\rightarrow_{\tau})_!\colon\D^{[n]}\to\D^{\cube{n}}$ as the subderivator consisting of those objects $x\in\D^{\cube{n}}$ such that the counit 
\[
(\rightarrow_{\tau})_!(\rightarrow_{\tau})^*x\toiso x
\]
is invertible. Since $\rightarrow_{\tau} \dashv p$ we reformulate this by contemplating the pastings
\[
\xymatrix{
\cube{n}\ar[r]^-p\ar[d]_-\id\drtwocell\omit{\varepsilon}&[n]\ar[r]^{\rightarrow_{\tau}}\ar[d]^{\rightarrow_{\tau}}\drtwocell\omit{\id}&\cube{n}\ar[d]^-\id\ar@{}[drr]|{=}&&\cube{n}\ar[r]^-{\rightarrow_{\tau} p}\ar[d]_-\id\drtwocell\omit{\varepsilon}&\cube{n}\ar[d]^-\id\\
\cube{n}\ar[r]_-\id&\cube{n}\ar[r]_-\id&\cube{n}&&\cube{n}\ar[r]_-\id&\cube{n}.
}
\]
Using the homotopy exactness of the square on the left, we conclude that $(\rightarrow_{\tau})_!(\rightarrow_{\tau})^\ast x\to x$ is invertible if and only if the canonical mate 
\[
\varepsilon^\ast\colon p^\ast(\rightarrow_{\tau})^\ast x\to x
\]
is invertible. But this is the case if and only if $\iota_i^\ast x$ is constant (i.e. in the essential image of $\pi_i^*$) for all $0\leq i\leq n$. We denote by $\D^{\cube{n},\kappa}\subseteq\D^{\cube{n}}$ the subderivator of all such objects $x$.
\item We claim that $\cof^{\underline{1}}\colon\D^{\cube{n}}\toiso\D^{\cube{n}}$ restricts to an equivalence of derivators
\begin{equation}\label{eq:dual-Sdot-II}
\D^{\cube{n}}_{\rightarrow}\toiso\D^{\cube{n},\kappa}
\end{equation}
and show this via induction on $n\geq 0$. For $n=0$ the morphism in question can be identified with $\id\colon\D\to\D$, and are therefore equivalences. Let now $n_0\geq 0$ fixed. We assume \eqref{eq:dual-Sdot-II} for all $0\leq n\leq n_0$. 

\noindent Let $x\in\D^{\cube{n_0+1}}_{\rightarrow}$. Then we observe
\begin{equation}\label{eq:dual-Sdot-III}
(\mathrm{d}_1^0)^*x\in\D^{\cube{n_0}}_{\rightarrow}\qquad\text{and}\qquad(\mathrm{d}_0^0)^*x\in\D^{\cube{n_0},\cube{n_0}_{0,n_0-1}}.
\end{equation}
We can regard $x$ as an object of $(\D^{\cube{n_0}})^{[1]}$ with underlying diagram
\[
(\mathrm{d}_1^0)^*x \to (\mathrm{d}_0^0)^*x.
\]
As a consequence, we can compute $\cof^{\underline{1}}(x)$ as the cofiber of an object of $(\D^{\cube{n_0}})^{[1]}$ with underlying diagram
\[
y_0=\cof^{\underline{1}}(\mathrm{d}_1^0)^*x \to y_1=\cof^{\underline{1}}(\mathrm{d}_0^0)^*x.
\]
By induction assumption, \cite[Lem.~10.2]{bg:cubical} and \eqref{eq:dual-Sdot-III} we conclude that $y_0\in\D^{\cube{n_0},\kappa}$ and that $y_1$ is constant, which in turn implies that $\cof^{\underline{1}}(x)\in\D^{\cube{n_0+1},\kappa}$.

\noindent Conversely, let $y\in\D^{\cube{n_0+1},\kappa}$. We observe that $(\mathrm{d}_1^0)^*y$ is constant and $(\mathrm{d}_y^0)^*y\in\D^{\cube{n_0+1},\kappa}$. We can regard $y$ as an object of $(\D^{\cube{n_0}})^{[1]}$ with underlying diagram
\[
(\mathrm{d}_1^0)^*y \to (\mathrm{d}_0^0)^*y.
\]
Similarly as above, we can compute $\fib^{\underline{1}}(y)$ as the fiber of an object of $(\D^{\cube{n_0}})^{[1]}$ with underlying diagram
\[
x_0=\fib^{\underline{1}}(\mathrm{d}_1^0)^*y \to x_1=\fib^{\underline{1}}(\mathrm{d}_0^0)^*y.
\]
By induction assumption, the assumption on $y$ and \cite[Lem.~10.2]{bg:cubical} we obtain $x_1\in\D^{\cube{n_0}}_{\rightarrow}$ and $x_0\in\D^{\cube{n_0},\cube{n_0}_{0,n_0-1}}$. Thus $\fib^{\underline{1}}(y)\in\D^{\cube{n_0+1}}_{\rightarrow}$, which completes the induction step.
\end{enumerate}
\end{proof}

\begin{cor}\label{thm:dual-Sdot}
Let $n\geq 0$ and $\D$ a stable derivator. Then there is an equivalence of derivators
\[
\Psi_n^{\square,ex}\colon\D^{[n]}\toiso\D^{\cube{n+1},ex}_{\rightarrow}.
\]
\end{cor}

\begin{proof}
We consider the diagram
\begin{equation}\label{eq:dual-Sdot-IV}
\xymatrix{
\D^{[n]} \ar[r]^{(\mathrm{d}_{n+1})_*} & \D^{[n+1],\infty}  \ar@{^{(}->}[d] \ar[r]^{(\rightarrow_{\tau})_!} & \D^{\cube{n+1},\kappa,\infty} \ar@{^{(}->}[d] \ar[r]^{\fib^{\underline{1}}} & \D^{\cube{n+1},ex}_{\rightarrow}  \ar@{^{(}->}[d]\\
& \D^{[n+1]}  \ar[r]^{(\rightarrow_{\tau})_!} & \D^{\cube{n+1},\kappa} \ar[r]^{\fib^{\underline{1}}} & \D^{\cube{n+1}}_{\rightarrow}.
}
\end{equation}
Since $\mathrm{d}_{n+1}$ is sieve, \cite[Prop.~3.6]{groth:ptstab} implies that the associated right Kan extension is an extension by zero morphism and 
\[
\mathrm{essim}((\mathrm{d}_{n+1})_*)=\D^{[n+1],\infty}.
\]
Hence the left morphism in the upper row of \eqref{eq:dual-Sdot-IV} is an equivalence. Moreover, it follows from the proof of \autoref{thm:dual-Sdot-sl} that both morphisms in the lower row are equivalences. The left of these equivalences restricts to $(\rightarrow_{\tau})_!\colon\D^{[n+1],\infty}\toiso\D^{\cube{n+1},\kappa,\infty}$ because of $L_{n+1}=\lbrace\infty\rbrace\subset\cube{n+1}$ and the right equivalence restricts to $\fib^{\underline{1}}\colon\D^{\cube{n+1},\kappa,\infty}\toiso\D^{\cube{n+1},ex}_{\rightarrow}$ by \autoref{prop:bicart-obstruction}.
\end{proof}

\begin{cor}\label{cor:Sdot-duality}
Let $n\geq 1$ and $\D$ a stable derivator. Then there is an equivalence of derivators
\[
\tilde{\Psi}_n\colon\D_{n,2}\toiso\D_{1,n+1}.
\]
\end{cor}

\begin{proof}
Let $\tilde{\Psi}_n$ be defined by the following chain of equivalences
\begin{equation}\label{eq:psi-tilde}
\D_{n,2}\toiso sl\D_{n,2}\toiso \D^{[n]}\xrightarrow{\Psi_n^{\square,ex}} \D^{\cube{n+1},ex}_{\rightarrow}\toiso do\D_{1,n+1}\toiso  \D_{1,n+1},
\end{equation}
where the first and fifth equivalence is \autoref{thm:slices}, the second equivalence is \autoref{egs:simplex-slice} (i), the third equivalence is \autoref{thm:dual-Sdot}, and the fourth equivalence is \autoref{prop:higher-cof}.
\end{proof}

\begin{thm}\label{thm:dual-Serre}
Let $n\geq 1$ and $\D$ a stable derivator. Then there are natural isomorphisms
\begin{enumerate}
\item $\xi^*\circ(\mathsf{s}_3^*)^{n+1}\circ\tilde{\Psi}_n\cong\xi^*\colon\D_{n,2}\to\D$ and
\item $\mathsf{s}_3^*\circ\tilde{\Psi}_n\cong\tilde{\Psi}_n\circ\mathsf{s}_3^*\colon\D_{n,2}\toiso\D_{1,n+1}$.
\end{enumerate}
\end{thm}

\begin{proof}
For the first statement, we observe that $\xi^*\colon\D_{n,2}\to\D$ corresponds to $0^*\colon\D^{[n]}\to\D$ under the first two equivalences in \eqref{eq:psi-tilde}, whereas the composition $\xi^*\circ(\mathsf{s}_3^*)^{n+1}\colon\D_{1,n+1}\to\D$ corresponds to $\infty^*\colon\D^{\cube{n+1},ex}_{\rightarrow}\to\D$ under last two equivalences in \eqref{eq:psi-tilde}. Hence the commutativity of the diagram
\[
\xymatrix{
\D^{[n]} \ar[drr]_{0^*} \ar[r]^{\cong} & \D^{[n+1],\infty}  \ar[dr]^{0^*} \ar[r]^{\cong} & \D^{\cube{n+1},\kappa,\infty} \ar[d]^{\emptyset^*} \ar[r]^{\cong} & \D^{\cube{n+1},ex}_{\rightarrow}  \ar[dl]^{\infty^*}\\
& & \D,
}
\]
where the top row is $\Psi_n^{\square,ex}$, completes the proof of (i).

\noindent For the second part we show the equivalent statement
\begin{equation}\label{eq:psi-serre}
(\mathsf{s}_3^*)^{-1}\circ\tilde{\Psi}_n\cong\tilde{\Psi}_n\circ(\mathsf{s}_3^*)^{-1}\colon\D_{1,n+1}\toiso\D_{2,n}.
\end{equation}
Let $q\colon[1]\times[n]\setminus\lbrace(1,n)\rbrace\to[n]$ be the functor defined by $q(0,i)=0$ and $q(1,i)=i+1$ and $\iota\colon[1]\times[n]\setminus\lbrace(1,n)\rbrace\to[1]\times[n]$ the natural inclusion.

\begin{figure}
\begin{displaymath}
\xymatrixcolsep{0.48cm}
\xymatrix{
\D_{n,2}\ar[rrrr]^{\mathsf{s}_3^{-1}}\ar[d]_{sl_{n,2}^*}&&&&\D_{n,2}\ar[d]^{sl_{n,2}^*}\\
sl\D_{n,2}\ar[rrrr]^{\mathsf{s}_3^{-1}}\ar[d]_{c_n^*}&&&&sl\D_{n,2}\ar[d]^{c_n^*}\\
\D^{[n]}\ar[d]_{(\mathrm{d}_{n+1})_*}\ar[r]^{q^*}&\D^{[1]\times[n]\setminus\lbrace(1,n)\rbrace}\ar[rr]^{\iota_*}&&\D^{[1]\times[n]}\ar[r]^{F\times\id}&\D^{[n]}\ar[d]^{(\mathrm{d}_{n+1})_*}\ar@/^1pc/[lldd]^{(\rightarrow_{\tau})_!=p_n^*}\\
\D^{[n+1],\infty}\ar[d]_{p_{n+1}^*=(\rightarrow_{\tau})_!}&&&&\D^{[n+1],\infty}\ar[d]^{p_{n+1}^*=(\rightarrow_{\tau})_!}\\
\D^{\cube{n+1},\kappa,\infty}\ar[rrruu]^{(\id\times(\rightarrow_{\tau}))^*}\ar[d]_{\fib^{\underline{1}}}\ar[rr]_{F\times\id}&&\D^{\cube{n},\kappa}\ar[d]_{\fib^{\underline{1}}}\ar@/^1pc/[rruu]^{(\rightarrow_{\tau})^*}&&\D^{\cube{n+1},\kappa,\infty}\ar[d]^{\fib^{\underline{1}}}\ar[ll]^{(\mathrm{d}_1^n)^*}\\
\D^{\cube{n+1},ex}_{\rightarrow}\ar[d]_{(\square_{\xi}^*)^{-1}}\ar[rr]_{(\mathrm{d}_1^0)^*}&&\D^{\cube{n}}_{\rightarrow}\ar[d]_{(\square_{\xi}^*)^{-1}}&&\D^{\cube{n+1},ex}_{\rightarrow}\ar[d]^{(\square_{\xi}^*)^{-1}}\ar[ll]^{(\mathrm{d}_0^n)^*}\\
do\D_{1,n+1}\ar[d]_{(do_{1,n+1}^*)^{-1}}\ar[rr]_{sd^*}&&sl\D_{1,n+1}&&do\D_{1,n+1}\ar[d]^{(do_{1,n+1}^*)^{-1}}\ar[ll]^{\tilde{sd}^*}\\
\D_{1,n+1}\ar[rrrr]_{\mathsf{s}_3^{-1}}&&&&\D_{1,n+1}
}
\end{displaymath}
\caption{The commutativity of the diagram above implies statement (ii)}
\label{fig:proof1}
\end{figure}

We consider \autoref{fig:proof1} and observe that the outer vertical compositions are by definition $\tilde{\Psi}_n$. Moreover, the top cell commutes by \autoref{rmk:slices} (i), the two squares involving $\fib^{\underline{1}}$ in the vertical direction commute by \autoref{rmk:elementary}, the two squares involving $(\square_{\xi}^*)^{-1}$ commute since they are induced by restriction from inverse image squares associated to commutative squares of functors in $Cat$, and the bottom cell commutes by \autoref{rmk:slices} (iii). In the next step we show that also the cell directly below the top one commutes. Let $r\colon\cube{2}\times[n]\to\underline{\Lambda}_{n+1,1}$ be the functor adjoint to the functor $[n]\to(\underline{\Lambda}_{n+1,1})^{\cube{2}}$ which maps $i$ to the square
\[
\xymatrix{
(-n+i,1)\ar[r]\ar[d] & (-n+i,i+2)\ar[d]\\
(0,1) \ar[r] & (0,i+2).
}
\]
Since the squares above are concatinations of elementary subsquares of $\underline{\Lambda}_{n+1,1}$, we deduce from \autoref{cor:bicart-concat}, that the associated inverse image morphism restricts to $r^*\colon\D_{n,2}\to(\D^{[n]})^{\cube{2},ex}$ and to conclude we note
\begin{itemize}
\item $(0,0)^*r^*=sl_{n,2}^*(\mathsf{s}_3^*)^{-1}$,
\item $(0,1)^*r^*=0$,
\item $(\mathrm{d}_0^0)^*r^*=\iota_*q^*sl_{n,2}^*$.
\end{itemize}
Furthermore, since $\tilde{sd}^*\colon do\D_{1,n+1}\to sl\D_{1,n+1}$ is an equivalence by \autoref{rmk:slices} and the vertical morphisms in the two squares above are equivalences by \autoref{prop:higher-cof}, we conclude that also $(\mathrm{d}_0^n)^*\colon\D^{\cube{n+1},ex}_{\rightarrow}\to\D^{\cube{n}}_{\rightarrow}$ and $(\mathrm{d}_1^n)^*\colon\D^{\cube{n+1},\kappa,\infty}\to\D^{\cube{n},\kappa}$ are equivalences. In the following we construct the isomorphism \eqref{eq:psi-serre} as the composition of isomorphisms
\begin{align}
&\tilde{\Psi}_n\circ\mathsf{s}_3^{-1}\\
\cong&f_2\circ(\rightarrow_{\tau})_!\circ(\mathrm{d}_{n+1})_*\circ(F\times\id)\circ\iota_*\circ q^*\circ f_1\\
\cong&f_2\circ(\rightarrow_{\tau})_!\circ(\mathrm{d}_{n+1})_*\circ(F\times\id)\circ(\id\times(\rightarrow_{\tau}))^*\circ p_{n+1}^*\circ(\mathrm{d}_{n+1})_* \circ f_1\\
\cong&f_2\circ(\rightarrow_{\tau})_!\circ(\mathrm{d}_{n+1})_*\circ(\rightarrow_{\tau})^*\circ(F\times\id)\circ p_{n+1}^*\circ(\mathrm{d}_{n+1})_* \circ f_1\\
\cong&f_2\circ((\mathrm{d}_1^n)^*)^{-1}\circ(\rightarrow_{\tau})_!\circ(\rightarrow_{\tau})^*\circ(F\times\id)\circ p_{n+1}^*\circ(\mathrm{d}_{n+1})_* \circ f_1\\
\cong&f_2\circ((\mathrm{d}_1^n)^*)^{-1}\circ(F\times\id)\circ p_{n+1}^*\circ(\mathrm{d}_{n+1})_* \circ f_1\\
\cong&\mathsf{s}_3^{-1}\circ\tilde{\Psi}_n.
\end{align}
Here we have used the abbreviations $f_1$ and $f_2$ for the composition of the first two morphisms in the left column and the composition of the last three morphisms in the right column of \autoref{fig:proof1}, respectively. Moreover, the first and last isomorphism above follows from the commutativity of top two and bottom three rows in the diagram, respectively. For the remaining isomorphisms we consider the following.
\begin{itemize}
\item For the second isomorphism we show that the left triangle commutes. We denote by $q'$ the composition
\[
[1]\times[n]\xrightarrow{\id\times(\rightarrow_{\tau})}\cube{n+1}\xrightarrow{p_{n+1}}[n+1]
\] 
and observe that $q'(0,i)=0$ and $q'(1,i)=i+1$. It is sufficient to show that the square
\[
\xymatrix{
\D^{[n]}\ar[r]^{q^*}\ar[d]_{(\mathrm{d}_{n+1})_*}&\D^{[1]\times[n]\setminus\lbrace(1,n)\rbrace}\ar[d]^{\iota_*}\\
\D^{[n+1]}\ar[r]_{q'^*}&\D^{[1]\times[n]}
}
\]
commutes. The above description of $q'$ implies $(1,n)^*\circ q'^*\circ(\mathrm{d}_{n+1})_*=0$. Hence the essential image of $q'^*\circ(\mathrm{d}_{n+1})_*$ is contained in the essential image of $\iota_*$ (\cite[Prop.~1.23]{groth:ptstab}). Therefore, it is sufficient to show that $\iota^*\circ q'^*\circ(\mathrm{d}_{n+1})_*\cong q^*$ (\cite[Prop.~1.20]{groth:ptstab}). But this follows from $q'\circ\iota=\mathrm{d}_{n+1}\circ q$, since the counit of the adjunction $(\mathrm{d}_{n+1})^*\dashv(\mathrm{d}_{n+1})_*$ is invertible (again \cite[Prop.~1.20]{groth:ptstab}).
\item For the third isomorphisms we note that morphisms of stable derivators (in particular the inverse image morphisms $(\rightarrow_{\tau})^*$) commute with cocones.
\item For the fourth isomorphism we show that the right triangle commutes. For this we observe that $p_{n+1}\circ\mathrm{d}_1^n=\mathrm{d}_{n+1}\circ p_n$. Again the counit of the adjunction $(\mathrm{d}_{n+1})^*\dashv(\mathrm{d}_{n+1})_*$ induces the desired isomorphism.
\item The fifth isomorphism is induced by the mutually inverse equivalences $(\rightarrow_{\tau})_!$ and $(\rightarrow_{\tau})^*$.
\end{itemize} 
\end{proof}

\begin{rmk}
In fact the statement of \autoref{thm:dual-Serre} is a special case of \autoref{thm:main} which will be proven independently. We decided to give an explicit proof at this point nevertheless, since \autoref{thm:dual-Serre} is the central motivation for defining the duality morphisms \autoref{defn:duality}. Moreover, \autoref{thm:dual-Serre} we be the the main ingredient for the construction of derivator Toda brackets in \S\ref{sec:Toda}.
\end{rmk}

\begin{rmk}\label{rmk:psi-prime}
We also define the equivalence $\Psi'_n\colon\D_{n,2}\toiso\D_{1,n+1}$ as the composition
\begin{equation}\label{eq:psi-prime}
\D_{n,2}\toiso sl\D_{n,2}\toiso \D^{[n]}\xrightarrow{\Psi_n^{\square}} \D^{\cube{n}}_{\rightarrow}\toiso sl\D_{1,n+1}\toiso  \D_{1,n+1}.
\end{equation}
Then the proof of \autoref{thm:dual-Serre} yields the relation $\Psi'_n\cong\mathsf{s}_3^*\circ\tilde{\Psi}_n$.
\end{rmk}

Both statements of \autoref{thm:dual-Serre} together imply that the equivalences $\tilde{\Psi}_n$ respect the values at $\mathsf{s}_3^i(\xi)$ up to some shift in $i\in\mathbb{Z}$. Because of this, we redefine $\tilde{\Psi}_n$ such that the compatibility above holds without any shift.

\begin{defn}\label{defn:duality}
Let $n\geq 1$ and $\D$ a stable derivator. The \textbf{duality morphism for the} $n$\textbf{-simplex} is the equivalence of derivators
\[
\Psi_n:=(\mathsf{s}_3^*)^{n+1}\circ\tilde{\Psi}_n\colon\D_{n,2}\toiso\D_{1,n+1}.
\]
\end{defn}

\begin{cor}
Let $n\geq 1$, $i\in\mathbb{Z}$ and $\D$ a stable derivator. Then there are natural isomorphisms
\begin{enumerate}
\item $\mathsf{s}_3^*\circ\Psi_n\cong\Psi_n\circ\mathsf{s}_3^*\colon\D_{n,2}\toiso\D_{1,n+1}$ and
\item $\xi^*\circ(\mathsf{s}_3^*)^i\circ\tilde{\Psi}_n\cong\xi^*\circ(\mathsf{s}_3^*)^i\colon\D_{n,2}\to\D$.
\end{enumerate}
\end{cor}

\begin{proof}
Both statements follow from straight-forward computations using \autoref{thm:dual-Serre} and $\Psi_n=(\mathsf{s}_3^*)^{n+1}\circ\tilde{\Psi}_n$.
\end{proof}

\begin{egs}
We illustrate the effect of the constructions in the proofs of \autoref{thm:dual-Sdot} and \autoref{thm:dual-Serre} in the cases $n=1$ and $n=2$. Let $\D$ be a stable derivator.
\begin{enumerate}
\item [(n=1)] We identify $\D_{1,2}\cong\D^{[1]}$ and that under this equivalence $\mathsf{s}_3^*\colon\D_{1,2}\toiso\D_{1,2}$ corresponds to $\cof\colon\D^{[1]}\toiso\D^{[1]}$. Let $x\in\D^{[1]}$ be an object with underlying diagram $x_0\xrightarrow{f}x_1$. Then $(\mathrm{d}_2)_*(x)$ looks like $x_0\xrightarrow{f}x_1\to 0$. Hence we can compute the underlying diagram of $x'=p_2^*\circ(\mathrm{d}_2)_*(x)$ as
\[
\xymatrix{
x_0 \ar[r]^{\id} \ar[d]_{f} & x_0 \ar[d]\\
x_1 \ar[r] & 0.
}
\]
To determine $x''=\fib^{\underline{1}}(x')$, we extend all morphisms in $x'$ to the left to fiber sequences
\[
\xymatrix{
\Omega x_1 \ar[r]\ar[d] & Ff \ar[r]\ar[d] & x_0 \ar[d]^{\id}\\
0 \ar[r]\ar[d] & x_0 \ar[r]^{\id} \ar[d]^{f} & x_0 \ar[d]\\
x_1 \ar[r]^{\id} & x_1 \ar[r] & 0.
}
\]
Here we find $x''$ as the upper left square, and hence $\tilde{\Psi}_1(x)$ as the upper left horizontal morphism. Furthermore, the middle vertical sequence exhibits the right vertical morphism of $x''$ as $(\mathsf{s}_3^*)^{-1}x$, such that the bicartesianess of $x''$ implies the $\tilde{\Psi}_1(x)=(\mathsf{s}_3^*)^{-2}(x)$. Since all the previous constructions are natural with respect to $x\in\D^{[1]}$, we obtain the important relation
\[
\Psi_1\cong\id\colon\D^{1,2}\to\D^{1,2}.
\]
\item [(n=2)] Similar to the previous case we identify $\D_{2,2}\cong\D^{[2]}$ and consider $x\in\D^{[2]}$ with underlying diagram $x_0\xrightarrow{f}x_1\xrightarrow{g}x_2$. Hence the underlying diagram of $x'=p_3^*\circ(\mathrm{d}_3)_*(x)$ looks like
\[
\xymatrix{
& x_1 \ar[rr]\ar[dd] & & x_1 \ar[dd]\\
x_0 \ar[rr]\ar[dd]\ar[ur] & & x_0 \ar[dd] \ar[ur]\\
& x_2 \ar[rr] & & 0\\
x_0 \ar[rr]\ar[ur] & & x_0. \ar[ur]
}
\]
To compute  $x''=\fib^{\underline{1}}(x')$, we again have to extend all morphisms in $x'$ to fiber sequences. As a preparation for this computation we look at the canonical extension of $x$ to an object of $\D_{2,2}$, whose underlying diagram
\[
\xymatrixrowsep{0.7cm}
\xymatrixcolsep{0.7cm}
\xymatrix{
0 \ar[r] & \Omega^2 x_0 \ar[r] \ar[d] & \Omega^2 x_1 \ar[r] \ar[d] & \Omega^2 x_2 \ar[r] \ar[d] & 0 \ar[d]\\
& 0 \ar[r] & \Omega Ff \ar[r] \ar[d] & \Omega Fgf \ar[r] \ar[d] & \Omega x_0 \ar[r] \ar[d] & 0 \ar[d]\\
& & 0 \ar[r] & \Omega Fg \ar[r] \ar[d] & \Omega x_1 \ar[r] \ar[d] & Ff \ar[r] \ar[d] & 0 \ar[d]\\
& & & 0 \ar[r] & \Omega x_2 \ar[r] \ar[d] & Fgf \ar[r] \ar[d] & Fg \ar[r] \ar[d] & 0 \ar[d]\\
& & & & 0 \ar[r] & x_0 \ar[r] & x_1 \ar[r] & x_2 \ar[r] & 0 
}
\]
encodes all required iterated fibers. We obtain $x''$, which is also the fundamental domain of $\tilde{\Psi}_2(x)$, as the front upper left cube in 
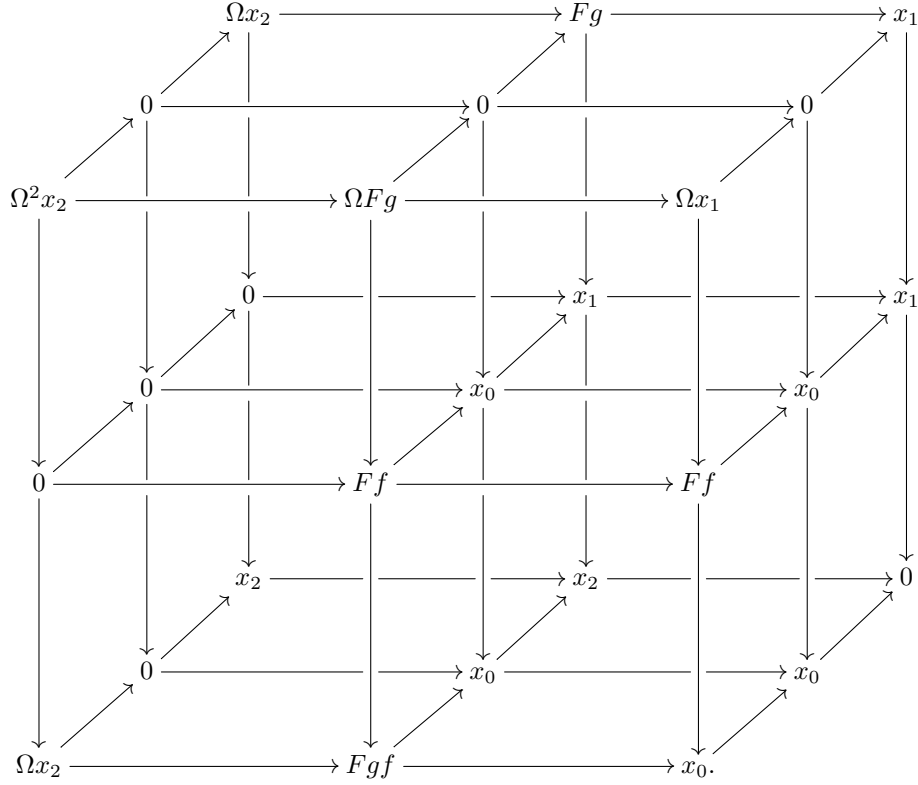
\[
\xymatrixrowsep{0.75cm}
\xymatrixcolsep{0.75cm}
\xymatrix{
& & \Omega x_2 \ar[rrr] \ar'[d]'[dd][ddd] & & & Fg \ar[rrr] \ar'[d]'[dd][ddd] & & & x_1 \ar[ddd]\\
& 0 \ar[ur] \ar[rrr] \ar'[d][ddd] & & & 0 \ar[ur] \ar[rrr] \ar'[d][ddd] & & & 0 \ar[ur] \ar[ddd]\\
\Omega^2x_2 \ar[ur] \ar[rrr] \ar[ddd] & & & \Omega Fg \ar[ur] \ar[rrr] \ar[ddd] & & & \Omega x_1 \ar[ur] \ar[ddd]\\
& & 0 \ar'[r]'[rr][rrr] \ar'[d]'[dd][ddd] & & & x_1 \ar'[r]'[rr][rrr] \ar'[d]'[dd][ddd] & & & x_1 \ar[ddd]\\
& 0 \ar[ur] \ar'[rr][rrr] \ar'[d][ddd] & & & x_0 \ar[ur] \ar'[rr][rrr] \ar'[d][ddd] & & & x_0 \ar[ur] \ar[ddd]\\
0 \ar[ur] \ar[rrr] \ar[ddd] & & & Ff \ar[ur] \ar[rrr] \ar[ddd] & & & Ff \ar[ur] \ar[ddd]\\
& & x_2 \ar'[r]'[rr][rrr] & & & x_2 \ar'[r]'[rr][rrr] & & & 0\\
& 0 \ar[ur] \ar'[rr][rrr] & & & x_0 \ar[ur] \ar'[rr][rrr] & & & x_0 \ar[ur]\\
\Omega x_2 \ar[ur] \ar[rrr] & & & Fgf \ar[ur] \ar[rrr] & & & x_0. \ar[ur]
}
\]
Moreover, the diagram above also allows us to check \autoref{thm:dual-Serre} (ii) in this case.
\begin{itemize}
\item The front square $(\mathrm{d}_1^0)^*(x'')$ of $x''$ is the fundamental slice of $\tilde{\Psi}_2(x)$, and the front $[2]\times[2]$-shaped face of the diagram exhibits $(\mathrm{d}_1^0)^*(x'')$ as 
\[
\fib^{\underline{1}}\circ(F\times\id_{\cube{2}})(x')\cong\fib^{\underline{1}}\circ p_2^*((\mathsf{s}_3^*)^{-1}(x))=\Psi^{\square}_2\circ(\mathsf{s}_3^*)^{-1}(x).
\]
The isomorphism above is in fact the key step of \autoref{thm:dual-Serre} (ii).
\item On the other hand the right square $(\mathrm{d}_0^2)^*(x'')$ of $x''$  is the fundamental slice of $\mathsf{s}_3^*\circ\tilde{\Psi}_2(x)$. But the middle slice of the $[2]\times[2]$-shaped diagram exhibits $(\mathrm{d}_0^2)^*(x'')$ as 
\[
\fib^{\underline{1}}\circ(\mathrm{d}_1^2)^*(x')=\fib^{\underline{1}}\circ p_2^*(x')=\Psi^{\square}_2(x).
\]
\item Hence the two points above together imply
\[
\mathsf{s}_3^*\circ\Psi^{\square}_2\circ(\mathsf{s}_3^*)^{-1}(x)\cong\Psi^{\square}_2(x).
\]
\end{itemize}
\end{enumerate}
\end{egs}

\begin{rmk}\label{rmk:dual-con}
The following result revisits the situation of \autoref{thm:dual-Sdot-sl} and will be used in \S\ref{sec:main}. More precisely, for a stable derivator $\D$ there are constructions completely dual to $\Psi^{\square}_n$ and $\Psi^{\square,ex}_n$ defined by the compositions
\[
\Psi^{\square\vee}_n\colon\D^{[n]}\xrightarrow{(\rightarrow_{\tau})_*}\D^{\cube{n},\kappa^{\vee}}\xrightarrow{\cof^{\underline{1}}}\D^{\cube{n}}_{\rightarrow}.
\]
and
\[
\Psi^{\square\vee,ex}_n\colon\D^{[n]}\xrightarrow{(\mathrm{d}_0)_!}\D^{[n+1],\emptyset}\xrightarrow{(\rightarrow_{\tau})_*}\D^{\cube{n+1},\kappa^{\vee},\emptyset}\xrightarrow{\cof^{\underline{1}}}\D^{\cube{n+1},ex}_{\rightarrow}.
\]
Here $\D^{\cube{n+1},\kappa^{\vee}}$ denotes the essential image of the right Kan extension morphism along $(\rightarrow_{\tau})\colon[n+1]\to\cube{n+1}$. Using exactly the dual arguments as before one shows that the compositions above consist of equivalences. We show that $\Psi^{\square,ex}_n$ and $\Psi^{\square\vee,ex}_n$ coincide up to shift, which is immediate from the following.
\end{rmk}

\begin{prop}\label{prop:kappa}
Let $n\geq 0$ and $\D$ a stable derivator. Then there is a natural isomorphism
\[
\Sigma\circ(\rightarrow_{\tau})_*\circ(\mathrm{d}_{0})_!\cong\cof^{\underline{1}}\circ(\rightarrow_{\tau})_!\circ(\mathrm{d}_{n+1})_*\colon\D^{[n]}\to\D^{\cube{n+1}}.
\]
\end{prop}

\begin{proof}
Since $\mathrm{d}_0$ and $\rightarrow_{\tau}$ are fully faithful functors it is sufficient to show 
\begin{itemize}
\item that the essential image of $\cof^{\underline{1}}\circ(\rightarrow_{\tau})_!\circ(\mathrm{d}_{n+1})_*$ is contained in the essential image of $(\rightarrow_{\tau})_*\circ(\mathrm{d}_{0})_!=\D^{\cube{n+1},\kappa^{\vee},\emptyset}$,
\item that there is an equivalence $G:=(\mathrm{d}_0)^*\circ(\rightarrow_{\tau})^*\circ\cof^{\underline{1}}\circ(\rightarrow_{\tau})_!\circ(\mathrm{d}_{n+1})_*\cong\Sigma$.
\end{itemize}
For the first point we observe that \autoref{rmk:elementary} implies
\[
\cof^{\underline{1}}\circ(\rightarrow_{\tau})_!\circ(\mathrm{d}_{n+1})_*=\Sigma^{n+1}\circ(\fib^{\underline{1}})^2\circ(\rightarrow_{\tau})_!\circ(\mathrm{d}_{n+1})_*
\]
and that by \autoref{rmk:dual-con} and \autoref{thm:dual-Sdot-sl}
\[
\D^{[n]}\xrightarrow{(\mathrm{d}_{n+1})_*}\D^{[n+1],\infty}\xrightarrow{(\rightarrow_{\tau})_!}\D^{\cube{n+1},\kappa,\infty}\xrightarrow{\fib^{\underline{1}}}\D^{\cube{n+1},ex}_{\rightarrow}\xrightarrow{\fib^{\underline{1}}}\D^{\cube{n+1},\kappa^{\vee},\emptyset}.
\]
For the second point we consider the functor $q:[n]\times[n]\to[n],(i,j)\mapsto\mathrm{min}\lbrace i,j\rbrace$ and observe that the functors $l,r\colon[n]\to[n]\times[n]$ defined by $l(i)=(i,n)$ and $r(i)=(n,i)$ are sections of $q$, in particular
\begin{equation}\label{eq:retract}
l^*\circ q^*\cong\id_{\D^{[n]}}\qquad\text{and}\qquad r^*\circ q^*\cong\id_{\D^{[n]}}.
\end{equation}
Since $G$ is a morphism of derivators, it commutes in particular with inverse images, hence
\[
l^*\circ(G\times\id^*)\circ q^*\cong(\id\times n)^*\circ(G\times\id^*)\circ q^*\cong G\circ(\id\times n)^*\circ q^*\cong G.
\]
By plugging in the definitions and using the 2-functoriality of inverse images we obtain
\begin{equation}\label{eq:eq1}
G\cong ((\rightarrow_{\tau})((-)+1),n)^*\circ (\cof^{\underline{1}}\times\id)\circ((\rightarrow_{\tau})_!\times\id)\circ((\mathrm{d}_{n+1})_*\times\id)\circ q^*.
\end{equation}
Consider the functors 
\begin{enumerate}
\item $\gamma_1\colon[n]\to\cube{n+1}\times[n],i\mapsto((\rightarrow_{\tau})(i+1),n)$,
\item $\gamma_2\colon[n]\to\cube{n+1}\times[n],i\mapsto((\rightarrow_{\tau})(i+1),i)$,
\item $\gamma_3\colon[n]\to\cube{n+1}\times[n],i\mapsto((\rightarrow_{\tau})(n+1),i)$,
\item $\gamma_4\colon[n]\to\cube{n+1}\times[n],i\mapsto((\rightarrow)(1),i)$,
\item $\gamma_5\colon[n]\to\cube{n+1}\times[n],i\mapsto((\rightarrow_{\tau})(n),i)$,
\end{enumerate}
and the unique natural transformations
\[
\alpha_1\colon\gamma_2\to\gamma_1,\qquad\alpha_2\colon\gamma_2\to\gamma_3\qquad\text{and}\qquad\alpha_3\colon\gamma_4\to\gamma_3.
\]
We claim that there are isomorphisms
\begin{align}
G\cong &\gamma_1^*\circ (\cof^{\underline{1}}\times\id)\circ((\rightarrow_{\tau})_!\times\id)\circ((\mathrm{d}_{n+1})_*\times\id)\circ q^*\\
\cong &\gamma_2^*\circ (\cof^{\underline{1}}\times\id)\circ((\rightarrow_{\tau})_!\times\id)\circ((\mathrm{d}_{n+1})_*\times\id)\circ q^*\\
\cong &\gamma_3^*\circ (\cof^{\underline{1}}\times\id)\circ((\rightarrow_{\tau})_!\times\id)\circ((\mathrm{d}_{n+1})_*\times\id)\circ q^*\\
\cong &\gamma_4^*\circ (\cof^{\underline{1}}\times\id)\circ((\rightarrow_{\tau})_!\times\id)\circ((\mathrm{d}_{n+1})_*\times\id)\circ q^*\\
\cong &\Sigma\circ\gamma_5^*\circ ((\rightarrow_{\tau})_!\times\id)\circ((\mathrm{d}_{n+1})_*\times\id)\circ q^*\\
\cong&\Sigma\circ r^*\circ q^*\\
\cong&\Sigma.
\end{align}
The single isomorphisms above are constructed as follows
\begin{itemize}
\item The first isomorphism is a reformulation of \eqref{eq:eq1}.
\item For the second isomorphism, we claim that 
\[
\tilde{\alpha_1}:=\alpha_1^*\circ (\cof^{\underline{1}}\times\id)\circ((\rightarrow_{\tau})_!\times\id)\circ((\mathrm{d}_{n+1})_*\times\id)\circ q^*
\]
is invertible. For this it is sufficient (using the axiom (Der 2)) to show that $i^*\circ\tilde{\alpha_1}$ is an isomorphism for all $i\in[n]$. In the following we fix $i\in[n]$. We consider the unique natural transformation $\beta_1\colon c_1\to c_2\colon\cube{i+1}\to\cube{n+1}\times[n]$ between the inclusions of $\cube{n+1}_{/(\rightarrow)(n-i)}$ at the coordinates $i$ and $n$, respectively. Hence we have
\[
i^*\circ\tilde{\alpha_1}=\tcof\circ\beta_1^*\circ((\rightarrow_{\tau})_!\times\id)\circ((\mathrm{d}_{n+1})_*\times\id)\circ q^*.
\]
Let $\tilde{\beta_1}:=\beta_1^*\circ((\rightarrow_{\tau})_!\times\id)\circ((\mathrm{d}_{n+1})_*\times\id)\circ q^*$ and $x\in\cube{n+1}_{/(\rightarrow)(n-i)}$. We claim that $x^*\circ\tilde{\beta_1}$ is an isomorphism. This would imply that $\tilde{\beta_1}$ is an isomorphism by (Der 2), and hence also that $\tilde{\alpha_1}$ is an isomorphism. If $x=\infty$ then $x^*\circ\tilde{\beta_1}$ is the identity on the zero object (because of the right Kan extension along $\mathrm{d}_{n+1}$). If $x\neq\infty$ then $p(x)\leq i$. Denoting by $\delta$ the natural transformation comparing the inclusions of $(p(x),i)$ and $(p(x),n)$ into $[n]\times[n]$. Using the 2-functoriality of $\D$, we conclude that $x^*\circ\tilde{\beta_1}=\delta^*\circ q^*$. But $q\circ\delta=\id_{p(x)}$. Invoking the 2-functoriality of $\D$ again $x^*\circ\tilde{\beta_1}$ is seen to be an isomorphism.
\item  For the third isomorphism, we claim that 
\[
\tilde{\alpha_2}:=\alpha_2^*\circ (\cof^{\underline{1}}\times\id)\circ((\rightarrow_{\tau})_!\times\id)\circ((\mathrm{d}_{n+1})_*\times\id)\circ q^*
\]
is invertible. Again it is sufficient (using the axiom (Der 2)) to show that $i^*\circ\tilde{\alpha_2}$ is an isomorphism for all $i\in[n]$. In the following we fix $i\in[n]$. Let $d\colon\cube{n-i}\to\cube{n+1}\times[n]$ be the inclusion of $\cube{n+1}_{/(\rightarrow_{\tau})(i+1)}\times\lbrace i\rbrace$. We claim that 
\[
D:=d^*\circ(\cof^{\underline{1}}\times\id)\circ((\rightarrow_{\tau})_!\times\id)\circ((\mathrm{d}_{n+1})_*\times\id)\circ q^*
\]
is constant. This would imply that $i^*\circ\tilde{\alpha_2}$, which is the diagonal in the underlying diagram of $D$, is an isomorphism. By \cite[Cor.~9.8]{bg:cubical} it is sufficient to show that for all $i+1\leq j\leq n$ the cocone $F^j\circ D$ in the direction of the $j$th coordinate vanishes. For this we consider a subset $M\subseteq\lbrace i+1,\cdots,n\rbrace\setminus\lbrace j\rbrace$. We have to show that 
\[
M^*\circ F^j\circ D\cong 0
\]
for all such M. Let $M^{\vee}=\mathbf{n+1}\setminus (M\cup\lbrace j\rbrace)$. We use the notation
\[
\mathrm{d}_{\varepsilon}^N:=\prod_{n\in N}\mathrm{d}_{\varepsilon}^n
\]
for $N\subseteq\mathbf{n+1}$ and $\varepsilon\in\lbrace 0,1\rbrace$. Then by \autoref{rmk:elementary}
\begin{align}
&M^*\circ F^j\circ D\\
\cong&(((\mathrm{d}_1^M)^*\times(\mathrm{d}_0^{M^{\vee}})^*\times F^j)\times i^*)\circ(\cof^{\underline{1}}\times\id)\circ((\rightarrow_{\tau})_!\times\id)\circ((\mathrm{d}_{n+1})_*\times\id)\circ q^*\\
\cong&((\mathrm{d}_0^M)^*\times C^{M^{\vee}}\times(\mathrm{d}_1^j)^*)\circ(\rightarrow_{\tau})_!\circ(\mathrm{d}_{n+1})_*\circ(\id\times i)^*\circ q^*
\end{align}
is exihibited as a total cofiber of the cube 
\[
E:=((\mathrm{d}_0^M)^*\times \id^{M^{\vee}}\times(\mathrm{d}_1^j)^*)\circ(\rightarrow_{\tau})_!\circ(\mathrm{d}_{n+1})_*\circ(\id\times i)^*\circ q^*
\]
parametrized by the coordinates $M^{\vee}$. We observe that $i\in M^{\vee}$ and claim that $C^i\circ E\cong 0$. This would imply by \cite[Cor.~9.8]{bg:cubical} that $\tcof\circ E\cong 0$ and hence that $\tilde{\alpha_2}$ is invertible. We consider $M'\subseteq M^{\vee}\setminus\lbrace i\rbrace$ and compute $M'^*\circ C^i\circ E$. We observe that 
\[
M'^*\circ C^i\circ E\cong C\circ t_1^*\circ (\rightarrow_{\tau})_!\circ(\mathrm{d}_{n+1})_*\circ(\id\times i)^*\circ q^*
\]
for some map $t_1\colon[1]\to\cube{n+1}$ satisfying the following properties
\begin{enumerate}
\item $t_1(0)_i=0$,
\item $t_1(1)_i=1$,
\item $t_1(0)_l=t_1(1)_l$ for $l\in\mathbf{n+1}\setminus\lbrace i\rbrace$,
\item $t_1(0)_j=t_1(1)_j=0$.
\end{enumerate}
From the fourth property we deduce that $t_1(1)\neq\infty$. Therefore, we obtain an isomorphism
\[
t_1^*\circ (\rightarrow_{\tau})_!\circ(\mathrm{d}_{n+1})_*\circ(\id\times i)^*\circ q^*\cong t_2^*\circ q^*
\]
where $t_2\colon[1]\to[n]\times[n]$ is the map classifying $(p(t_1(0)),i)\to(p(t_1(1)),i)$. It is sufficient to show that $q\circ t_2$ is constant, but this holds by construction. More precisely, property (iii) above implies that if $p(t_1(0))\neq p(t_1(1))$ then
\[
p(t_1(1))>p(t_1(0))=i.
\]
\item For the fourth isomorphism, we claim that 
\[
\tilde{\alpha_3}:=\alpha_3^*\circ (\cof^{\underline{1}}\times\id)\circ((\rightarrow_{\tau})_!\times\id)\circ((\mathrm{d}_{n+1})_*\times\id)\circ q^*
\]
is invertible. Again it is sufficient (using the axiom (Der 2)) to show that $i^*\circ\tilde{\alpha_3}$ is an isomorphism for all $i\in[n]$. In the following we fix $i\in[n]$. Let $d'\colon\cube{n}\to\cube{n+1}\times[n]$ be the inclusion of $\cube{n+1}_{/(\rightarrow)(1)}\times\lbrace i\rbrace$. Then 
\begin{align}
D':=&d'^*\circ(\cof^{\underline{1}}\times\id)\circ((\rightarrow_{\tau})_!\times\id)\circ((\mathrm{d}_{n+1})_*\times\id)\circ q^*\\
\cong&(\mathrm{d}_0^{n+1})^*\circ\cof^{\underline{1}}\circ(\rightarrow_{\tau})_!\circ(\mathrm{d}_{n+1})_*\circ(\id\times i)^*\circ q^*
\end{align}
is constant (using exactly the dual argument as in step (ii) of the proof of \autoref{thm:dual-Sdot-sl}), because the essential image of
\[
\cof^{\underline{1}}\circ(\rightarrow_{\tau})_!(\mathrm{d}_{n+1})_*\circ(\id\times i)^*\circ q^*
\]
is by the first part of this proof contained in $\D^{\cube{n+1},\kappa^{\vee},\emptyset}$. In particular, $i^*\circ\tilde{\alpha_3}$, which is the diagonal in the underlying diagram of $D'$, is an isomorphism.
\item We use the notation $\D':=\D^{[n]}$. Moreover, let $u\colon[1]\to\cube{n+1}$ be the map classifying $(\rightarrow_{\tau})(n)\to\infty$. We obtain by \autoref{rmk:elementary} an isomorphism
\[
C\circ u^*\cong ((\rightarrow)(1))^*\circ\cof^{\underline{1}}\colon\D'^{\cube{n+1}}\to\D'.
\]
We invoke again \autoref{rmk:elementary} to identify the restriction of this isomorphism to $\D'^{\cube{n+1},\infty}$
\[
\Sigma\circ((\rightarrow_{\tau})(n))^*\cong C\circ u^*\cong ((\rightarrow)(1))^*\circ\cof^{\underline{1}}\colon\D'^{\cube{n+1},\infty}\to\D'.
\]
By using the canonical identification $\D'^{\cube{n+1}}\cong\D^{\cube{n+1}\times[n]}$, we obtain the fifth isomorphism by appropriate restriction.
\item The sixth isomorphism is induced by $(p\times\id)\circ\gamma_5=r$.
\item The seventh isomorphism is \eqref{eq:retract}.
\end{itemize}
\end{proof}

\section{The symmetry theorem - The general case}
\label{sec:main}

The main objective of this section is the construction of well-behaved equivalences $\D_{n,k}\to\D_{k-1,n+1}$ (\autoref{thm:main}). For this we recall from \autoref{prop:square_I} the vanishing of the composition $(d^h)^*\circ d^v_!\colon\D_{n,k+1}\to\D_{n+1,k}$ induced by the sieve $d^h$ and the (almost) complementary cosieve $d^v$. By considering also the left and right adjoints, one can extend the composition above to a recollement
\[
\xymatrix{
\D_{n,k+1} \ar@{<-}@<1.5mm>[r] \ar@{<-}@<-1.5mm>[r] \ar[r] & \D_{n+1,k+1} \ar@{<-}@<1.5mm>[r] \ar@{<-}@<-1.5mm>[r] \ar[r] & \D_{n+1,k}.
}
\] 
The central part of the proof of \autoref{thm:main} can be summarized as follows.
\begin{enumerate}
\item The first step consists of the construction of morphisms of recollements
\begin{equation}\label{eq:central}
\xymatrix{
\D_{n,k+1} \ar@{<-}@<1.5mm>[r] \ar@{<-}@<-1.5mm>[r] \ar[r]\ar[d] & \D_{n+1,k+1} \ar@{<-}@<1.5mm>[r] \ar@{<-}@<-1.5mm>[r] \ar[r] \ar[d]& \D_{n+1,k} \ar[d]\\
\D_{k,n+1} \ar@{<-}@<1.5mm>[r] \ar@{<-}@<-1.5mm>[r] \ar[r] & \D_{k,n+2} \ar@{<-}@<1.5mm>[r] \ar@{<-}@<-1.5mm>[r] \ar[r] & \D_{k-1,n+2}
}
\end{equation}
for $n\geq 1,k\geq 2$, such that the left (resp. right) vertical morphisms are in the case $n=1$ (resp. $k=2$) the equivalences induced by \autoref{thm:dual-Sdot-sl} (resp. \autoref{rmk:dual-con}). 
\item In the second step we show that, if the two outer vertical morphisms of \eqref{eq:central} are equivalences then so is the inner vertical morphism.
\end{enumerate}
The second step will rely on the following alternative description of the middle term in a recollement induced by a sieve-cosieve-pair. 

\begin{lem}\label{lem:recollement}
Let $\D$ a stable derivator and $u\colon A\to B$ be a sieve with complementary cosieve $v\colon A' \to B$. Let $\D^{B\times [1]}_{A,A'}$ be the full subderivator on those objects $x$, such that $(\id\times\mathrm{d}_1)^*(x)\in\D^{B,A'}$ and $(\id\times\mathrm{d}_0)^*(x)\in\D^{B,A}$. Then there is an equivalence of derivators
\[
\D^B\toiso\D^{B\times [1]}_{A,A'}.
\]
\end{lem}

\begin{proof}
Let $C=(B\times [1])\setminus(A'\times\lbrace 1\rbrace)$, $i\colon B\to C$ the inclusion of $B\times \lbrace 0\rbrace$ and $p\colon C\to B$ the projection. Then $i \dashv p$ and the essential image of $i_!=p*$ consists precisely of those objects on which the counit of $i_!\dashv i^*$ is invertible. This latter condition holds exactly when $(a\times\id)^*$ is constant for all $a\in A$ and we denote by $\D^C_0\subseteq\D^C$ the subderivator of all those objects. Next, we consider the inclusion $j\colon C\to B\times [1]$, which is clearly a sieve. We define $C'=(B\times [1])\setminus C$ invoke \cite[Prop.~1.23]{groth:ptstab} for the equivalence $j_*\colon\D^C\toiso\D^{B\times [1],C'}$. We denote by $\D^{B \times[1]}_0$ the essential image of $\D^C_0$ under this equivalence (which consists by \cite[Prop.~1.23]{groth:ptstab} precisely of those objects such that the restriction to is in $\D^C_0$ and the restriction to $C'$ vanishes). We claim the equivalence $\id\times\cof\colon\D^{B\times [1]}\toiso\D^{B\times [1]}$ restricts to an equivalence $\D^{B \times[1]}_0\toiso\D^{B\times [1]}_{A,A'}$. To see this, we invoke \cite[Lem.~8.19]{bg:cubical} which shows that the vanishing on $C'$ exactly corresponds to $(\id\times\mathrm{d}_1)^*(x)\in\D^{B,A'}$ and \cite[Cor.~8.6]{bg:cubical} which shows that the constantness condition for objects in $A$ exactly corresponds to $(\id\times\mathrm{d}_0)^*(x)\in\D^{B,A}$. Summing up, there is a chain of equivalences
\[
\D^B\xrightarrow{i_!}\D^C_0\xrightarrow{j_*}\D^{B\times [1]}_0\xrightarrow{\id\times\cof}\D^{B\times [1]}_{A,A'}.
\]
\end{proof}

\begin{rmk}\label{rmk:inverse}
We observe that the inverse of the equivalence in \autoref{lem:recollement} is induced by the cocone morphism
\[
\id^*\times F\colon\D^{B\times[1]}\to\D^B.
\]
\end{rmk}

\begin{con}\label{con:sl-cube}
Recall the cubical slice $Sl^{\square}_{n,k}=\mathbb{Z}^k_{\xi/\mathsf{s}_3^{k-1}(\xi)}$ (c.f. \autoref{rmk:slnk-hat}). We consider the isomorphism of categories
\[
c_{n,k}\colon[n]^{k-1}\to Sl^{\square}_{n,k},(f_0,\cdots,f_{k-2})\mapsto(0,f_0+1,\cdots,f_{k-2}+k+1).
\]
Let $\D_0^{[n]^{k-1}}$ denote the essential image of $c_{n,k}^*\colon sl\hat{\D}_{n,k}\to\D^{[n]^{k-1}}$. Furthermore, we recall the functors
\begin{itemize}
\item $d^v\colon Sl_{n,k}\to Sl_{n+1,k},(0,f_1,\cdots,f_{k-1})\mapsto(0,f_0+1,\cdots,f_{k-1}+1)$,
\item $d^h\colon Sl_{n,k}\to Sl_{n,k+1},(0,f_1,\cdots,f_{k-1})\mapsto(0,1,f_0+1,\cdots,f_{k-1}+1)$.
\end{itemize}
These assignments extend formally to functors $d^v\colon Sl^{\square}_{n,k}\to Sl^{\square}_{n+1,k}$ and $d^h\colon Sl^{\square}_{n,k}\to Sl^{\square}_{n,k+1}$ such that there are (strictly) commutative diagrams
\[
\resizebox{.95\hsize}{!}{$
\xymatrix{
Sl_{n,k}\ar[r]^{\mathsf{j}}\ar[d]_{d^v}&Sl^{\square}_{n,k}\ar[r]^{c_{n,k}}\ar[d]_{d^v}&[n]^{k-1}\ar[d]^{\mathrm{d}_0^{k-1}} &Sl_{n,k}\ar[r]^{\mathsf{j}}\ar[d]_{d^h}&Sl^{\square}_{n,k}\ar[r]^{c_{n,k}}\ar[d]_{d^h}&[n]^{k-1}\ar[d]^{\emptyset\times\id_{[n]^{k-1}}}\\
Sl_{n+1,k}\ar[r]^{\mathsf{j}}&Sl^{\square}_{n+1,k}\ar[r]^{c_{n+1,k}}&[n+1]^{k-1} &Sl_{n,k+1}\ar[r]^{\mathsf{j}}&Sl^{\square}_{n,k+1}\ar[r]^{c_{n,k+1}}&[n]^{k}.
}$}
\]
By passing to inverse images and suitable restrictions we obtain strictly commutative squares of derivators
\begin{equation}\label{eq:dvdh}
\xymatrix{
sl\D_{n,k}&\D_0^{[n]^{k-1}}\ar[l] &sl\D_{n,k}&\D_0^{[n]^{k-1}}\ar[l]\\
sl\D_{n+1,k}\ar[u]^{(d^v)^*}&\D_0^{[n+1]^{k-1}}\ar[u]_{(\mathrm{d}_0^{k-1})^{*}}\ar[l], &sl\D_{n,k+1}\ar[u]^{(d^h)^*}&\D_0^{[n]^{k}}\ar[u]_{(\emptyset\times\id_{[n]^{k-1}})^{*}}\ar[l].
}
\end{equation}
Using the same argument as in \autoref{prop:the-sieve} and \autoref{prop:dh-ext-zero} we see that the adjunctions 
\[
\resizebox{.95\hsize}{!}{$
(\mathrm{d}_0^{k-1})_{!}\colon\D^{[n]^{k-1}}\rightleftarrows\D^{[n+1]^{k-1}}\colon(\mathrm{d}_0^{k-1})^{*}\text{ and }(\emptyset\times\id_{[n]^{k-1}})^{*}\colon\D^{[n]^{k}}\rightleftarrows\D^{[n]^{k-1}}\colon(\emptyset\times\id_{[n]^{k-1}})_{*}
$}
\]
restrict to adjunctions
\[
\resizebox{.95\hsize}{!}{$
(\mathrm{d}_0^{k-1})_{!}\colon\D_0^{[n]^{k-1}}\rightleftarrows\D_0^{[n+1]^{k-1}}\colon(\mathrm{d}_0^{k-1})^{*}\text{ and }(\emptyset\times\id_{[n]^{k-1}})^{*}\colon\D_0^{[n]^{k}}\rightleftarrows\D_0^{[n]^{k-1}}\colon(\emptyset\times\id_{[n]^{k-1}})_{*}.
$}
\]
The analogue of this is not true for $(\mathrm{d}_0^{k-1})_{*}$ and $(\emptyset\times\id_{[n]^{k-1}})_{!}$. Instead we invoke \autoref{lem:adjoint-lem} and its dual to see the following.
\begin{itemize}
\item The right adjoint of $(\mathrm{d}_0^{k-1})^{*}\colon\D_0^{[n+1]^{k-1}}\to\D_0^{[n]^{k-1}}$ is given by 
\[
(\mathrm{d}_0)_{*}\times(\mathrm{d}_0^{k-2})_{!}\colon\D_0^{[n]^{k-1}}\to\D_0^{[n+1]^{k-1}}.
\]
\item The essential image of the left adjoint of $(\emptyset\times\id_{[n]^{k-1}})^*\colon\D_0^{[n]^{k}}\to\D_0^{[n]^{k-1}}$ is given by those objects in $\D_0^{[n]^{k}}$ satisfying property (P1) at those objects in $[n]^{k}$ which are increasing sequences in $[n-1]^k\subset[n]^k$ (c.f. \autoref{prop:dh-domain}).
\end{itemize}
\end{con}

Next, we apply \autoref{lem:recollement} to the situation of \autoref{con:sl-cube}.

\begin{cor}\label{cor:recollement}
Let $n\geq 1,k\geq2$ and $\D$ a stable derivator. Let $(sl\hat{\D}_{n+1,k+1})^{[1]}_{A,A'}$ be the full subderivator of $\D^{(Sl^{\square}_{n+1,k+1}\times[1])}$ on those objects $x$, such that $(\id\times\mathrm{d}_1)^*(x)$ is contained in the essential image of $d^h_*\colon sl\hat{\D}_{n+1,k}\to sl\hat{\D}_{n+1,k+1}$ and $(\id\times\mathrm{d}_0)^*(x)$ is contained in the essential image of $d^v_!\colon sl\hat{\D}_{n,k+1}\to sl\hat{\D}_{n+1,k+1}$. Then there is an equivalence of derivators
\[
sl\hat{\D}_{n+1,k+1}\cong(sl\hat{\D}_{n+1,k+1})^{[1]}_{A,A'}
\]
\end{cor}

\begin{proof}
We apply \autoref{lem:recollement} to the case $B=Sl^{\square}_{n+1,k+1}$ and the sieve $u=d^h\colon A=Sl^{\square}_{n+1,k}$ with complementary cosieve $v\colon A'\to B$ and conclude by identifying the essential image of the composition
\[
sl\hat{\D}_{n+1,k+1}\to\D^{Sl^{\square}_{n+1,k+1}}\toiso\D^{(Sl^{\square}_{n+1,k+1}\times[1])}_{A,A'},
\]
where the equivalence is induced by \autoref{lem:recollement}, with $(sl\hat{\D}_{n+1,k+1})^{[1]}_{A,A'}$.
\end{proof}

We note that there is also a dual picture. For the sieve
\[
d^{v\vee}=\underline{\Lambda}(\Lambda_{k},\mathrm{d}_{n+k+1})\vert_{Sl_{n,k+1}}\colon Sl_{n,k+1}\to Sl_{n+1,k+1},(0,f_1,\cdots,f_k)\mapsto(0,f_1,\cdots,f_k)
\]
we can consider the almost complementary cosieve
\[
d^{h\vee}\colon Sl_{n+1,k}\to Sl_{n+1,k+1}, (0,f_1,\cdots,f_{k-1})\mapsto(0,f_1,\cdots,f_{k-1},n+k+1).
\]
\autoref{eg:ver-fd} yields the identification $(d^{v\vee})^*\cong\mathsf{d}^v[2n]$. In the following we establish a similar identification for $(d^{h\vee})^*$.

\begin{prop}\label{prop:dhvee}
Let $n \geq 1, k\geq 2$ and $\D$ a stable derivator. Then there is a natural isomorphism
\[
(d^{h\vee})^*\cong\Sigma^n\circ\mathsf{d}^h[2(n+k)]\colon sl\D_{n,k+1} \to sl\D_{n,k}.
\]
\end{prop}

\begin{proof}
We note that $d^{h\vee}=\mathsf{s}_3\circ d^h\colon Sl_{n,k}\to \underline{\Lambda}(\Lambda_{k},\Lambda_{n+k})$, which gives rise to first isomorphism below
\begin{align}
(d^{h\vee})^*\cong & \mathsf{d}^h\circ\mathsf{s}_3^*\\
\cong & (\mathsf{s}_3^*)^{-n-k}\circ \mathsf{d}^h[2(n+k)]\circ (\mathsf{s}_3^*)^{n+k}\circ\mathsf{s}_3^*\\
\cong & \Sigma^{-n(k-1)}\circ \mathsf{d}^h[2(n+k)]\circ \Sigma^{nk}\\
\cong& \Sigma^n\circ \mathsf{d}^h[2(n+k)].
\end{align}
The second isomorphism is \autoref{cor:fracCY}, and the third is \autoref{cor:gsd-dual}.
\end{proof}

\begin{rmk}\label{rmk:dual-recoll}
We observe, the pair of morphisms $(d^{v\vee},d^{h\vee})$ behaves completely dually to $(d^h,d^v)$. 
\begin{enumerate}
\item The commutativity of the squares
\begin{equation}\label{eq:dvdhvee}
\xymatrix{
sl\D_{n,k}&\D_0^{[n]^{k-1}}\ar[l] &sl\D_{n,k}&\D_0^{[n]^{k-1}}\ar[l]\\
sl\D_{n+1,k}\ar[u]^{(d^{v\vee})^*}&\D_0^{[n+1]^{k-1}}\ar[u]_{(\mathrm{d}_n^{k-1})^{*}}\ar[l], &sl\D_{n,k+1}\ar[u]^{(d^{h\vee})^*}&\D_0^{[n]^{k}}\ar[u]_{(\id_{[n]^{k-1}}\times\infty)^{*}}\ar[l].
}
\end{equation}
is dual to \eqref{eq:dvdh}.
\item There is an equivalence dual to the one of \autoref{cor:recollement}. For $n\geq 1,k\geq2$ and $(sl\hat{\D}_{n+1,k+1})^{[1]}_{A^{\vee},A^{\vee'}}$ the full subderivator of $\D^{(Sl^{\square}_{n+1,k+1}\times[1])}$ on those objects $x$, such that $(\id\times\mathrm{d}_1)^*(x)$ is contained in the essential image of $d^{v\vee}_*\colon sl\hat{\D}_{n+1,k}\to sl\hat{\D}_{n+1,k+1}$ and $(\id\times\mathrm{d}_0)^*(x)$ is contained in the essential image of $d^{h\vee}_!\colon sl\hat{\D}_{n,k+1}\to sl\hat{\D}_{n+1,k+1}$. Then there is an equivalence of derivators
\[
sl\hat{\D}_{n+1,k+1}\cong(sl\hat{\D}_{n+1,k+1})^{[1]}_{A^{\vee},A^{\vee'}}.
\]
\end{enumerate}
\end{rmk}

\begin{thm}\label{thm:main}
Let $n\geq 1,k\geq 2$, $a\in\mathbb{Z}$ and $\D$ a stable derivator. Then the are equivalences of derivators
\[
\Phi_{n,k}\colon\D_{n,k}\toiso\D_{k-1,n+1},
\]
such that the following properties are satisfied

\begin{enumerate}
\item $\Phi_{n,k}\circ\Phi_{k-1,n+1}\cong\id$,
\item $\mathsf{s}_3^*\circ\Phi_{n,k}\cong\Phi_{n,k}\circ\mathsf{s}_3^*$,
\item $\mathsf{d}_{k-1,n+1}^v[a]\circ\Phi_{n,k+1}\cong\Phi_{n,k}\circ\mathsf{d}_{n,k}^h[a]$ for $a\in\mathbb{Z}$ even,\newline $\mathsf{d}_{k-1,n+1}^v[a]\circ\Phi_{n,k}\cong\Phi_{n+1,k}\circ\mathsf{d}_{n,k}^h[a]$ for $a\in\mathbb{Z}$ odd,
\item $\mathsf{d}_{k-1,n+1}^h[a]\circ\Phi_{n+1,k}\cong\Phi_{n,k}\circ\mathsf{d}_{n,k}^v[a]$ for $a\in\mathbb{Z}$ even, \newline $\mathsf{d}_{k-1,n+1}^h[a]\circ\Phi_{n,k}\cong\Phi_{n+1,k}\circ\mathsf{d}_{n,k}^v[a]$ for $a\in\mathbb{Z}$ odd,
\item $\xi^*\Phi_{n,k}\cong\xi^*$,
\item $\Phi_{n,2}=\Psi_n$.
\end{enumerate}

\end{thm}

\begin{proof}
By \autoref{thm:slices}, \autoref{prop:dnk-hat} and \autoref{rmk:slnk-hat} it is sufficient to show the corresponding statements for the derivators $sl\hat{\D}_{n,k}$.

\noindent We define the map $ad_{n,k}\colon[k-1]^n\to(\cube{n})^{k-1}$ by the assignment
\[
(i_0\cdots,i_{n-1})\mapsto(\lbrace j\in\mathbf{n}\vert 1\leq i_j\rbrace,\lbrace j\in\mathbf{n}\vert 2\leq i_j\rbrace,\cdots,\lbrace j\in\mathbf{n}\vert k-1\leq i_j\rbrace).
\]
We define the morphism of derivators $\Psi_{n,k}\colon sl\hat{\D}_{n,k}\to sl\hat{\D}_{k-1,n+1}$ as the composition
\[
sl\hat{\D}_{n,k}\xrightarrow{c_{n,k}^*}\D_0^{[n]^{k-1}}\xrightarrow{(\Psi_n^{\square})^{k-1}}\D^{(\cube{n})^{k-1}}\xrightarrow{ad_{n,k}^*}\D_0^{[k-1]^n}\xrightarrow{(c_{k-1,n+1}^{-1})^*}sl\hat{\D}_{k-1,n+1}.
\]
First we show that $\Psi_{n,k}$ is well defined. Consider a non-injective object $f=(0,f_1,\cdots,f_n)\in Sl^{\square}_{k-1,n+1}$ and  let $g=(g_0,\cdots,g_{n-1}):=c_{k-1,n+1}^{-1}(f)\in[k-1]^n$. The non-injectivity of $f$ implies that there is $0\leq i\leq n-2$ such that $j=g_i>g_{i+1}$. As a consequence, the $j$th coordinate of $ad_{n,k}(g)$ is not in the image of the standard maximal path $\rightarrow\colon[n]\to\cube{n}$. Therefore, 
\[
f^*\circ\Psi_{n,k}=(ad_{n,k}(g))^*\circ(\Psi_n^{\square})^{k-1}\circ c_{n,k}^*=0
\]
holds by construction of $\Psi_n^{\square}$. This establishes property (P2) for $f$ on objects in the essential image of $\Psi_{n,k}$.

\noindent We proceed in the following steps.
\begin{enumerate}
\item In this part we establish an alternative description of $\Psi^{\square}_{n,k}:=ad_{n,k}^*\circ(\Psi_n^{\square})^{k-1}$ and construct adjoint morphisms. For this we consider the shuffle permutation $sh_{n,k}\colon(\cube{n})^k\to(\cube{k})^n$ defined by the assignment
\[
((l_0^0,\cdots,l_{n-1}^0),\cdots,(l_0^{k-1},\cdots,l_{n-1}^{k-1}))\mapsto((l_0^0,\cdots,l_0^{k-1}),\cdots,(l_{n-1}^0,\cdots,l_{n-1}^{k-1}))
\]
and note that $(sh_{n,k})^{-1}=sh_{k,n}$. We claim that $sh_{n,k-1}\circ\mathrm{ad}_{n,k}=(\rightarrow_{\tau})^n$.
To establish the claim we consider $i=(i_0,\cdots,i_{n-1})\in[k-1]^n$ and let $sh_{n,k-1}\circ\mathrm{ad}_{n,k}(i)$ be of the following form
\[
((l_0^0,\cdots,l_{k-2}^0),\cdots,(l_0^{n-1},\cdots,l_{k-1}^{n-1}))\in(\cube{k-1})^n.
\]
Consider $j\in\mathbf{n}$. Then by the definition of $\mathrm{ad}_{n,k}$ we have $l^j_{\alpha}=1$ if and only if $\alpha-1\leq i_j$. Hence $(sh_{n,k-1}\circ\mathrm{ad}_{n,k}(i))_j=(l_0^j,\cdots,l_{k-2}^j)=(\rightarrow_{\tau})(i_j)$ as claimed.

As a consequence, we obtain the first isomorphism in
\begin{align}
\Psi^{\square}_{n,k}=&\mathrm{ad}_{n,k}^*\circ(\Psi_n^{\square})^{k-1}\\
\cong&((\rightarrow_{\tau})^n)^*\circ sh_{k-1,n}^*\circ(\Psi_n^{\square})^{k-1}\label{eq:psink}\\
\cong&((\rightarrow_{\tau})^n)^*\circ sh_{k-1,n}^*\circ\fib^{\underline{1}}\circ((\rightarrow_{\tau})^{k-1})_!.
\end{align}
For the second isomorphism above we use that $(\rightarrow_{\tau})_!$ is a morphism of derivators, and therefore commutes with fibers in unrelated coordinates. In particular, we see that $\Psi^{\square}_{n,k}$ is left adjoint to
\[
\Psi^{\square\vee}_{n,k}:=((\rightarrow_{\tau})^{k-1})^*\circ\cof^{\underline{1}}\circ sh_{n,k-1}^*\circ((\rightarrow_{\tau})^{n})_*.
\]
\item We observe that $sh_{n,1}\colon\cube{n}\to\cube{n}$ and $(\rightarrow_\tau)\colon[1]\to[1]$ are both identities. Hence we obtain in the special cases $k=2$ and $n=1$ the identifications
\begin{equation}\label{eq:main-boundary}
\Psi_{n,2}^{\square}\cong\Psi_n^{\square}\qquad\text{and}\qquad\Psi^{\square\vee}_{1,k}\cong\Psi_{k-1}^{\square\vee},
\end{equation}
respectively. These are equivalences of derivators by \autoref{thm:dual-Sdot-sl} and \autoref{rmk:dual-con}.
\item In this step we show that $\Psi_{n,k}$ is an equivalence of derivators under the assumption that the diagrams
\begin{equation}\label{eq:main-comm-I}
\xymatrix{
sl\hat{\D}_{n+1,k}\ar[r]^{(d^h)_*}\ar[d]_{\Psi_{n+1,k}}&sl\hat{\D}_{n+1,k+1}\ar[d]_{\Psi_{n+1,k+1}}&sl\hat{\D}_{n,k+1}\ar[l]_{(d^v)_!}\ar[d]^{\Psi_{n,k+1}\circ\Omega^k}\\
sl\hat{\D}_{k-1,n+2}\ar[r]^{(d^v)_*}\ar[d]_{\Psi_{n+1,k}^{\vee}}&sl\hat{\D}_{k,n+2}\ar[d]_{\Psi_{n+1,k+1}^{\vee}}&sl\hat{\D}_{k,n+1}\ar[l]_{(d^h)_*}\ar[d]^{\Sigma^k\circ\Psi_{n,k+1}^{\vee}}\\
sl\hat{\D}_{n+1,k}\ar[r]^{(d^h)_*}&sl\hat{\D}_{n+1,k+1}&sl\hat{\D}_{n,k+1}\ar[l]_{(d^v)_!}
}
\end{equation}
and
\begin{equation}\label{eq:main-comm-II}
\xymatrix{
sl\hat{\D}_{k,n+1}\ar[r]^{(d^h)_*}\ar[d]_{\Sigma^k\circ\Psi_{n,k+1}^{\vee}}&sl\hat{\D}_{k,n+2}\ar[d]_{\Psi_{n+1,k+1}^{\vee}}&sl\hat{\D}_{k-1,n+2}\ar[l]_{(d^v)_!}\ar[d]^{\Psi_{n+1,k}^{\vee}}\\
sl\hat{\D}_{n,k+1}\ar[r]^{(d^v)_!}\ar[d]_{\Psi_{n,k+1}\circ\Omega^k}&sl\hat{\D}_{n+1,k+1}\ar[d]_{\Psi_{n+1,k+1}}&sl\hat{\D}_{n+1,k}\ar[l]_{(d^h)_!}\ar[d]^{\Psi_{n+1,k}}\\
sl\hat{\D}_{k,n+1}\ar[r]^{(d^h)_*}&sl\hat{\D}_{k,n+2}&sl\hat{\D}_{k-1,n+2}\ar[l]_{(d^v)_!}
}
\end{equation}
commute up to natural isomorphism. Our notation is slightly abusive, since the morphisms $(d^v)_*$ and $(d^h)_!$ are not the restrictions of the respective Kan extensions (c.f. \autoref{con:sl-cube}). The commutativity of \eqref{eq:main-comm-I} and \eqref{eq:main-comm-II} will be established in steps (iv)-(ix). We can assume by induction that the morphisms $\Psi_{n',k'}$ are equivalences whenever $n'+k'\leq n+k+1$. We consider the unit $\eta_{n+1,k+1}$ of the adjunction $\Psi_{n+1,k+1}\dashv\Psi^{\vee}_{n+1,k+1}$ and use that commutativity of \eqref{eq:main-comm-I} to show that it is an isomorphism. For this we consider the diagram
\begin{equation}
\xymatrixcolsep{0.48cm}
\xymatrix{
(sl\hat{\D}_{n+1,k+1})^{[1]}_{A,A'}\ar@{<-}@<0.5mm>[rrr]_{(\Psi_{n+1,k+1})^{[1]}}\ar[d]_{F}^{\sim}&&&(sl\hat{\D}_{k,n+2})^{[1]}_{0}\ar[d]_{F}\ar@{<-}@<0.5mm>[lll]_{(\Psi^{\vee}_{n+1,k+1})^{[1]}}\ar[rrr]^{((\mathsf{s}^*_3)^{n+k})^{[1]}}&&&(sl\hat{\D}_{k,n+2})^{[1]}_{A^{\vee},A^{\vee'}}\ar[d]_{F}^{\sim}\\
sl\hat{\D}_{n+1,k+1}\ar@{<-}@<0.5mm>[rrr]_{\Psi_{n+1,k+1}}&&&sl\hat{\D}_{k,n+2}\ar@{<-}@<0.5mm>[lll]_{\Psi^{\vee}_{n+1,k+1}}\ar[rrr]^{(\mathsf{s}^*_3)^{n+k}}&&&sl\hat{\D}_{k,n+2},
}
\end{equation}
where $(sl\hat{\D}_{k,n+2})^{[1]}_{0}:=\mathrm{essim}((\mathsf{s}_3^*)^{-(n+k)}\colon(sl\hat{\D}_{k,n+2})^{[1]}_{A^{\vee},A^{\vee'}}\to(sl\hat{\D}_{k,n+2})^{[1]})$, and observe the following.
\begin{itemize}
\item The left vertical morphism is an equivalence by \autoref{cor:recollement} and \autoref{rmk:inverse}.
\item The right vertical morphism is an equivalence by \autoref{rmk:dual-recoll} (ii) and \autoref{rmk:inverse}.
\item The two right horzontal morphisms are equivalences, and the right square is commutative.
\item By the two points above also the middle vertical morphism is an equivalence.
\item By \autoref{cor:gsdver} there is an isomorphism
\begin{align}
d^v_*=&\mathsf{d}^v[-2n-1]\\
\cong&(\mathsf{s}_3^*)^{-(n+k)}\circ\mathsf{d}^v[2k-1]\circ(\mathsf{s}_3^*)^{n+k}\\
=&(\mathsf{s}_3^*)^{-(n+k)}\circ d^{v\vee}_*\circ(\mathsf{s}_3^*)^{n+k}\colon\D_{k-1,n+2}\to\D_{k,n+2},
\end{align}
where the two equalities are \autoref{eg:ver-fd}.
\item There are isomorphisms
\begin{align}
\Omega^{k}\circ d^h_*=&\Omega^{k}\circ\mathsf{d}^h[1]\\
\cong&(\mathsf{s}_3^*)^{-(n+k)}\circ\Omega^{k}\circ\mathsf{d}^h[2(n+k)+1]\circ(\mathsf{s}_3^*)^{n+k}\\
\cong&(\mathsf{s}_3^*)^{-(n+k)}\circ d^{h\vee}_!\circ(\mathsf{s}_3^*)^{n+k}\colon\D_{k,n+1}\to\D_{k,n+2},
\end{align}
where the first isomorphism is \autoref{cor:gsd-dual} and the second isomorphism is \autoref{prop:dhvee}.
\item The two points above imply that the two upper left horizontal morphisms are well-defined, and hence form an adjunctions.
\end{itemize}
The first, fourth and last point above together imply that $\eta_{n+1,k+1}$ is an isomorphism if and only if the unit $\eta'$ of the upper adjunction in the diagram above is an isomorphism. The commutativity of the left half of \eqref{eq:main-comm-I} implies that
\begin{equation}\label{eq:main-unit-I}
(d^h\times\mathrm{d}_1)^*\circ\eta'\cong\eta_{n+1,k}\qquad\text{and}\qquad(d_v\times\mathrm{d}_1)^*\circ\eta'\cong\id_0
\end{equation}
and the commutativity of the right half of \eqref{eq:main-comm-I} implies that
\begin{equation}\label{eq:main-unit-II}
(d^v\times\mathrm{d}_0)^*\circ\eta'\cong\Sigma\circ\eta_{n,k+1}\qquad\text{and}\qquad(d^h\times\mathrm{d}_0)^*\circ\eta'\cong\id_0.
\end{equation}
Furthermore, the property (P2) for non-injective objects implies that 
\begin{equation}\label{eq:main-unit-III}
(x\times\id_{[1]})^*\circ\eta'\cong\id_0
\end{equation}
for all $x\in Sl^{\square}_{n+1,k+1}\setminus(d^h(Sl^{\square}_{n+1,k})\cup d^v(Sl^{\square}_{n,k+1}))$. The axiom (Der2) together with \eqref{eq:main-unit-I}, \eqref{eq:main-unit-II} and \eqref{eq:main-unit-III} and the induction assumption implies that $\eta'$, and hence also $\eta_{n+1,k+1}$ is invertible. We use anologous arguments applied to \eqref{eq:main-comm-II} to show that also the counit of $\Psi_{n+1,k+1}\dashv\Psi^{\vee}_{n+1,k+1}$ is invertible. Consequently, the morphisms $\Psi_{n',k'}$ are equivalences for $n'+k'\leq n+k+2$ and $n'\geq 2$, $k'\geq 3$. We note that step (ii) takes care of the remaining cases $n'=1$ and $k'=2$ and therefore finishes the induction step.
\item We claim that the diagrams
\[
\xymatrix{
(\cube{n})^{k-1}\ar[r]^{sh_{n,k-1}}\ar[d]_{(\mathrm{d}_0^{n+1})^{k-1}}&(\cube{k-1})^{n}\ar[d]^{(\id_{\cube{k-1}})^n\times\infty}&&(\cube{n})^{k-1}\ar[r]^{sh_{n,k-1}}\ar[d]_{(\mathrm{d}_1^{0})^{k-1}}&(\cube{k-1})^{n}\ar[d]^{\emptyset\times(\id_{\cube{k-1}})^n}\\
(\cube{n+1})^{k-1}\ar[r]_{sh_{n+1,k-1}}&(\cube{k-1})^{n+1},&&(\cube{n+1})^{k-1}\ar[r]_{sh_{n+1,k-1}}&(\cube{k-1})^{n+1}
}
\]
commute. Consider $l=((l_0^0,\cdots,l_{n-1}^0),\cdots,(l_0^{k-2},\cdots,l_{n-1}^{k-2}))\in(\cube{n})^{k-1}$ then the maps of the left square above operate on $l$ as follows
\[
\resizebox{.95\hsize}{!}{$
\xymatrix{
((l_0^0,\cdots,l_{n-1}^0),\cdots,(l_0^{k-2},\cdots,l_{n-1}^{k-2}))\ar@{|->}[r]\ar@{|->}[d]&((l_0^0,\cdots,l_0^{k-2}),\cdots,(l_{n-1}^0,\cdots,l_{n-1}^{k-2}))\ar@{|->}[d]\\
((l_0^0,\cdots,l_{n-1}^0,1),\cdots,(l_0^{k-2},\cdots,l_{n-1}^{k-2},1))\ar@{|->}[r]&((l_0^0,\cdots,l_0^{k-2}),\cdots,(l_{n-1}^0,\cdots,l_{n-1}^{k-2}),(1,\cdots,1)).
}$}
\]
Dually, for the right square we have
\[
\resizebox{.95\hsize}{!}{$
\xymatrix{
((l_0^0,\cdots,l_{n-1}^0),\cdots,(l_0^{k-2},\cdots,l_{n-1}^{k-2}))\ar@{|->}[r]\ar@{|->}[d]&((l_0^0,\cdots,l_0^{k-2}),\cdots,(l_{n-1}^0,\cdots,l_{n-1}^{k-2}))\ar@{|->}[d]\\
((0,l_0^0,\cdots,l_{n-1}^0),\cdots,(0,l_0^{k-2},\cdots,l_{n-1}^{k-2}))\ar@{|->}[r]&((0,\cdots,0),(l_0^0,\cdots,l_0^{k-2}),\cdots,(l_{n-1}^0,\cdots,l_{n-1}^{k-2})).
}$}
\]
\item The commutativity of the right square in step (iv) implies the commutativity of the third square of the left and the fourth square on the right in the diagrams below
\[
\xymatrix{
sl\hat{\D}_{n,k+1}\ar[r]_{(d^h)^*}\ar[d]_{c_{n,k+1}^*}&sl\hat{\D}_{n,k}\ar[d]^{c_{n,k}^*} &&sl\hat{\D}_{k,n+1}\ar[r]_{(d^v)^*}\ar[d]_{c_{k,n+1}^*}&sl\hat{\D}_{k-1,n+1}\ar[d]^{c_{k-1,n+1}^*}\\
\D_0^{[n]^{k}}\ar[r]_{(\emptyset\times(\id_{[n]})^{k-1})^*}\ar[d]_{((\rightarrow_{\tau})^{k})_!}&\D_0^{[n]^{k-1}}\ar[d]^{((\rightarrow_{\tau})^{k-1})_!} &&\D_0^{[k]^{n}}\ar[r]_{((\mathrm{d}_0)^n)^*}\ar[d]_{((\rightarrow_{\tau})^{n})_*}&\D_0^{[k-1]^{n}}\ar[d]^{((\rightarrow_{\tau})^{n})_*}\\
\D^{(\cube{n})^{k}}\ar[r]_{(\emptyset\times(\id_{\cube{n}})^{k})^*}\ar[d]_{sh_{k,n}^*}&\D^{(\cube{n})^{k-1}}\ar[d]^{sh_{k-1,n}^*} &&\D^{(\cube{k})^{n}}\ar[r]_{((\mathrm{d}_0^{0})^{n})^*}\ar[d]_{\cof^{\underline{1}}}&\D^{(\cube{k-1})^{n}}\ar[d]^{\cof^{\underline{1}}}\\
\D^{(\cube{k})^{n}}\ar[r]_{((\mathrm{d}_1^{0})^{n})^*}\ar[d]_{\fib^{\underline{1}}}&\D^{(\cube{k-1})^{n}}\ar[d]^{\fib^{\underline{1}}} &&\D^{(\cube{k})^{n}}\ar[r]_{((\mathrm{d}_1^{0})^{n})^*}\ar[d]_{sh_{n,k}^*}&\D^{(\cube{k-1})^{n}}\ar[d]^{sh_{n,k-1}^*}\\
\D^{(\cube{k})^{n}}\ar[r]_{((\mathrm{d}_0^{0})^{n})^*}\ar[d]_{((\rightarrow_{\tau})^n)^*}&\D^{(\cube{k-1})^{n}}\ar[d]^{((\rightarrow_{\tau})^n)^*} &&\D^{(\cube{n})^{k}}\ar[r]_{(\emptyset\times(\id_{\cube{n}})^{k})^*}\ar[d]_{((\rightarrow_{\tau})^{k})^*}&\D^{(\cube{n})^{k-1}}\ar[d]^{((\rightarrow_{\tau})^{k-1})^*}\\
\D_0^{[k]^{n}}\ar[r]_{((\mathrm{d}_0)^n)^*}\ar[d]_{(c_{k,n+1}^{-1})^*}&\D_0^{[k-1]^{n}}\ar[d]^{(c_{k-1,n+2}^{-1})^*} &&\D_0^{[n]^{k}}\ar[r]_{(\emptyset\times(\id_{[n]})^{k-1})^*}\ar[d]_{(c_{n,k+1}^{-1})^*}&\D_0^{[n]^{k-1}}\ar[d]^{(c_{n,k}^{-1})^*}\\
sl\hat{\D}_{k,n+1}\ar[r]_{(d^v)^*}&sl\hat{\D}_{k-1,n+1}, &&sl\hat{\D}_{n,k+1}\ar[r]_{(d^h)^*}&sl\hat{\D}_{n,k}.\\
}
\]
The first and sixth squares follow in both cases from \autoref{con:sl-cube}. For the fourth square on the left and the third square on the right we invoke \autoref{rmk:elementary}. The second and fifth squares are induced by inverse images from commutative squares in $Cat$. In the case of the second squares we additionally use that $(\rightarrow_{\tau})$ has a left and a right adjoint. Moreover, we use step (i) and the fact that $\fib^{\underline{1}}$ preserves permutations of coordinates to identify the columns with $\Psi_{n,k+1},\Psi_{n,k}$ and $\Psi^{\vee}_{n,k+1}$,$\Psi^{\vee}_{n,k}$, respectively. Hence we have isomorphisms
\begin{equation}\label{eq:main-square-I}
\Psi_{n,k}\circ(d^h)^*\cong(d^v)^*\circ\Psi_{n,k+1}\qquad\text{and}\qquad\Psi^{\vee}_{n,k}\circ(d^v)^*\cong(d^h)^*\circ\Psi^{\vee}_{n,k+1}
\end{equation}
By passing to the right (resp. left) adjoints in \eqref{eq:main-square-I} we obtain
\begin{equation}\label{eq:main-square-II}
(d^h)_*\circ\Psi^{\vee}_{n,k}\cong\Psi^{\vee}_{n,k+1}\circ(d^v)_*\qquad\text{and}\qquad(d^v)_!\circ\Psi_{n,k}\cong\Psi_{n,k+1}\circ(d^h)_!,
\end{equation}
which exactly yields the commutativity of the lower left (resp. lower right) square of \eqref{eq:main-comm-I} (resp. \eqref{eq:main-comm-II}).
\item Consider the upper left square of \eqref{eq:main-comm-I}
\begin{equation}\label{eq:square-ul}
\xymatrix{
sl\hat{\D}_{n,k}\ar[r]^{(d^h)_*}\ar[d]_{\Psi_{n,k}}&sl\hat{\D}_{n,k+1}\ar[d]_{\Psi_{n,k+1}}\\
sl\hat{\D}_{k-1,n+1}\ar[r]^{(d^v)_*}&sl\hat{\D}_{k,n+1}.
}
\end{equation}
The isomorphism \eqref{eq:main-square-I} implies that
\[
(d^v)^*\circ\Psi_{n,k+1}\circ(d^h)_*\cong\Psi_{n,k}\circ(d^h)^*\circ(d^h)_*\cong\Psi_{n,k}.
\]
Hence, for the commutativity of \eqref{eq:square-ul} it is sufficient to show that the essential image of $\Psi_{n,k+1}\circ(d^h)_*$ is contained in the essential image of $(d^v)_*$. We restrict to the situation of $\Psi^{\square}_{n,k}$ by composing with the equivalences $c_{n,k}$. In particular,  \eqref{eq:square-ul} becomes
\begin{equation}\label{eq:square-uln}
\xymatrix{
\D_0^{[n]^{k-1}}\ar[rr]^{(\emptyset\times\id_{[n]}^{k-1})_*}\ar[d]_{\Psi^{\square}_{n,k}}&&\D_0^{[n]^{k}}\ar[d]_{\Psi^{\square}_{n,k+1}}\\
\D_0^{[k-1]^{n}}\ar[rr]^{(\mathrm{d}_0)_*\times(\mathrm{d}_0^{n-1})_!}&&\D_0^{[k]^{n}}.
}
\end{equation}
by \autoref{con:sl-cube}. Under these identification, we observe that the condition 'contained in the essential image of $(d^v)_*$' is equivalent to the constantness of the restrictions along the 1-simplices classifying the maps 
\[
(0,i_1,\cdots,i_{n-1})\to(1,i_1,\cdots,i_{n-1})
\]
for all $(i_1,\cdots,i_{n-1})\in [k]^{n-1}$. By precomposing with $((\rightarrow_{\tau})^n)^*$, the above constantness conditions follow form the constantness of the restrictions along the 1-simplices classifying the maps
\[
\resizebox{.95\hsize}{!}{$
((0,\cdots,0),(\rightarrow_{\tau})(i_1),\cdots,(\rightarrow_{\tau})(i_{n-1}))\to((1,0,\cdots,0),(\rightarrow_{\tau})(i_1),\cdots,(\rightarrow_{\tau})(i_{n-1})).
$}
\]
Let $(i_1,\cdots,i_{n-1})\in[k]^{n-1}$, and $s\colon[1]\to\cube{k^n}$ the map described above. Then we compute the cofiber
\begin{align}
&C\circ s^*\circ\fib^{\underline{1}}\circ sh_{k,n}^*\circ((\rightarrow_{\tau})^k)_!\circ(\emptyset\times\id_{[n]}^{k-1})_*\\
\cong&\tfib\circ\iota_1^*\circ sh_{k,n}^*\circ((\rightarrow_{\tau})^k)_!\circ(\emptyset\times\id_{[n]}^{k-1})_*\\
\cong&\tfib\circ\iota_2^*\circ((\rightarrow_{\tau})^k)_!\circ(\emptyset\times\id_{[n]}^{k-1})_*
\end{align}
by \autoref{rmk:elementary} (first isomorphism) and the functoriality of inverse images (second isomorphism), and where 
\[
\iota_1\colon((\cube{k})^n)_{s(1)/}\to(\cube{k})^n \qquad\text{and}\qquad\iota_2\colon((\cube{n})^k)_{sh_{k,n}\circ s(1)/}\to(\cube{n})^k
\]
denote the natural inclusions of undercategories. We claim that the latter total cofiber is compute over a constant cube with value zero. For this we consider $j=(j_0,\cdots,j_{k-1})\in((\cube{n})^k)_{sh_{k,n}\circ s(1)/}$. Since the functor $(\rightarrow_{\tau})$ admits a right adjoint $p$ (c.f. proof of \autoref{thm:dual-Sdot-sl}) we have $((\rightarrow_{\tau})^k)_!\cong(p^k)^*$. Using that the first index of $j_0$ is 1, we see that $p(j)$ is contained in the complement of the image of $\emptyset\times\id_{[n]}$. Hence, 
\[
p^*\circ((\rightarrow_{\tau})^k)_!\circ(\emptyset\times\id_{[n]}^{k-1})_*\cong p^*\circ(p^k)^*\circ(\emptyset\times\id_{[n]}^{k-1})_*\cong 0.
\]
\item Next, we show the commutativity of the upper right square 
\[
\xymatrix{
sl\hat{\D}_{n+1,k}\ar[d]_{\Psi_{n+1,k}}&sl\hat{\D}_{n,k}\ar[l]_{(d^v)_!}\ar[d]^{\Psi_{n,k}\circ\Omega^{k-1}}\\
sl\hat{\D}_{k-1,n+2}&sl\hat{\D}_{k-1,n+1}\ar[l]_{(d^h)_*}
}
\]
of \eqref{eq:main-comm-I} (which coincides with the lower left square of \eqref{eq:main-comm-II}). We invoke \eqref{eq:main-square-I} for 
\[
(d^v)^*\circ\Psi_{n+1,k}\circ(d^v)_!\cong\Psi_{n+1,k-1}\circ(d^h)^*\circ(d^v)_!\cong 0.
\]
Hence the essential image of $\Psi_{n+1,k}\circ(d^v)_!$ is contained in the essential image of $(d^h)_*\circ\Psi_{n,k}\circ\Omega^{k-1}$. Therefore it is sufficient to show that
\begin{equation}\label{eq:main-twisted}
(d^h)^*\circ\Psi_{n+1,k}\circ(d^v)_!\cong\Psi_{n,k}\circ\Omega^{k-1}.
\end{equation}
By composition with inverse images of the form $c_{n,k}^*$ this in turn can be reformulated as the commutativity of
\[
\xymatrix{
\D_0^{[n+1]^{k-1}}\ar[d]_{\Psi_{n+1,k}^{\square}}&\D_0^{[n]^{k-1}}\ar[d]^{\Psi_{n,k}^{\square}\circ\Omega^{k-1}}\ar[l]_{(d_0)^{k-1}_!}\\
\D_0^{[k-1]^{n+1}}\ar[r]_{(\emptyset\times\id_{[k-1]}^n)^*}&\D_0^{[k-1]^n}.
}
\]
We compute the composition through the left hand side 
\begin{align}
&(\emptyset\times\id_{[k-1]}^n)^*\circ\Psi_{n+1,k}^{\square}\circ(d_0)^{k-1}_!\\
=&(\emptyset\times\id_{[k-1]}^n)^*\circ((\rightarrow_{\tau})^{n+1})^*\circ(sh_{k-1,n+1})^*\circ\fib^{\underline{1}}\circ(\rightarrow_{\tau})^{k-1}_!\circ(d_0)^{k-1}_!\\
\cong&((\rightarrow_{\tau})^{n})^*\circ(sh_{k-1,n})^*\circ((\mathrm{d}_1^0)^{k-1})^*\circ\fib^{\underline{1}}\circ(\rightarrow_{\tau})^{k-1}_!\circ(d_0)^{k-1}_!\\
\cong&((\rightarrow_{\tau})^{n})^*\circ(sh_{k-1,n})^*\circ((\mathrm{d}_1^0)^{k-1})^*\circ\fib^{\underline{1}}\circ(\mathrm{d}_0^0)^{k-1}_!\circ(\rightarrow_{\tau})^{k-1}_!\\
\cong&((\rightarrow_{\tau})^{n})^*\circ(sh_{k-1,n})^*\circ\fib^{\underline{1}}\circ(\rightarrow_{\tau})^{k-1}_!\circ\Omega^{k-1}\\
=&\Psi^{\square}_{n,k}\circ\Omega^{k-1}.
\end{align}
The equalities are the definition of $\Psi^{\square}_{n,k}$, the first two isomorphisms follow from the 2-functoriality of inverse images and the pseudofunctoriality of left Kan extensions, respectively. The third isomorphism is the $(\cube{n})^{k-1}$-parametrized version of 
\[
\mathrm{d}_1^*\circ\fib\circ(\mathrm{d}_0)_!\cong\Omega.
\]
\item For the commutativity of the lower right square
\[
\xymatrix{
sl\hat{\D}_{k-1,n+2}\ar[d]_{\Psi_{n+1,k}^{\vee}}&sl\hat{\D}_{k-1,n+1}\ar[l]_{(d^h)_*}\ar[d]^{\Sigma^{k-1}\circ\Psi_{n,k}^{\vee}}\\
sl\hat{\D}_{n+1,k}&sl\hat{\D}_{n,k}\ar[l]_{(d^v)_!}
}
\]
of \eqref{eq:main-comm-I} (which coincides with the upper left square of \eqref{eq:main-comm-II}) we invoke \eqref{eq:main-square-I} for
\[
(d^h)^*\circ\Psi^{\vee}_{n+1,k}\circ(d^h)_*\cong\Psi^{\vee}_{n+1,k-1}\circ(d^v)^*\circ(d^h)_*\cong 0.
\]
Hence the essential image of $\Psi^{\vee}_{n+1,k}\circ(d^h)_*$ is contained in the essential image of $(d^v)_!\circ\Psi^{\vee}_{n,k}\circ\Sigma^{k-1}$. Therefore it is sufficient to show that
\[
(d^v)^*\circ\Psi^{\vee}_{n+1,k}\circ(d^h)_*\cong\Psi^{\vee}_{n,k}\circ\Sigma^{k-1}.
\]
But this is obtained as the right adjoint isomorphism of \eqref{eq:main-twisted}, which was shown in the previous step.
\item For the verification of the assumptions of step (iii), it remains to show the commutativity of the upper right square
\begin{equation}\label{eq:square-ur}
\xymatrix{
sl\hat{\D}_{k,n+1}\ar[d]_{\Psi_{n,k+1}^{\vee}}&sl\hat{\D}_{k-1,n+1}\ar[l]_{(d^v)_!}\ar[d]^{\Psi_{n,k}^{\vee}}\\
sl\hat{\D}_{n,k+1}&sl\hat{\D}_{n,k}\ar[l]_{(d^h)_!}
}
\end{equation}
of \eqref{eq:main-comm-II}. The strategy for this step is similar to step (vi). We apply \eqref{eq:main-square-I} for 
\[
(d^h)^*\circ\Psi^{\vee}_{n,k+1}\circ(d^v)_!\cong\Psi^{\vee}_{n,k}\circ(d^v)^*\circ(d^v)_!\cong\Psi^{\vee}_{n,k}.
\]
Hence for the commutativity of \eqref{eq:square-ur} it is sufficient to show that the essential image of $\Psi^{\vee}_{n,k+1}\circ(d^v)_!$ is contained in the essential image of $(d^h)_!$. By \autoref{con:sl-cube}, it suffices to show that for every injective $\tilde{i}\in Sl_{n-1,k+1}$ the elementary subcube
\begin{equation}\label{eq:main-bicart}
\square_{\tilde{i}}^*\circ\Psi^{\vee}_{n,k+1}\circ(d^v)_!
\end{equation}
is bicartesian. We use that the injective elements in $Sl_{n-1,k+1}$  correspond under $c_{n,k+1}$ to increasing sequences in $[n-1]^{k}$. We denote $c_{n,k+1}(\tilde{i})=i=(i_0,\cdots,i_{k-1})$. By composition with the equivalence, the bicartesianess of \eqref{eq:main-bicart} is seen to be equivalent to the biartesianess of 
\begin{equation}\label{eq:main-bicart-II}
F_i:=\square_i^*\circ\Psi^{\square\vee}_{n,k+1}\circ(\mathrm{d}_0)^n_!.
\end{equation}
Using the functoriality of inverse images for the composition $(\rightarrow_{\tau})^k\circ\square_i$ and \autoref{rmk:elementary}, we obtain
\[
F_i\cong\prod_{j=0}^{k-1}(\mathrm{d}_0^0\times\cdots\times\mathrm{d}_0^{i_{j-1}}\times\id\times\mathrm{d}_1^{i_{j+1}}\times\cdots\times\mathrm{d}_1^{n-1})^*\circ\cof^{\underline{1}}\circ sh_{n,k}^*\circ(\rightarrow_{\tau})_*^n\circ (\mathrm{d}_0)^n_!.
\]
It is sufficient to show that $\tfib\circ F_i\cong 0$. To establish this, we observe
\begin{align}
&\tfib\circ F_i\\
\cong&\prod_{j=0}^{k-1}(\mathrm{d}_0^0\times\cdots\times\mathrm{d}_0^{i_{j-1}}\times F\times\mathrm{d}_1^{i_{j+1}}\times\cdots\times\mathrm{d}_1^{n-1})^*\circ\cof^{\underline{1}}\circ sh_{n,k}^*\circ(\rightarrow_{\tau})_*^n\circ (\mathrm{d}_0)^n_!\\
\cong&\prod_{j=0}^{k-1}(C\times\cdots\times C\times \mathrm{d}_1^{i_j}\times\mathrm{d}_0^{i_{j+1}}\times\cdots\times\mathrm{d}_0^{n-1})^*\circ sh_{n,k}^*\circ(\rightarrow_{\tau})_*^n\circ (\mathrm{d}_0)^n_!\\
\cong& \tcof\circ\prod_{j=0}^{k-1}(\id\times\cdots\times\id\times \mathrm{d}_1^{i_j}\times\mathrm{d}_0^{i_{j+1}}\times\cdots\times\mathrm{d}_0^{n-1})^*\circ sh_{n,k}^*\circ(\rightarrow_{\tau})_*^n\circ (\mathrm{d}_0)^n_!.
\end{align}
We claim that the latter cofiber is computed over a constant cube with value $0$. The value at the initial vertex of this cube vanishes, because of the precomposition with the extension-by-zero morphism $(\mathrm{d}_0)^{n}_!$. For every other vertex
\[
\delta=((\delta_{0,0},\cdots,\delta_{0,i_0-1},0,1,\cdots,1),\cdots,(\delta_{k-1,0},\cdots,\delta_{0,i_{k-1}-1},0,1,\cdots,1))
\]
with $\delta_{j,i}\in\lbrace 0,1 \rbrace$ there exists $(j,i)$ such that $\delta_{j,i}=1$. We consider $j_m\in\mathbf{k}$ maximal among those $j$ such that $\delta_{j,i}=1$ for some $i\in\lbrace 0, \cdots, i_{j}-1\rbrace$, and $i_m\in\lbrace 0, \cdots, i_{j_m}-1\rbrace$ maximal among those $i$ such that $\delta_{j_m,i}=1$. Consequently, $\delta^*\circ sh_{n,k}^*\circ(\rightarrow_{\tau})_*^n$ is by construction of $(\rightarrow_{\tau})_*$ exhibited as  $l^*$ for some $l=(l_0\cdots,l_{n-1})\in[k]^{n}$ with $l_{i_m}=j_m$ and $l_{i_m+1}<j_m$. This yields the desired vanishing. 
\item Using \autoref{rmk:dh-welldef} and \autoref{eg:ver-fd} we can reformulate \eqref{eq:main-square-I} as 
\[
\Psi_{n,k}\circ\mathsf{d}^h\cong\mathsf{d}^v[-2n]\circ\Psi_{n,k+1}.
\]
Since $\Psi_{n,k}$ is an equivalence with inverse $\Psi^{\vee}_{n,k}$ (step (iii)), we can pass to iterated adjoints on both sides of the above isomorphism, to obtain for $a\in\mathbb{Z}$
\begin{equation}\label{eq:main-face-I}
\Psi_{n,k}\circ\mathsf{d}^h[a]\cong\mathsf{d}^v[-2n+a]\circ\Psi_{n,k+1}
\end{equation}
and
\begin{equation}\label{eq:main-face-dual}
\mathsf{d}^h[a]\circ\Psi_{n,k}^{\vee}\cong\Psi_{n,k+1}^{\vee}\circ\mathsf{d}^v[-2n+a]
\end{equation}
for $a$ even, respectively odd. In the odd case, precomposing with $\Psi_{n,k}$ and postcomposing with $\Psi_{n,k+1}$ leads to
\begin{equation}\label{eq:main-face-II}
\Psi_{n,k+1}\circ\mathsf{d}^h[a]\cong\mathsf{d}^v[-2n+a]\circ\Psi_{n,k}.
\end{equation}
In particular, we have the following composition of isomorphisms
\begin{align}
\Psi_{n,k}\circ\mathsf{s}_3^*\cong&\Psi_{n,k}\circ\mathsf{d}^h\circ\mathsf{d}^h[-1]\circ\mathsf{s}_3^*\\
\cong&\Psi_{n,k}\circ\mathsf{d}^h\circ(\mathsf{s}_3^*)^{n+k+1}\circ\mathsf{d}^h[-2(n+k)-3]\circ(\mathsf{s}_3^*)^{-(n+k+1)}\circ\mathsf{s}_3^*\\
\cong&\Psi_{n,k}\circ\mathsf{d}^h\circ\Sigma^{nk}\circ\mathsf{d}^h[-2(n+k)-3]\circ\Omega^{n(k-1)}\\
\cong&\mathsf{d}^v[-2n]\circ\Psi_{n,k+1}\circ\mathsf{d}^h[-2(n+k)-3]\circ\Sigma^{n}\\
\cong&\mathsf{d}^v[-2n]\circ\mathsf{d}^v[-4n-2k-3]\circ\Psi_{n,k}\circ\Sigma^{n}\label{eq:main-sym}\\
\cong&\mathsf{d}^v[-2n]\circ(\mathsf{s}_3^*)^{-(n+k+1)}\circ\mathsf{d}^v[-2n-1]\circ(\mathsf{s}_3^*)^{n+k+1}\circ\Psi_{n,k}\circ\Sigma^{n}\\
\cong&\mathsf{d}^v[-2n]\circ\Omega^{kn}\circ\mathsf{d}^v[-2n-1]\circ\Sigma^{(k-1)n}\circ\mathsf{s}_3^*\circ\Psi_{n,k}\circ\Sigma^{n}\\
\cong&\mathsf{d}^v[-2n]\circ\mathsf{d}^v[-2n-1]\circ\mathsf{s}_3^*\circ\Psi_{n,k}\\
\cong&\mathsf{s}_3^*\circ\Psi_{n,k},
\end{align}
where the single isomorphisms are induced by, first, the unit of the adjunction $\mathsf{d}^h[-1]\dashv\mathsf{d}^h$ (\autoref{prop:preparations-h}), second, \autoref{cor:gsd-dual}, third, \autoref{cor:fracCY}, fourth, \eqref{eq:main-face-I}, fifth, \eqref{eq:main-face-II}, sixth, \autoref{cor:gsdver}, seventh, \autoref{cor:fracCY} and ninth, the inverse of the unit of the adjunction $\mathsf{d}^v[-2n-1]\dashv\mathsf{d}^v[-2n]$. We have also used that all morphisms of derivators commute with $\Sigma$ and $\Omega$.
\item We define $\Phi_{n,k}=(\mathsf{s_3}^*)^n\circ\Psi_{n,k}\colon\D_{n,k}\toiso\D_{k-1,n+1}$, which is an equivalence by the previous steps. Hence it remains to prove part (i) to (vi) of the theorem. We use \eqref{eq:main-sym} to conclude
\[
\Phi_{n,k}\circ\mathsf{s}_3^*=(\mathsf{s_3}^*)^n\circ\Psi_{n,k}\circ\mathsf{s}_3^*\cong\mathsf{s}_3^*\circ(\mathsf{s_3}^*)^n\circ\Psi_{n,k}\cong\mathsf{s}_3^*\circ\Phi_{n,k}.
\]
and hence part (ii). For part (vi) of the theorem we use step (ii) for 
\[
\Phi_{n,2}=(\mathsf{s_3}^*)^n\circ\Psi_{n,2}\cong(\mathsf{s_3}^*)^n\circ\Psi'_n=\Psi_n.
\]
\item For part (iii) of the theorem, we consider
\begin{align}
\Phi_{n,k}\circ\mathsf{d}^h=&(\mathsf{s}_3^*)^n\circ\Psi_{n,k}\circ\mathsf{d}^h\\
\cong&(\mathsf{s}_3^*)^n\circ\mathsf{d}^v[-2n]\circ\Psi_{n,k+1}\\
\cong&\mathsf{d}^v\circ(\mathsf{s}_3^*)^n\circ\Psi_{n,k+1}\\
=&\mathsf{d}^v\circ\Phi_{n,k+1},
\end{align}
where the first isomorphism is \eqref{eq:main-face-I} and second isomorphism is \autoref{cor:gsdver}. We obtain part (iii) of the theorem in the case $a\in\mathbb{Z}$ even, by passing to iterated adjoint isomorphisms (which exist by \autoref{cor:gsdver}, \autoref{cor:gsd-dual} and step (iii)). For $a\in\mathbb{Z}$ odd apply the analogous argument building on \eqref{eq:main-face-dual}.
\item 
\begin{figure}
\resizebox{.9\linewidth}{!}{
  \begin{minipage}{\linewidth}
  \begin{align}
&\Phi_{n,k}^{-1}\\
=&\Psi_{n,k}^{\vee}\circ(\mathsf{s}_3^*)^{-n}\\
\cong&(\mathsf{s}_3^*)^{-n}\circ\Psi_{n,k}^{\vee}\\
\cong&(\mathsf{s}_3^*)^{-n}\circ (d^h)^*\circ d^h_!\circ\Psi_{n,k}^{\vee}\\
\cong&(\mathsf{s}_3^*)^{-n}\circ (d^h)^*\circ\Psi_{n,k+1}^{\vee}\circ d^v_!\\
=&(\mathsf{s}_3^*)^{-n}\circ (d^h)^*\circ (c_{n,k+1}^{-1})^* \circ((\rightarrow_{\tau})^{k})^*\circ\cof^{\underline{1}}\circ sh_{n,k}^*\circ((\rightarrow_{\tau})^{n})_*\circ c_{k,n+1}^*\circ d^v_!\\
\cong&(\mathsf{s}_3^*)^{-n}\circ (d^h)^*\circ (c_{n,k+1}^{-1})^* \circ((\rightarrow_{\tau})^{k})^*\circ(\cof^{\underline{1}})^2\circ sh_{n,k}^*\circ((\rightarrow_{\tau})^{n})_!\circ c_{k,n+1}^*\circ d^{v\vee}_*\circ\Omega^{n}\\
\cong&(\mathsf{s}_3^*)^{-n}\circ (d^h)^*\circ (c_{n,k+1}^{-1})^* \circ((\rightarrow_{\tau})^{k})^*\circ\fib^{\underline{1}}\circ sh_{n,k}^*\circ((\rightarrow_{\tau})^{n})_!\circ c_{k,n+1}^*\circ d^{v\vee}_*\circ\Sigma^{(k-1)n}\\
=&(\mathsf{s}_3^*)^{-n}\circ (d^h)^*\circ \Psi_{k,n+1}\circ d^{v\vee}_*\circ\Sigma^{(k-1)n}\\
\cong&\Sigma^{(k-1)n}\circ(\mathsf{s}_3^*)^{-n}\circ (d^h)^*\circ \Psi_{k,n+1}\circ d^{v\vee}_*\\
\cong&(\mathsf{s}_3^*)^{k}\circ (d^h)^*\circ \Psi_{k,n+1}\circ d^{v\vee}_*\\
\cong&(\mathsf{s}_3^*)^{k-1}\circ (d^{h\vee})^*\circ \Psi_{k,n+1}\circ d^{v\vee}_*\\
\cong&(\mathsf{s}_3^*)^{k-1}\circ \Psi_{k-1,n+1}\circ (d^{v\vee})^*\circ d^{v\vee}_*\\
\cong&(\mathsf{s}_3^*)^{k-1}\circ \Psi_{k-1,n+1}\\
=&\Phi_{k-1,n+1}.
\end{align}
  \end{minipage}
}
\caption{The proof of part (i).}
\label{fig:mainthm}
\end{figure}
In this step we first show that there are isomorphisms
\begin{equation}\label{eq:main-assumption}
\fib^{\underline{1}}\circ(\rightarrow_{\tau})^n_*\circ c_{k,n+1}^*\circ d^v_!\cong\Omega^{n}\circ(\rightarrow_{\tau})^n_!\circ c_{k,n+1}^*\circ d^{v\vee}_*\colon sl\hat{\D}_{k-1,n+1}\to\D^{(\cube{k})^n}
\end{equation}
as follows
\begin{align}
&\fib^{\underline{1}}\circ(\rightarrow_{\tau})^n_*\circ c_{k,n+1}^*\circ d^v_!\\
\cong&\fib^{\underline{1}}\circ(\rightarrow_{\tau})^n_*\circ ((\mathrm{d}_0)^n)_! \circ c_{k-1,n+1}^*\\
\cong&(\fib^{\underline{1}}\circ(\rightarrow_{\tau})_*\circ (\mathrm{d}_0)_!)^n \circ c_{k-1,n+1}^*\\
\cong&(\Omega\circ(\rightarrow_{\tau})_!\circ (\mathrm{d}_{k+1})_*)^n \circ c_{k-1,n+1}^*\\
\cong&\Omega^n\circ(\rightarrow_{\tau})^n_!\circ ((\mathrm{d}_{k+1})_*)^n \circ c_{k-1,n+1}^*\\
\cong&\Omega^n\circ(\rightarrow_{\tau})^n_! \circ c_{k,n+1}^*\circ d^{v\vee}_*.
\end{align}
Here, the first and fifth isomorphisms are \autoref{con:sl-cube}, the second and fourth isomorphisms follow from the compatibility with products, and the third isomorphism is \autoref{prop:kappa}.
The isomorphism \eqref{eq:main-assumption} is the key ingredient for part (i) of the theorem, i.e. that $\Phi_{k-1,n+1}\cong\Phi_{n,k}^{-1}$. This is established by \autoref{fig:mainthm}, where
\begin{itemize}
\item The first isomorphism is \eqref{eq:main-sym}, 
\item the second one is the unit of the adjunction $d^h_!\dashv(d^h)^*$ (\autoref{prop:preparations-h}),
\item the third one is \eqref{eq:square-ur}, 
\item the fourth one is \eqref{eq:main-assumption}, 
\item the fifth one is \autoref{rmk:elementary},
\item the sixth one is induced by exactness of morphisms of stable derivators,
\item the seventh one is \autoref{cor:fracCY}, 
\item the eighth one is \autoref{prop:dhvee},
\item the ninth one is induced by exactly the dual construction of step (v),
\item and the tenth one is the counit of the adjunction $(d^{v\vee})^*\dashv d^{v\vee}_*$.
\end{itemize}

\item We obtain part (iv) of the theorem by passing to adjoint isomorphisms of part (iii) and using part (i).
\item Finally, for part (v) we consider the diagram
\[
\xymatrix{
\D_0^{[n]^{k-1}}\ar[r]^{((\rightarrow_{\tau})^{k-1})_!}\ar[rrd]_{\emptyset^*}&\D^{(\cube{n})^{k-1}}\ar[r]^{\fib^{\underline{1}}}\ar[rd]^{\emptyset^*}&\D^{(\cube{n})^{k-1}}\ar[r]^{sh_{k-1,n}^*}\ar[d]_{\infty^*}&\D^{(\cube{k-1})^{n}}\ar[r]^{((\rightarrow_{\tau})^n)^*}\ar[dl]_{\infty^*}&\D_0^{[k-1]^n}\ar[d]^{(c_{k-1,n+1}^{-1})^*}\ar[dll]^{\infty^*}\\
sl\hat{\D}_{n,k}\ar[u]^{c^*_{n,k}}\ar[rr]_{\xi^*}&&\D&&sl\hat{\D}_{k-1,n+1},\ar[ll]^{\xi^*\circ(\mathsf{s}_3^*)^n}
}
\]
where the composition through the top is $\Psi_{n,k}$. All triangles, except the third one, commute by the 2-functoriality of $\D$ (for the second triangle we additionally use $((\rightarrow_{\tau})^{k-1})_!=(p^k)^*$). For the commutativity of the third triangle we invoke \autoref{rmk:elementary}. Hence we obtain
\[
\xi^*\cong\xi^*\circ(\mathsf{s}_3^*)^n\circ\Psi_{n,k}=\xi^*\circ\Phi_{n,k}.
\]
\end{enumerate}
\end{proof}

\begin{rmk}
\begin{enumerate}
\item We have rarely used the first definition of $\Psi^{\square}_{n,k}$, however using this definition we can describe the underlying diagram of $\Psi^{\square}_{n,k}$. For this we note that it follows from the proof of \autoref{prop:filtered} that
\begin{enumerate}
\item $\rightarrow(n)^*\circ\Psi_n^{\square}\cong 0^*$,
\item $\rightarrow(i)^*\circ\Psi_n^{\square}\cong\Omega^{n-i-1}\circ F\circ[n-i-1,n-i]^*$ for $0\leq i\leq n-1$.
\end{enumerate}
We have already seen in the first part proof of \autoref{thm:main} that $j^*\circ\Psi_{n,k}^{\square}=0$ whenever $j=(j_0,\cdots,j_{n-1})$ is not a non-decreasing sequence in $[k-1]^n$. Therefore it is sufficient to determine $j^*\circ\Psi_{n,k}^{\square}$ for non-decreasing sequences $j$. In this case we can regard $j$ as a functor $[n-1]\to[k-1]$ and
\begin{equation}\label{eq:ad}
ad_{n,k}(j)=(\rightarrow(n-\mathrm{min}\lbrace i\vert 1\leq j(i)\rbrace),\cdots,\rightarrow(n-\mathrm{min}\lbrace i\vert k-1\leq j(i)\rbrace)).
\end{equation}
In particular, we obtain a factorization
\[
ad_{n,k}\vert_{[k-1]^{[n-1]}}\colon[k-1]^{[n-1]}\xrightarrow{m_{n,k}}[n]^{[k-2]}\xrightarrow{(\rightarrow)^{[k-2]}}(\cube{n})^{[k-2]},
\]
where $[k-1]^{[n-1]}\xrightarrow{m_{n,k}}[n]^{[k-2]},l\mapsto\mathrm{min}\lbrace i\vert l+1\leq j(i)\rbrace$. We conclude that
\[
j^*\circ\Psi_{n,k}^{\square}\cong \Omega^{j'}\circ\tfib\circ\square_{m_{n,k}(j)},
\]
where the total fiber is computed of the subcube of maximal dimension ending in $m_{n,k}(j)$ and $j'=\sum_{l\vert m_{n,k}(j)_l\geq 1}m_{n,k}-1$. We will make this more explicit in a specific example (\autoref{eg:D34}). Furthermore, we observe that \eqref{eq:ad} relates $ad_{n,k}$ to $\mathrm{ad}_{n,k}$ (c.f. \autoref{adjunction-duality}) and hence justifies the notation.
\item The statements of \autoref{lem:recollement} and \autoref{cor:recollement} mimic the first steps of a general result concerning recollements in the context of $\infty$-categories \cite{BG-recollement}. We emphasize that the proofs of the main results of \emph{loc.~cit.} do not generalize to the context of stable derivators. In particular, we highly expect that the analogue of step (iii) of the proof of \autoref{thm:main} will be significantly simpler in the context of stable $\infty$-categories.

\noindent We refer to \cite{DJW-Sdot} for a precise definition of the $\infty$-categorical analogues of the derivators $\D_{n,k}$ and in particular section 2.4 of \emph{loc.~cit.} for a discussion of the relevant ladders of recollents in this setting.
\item It is possible to extend the statement of \autoref{thm:main} by considering some additional boundary cases. More precisely, we define for a stable derivator $\D$, $n,k\in\mathbb{Z}, k\leq 1, n+k\geq 1$
\begin{itemize}
\item $\D_{n,k}=0$ for $k\leq 0$ and $\D_{n,k}=\D$ for $k=1$,
\item $\mathsf{d}^v=\id$,
\item $\mathsf{d}^h\colon\D_{n,2}\to\D_{n,1}$ the morphism induced by the inverse image of $\emptyset\colon\bbone\to Sl_{n,2}$,
\item $\Phi_{n,k}=\id$.
Then it is straight forward to check that \autoref{thm:main} holds for $n+k\geq 1$.
\end{itemize}
\item Let $\#(Sl_{n,k})$ denote the cardinality of the set of injective objects in $Sl_{n,k}$. We invoke \autoref{egs:simplex-slice} and \autoref{prop:higher-cof} to see that $\#(Sl_{n,2})=n+1$ and $\#(Sl_{1,n+1})=n+1$, respectively. Moreover, the functors $d^v$ and $d^h$ show that 
\[
\#(Sl_{n+1,k+1})=\#(Sl_{n+1,k})+\#(Sl_{n,k+1}).
\]
Hence, by induction we conclude that $\#(Sl_{n,k})=\binom{n+k-1}{k-1}=\binom{n+k-1}{n}$. Therefore, the derivators $\D_{n,k}$ can be regarded as a categorification of Pascal's triangle, and the equivalences $\Phi_{n,k}$ as a categorification of the symmetry of Pascal's triangle.
\end{enumerate}
\end{rmk}

\begin{eg}\label{eg:D34}
Consider an object $X\in sl\D_{3,4}$ such that the underlying diagram of $X$ is of the form
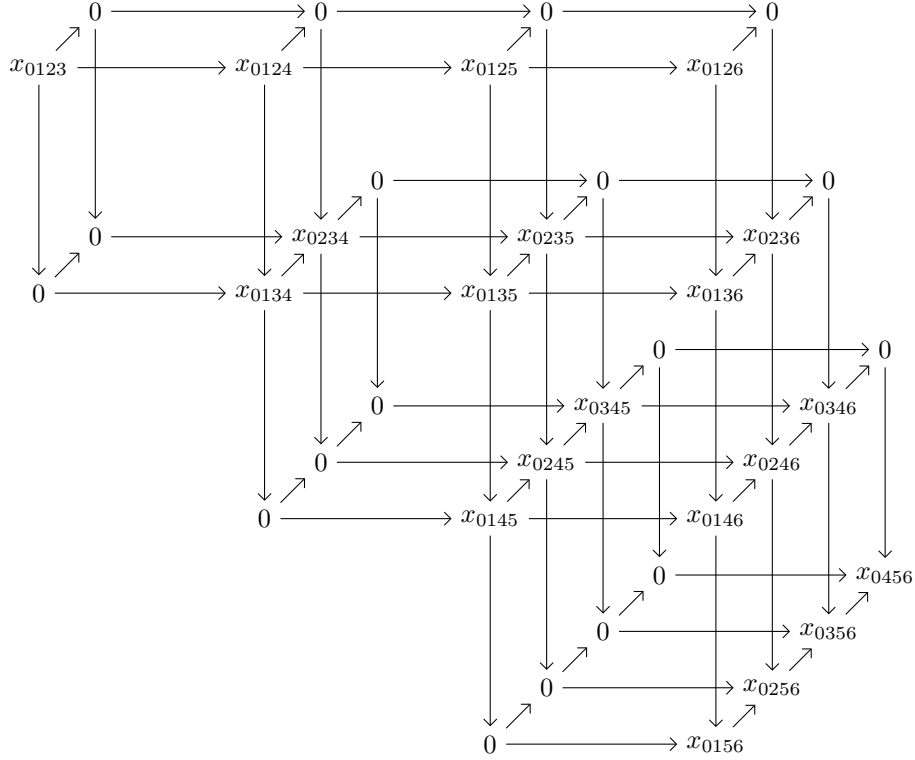
\[
\begin{tikzpicture}[scale=0.75]
\node (A) at (1,13) {$0$};
\node (B) at (5,13) {$0$};
\node (C) at (9,13) {$0$};
\node (D) at (13,13) {$0$};
\node (E) at (0,12) {$x_{0123}$};
\node (F) at (4,12) {$x_{0124}$};
\node (G) at (8,12) {$x_{0125}$};
\node (H) at (12,12) {$x_{0126}$};
\node (I) at (6,10) {$0$};
\node (J) at (10,10) {$0$};
\node (K) at (14,10) {$0$};
\node (L) at (1,9) {$0$};
\node (M) at (5,9) {$x_{0234}$};
\node (N) at (9,9) {$x_{0235}$};
\node (O) at (13,9) {$x_{0236}$};
\node (P) at (0,8) {$0$};
\node (Q) at (4,8) {$x_{0134}$};
\node (R) at (8,8) {$x_{0135}$};
\node (S) at (12,8) {$x_{0136}$};
\node (T) at (11,7) {$0$};
\node (U) at (15,7) {$0$};
\node (V) at (6,6) {$0$};
\node (W) at (10,6) {$x_{0345}$};
\node (X) at (14,6) {$x_{0346}$};
\node (Y) at (5,5) {$0$};
\node (Z) at (9,5) {$x_{0245}$};
\node (AA) at (13,5) {$x_{0246}$};
\node (AB) at (4,4) {$0$};
\node (AC) at (8,4) {$x_{0145}$};
\node (AD) at (12,4) {$x_{0146}$};
\node (AE) at (11,3) {$0$};
\node (AF) at (15,3) {$x_{0456}$};
\node (AG) at (10,2) {$0$};
\node (AH) at (14,2) {$x_{0356}$};
\node (AI) at (9,1) {$0$};
\node (AJ) at (13,1) {$x_{0256}$};
\node (AK) at (8,0) {$0$};
\node (AL) at (12,0) {$x_{0156}$};
\path[->,font=\scriptsize,>=angle 90]
(A) edge (B)
(B) edge (C)
(C) edge (D)
(E) edge (F)
(F) edge (G)
(G) edge (H)
(I) edge (J)
(J) edge (K)
(L) edge (M)
(M) edge (N)
(N) edge (O)
(P) edge (Q)
(Q) edge (R)
(R) edge (S)
(T) edge (U)
(V) edge (W)
(W) edge (X)
(Y) edge (Z)
(Z) edge (AA)
(AB) edge (AC)
(AC) edge (AD)
(AE) edge (AF)
(AG) edge (AH)
(AI) edge (AJ)
(AK) edge (AL)
(E) edge (P)
(A) edge (L)
(F) edge (Q)
(Q) edge (AB)
(B) edge (M)
(M) edge (Y)
(I) edge (V)
(G) edge (R)
(R) edge (AC)
(AC) edge (AK)
(C) edge (N)
(N) edge (Z)
(Z) edge (AI)
(J) edge (W)
(W) edge (AG)
(T) edge (AE)
(H) edge (S)
(S) edge (AD)
(AD) edge (AL)
(D) edge (O)
(O) edge (AA)
(AA) edge (AJ)
(K) edge (X)
(X) edge (AH)
(U) edge (AF)
(E) edge (A)
(P) edge (L)
(F) edge (B)
(Q) edge (M)
(M) edge (I)
(AB) edge (Y)
(Y) edge (V)
(G) edge (C)
(R) edge (N)
(N) edge (J)
(AC) edge (Z)
(Z) edge (W)
(W) edge (T)
(AK) edge (AI)
(AI) edge (AG)
(AG) edge (AE)
(H) edge (D)
(S) edge (O)
(O) edge (K)
(AD) edge (AA)
(AA) edge (X)
(X) edge (U)
(AL) edge (AJ)
(AJ) edge (AH)
(AH) edge (AF);
\end{tikzpicture}
\]
Let $f=(f_0,f_1,f_2,f_3)\in Sl_{3,4}$ and $i=\lbrace i_0,\cdots,i_j\rbrace\subset\mathbf{4}$. Consider the cube
\[
\cube{i_0\cdots i_j}_{f}\colon\cube{\#(i)}\to Sl_{3,4},(\delta_0,\cdots,\delta_j)\mapsto f+\sum_{l=0}^j\delta_le_{i_l},
\]
where $e_i$ denotes the $i$th basis vector. Moreover, we define
\[
F^{i_0\cdots i_j}_{f}=\tfib\circ(\cube{i_0\cdots i_j}_{f})^*(X).
\]
Then the underlying diagram of the object $\Psi_{3,4}(X)$ can be described as
\[
\begin{tikzpicture}[scale=0.75]
\node (A) at (1,13) {$0$};
\node (B) at (5,13) {$0$};
\node (C) at (9,13) {$0$};
\node (D) at (13,13) {$0$};
\node (E) at (0,12) {$\Omega^6F^{123}_{0345}$};
\node (F) at (4,12) {$\Omega^5F^{123}_{0245}$};
\node (G) at (8,12) {$\Omega^4F^{123}_{0235}$};
\node (H) at (12,12) {$\Omega^3F^{123}_{0234}$};
\node (I) at (6,10) {$0$};
\node (J) at (10,10) {$0$};
\node (K) at (14,10) {$0$};
\node (L) at (1,9) {$0$};
\node (M) at (5,9) {$\Omega^4F^{23}_{0145}$};
\node (N) at (9,9) {$\Omega^3F^{23}_{0135}$};
\node (O) at (13,9) {$\Omega^2F^{23}_{0134}$};
\node (P) at (0,8) {$0$};
\node (Q) at (4,8) {$\Omega^4F^{123}_{0145}$};
\node (R) at (8,8) {$\Omega^3F^{123}_{0135}$};
\node (S) at (12,8) {$\Omega^2F^{123}_{0134}$};
\node (T) at (11,7) {$0$};
\node (U) at (15,7) {$0$};
\node (V) at (6,6) {$0$};
\node (W) at (10,6) {$\Omega^2F^3_{0125}$};
\node (X) at (14,6) {$\Omega F^3_{0124}$};
\node (Y) at (5,5) {$0$};
\node (Z) at (9,5) {$\Omega^2F^{23}_{0125}$};
\node (AA) at (13,5) {$\Omega F^{23}_{0124}$};
\node (AB) at (4,4) {$0$};
\node (AC) at (8,4) {$\Omega^2F^{123}_{0125}$};
\node (AD) at (12,4) {$\Omega F^{123}_{0124}$};
\node (AE) at (11,3) {$0$};
\node (AF) at (15,3) {$x_{0123}$};
\node (AG) at (10,2) {$0$};
\node (AH) at (14,2) {$F^3_{0123}$};
\node (AI) at (9,1) {$0$};
\node (AJ) at (13,1) {$F^{23}_{0123}$};
\node (AK) at (8,0) {$0$};
\node (AL) at (12,0) {$F^{123}_{0123}$};
\path[->,font=\scriptsize,>=angle 90]
(A) edge (B)
(B) edge (C)
(C) edge (D)
(E) edge (F)
(F) edge (G)
(G) edge (H)
(I) edge (J)
(J) edge (K)
(L) edge (M)
(M) edge (N)
(N) edge (O)
(P) edge (Q)
(Q) edge (R)
(R) edge (S)
(T) edge (U)
(V) edge (W)
(W) edge (X)
(Y) edge (Z)
(Z) edge (AA)
(AB) edge (AC)
(AC) edge (AD)
(AE) edge (AF)
(AG) edge (AH)
(AI) edge (AJ)
(AK) edge (AL)
(E) edge (P)
(A) edge (L)
(F) edge (Q)
(Q) edge (AB)
(B) edge (M)
(M) edge (Y)
(I) edge (V)
(G) edge (R)
(R) edge (AC)
(AC) edge (AK)
(C) edge (N)
(N) edge (Z)
(Z) edge (AI)
(J) edge (W)
(W) edge (AG)
(T) edge (AE)
(H) edge (S)
(S) edge (AD)
(AD) edge (AL)
(D) edge (O)
(O) edge (AA)
(AA) edge (AJ)
(K) edge (X)
(X) edge (AH)
(U) edge (AF)
(E) edge (A)
(P) edge (L)
(F) edge (B)
(Q) edge (M)
(M) edge (I)
(AB) edge (Y)
(Y) edge (V)
(G) edge (C)
(R) edge (N)
(N) edge (J)
(AC) edge (Z)
(Z) edge (W)
(W) edge (T)
(AK) edge (AI)
(AI) edge (AG)
(AG) edge (AE)
(H) edge (D)
(S) edge (O)
(O) edge (K)
(AD) edge (AA)
(AA) edge (X)
(X) edge (U)
(AL) edge (AJ)
(AJ) edge (AH)
(AH) edge (AF);
\end{tikzpicture}
\]
\end{eg}

\begin{cor}\label{cor:contra-S}
Let $\D$ be a stable derivator and $k\geq 2$. Then there is a pseudofunctor
\[
\mathsf{S}^{\bullet}_{(k-1)}(\D)\colon\underline{\Lambda}^{op}\to Der
\]
satisfying the following properties
\begin{enumerate}
\item $\mathsf{S}^{\bullet}_{(k-1)}(\D)(\Lambda_m)=\D_{k-1,m-k+2}$,
\item $\mathsf{S}^{\bullet}_{(k-1)}(\D)(\mathrm{d}_i\colon\Lambda_m\to\Lambda_{m+1})\cong\mathsf{d}^h[2(i-k)]\colon\D_{k-1,m-k+3}\to\D_{k-1,m-k+2}$,
\item $\mathsf{S}^{\bullet}_{(k-1)}(\D)(\mathrm{s}_i\colon\Lambda_{m+1}\to\Lambda_{m})\cong\mathsf{d}^h[2(i-k)+1]\colon\D_{k-1,m-k+2}\to\D_{k-1,m-k+3}$,
\item there is a pseudonatural equivalence $\mathsf{S}_{\bullet}^{(k-1)}\to\mathsf{S}^{\bullet}_{(k-1)}$.
\end{enumerate}
\end{cor}

\begin{proof}
We apply \autoref{prop:conjugation} to the 2-functor $\mathsf{S}_{\bullet}^{(k-1)}$ and the set of equivalences $\mathsf{S}_{\Lambda_m}=\Phi_{m-k+1,k}\colon\D_{m-k+1,k}\to\D_{k-1,m-k+2}$ (\autoref{thm:main}). The equivalences on the values of 1-morphisms follow from \autoref{eg:ver-fd} and \autoref{thm:main} (iii) and (iv).
\end{proof}

In particular, we conclude the the generalized horizontal face and degeneracy morphisms satisfy the simplicial relations in the same way as the vertical structure morphisms. We can apply this to \autoref{cor:dh-comp}. In the following the functor $\tau\colon\Delta^{op}\times\Delta^{op}\to\Delta^{op}\times\Delta^{op}$ denotes the interchange of factors.

\begin{cor}\label{cor:bisimplicial}
Let $\D$ be a stable derivator. Then there is a pseudofunctor
\[
\mathsf{S}_{\bullet,\bullet}(\D)\colon\Delta^{op}\times\Delta^{op}\to Der
\]
satisfying the following properties
\begin{enumerate}
\item $\mathsf{S}_{\bullet,\bullet}(\D)(\Delta_n,\Delta_k)=\D_{n+1,k+2}$,
\item $\mathsf{S}_{\bullet,\bullet}(\D)(\mathrm{d}_i,\id)\cong\mathsf{d}^v[2i]$,
\item $\mathsf{S}_{\bullet,\bullet}(\D)(\mathrm{s}_i,\id)\cong\mathsf{d}^v[2i+1]$,
\item $\mathsf{S}_{\bullet,\bullet}(\D)(\id,\mathrm{d}_i)\cong\mathsf{d}^h[2i]$,
\item $\mathsf{S}_{\bullet,\bullet}(\D)(\id,\mathrm{s}_i)\cong\mathsf{d}^h[2i+1]$,
\item there is a peudonatural equivalence $\mathsf{S}_{\bullet,\bullet}(\D)\cong\mathsf{S}_{\bullet,\bullet}(\D)\circ\tau$.
\end{enumerate}
\end{cor}

\begin{proof}
By \autoref{cor:dh-comp} there are morphisms in $PsFun(\Delta^{op},Der)$ defined locally by $\mathsf{d}_{n,k+1}^h[a]$. These assemble into a pseudofunctor $\Delta^{op}\to PsFun(\Delta^{op},Der)$ by \autoref{cor:contra-S}.
\end{proof}

\begin{rmk}
Let $\D$ be a stable derivator. Then we call the structure defined by the derivators $\D_{n,k}$, the morphisms $\mathsf{s}_3, \mathsf{d}^v[a], \mathsf{d}^h[a]$ and $\Phi_{n,k}$, and the 2-morphisms defined by \autoref{thm:full-square} and \autoref{prop:dual-square} the bivariant parasimplicial $\mathsf{S}_{\bullet}$-con-struction.
\begin{enumerate}
\item It is clear that the bisimplicial object from \autoref{cor:bisimplicial} only provides a very coarse approximation of the bivariant parasimplicial $\mathsf{S}_{\bullet}$-construction, since we have discarded a lot of the structure morphisms.
\item In fact, our results be regarded as a first step toward a description of the bivariant parasimplicial $\mathsf{S}_{\bullet}$-construction as a derivator-valued presheaf on a 2-category $\tilde{\Lambda}$, which can be described as sub-2-category of the 2-category of pseudofunctors $PsFun(Der^{st},Der^{st})$ with
\begin{itemize}
\item objects, the pseudofunctors $\D\mapsto\D_{n,k}$ for $n+k\geq 1$,
\item morphisms, compositions of (elementary) pseudonatural transformations of the form $\mathsf{S}_{\bullet}^{(k)}(f)$ and $\mathsf{S}^{\bullet}_{(k)}(g)$ for $f,g$ morphisms in $\underline{\Lambda}$,
\item 2-morphisms, compositions of modifications between elementary 1-mor- phisms induced by 2-morphisms in $\underline{\Lambda}$ and the isomorphisms of \autoref{thm:full-square} and \autoref{prop:dual-square}.
\end{itemize}
\item From the above definition the morphism categories of $\tilde{\Lambda}$ are in general hard to understand. It should be an interesting problem to describe the 2-category $\tilde{\Lambda}$ purely combinatorially, since we expect this 2-category to encode further structures relevant for stable homotopy theory and representation theory in a systematic way.
\end{enumerate}
\end{rmk}

\section{Higher Toda brackets for derivators}
\label{sec:Toda}

In this section we discuss a first application of \autoref{thm:dual-Sdot} and \autoref{thm:dual-Serre} concerning higher Toda brackets. Recall that Toda brackets are operations defined on certain strings of composable morphisms in the homotopy category of a stable model category. However, in general Toda brackets are not defined for all such strings, and if they are, there is often a set of different values, i.e. they are only defined up to some indeterminacy. We will show that for a strong stable derivator $\D$ there is a functorial construction lifting the higher Toda brackets. The following definitions of filtered objects and Toda brackets are based on \cite[Appendix]{Shipley-rational} and \cite[3.3]{Sagave-universal}. To establish a relation to the theory of triangulated categories, we recall the notion of a strong stable derivator.

\begin{defn}
A stable derivator $\D$ is called \textbf{strong} if for every finite free category $A$ the underlying diagram functor
\[
\mathsf{dia}_A\colon\D(A)\to\D(\bbone)^A
\]
is an epivalence of categories (i.e. is full and essentially surjective).
\end{defn}

\begin{thm}
Let $\D$ be a strong stable derivator and $A\in Cat$. Then there is a canonical triangulation on $\D(A)$ defined by the suspension functor $\Sigma$ and the class of distinguished triangles which are isomorphic to underlying diagrams of 1-cofiber sequences.
\end{thm}

\begin{proof}
This is due to Maltsiniotis \cite{maltsiniotis:triang}, a published proof can be found in \cite[Thm.~4.16]{groth:ptstab}, the ideas go back at least to \cite{franke:adams}.
\end{proof}

\begin{defn}
Let $\mathcal{T}$ be a triangulated category and 
\[
x_{n-1}\xrightarrow{u_{n-1}}x_{n-2}\xrightarrow{u_{n-2}}\cdots\xrightarrow{u_{1}}x_0
\]
be an $(n-1)$-simplex in $\mathcal{T}$. An $n$\textbf{-filtered object} $y\in[u_1,\cdots,u_{n-1}]$ consists of an $n$-simplex
\[
y_0\xrightarrow{v_1}y_1\xrightarrow{v_2}\cdots\xrightarrow{v_n}y_n
\]
in $\mathcal{T}$, such that $y_0=0$, $y_n=y$, and choices of distinguished triangles
\begin{equation}\label{eq:filtration}
y_j\xrightarrow{v_{j+1}}y_{j+1}\xrightarrow{r_{j+1}}\Sigma^jx_j\xrightarrow{q_j}\Sigma y_j
\end{equation}
such that $\Sigma r_j\circ q_j=\Sigma^ju_j$. Moreover, the map $x_0=y_1\to y$ is denoted by $\sigma_y$.
\end{defn}

\begin{defn}\label{defn:Toda}
Let $\mathcal{T}$ be a triangulated category and 
\[
x_{n}\xrightarrow{u_{n}}x_{n-1}\xrightarrow{u_{n-1}}\cdots\xrightarrow{u_{1}}x_0
\]
an $n$-simplex in $\mathcal{T}$. A map $\Sigma^{n-2}x_n\xrightarrow{\gamma}x_0$ lies in the $n$\textbf{-fold Toda bracket} of the above sequence, if there is an $(n-1)$-filtered object $y\in[u_2,\cdots,u_{n-1}]$ and a decomposition  $\gamma\colon\Sigma^{n-2}x_n\xrightarrow{\gamma_n}y\xrightarrow{\gamma_0}x_0$ such that there is a commutative diagram
\[
\xymatrix{
&x_1\ar[d]_{\sigma_y}\ar[rd]^{u_1}\\
\Sigma^{n-2}x_n\ar[r]^{\gamma_n}\ar[rd]_{\Sigma^{n-2}u_n}&y\ar[r]_{\gamma_0}\ar[d]^{r_{n-1}}&x_0\\
&\Sigma^{n-2}x_{n-1}.
}
\]
\end{defn}

\begin{prop}\label{prop:filtered}
Let $n\geq 1$ and $\D$ a strong stable derivator. Let $X\in\D^{\cube{n}}_{\rightarrow}$. Then the underlying diagram of $(\mathrm{d}_0)_!\circ(\Psi_n^{\square})^{-1}(X)$ is canonically an $(n+1)$-filtered object of the underlying diagram of $\rightarrow^*(X)$.
\end{prop}

\begin{proof}
Let $X\in\D^{\cube{n}}_{\rightarrow}$. It follows from the proof of \autoref{thm:slices} and \autoref{prop:higher-cof} that the lower square in the diagram
\[
\xymatrix{
\D_{n,2}\ar[r]^{\sim}\ar[d]_{\Psi'_n}&sl\D_{n,2}\ar[r]^{\sim}&\D^{[n]}\ar[d]^{\Psi_n^{\square}}\\
\D_{1,n+1}\ar[rd]_{\sim}\ar[r]^{\sim}&sl\D_{1,n+1}\ar[r]^{\square_{\xi}^*}\ar[d]_{(sd^*)^{-1}}&\D^{\cube{n}}_{\rightarrow}\ar[d]^{\beta_!\circ\alpha_*}\\
&do\D_{1,n+1}\ar[r]_{\square_{\xi}^*}&\D^{\cube{n+1},ex}_{\rightarrow},
}
\]
where the vertical maps on the right are those from \autoref{prop:tcof-new}, commutes. The upper cell commutes by definition, and for the triangle we invoke \autoref{rmk:slices}. Let 
\begin{itemize}
\item $Y\in\D^{[n]}$,
\item $Z\in\D_{n,2}$,
\item $V\in\D^{\cube{n+1},ex}_{\rightarrow}$,
\item $W\in\D_{1,n+1}$
\end{itemize}
be objects corresponding to $X$. Furthermore, for a poset $A$ with $a\leq b\in A$ we use the notation $[a,b]\colon[1]\to A$ for the functor $0\mapsto a$, $1\mapsto b$. Then \autoref{prop:tcof-new} induces the second isomorphism in
\begin{align}
&[\mathsf{s}_3^n(\xi),\mathsf{s}_3^{n+1}(\xi)]^*(W)\\
\cong& [\rightarrow(n),\rightarrow(n+1)]^*(V)\\
\cong&[\emptyset,\infty]^*\circ\cof^{\underline{1}}(X)\label{eq:filt}\\
\cong&[0,n]^*(Y)\\
\cong&[(0,1),(0,n+1)]^*(Z).
\end{align}
On the other hand, we have for $0\leq i\leq n-1$
\begin{align}
&[\rightarrow(i),\rightarrow(i+1)]^*(X)\\
\cong &[\mathsf{s}_3^i(\xi),\mathsf{s}_3^{i+1}(\xi)]^*(W)\\
\cong &[\mathsf{s}_3^n(\xi),\mathsf{s}_3^{n+1}(\xi)]^*((\mathsf{s}_3^{i-n})^*(W))\\
\cong &[(0,1),(0,n+1)]^*((\Psi'_n\circ\mathsf{s}_3^{i-n})^*(Z))\label{eq:filtr2}\\
\cong &[(0,1),(0,n+1)]^*((\mathsf{s}_3^{i-n})^*(Z))\\
\cong &\Omega^{n-i}\circ[(0,1),(0,n+1)]^*((\mathsf{s}_1^{n-i})^*(Z))\\
\cong &\Omega^{n-i}\circ[(n-i,n-i+1),(n-i,2n-i+1)]^*(Z),
\end{align}
where the third isomorphism is \eqref{eq:filt}, the fourth isomorphisms is \autoref{thm:dual-Serre} and the fifth isomorphism follows from \autoref{cor:dnk-sigma} and the definition of $\mathsf{s}_3$. 

\noindent We observe that the underlying diagram of $(\mathrm{d}_0)_!(Y)$ can be identified the following restriction of the underlying diagram of $Z$
\[
(0,0)^*(Z)\xrightarrow{v_1}(0,1)^*(Z)\xrightarrow{v_2}\cdots\xrightarrow{v_{n+1}}(0,n+1)^*(Z).
\]
In the next step, we consider restrictions of the underlying diagram of $Z$ along inclusions of the form
\[
\xymatrix{
(0,i)\ar[r]\ar[d]&(0,i+1)\ar[r]\ar[d]&(0,n+2)\ar[d]\\
(i,i)\ar[r]&(i,i+1)\ar[r]\ar[d]&(i,n+2)\ar[r]\ar[d]&(i,n+i+2)\ar[d]\\
&(i+1,i+1)\ar[r]&(i+1,n+2)\ar[r]&(i+1,n+i+2)
}
\]
We invoke the properties (P1) and (P2) for objects in $\D_{1,n+1}$ and \autoref{cor:bicart-concat} to conclude that we obtain distinguished triangles
\[
\xymatrix{
(0,i)^*(Z)\ar[r]^{v_{i+1}}\ar[d]&(0,i+1)^*(Z)\ar[r]\ar[d]^{r_{i+1}}&0\ar[d]\\
0\ar[r]&\Sigma^ix_i\ar[r]^{q_i}\ar[d]&\Sigma(0,i)^*(Z)\ar[r]\ar[d]^{-\Sigma v_{i+1}}&0\ar[d]\\
&0\ar[r]&\Sigma(0,i+1)^*(Z)\ar[r]^{-\Sigma r_{i+1}}&\Sigma^{i+1}x_i,
}
\]
where we have used \eqref{eq:filtr2} to identify $(i,i+1)^*(Z)\cong\Sigma^i\circ\rightarrow(n-i)^*(X)=\Sigma^ix_i$. Finally, we use the restriction of the underlying diagram along 
\[
\xymatrixcolsep{0.66cm}
\xymatrix{
(0,i-1)\ar[r]\ar[d]&(0,i)\ar[r]\ar[d]&(0,i+1)\ar[r]\ar[dd]&(0,n+2)\ar[d]\\
(i-1,i-1)\ar[r]&(i-1,i)\ar[rr]\ar[d]&&(i-1,n+2)\ar[d]\ar[r]&(i-1,n+i+1)\ar[d]\\
&(i,i)\ar[r]&(i,i+1)\ar[r]&(i,n+2)\ar[r]&(i,n+i+1),
}
\]
which consequently gives rise to the diagram
\[
\xymatrixcolsep{0.67cm}
\xymatrix{
(0,i-1)^*(Z)\ar[r]^{v_i}\ar[d]&(0,i)^*(Z)\ar[r]^{v_{i-1}}\ar[d]^{r_i}&(0,i+1)^*(Z)\ar[r]\ar[dd]^(.35){r_{i+1}}&0\ar[d]\\
0\ar[r]&\Sigma^{i-1}x_{i-1}\ar[rr]^(.35){q_{i-1}}\ar[d]&&\Sigma(0,i-1)^*(Z)\ar[d]^{-\Sigma v_i}\ar[r]&0\ar[d]\\
&0\ar[r]&\Sigma^{i-1}x_{i-1}\ar[r]^{q_i}&\Sigma(0,i)^*(Z)\ar[r]^{-\Sigma r_i}&\Sigma^ix_{i-1}.
}
\]
We invoke \eqref{eq:filtr2} again to identify $-\Sigma r_i\circ q_i=\Sigma_iu_i$. For $i$ even, we replace the distinguished triangles
\[
(0,i)^*(Z)\xrightarrow{v_{i+1}}(0,i+1)^*(Z)\xrightarrow{r_{i+1}}\Sigma^ix_i\xrightarrow{q_i}\Sigma(0,i)^*(Z)
\]
with the isomorphic, and hence also distinguished triangles
\[
(0,i)^*(Z)\xrightarrow{v_{i+1}}(0,i+1)^*(Z)\xrightarrow{-r_{i+1}}\Sigma^ix_i\xrightarrow{-q_i}\Sigma(0,i)^*(Z),
\]
and conclude that we indeed have constructed a filtration.
\end{proof}

\begin{eg}
We explain the procedure, in the case $n=3$. Let $X\in\D^{\cube{n}}_{\rightarrow}$ with underlying diagram
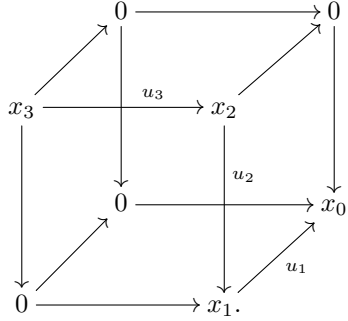
\[
\xymatrix{
&0\ar[rr]\ar[dd]&&0\ar[dd]\\
x_3\ar[ur]\ar[rr]^(.65){u_3}\ar[dd]&&x_2\ar[ur]\ar[dd]^(.35){u_2}\\
&0\ar[rr]&&x_0\\
0\ar[ur]\ar[rr]&&x_1\ar[ur]_{u_1}.
}
\]
Then $(\Psi_3^{\square})^{-1}(X)$ is an object in $\D^{[3]}$ with underlying diagram
\begin{equation}\label{eq:filtr-delta3}
y_1\xrightarrow{v_2}y_2\xrightarrow{v_3}y_3\xrightarrow{v_4}y_4.
\end{equation}
Using \autoref{eg:explicit}, we can extend $(\Psi_3^{\square})^{-1}(X)$ to an object of $\D_{3,2}$ with underlying diagram
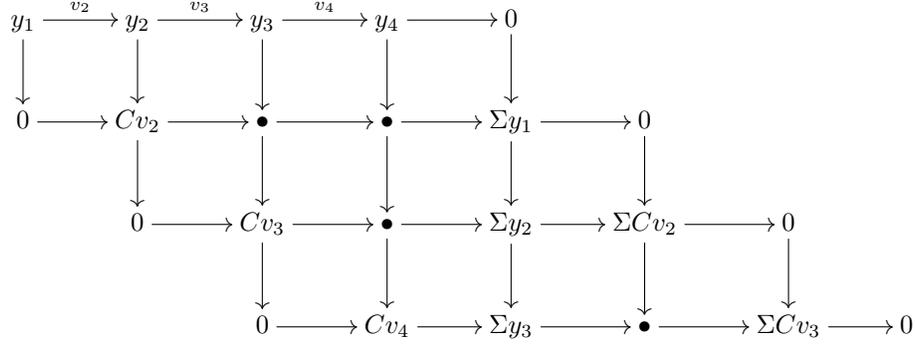
\[
\xymatrix{
y_1\ar[r]^{v_2}\ar[d]&y_2\ar[r]^{v_3}\ar[d]&y_3\ar[r]^{v_4}\ar[d]&y_4\ar[r]\ar[d]&0\ar[d]\\
0\ar[r]&Cv_2\ar[r]\ar[d]&\bullet\ar[r]\ar[d]&\bullet\ar[r]\ar[d]&\Sigma y_1\ar[r]\ar[d]&0\ar[d]\\
&0\ar[r]&Cv_3\ar[r]\ar[d]&\bullet\ar[r]\ar[d]&\Sigma y_2\ar[r]\ar[d]&\Sigma Cv_2\ar[r]\ar[d]&0\ar[d]\\
&&0\ar[r]&Cv_4\ar[r]&\Sigma y_3\ar[r]&\bullet\ar[r]&\Sigma Cv_3\ar[r]&0\\
}
\]
We invoke \eqref{eq:filtr2} to identify the composition of the central three morphisms in the $(i+1)$-st line with $\Sigma^iu_u\colon\Sigma^ix_i\to\Sigma^ix_{i-1}$. Furthermore, we indicate the distinguished triangles, which exhibit \eqref{eq:filtr-delta3} as a filtration.
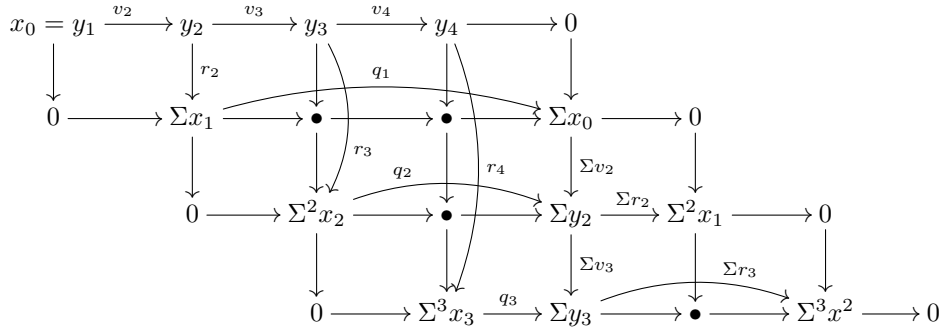
\[
\xymatrixrowsep{0.75cm}
\xymatrixcolsep{0.75cm}
\xymatrix{
x_0=y_1\ar[r]^{v_2}\ar[d]&y_2\ar[r]^{v_3}\ar[d]^{r_2}&y_3\ar@/^1pc/[dd]^(.65){r_3}\ar[r]^{v_4}\ar[d]&y_4\ar@/^1pc/[ddd]^{r_4}\ar[r]\ar[d]&0\ar[d]\\
0\ar[r]&\Sigma x_1\ar@/^1pc/[rrr]^{q_1}\ar[r]\ar[d]&\bullet\ar[r]\ar[d]&\bullet\ar[r]\ar[d]&\Sigma x_0\ar[r]\ar[d]^{\Sigma v_2}&0\ar[d]\\
&0\ar[r]&\Sigma^2 x_2\ar@/^1pc/[rr]^(.35){q_2}\ar[r]\ar[d]&\bullet\ar[r]\ar[d]&\Sigma y_2\ar[r]^{\Sigma r_2}\ar[d]^{\Sigma v_3}&\Sigma^2 x_1\ar[r]\ar[d]&0\ar[d]\\
&&0\ar[r]&\Sigma^3 x_3\ar[r]^{q_3}&\Sigma y_3\ar[r]\ar@/^1pc/[rr]^(.65){\Sigma r_3}&\bullet\ar[r]&\Sigma^3 x^2\ar[r]&0\\
}
\]
\end{eg}

\begin{defn}
Let $n\geq 3$ and $\D$ a stable derivator. 
\begin{enumerate}
\item Let $\rightarrow^t_n\colon[n]\to\cube{n-2}\times[2]$ be the functor defined by 
\begin{itemize}
\item $\rightarrow^t_n(0)=(\rightarrow(0),0)$,
\item $\rightarrow^t_n(i)=(\rightarrow(i-1),1)$ for $1\leq i\leq n-1$,
\item $\rightarrow^t_n(n)=(\rightarrow(n-2),2)$.
\end{itemize}
Let $\D^{\mathsf{T}_n}$ be the full subderivator of $\D^{\cube{n-2}\times[2]}$ spanned by those objects $X$ such that $M^*X=0$ for all $M\in\cube{n-2}\times[2]$ such that $M$ is not in the image of $\rightarrow^t_n$. The derivator $\D^{\mathsf{T}_n}$ is called the \textbf{derivator of} $n$\textbf{-fold Toda bracket data}.
\item Let $\mathrm{d}_1^{n-2}\colon\cube{n-2}\to\cube{n-1}_{0,n-2}$ be the inclusion of the 0-face with respect to the $(n-2)$-nd coordinate. Moreover, let $e\colon[1]\to\cube{n-1}\times[2]$ be the functor classifying $(\rightarrow(n-1),0)\to(\rightarrow(n-1),2)$.
Then the $n$\textbf{-fold Toda bracket morphism} for $\D$ is defined by the composition
\[
\mathsf{Toda}_n:=e^*\circ(\iota_{0,n-2}\times\id_{[2]})_!\circ(\mathrm{d}_1^{n-2}\times\id_{[2]})_*\colon\D^{\mathsf{T}_n}\to\D^{[1]}.
\]
\end{enumerate}
\end{defn}

\begin{thm}\label{thm:toda}
Let $n\geq 3$ and $\D$ a strong stable derivator. Let $X\in\D^{\mathsf{T}_n}(\bbone)$. Then the underlying diagram of $\mathsf{Toda}_n(X)$ lies in the $n$-fold Toda bracket of the underlying diagram of $(\rightarrow^t_n)^*(X)$.
\end{thm}

\begin{proof}
We define $X^1=(\id_{\cube{n-2}}\times 1)^*(X)\in\D^{\cube{n-2}}_{\rightarrow}$. Let the underlying diagram of $(\rightarrow^t_n)^*(X)$ be of the form
\[
x_{n}\xrightarrow{u_{n}}x_{n-1}\xrightarrow{u_{n-1}}\cdots\xrightarrow{u_{1}}x_0
\]
Consider the factorization $\cube{n-1}\xrightarrow{\alpha}(\cube{n-2}\times[2])\setminus\infty\xrightarrow{\beta}\cube{n-2}\times[2]$ of the inclusion $\id_{\cube{n-1}}\times\mathrm{d}_2$. Next, we consider $T:=\beta_!\circ\alpha_*\circ\mathsf{Toda}_n(X)$ pass to the inverse image defined by the following diagram in $\cube{n-2}\times[2]\times[2]$
\begin{equation}\label{eq:Toda}
\xymatrix{
&(\rightarrow(n-3),1,1)\ar[r]\ar[d]&(\rightarrow(n-3),1,2)\ar[d]\\
(\rightarrow(n-2),1,0)\ar[r]\ar[d]&(\rightarrow(n-2),1,1)\ar[r]\ar[d]&(\rightarrow(n-2),1,2)\\
(\rightarrow(n-2),2,0)\ar[r]&(\rightarrow(n-2),2,1)
}
\end{equation}
Since $\alpha$ and $\beta$ are fully faithful, we invoke \cite[Prop.~1.20]{groth:ptstab} to conclude that the underlying diagram of $\mathsf{Toda}_n(X)$ is obtained by restricting further to the middle row in \eqref{eq:Toda}. We claim that the inverse image of $T$ along \eqref{eq:Toda} is of the following form
\begin{equation}\label{eq:Toda-II}
\xymatrix{
&x_1\ar[r]^{u_1}\ar[d]&x_0\ar[d]^{\id}\\
\Sigma^{n-2}x_n\ar[r]\ar[d]_{\id}&y_{n-1}\ar[r]\ar[d]&x_0\\
\Sigma^{n-2}x_n\ar[r]_{\Sigma^{n-2}u_n}&\Sigma^{n-2}x_{n-1}
}
\end{equation}
To establish this, we use that the subcubes $\square_{i,j}^*(T)$ for 
\[
\square_{i,j}\colon\cube{n-1}\xrightarrow{(\id_{\cube{n-2}}\times\mathrm{d}_i)}\cube{n-2}\times[2]\xrightarrow{(\id_{\cube{n-2}\times[2]}\times j)}\cube{n-2}\times[2]\times[2]
\]
are bicartesian for $i,j\in\lbrace 0,1,2\rbrace$ by \autoref{prop:frankes-lemma} and \autoref{cor:bicart-concat}.  As a consequence, we conclude the following identifications.
\begin{itemize}
\item The bicatesianess of $\square_{2,0}^*(T)$ implies the relation $(\rightarrow(n-2),1,0)^*(T)\cong\Sigma^{n-2}x_n$ since (using fully faithfulness of Kan extensions again) the restriction of this cube to $\cube{n-1}_{0,n-2}$ is concentrated at the initial vertex with value $x_n$.
\item Building on this, and using that $\alpha_!$ is an extension-by-zero morphism \cite[Prop.~1.23]{groth:ptstab}, we invoke \cite[Thm.~8.11]{gst:basic} applied to the bicartesian cube $\square_{0,0}^*(T)$ to see that the left vertical morphism in \eqref{eq:Toda-II} is the identity.
\item The same argument applied to the bicartesian cube $\square_{2,2}^*(T)$ yields the identity on the right vertical morphism in \eqref{eq:Toda-II}.
\item The identification of the top row is \cite[Prop.~1.20]{groth:ptstab} applied to the four Kan extensions in the construction of $T$ (which was used implicitly before).
\item For the identification of the bottom row, we use additionally the bicartesianess of the cubes $\square_{1,0}^*(T)$ and $\square_{1,1}^*(T)$, and the observation, that the restrictions of these cubes to $\cube{n-1}_{0,n-2}$ are concentrated at the initial vertex.
\end{itemize}
Finally, we use the bicartesianess of $\square_{2,1}^*(T)$, and invoke \autoref{prop:tcof-new} to identify the top map in the central column of \eqref{eq:Toda-II} with the canonical map from $\infty^*(X^1)\to\tcof(X^1)$. Finally, \autoref{prop:filtered} yields that this map admits a filtration $y=y_{n-1}=\tcof{X^1}$ such that $y_{n-1}\in[u_2,\cdots,u_n-1]$ and use the dual of \autoref{prop:tcof-new} combined with the bicartesianess of $\square_{0,1}^*(T)$ to identify the remaining map in \eqref{eq:Toda-II} (which is a suspension of the canonical map $\tfib(X^1)\to\emptyset^*(X^1)$) with the corresponding cofiber in the filtration.
\end{proof}

\begin{eg}\label{eg:Toda-4}
We illustrate the construction of Toda bracket morphisms in the case $n=4$.
\begin{enumerate}
\item Let $X\in\D^{\mathsf{T}_4}$ with underlying diagram
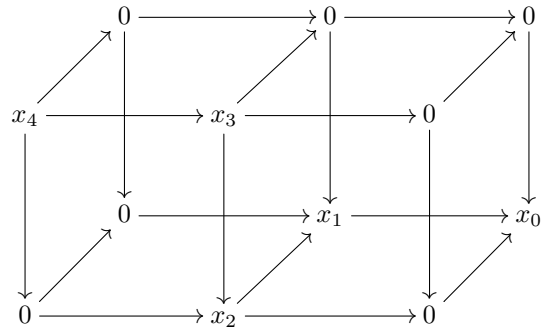
\[
\xymatrix{
&0\ar[rr]\ar[dd]&&0\ar[rr]\ar[dd]&&0\ar[dd]\\
x_4\ar[ur]\ar[rr]\ar[dd]&&x_3\ar[ur]\ar[rr]\ar[dd]&&0\ar[ur]\ar[dd]\\
&0\ar[rr]&&x_1\ar[rr]&&x_0\\
0\ar[ur]\ar[rr]&&x_2\ar[ur]\ar[rr]&&0\ar[ur]\\
}
\]
\item We extend this diagram via the right Kan extension $(\mathrm{d}_1^{2}\times[2])_*$ (dashed arrows), which is an extension-by-zero morphism, and the left Kan extension $(\iota_{0,2}\times[2])_!$ (dotted arrows), which 'adds bicartesian cubes'. Moreover, we have omitted all arrows in the direction of the last coordinate (which was displayed in the horizontal direction in the diagram above), with the exception of the morphisms which give rise to the Toda bracket. These are the curved arrows in the diagram below.
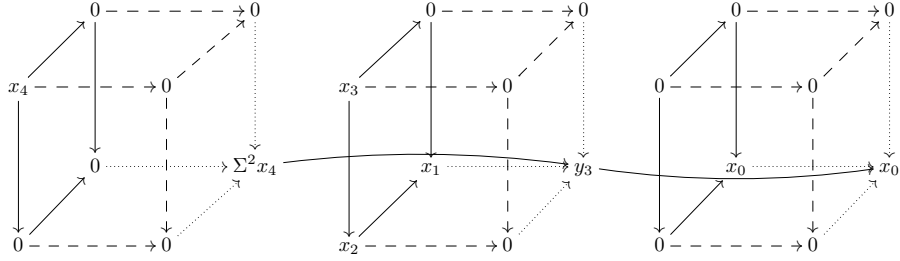
\[
\resizebox{.95\hsize}{!}{$
\xymatrix{
&0\ar@{-->}[rr]\ar[dd]&&0\ar@{.>}[dd]&&0\ar@{-->}[rr]\ar[dd]&&0\ar@{.>}[dd]&&0\ar@{-->}[rr]\ar[dd]&&0\ar@{.>}[dd]\\
x_4\ar[ur]\ar@{-->}[rr]\ar[dd]&&0\ar@{-->}[ur]\ar@{-->}[dd]&&x_3\ar[ur]\ar@{-->}[rr]\ar[dd]&&0\ar@{-->}[ur]\ar@{-->}[dd]&&0\ar[ur]\ar@{-->}[rr]\ar[dd]&&0\ar@{-->}[ur]\ar@{-->}[dd]\\
&0\ar@{.>}[rr]&&\Sigma^2 x_4\ar@/^/[rrrr]&&x_1\ar@{.>}[rr]&&y_3\ar@/_/[rrrr]&&x_0\ar@{.>}[rr]&&x_0\\
0\ar[ur]\ar@{-->}[rr]&&0\ar@{.>}[ur]&&x_2\ar[ur]\ar@{-->}[rr]&&0\ar@{.>}[ur]&&0\ar[ur]\ar@{-->}[rr]&&0\ar@{.>}[ur]\\
}$}
\]
\item In the next step, we indicate the construction of the object $T$ from the proof of \autoref{thm:toda}. Similar to the diagram above, we visualize the effect of the right Kan extension $\alpha_*$, which is an extension-by-zero morphism, with dashed arrows, and the effect of the left Kan extension $\beta_!$, which 'adds bicartesian cubes, with dotted arrows. Finally, we display those morphisms in the direction of the last coordinate, which give rise to \eqref{eq:Toda-II}, by curved arrows
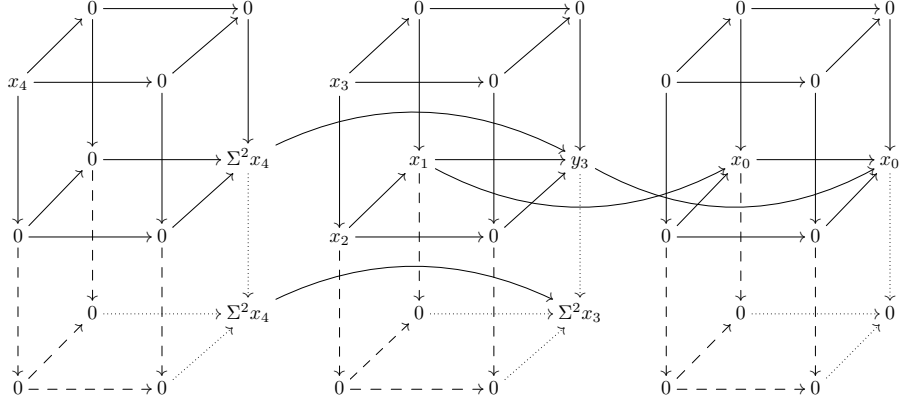
\[
\resizebox{.95\hsize}{!}{$
\xymatrix{
&0\ar[rr]\ar[dd]&&0\ar[dd]&&0\ar[rr]\ar[dd]&&0\ar[dd]&&0\ar[rr]\ar[dd]&&0\ar[dd]\\
x_4\ar[ur]\ar[rr]\ar[dd]&&0\ar[ur]\ar[dd]&&x_3\ar[ur]\ar[rr]\ar[dd]&&0\ar[ur]\ar[dd]&&0\ar[ur]\ar[rr]\ar[dd]&&0\ar[ur]\ar[dd]\\
&0\ar[rr]\ar@{-->}[dd]&&\Sigma^2 x_4\ar@/^2pc/[rrrr]\ar@{.>}[dd]&&x_1\ar@/_2pc/[rrrr]\ar[rr]\ar@{-->}[dd]&&y_3\ar@/_2pc/[rrrr]\ar@{.>}[dd]&&x_0\ar[rr]\ar@{-->}[dd]&&x_0\ar@{.>}[dd]\\
0\ar[ur]\ar[rr]\ar@{-->}[dd]&&0\ar[ur]\ar@{-->}[dd]&&x_2\ar[ur]\ar[rr]\ar@{-->}[dd]&&0\ar[ur]\ar@{-->}[dd]&&0\ar[ur]\ar[rr]\ar@{-->}[dd]&&0\ar[ur]\ar@{-->}[dd]\\
&0\ar@{.>}[rr]&&\Sigma^2 x_4\ar@/^2pc/[rrrr]&&0\ar@{.>}[rr]&&\Sigma^2 x_3&&0\ar@{.>}[rr]&&0\\
0\ar@{-->}[ur]\ar@{-->}[rr]&&0\ar@{.>}[ur]&&0\ar@{-->}[ur]\ar@{-->}[rr]&&0\ar@{.>}[ur]&&0\ar@{-->}[ur]\ar@{-->}[rr]&&0\ar@{.>}[ur]\\
}$}
\]

\end{enumerate}
\end{eg}

We conclude this section by providing an alternative construction of the derivator Toda bracket.
For this let $\alpha\colon[1]\times\cube{n-2}\times[2]\setminus\lbrace(0,\infty,1),(0,\infty,2)\rbrace\to\cube{n-2}\times[2]$ be the functor induced by $\mathrm{s}_0\times\id\times\id$ and $\beta\colon[1]\times\cube{n-2}\times[2]\setminus\lbrace{(0,\infty,1),(0,\infty,2)}\to[1]\times\cube{n-2}\times[2]$ be the inclusion. Finally, let $\gamma=\id\times\infty\times\infty\colon[1]\to[1]\times\cube{n-2}\times[2]$.

\begin{cor}\label{cor:Toda-alt}
Let $n\geq 3$ and $\D$ a strong stable derivator. Let $X\in\D^{\mathsf{T}_n}(\bbone)$. Then the underlying diagram of $\gamma^*\circ\beta_!\circ\alpha^*(X)$ lies in the $n$-fold Toda bracket of the underlying diagram of $(\rightarrow_n^t)^*(X)$
\end{cor}

\begin{proof}
By construction, $\beta_!\circ\alpha^*(X)$ can be considered as an object in $(\D^{[1]})^{\mathsf{T}_n}$. We consider the object
\[
Y=(\id\times\iota_{0,n-2}\times\id)_!\circ(\id\times\mathrm{d}_1^{n-2}\times\id)_*\circ\beta_!\circ\alpha^*(X)\in\D^{[1]\times\cube{n-1}\times[2]}(\bbone).
\]
Let $\delta=\id\times\infty\times\mathrm{d}_1\colon[1]\times[1]\to[1]\times\cube{n-1}\times[2]$. Then,
\begin{itemize}
\item $(\mathrm{d}_0\times\id)^*\circ\delta^*(Y)\cong\mathsf{Toda}_n(X)$,
\item $(\id\times\mathrm{d}_0)^*\circ\delta^*(Y)\cong\gamma^*\circ\beta_!\circ\alpha^*(X)$.
\end{itemize}
Hence, it is sufficient to show that $(\mathrm{d}_1\times\id)^*\circ\delta^*(Y)$ and $(\id\times\mathrm{d}_1)^*\circ\delta^*(Y)$ are constant.
\begin{enumerate}
\item It follows from \autoref{prop:frankes-lemma} and \autoref{cor:bicart-concat} that the $(n-1)$-cubes $(\mathrm{d}_0\times\id\times\mathrm{d}_i)^*\circ\beta_!\circ\alpha^*(X)$ are bicartesian for $0\leq i\leq 2$. In particular, $Z=(\mathrm{d}_0\times\id\times\mathrm{d}_1)^*\circ\beta_!\circ\alpha^*(X)$ is bicartesian and by construction $\iota_{1,n-2}^*(Z)=0$. Furthermore, $\emptyset^*(Z)=\emptyset^*(X)$. We conclude, that 
\[
C\circ(\mathrm{d}_1\times\id)^*\circ\delta^*(Y)=\tcof(Z)=0.
\]
Hence $(\mathrm{d}_1\times\id)^*\circ\delta^*(Y)$ is constant.
\item On the other hand we compute
\begin{align}
&(\id\times\mathrm{d}_1)^*\circ\delta^*(Y)\\
\cong &(\id\times\infty\times 0)^*\circ(\id\times\iota_{0,n-2}\times\id)_!\circ(\id\times\mathrm{d}_1^{n-2}\times\id)_*\circ\beta_!\circ\alpha^*(X)\\
\cong &(\id\times\infty)^*\circ(\id\times\iota_{0,n-2})_!\circ(\id\times\mathrm{d}_1^{n-2})_*\circ(\id\times\id\times 0)^*\circ\beta_!\circ\alpha^*(X)\\
\cong &(\id\times\infty)^*\circ(\id\times\iota_{0,n-2})_!\circ(\id\times\mathrm{d}_1^{n-2})_*\circ(\mathrm{s}_0\times\id)^*\circ(\id\times 0)^*(X)\\
\cong &\mathrm{s}_0^*\circ\infty^*\circ(\iota_{0,n-2})_!\circ(\mathrm{d}_1^{n-2})_*\circ(\id\times 0)^*(X),
\end{align}
which yields the constantness of $(\id\times\mathrm{d}_1)^*\circ\delta^*(Y)$.
\end{enumerate}
\end{proof}

\begin{rmk}
We illustrate also the construction of \autoref{cor:Toda-alt} in the case of 4-fold Toda brackets. Recall the object $\in\D^{\mathsf{T}_4}$ from \autoref{eg:Toda-4}. Then $\alpha^*(X)$ is obtained from $(\mathrm{s}_0\times\id)^*(X)$ by discarding two vertices. In the following diagram $X$ corresponds to the squares in the back
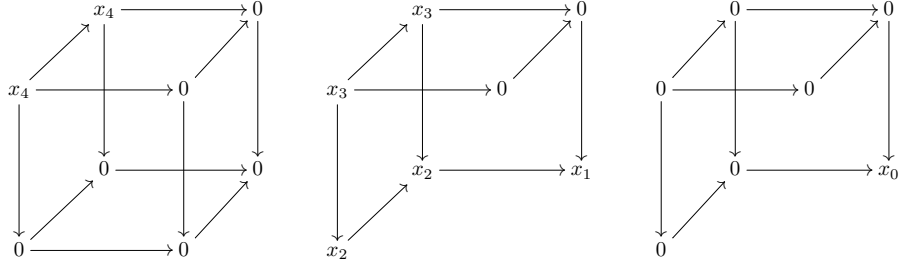
\[
\resizebox{.95\hsize}{!}{$
\xymatrix{
&x_4\ar[rr]\ar[dd]&&0\ar[dd]&&x_3\ar[rr]\ar[dd]&&0\ar[dd]&&0\ar[rr]\ar[dd]&&0\ar[dd]\\
x_4\ar[ur]\ar[rr]\ar[dd]&&0\ar[ur]\ar[dd]&&x_3\ar[ur]\ar[rr]\ar[dd]&&0\ar[ur]&&0\ar[ur]\ar[rr]\ar[dd]&&0\ar[ur]\\
&0\ar[rr]&&0&&x_2\ar[rr]&&x_1&&0\ar[rr]&&x_0\\
0\ar[ur]\ar[rr]&&0\ar[ur]&&x_2\ar[ur]&&&&0\ar[ur].
}$}
\]
Moreover, the object $\beta\circ\alpha^*(X)$ is obtained from the above by completing the front part of the diagram to a concatination of to bicartesian 3-cubes, as indicated by the dashed arrows, together with the induced maps to the back part of the diagram, as indicated by the dotted arrows below
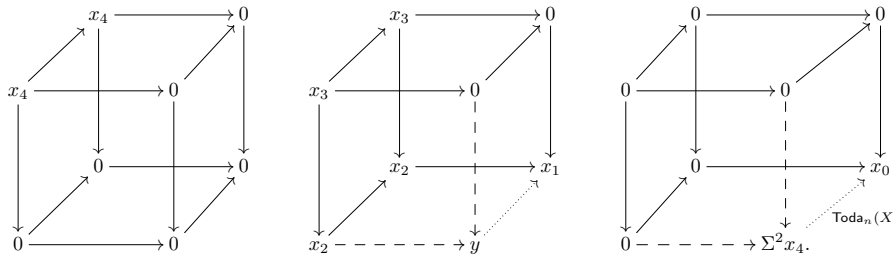
\[
\resizebox{.95\hsize}{!}{$
\xymatrix{
&x_4\ar[rr]\ar[dd]&&0\ar[dd]&&x_3\ar[rr]\ar[dd]&&0\ar[dd]&&0\ar[rr]\ar[dd]&&0\ar[dd]\\
x_4\ar[ur]\ar[rr]\ar[dd]&&0\ar[ur]\ar[dd]&&x_3\ar[ur]\ar[rr]\ar[dd]&&0\ar[ur]\ar@{-->}[dd]&&0\ar[ur]\ar[rr]\ar[dd]&&0\ar[ur]\ar@{-->}[dd]\\
&0\ar[rr]&&0&&x_2\ar[rr]&&x_1&&0\ar[rr]&&x_0\\
0\ar[ur]\ar[rr]&&0\ar[ur]&&x_2\ar@{-->}[rr]\ar[ur]&&y\ar@{.>}[ur]&&0\ar[ur]\ar@{-->}[rr]&&\Sigma^2x_4.\ar@{.>}[ur]_{\mathsf{Toda}_n(X)}
}$}
\]
The derivator Toda bracket of $X$ is now obtained as the map in the lower right of the diagram above.
We observe, that the front part of the diagram is a 3-cofiber sequence. Hence the cone of the n-fold Toda bracket provides a measure how far away a Toda bracket datum is from being an $(n-1)$-cofiber sequence.
\end{rmk}

\section{Higher triangulations}
\label{sec:triangles}

In this section we analyze the question, in which way the structure of the bivariant $\mathsf{S}_{\bullet}$-construction can be used to construct higher analogues of (strong) triangulations on a stable derivator $\D$. We begin by recalling the construction of strong triangulations using the covariant $\mathsf{S}_{\bullet}$-construction following Groth--{\v S}{\v t}ov{\'\i}{\v c}ek. The key step in this construction is based on the epivalence of the underlying diagram functor $\D([n])\to\D(\bbone)^{[n]}$. However, this property fails drastically if one replaces $[n]$ by $[n]^k$ for $k\geq 2$. Therefore, a satisfying axiomatization of the calculus of higher cofiber sequences becomes impossible on the level of underlying homotopy categories. Instead of this we describe how the results of the previous sections lead to a generalization of strong triangulations for coherent diagrams. Finally, we indicate some relations to the notion of $n$-angulated categories \cite{GKO-angulated} and cluster tilting theory.

\begin{defn}\label{defn:triangle}
Let $T$ be an additive category and $\Sigma\colon T\toiso T$ be an automorphism and $n\geq 1$. An $n$-\textbf{triangle} $(F,\phi)$ in $T$ consists of
\begin{enumerate}
\item a functor $F\colon\underline{\Lambda}_{n+1,1}\to T$ such that $F(f)\cong 0$ for $f\in\underline{\Lambda}_{n+1,1}$ non-injective.
\item a natural isomorphism $\phi\colon F\circ\mathsf{s}_2\toiso\Sigma\circ F$.
\end{enumerate}
A \textbf{morphism} of $n$-triangles $\psi\colon(F_1,\phi_1)\to(F_2,\phi_2)$ is a natural transformation $\psi\colon F_1\to F_2$ such that $(\Sigma\circ\psi)\circ\phi_1=\phi_2\circ(\psi\circ\mathsf{s}_2)$.
\end{defn}

\begin{egs}
Let $m,n\geq 1$, $T$ be an additive category and $\Sigma\colon T\toiso T$ be an automorphism.
\begin{enumerate}
\item Let $(F,\phi)$ be an $n$-triangle in $T$ and $\alpha\colon\Lambda_{m+1}\to\Lambda_{n+1}$ be a morphism of parasimplices. Then $\alpha^*(F,\phi):=(F\circ\alpha_*,\phi\circ\alpha_*)$ is an $m$-triangle in $T$. In particular, there is an $n$-triangle $\mathsf{s}_1^*(F,\phi):=\mathsf{t}^*(F,\phi)$.
\item Let $(F,\phi)$ be an $n$-triangle in $T$. Then $\mathsf{s}_2^*(F,\phi):=(F\circ\mathsf{s}_2,-\phi\circ\mathsf{s}_2)$ is an $n$-triangle in $T$.
\item Let $\D$ be a stable derivator and $X\in\D_{n,2}(\bbone)\subset\D(\underline{\Lambda}_{n+1,1})$. By \autoref{cor:dnk-sigma} we obtain an isomorphism $\mathsf{s}_2^*(X)\cong\Sigma(X)$. We obtain an $n$-triangle $(F_XX,\phi_X)$ in $\D(\bbone)$, where $F_X=\mathsf{dia}_{\underline{\Lambda}_{n+1,1}}(X)$ is the underlying diagram of $X$ and $\phi_X$ is the composition
\[
F_X\circ\mathsf{s}_2\toiso\mathsf{dia}_{\underline{\Lambda}_{n+1,1}}(\mathsf{s}_2^*(X))\toiso\mathsf{dia}_{\underline{\Lambda}_{n+1,1}}(\Sigma(X))\cong\Sigma\circ F_X.
\]
\end{enumerate}
\end{egs}

\begin{defn}
Let $T$ be an additive category and $\Sigma\colon T\toiso T$ be an automorphism, $n\geq 1$ and $(F,\phi)$ an $n$-triangle in $T$. The \textbf{base} of $(F,\phi)$ is the functor $b(F,\phi):=F\circ sl_{n,2}\colon[n]\cong Sl_{n,2}\to T$.
\end{defn}

\begin{defn}\label{defn:str-tria}
Let $T$ be an additive category and $\Sigma\colon T\toiso T$ be an automorphism. A \textbf{strong triangulation} on $(A,\Sigma)$ consists classes $T_n$ of $n$-triangles in $T$, which are closed under isomorphisms of $n$-triangles, for $n\geq 2$ such that
\begin{itemize}
\item[(Ex)] 
\begin{enumerate}
\item every functor $[n]\to T$ is the base of an $n$-triangle in $T_n$,
\item for every $(F,\phi)\in T_n$ and every morphism $\alpha\colon[m+1]\to[n+1]$ in $\Delta$, the $m$-triangle $\mathsf{i}(\alpha)^*(F,\phi)$ is in $T_m$,
\end{enumerate}
\item[(wF)] for $n$-triangles $(F_1,\phi_1),(F_2,\phi_2)\in T_n$ and a natural transformation between the bases $\psi\colon b(F_1,\phi_1)\to b(F_2,\phi_2)$ there is a morphism $\tilde{\psi}\colon(F_1,\phi_1)\to(F_2,\phi_2)$ such that $b(\tilde{\psi})=\psi$,
\item[(Rot)] for $(F,\phi)\in T_n$ also the triangles $\mathsf{s}_1^*(F,\phi)$ and $\mathsf{s}_2^*(F,\phi)$ are in $T_n$.
\end{itemize}
\end{defn}

\begin{thm}\label{thm:str-tria}
Let $\D$ be a strong, stable derivator, then the suspension functor $\Sigma\colon\D(\bbone)\toiso\D(\bbone)$ and for $n\geq 2$ the classes $T_n$ consisting of those triangles, which are isomorphic to $(F_X,\phi_X)$ for some $X\in\D_{n,2}$ define a strong triangulation on $\D(\bbone)$.
\end{thm}

\begin{proof}
This is \cite[Thm.~13.6]{gst:Dynkin-A}.
\end{proof}

\begin{rmk}\label{rmk:str-tria}
We have ordered the axiom of a strong triangulation in a slightly different way then in \cite{gst:Dynkin-A}. The reason for this is the observation that the axioms fall into a systematic pattern, which can be described as follows.
\begin{enumerate}
\item The axiom (Ex) is an existence axiom. Consider the the case $n=1$. The datum of a 1-triangle is is exactly the datum of a triangle, i.e. a sequence
\begin{equation}\label{eq:triangle}
\cdots\to x\to y\to z\to \Sigma x\to\Sigma y\to\cdots,
\end{equation}
and the first part of the axiom (Ex) ensures the existence of 1-triangles in $T_1$ extending arbitrary morphisms $x\to y$ in $T$. In this case we call the object $z$ in \eqref{eq:triangle} a cone of $x\to y$. Similarly, in the case $n=2$, the first part of the axiom (Ex) yields for a diagram
\[
x_1\to x_2\to x_3
\]
in $T$ the existence of a 2-triangle in $T_2$ of the form
\begin{equation}\label{eq:2-triangle}
\xymatrixrowsep{0.6cm}
\xymatrixcolsep{0.6cm}
\xymatrix{
\ddots\ar[d] &\ddots\ar[d] &\ddots\ar[d] &\ddots\ar[d]\\
0\ar[r]&x_1\ar[r]\ar[d]&x_2\ar[r]\ar[d]&x_3\ar[r]\ar[d]&0\ar[d]\\
&0\ar[r]&c_{12}\ar[r]\ar[d]&c_{13}\ar[r]\ar[d]&\Sigma x_1\ar[r]\ar[d]&0\ar[d]\\
&&0\ar[r]&c_{23}\ar[r]\ar[d]&\Sigma x_2\ar[r]\ar[d]&\Sigma c_{12}\ar[r]\ar[d]&0\ar[d]\\
&&&0\ar[r]&\Sigma x_3\ar[r]\ar[d]&\Sigma c_{13}\ar[r]\ar[d]&\Sigma c_{23}\ar[r]\ar[d]&0\ar[d]\\
&&&&0\ar[r]&\Sigma^2 x_1\ar[r]\ar[d]&\Sigma^2 x_2\ar[r]\ar[d]&\Sigma^2 x_3\ar[r]\ar[d]&0\ar[d]\\
&&&&&\ddots&\ddots&\ddots&\ddots,
}
\end{equation}
where for $i,j\in\lbrace 1,2,3\rbrace$ the object $c_{ij}$ is cone of $x_i\to x_j$. The second part of the axiom (Ex) ensures that the restrictions along the face morphisms $\mathrm{d}_0,\mathrm{d}_1,\mathrm{d}_2$ and $\mathrm{d}_3$ are in $T_1$. By unraveling the definitions (c.f. \autoref{eg:explicit}) we see that the restrictions along the latter three faces correspond to 1-triangles with bases $x_2\to x_3$, $x_1\to x_3$ and $x_1\to x_2$. The remaining face morphism $\mathrm{d}_0$ however yields a triangle of the form
\[
\cdots\to c_{12}\to c_{13}\to c_{23}\to\Sigma c_{12}\cdots.
\]
Moreover, we can reformulate \eqref{eq:2-triangle} as follows.
\begin{itemize}
\item Discarding all non-injective objects
\item Writing all horizontal and vertical compositions as simplices
\item Passing to $\mathsf{s}_2$-orbits and writing morphisms of the $x\to\Sigma y$ as $x\xrightarrow{+}y$
\end{itemize}
Hence we obtain a diagram
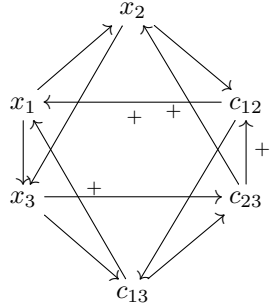
\[
\xymatrix{
&x_2\ar[dr]\ar[ddl]\\
x_1\ar[ur]\ar[d]&&c_{12}\ar[ll]^+\ar[ddl]\\
x_3\ar[rr]\ar[dr]&&c_{23}\ar[u]_+\ar[uul]^+\\
&c_{13}\ar[uul]_+\ar[ur],
}
\]
where the upper front and back and lower left and right triangles are the restrictions to the four faces above, and the remaining triangles commute. In particular, we conclude that the standard octahedral axiom for triangulated categories follows from (Ex) for 2-triangles. On the other hand we see that 2-triangles can be regarded as octahedral diagrams. The effect of the axiom (Ex) for $n\geq 3$ can be described in an analogous way.
\begin{itemize}
\item Every $n$-simplex $X$ in $T$ extends to some $n$-triangle $(F,\phi)$,
\item the $(n-1)$-triangle $\mathrm{d}_0^*(F,\phi)$ relates the cones of all edges of $bX$,
\item $(F,\phi)$ can be rewritten as a pasting of $n$-cells of the form of an $n$-simplex or an $(n-1)$-triangle, and there are $n+2$ of each of these types of cells.
\end{itemize}
Hence we can think of the axiom (Ex) for $n$ large as a generalized octahedral axiom (although the diagrams mentioned in the last point above are not of the shape of an orthoplex in general, as we will see below). For example the diagram in the case $n=3$ for a 3-simplex $x_1\to x_2\to x_3\to x_4$ is described in \autoref{fig:octa}
\begin{figure}
\begin{displaymath}
\begin{tikzpicture}[scale=0.75]
\coordinate[label=right:$x_1$] (A) at (0:8cm);
\coordinate[label=90:$x_4$] (B) at (72:8cm);
\coordinate[label=left:$c_{34}$] (C) at (144:8cm);
\coordinate[label=left:$c_{23}$] (D) at (216:8cm);
\coordinate[label=270:$c_{12}.$] (E) at (288:8cm);
\coordinate (F) at (36:1.5cm);
\coordinate (G) at (108:1.5cm);
\coordinate (H) at (180:1.5cm);
\coordinate (I) at (252:1.5cm);
\coordinate (J) at (324:1.5cm);
\coordinate[label=right:$c_{14}$] (K) at (36:5cm);
\coordinate[label=90:$x_3$] (L) at (108:5cm);
\coordinate[label=left:$c_{24}$] (M) at (180:5cm);
\coordinate[label=270:$c_{13}$] (N) at (252:5cm);
\coordinate[label=right:$x_2$] (O) at (324:5cm);
\draw (A) -- (B) node[midway,above,sloped]{$\leftarrow$};
\draw (B) -- (C) node[midway,above,sloped]{$\leftarrow$};
\draw (C) -- (D) node[midway,below,sloped]{$\rightarrow$};
\draw (D) -- (E) node[midway,below,sloped]{$\rightarrow$};
\draw (E) -- (A) node[midway,below,sloped]{$\rightarrow$};
\draw (F) -- (H) node[midway,above,sloped]{$\leftarrow$};
\draw (H) -- (J) node[midway,below,sloped]{$\rightarrow$};
\draw (J) -- (G) node[midway,above,sloped]{$\leftarrow$};
\draw (G) -- (I) node[midway,below,sloped]{$\rightarrow$};
\draw (I) -- (F) node[midway,below,sloped]{$\rightarrow$};
\draw (A) -- (J) node[midway,above,sloped]{$\leftarrow$};
\draw (B) -- (F) node[midway,above,sloped]{$\leftarrow$};
\draw (C) -- (G) node[midway,below,sloped]{$\rightarrow$};
\draw (D) -- (H) node[midway,below,sloped]{$\rightarrow$};
\draw (E) -- (I) node[midway,above,sloped]{$\leftarrow$};
\draw (F) -- (A) node[midway,below,sloped]{$\rightarrow$};
\draw (G) -- (B) node[midway,below,sloped]{$\rightarrow$};
\draw (H) -- (C) node[midway,above,sloped]{$\leftarrow$};
\draw (I) -- (D) node[midway,above,sloped]{$\leftarrow$};
\draw (J) -- (E) node[midway,below,sloped]{$\rightarrow$};
\draw (A) -- (G) node[midway,above,sloped]{$\leftarrow$};
\draw (B) -- (H) node[midway,above,sloped]{$\leftarrow$};
\draw (C) -- (I) node[midway,below,sloped]{$\rightarrow$};
\draw (D) -- (J) node[midway,below,sloped]{$\rightarrow$};
\draw (E) -- (F) node[midway,above,sloped]{$\leftarrow$};
\draw (I) -- (A) node[midway,below,sloped]{$\rightarrow$};
\draw (J) -- (B) node[midway,below,sloped]{$\rightarrow$};
\draw (F) -- (C) node[midway,above,sloped]{$\leftarrow$};
\draw (G) -- (D) node[midway,above,sloped]{$\leftarrow$};
\draw (H) -- (E) node[midway,below,sloped]{$\rightarrow$};
\draw[dashed] (F) -- (K);
\draw[dashed] (G) -- (L);
\draw[dashed] (H) -- (M);
\draw[dashed] (I) -- (N);
\draw[dashed] (J) -- (O);
\fill (A) circle (2pt);
\fill (B) circle (2pt);
\fill (C) circle (2pt);
\fill (D) circle (2pt);
\fill (E) circle (2pt);
\fill (F) circle (2pt);
\fill (G) circle (2pt);
\fill (H) circle (2pt);
\fill (I) circle (2pt);
\fill (J) circle (2pt);

\end{tikzpicture}
\end{displaymath}
\caption{The octahedral diagram for a 3-simplex.}
\label{fig:octa}
\end{figure}
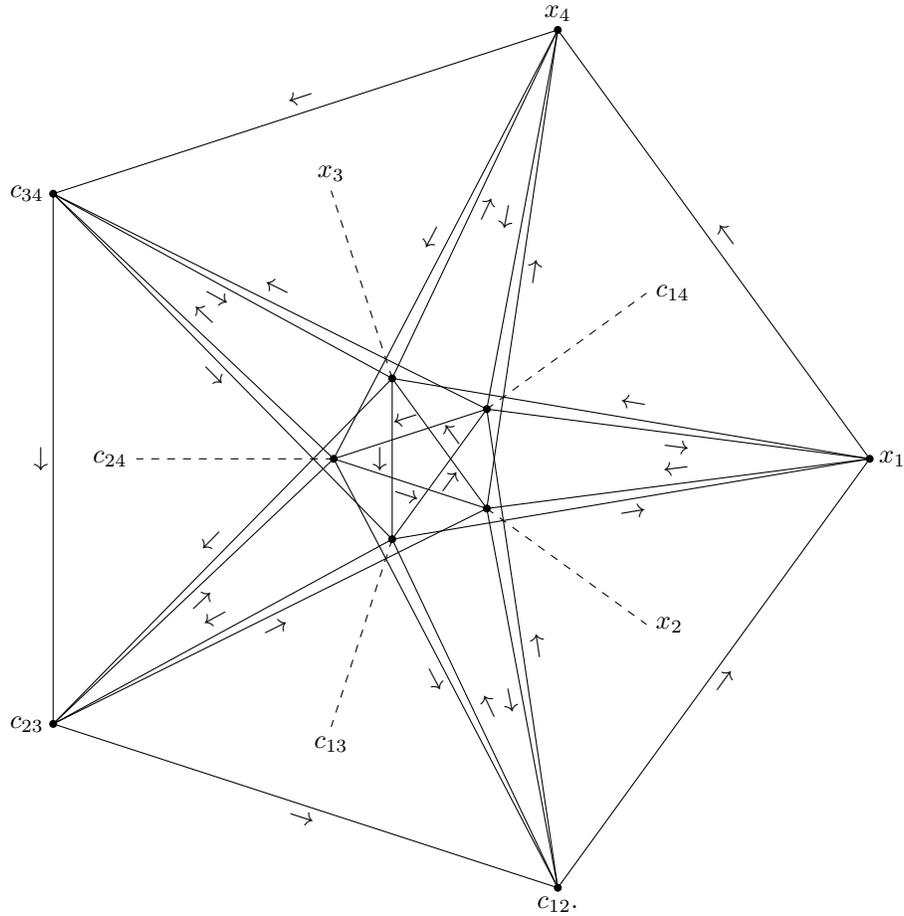
\item The axiom (wF) is the weak functoriality axiom. In the case $n=1$ this axiom states that for two 1-triangles $X_1$ and $X_2$ with cones $c_1$ and $c_2$, respectively, and a morphism $f\colon b(X_1)\to b(X_2)$, there is some morphism $f_c\colon c_1\to c_2$ such that $f$ and $f_c$ assemble into a morphism of 1-triangles. This generalizes in the expected way to $n\geq 2$. In this case we obtain two $n$-triangles $X_1$ and $X_2$ and a morphism $f\colon b(X_1)\to b(X_2)$ some morphism $f_c\colon \mathrm{d}_0^*(X_1)\to \mathrm{d}_0^*(X_2)$ of $(n-1)$-triangles such that $f$ and $f_c$ assemble into a morphism of $n$-triangles.
\item The axiom (Rot) is the rotation axiom. More precisely, for an $n$-triangle $X$ in $T_n$ the axiom (Rot) implies that
\[
\mathsf{s}_3^*(X):=(\mathsf{s}_1^{-1})^*(\mathsf{s}_2^*(X))
\]
is also in $T_n$. In the special cases, we have considered before, $\mathsf{s}_3^*(X)$ can be described explicitly (up to sign).
\begin{itemize}
\item For $n=1$, it corresponds to the shift to the right in \eqref{eq:triangle}.
\item For $n=2$, it corresponds to the rotation around the axis through $x_2$ and $c_{13}$ in the direction of the arrows orthogonal to this axis in \eqref{eq:2-triangle}
\item For $n=3$, it corresponds to the rotation around the center in the direction of the outer arrows in \autoref{fig:octa}.
\end{itemize}
In the general the rotation satisfies $(\mathsf{s}_3^*)^{n+2}\cong(\mathsf{s}_2^*)^n$.
\end{enumerate}
\end{rmk}

As a summary, we see that a strong triangulation can be regarded as an axiomatization of the calculus of cofiber sequences. We now turn to the question whether the calculus of $k$-cofiber sequences for $k\geq 2$ can be axiomatized in a similar way. Of course, a satisfying axiomatization should allow a sufficient amount of useful examples. Hence we ask more precisely, whether \autoref{thm:str-tria} can be generalized to $k$-cofiber sequences, at least for a sufficiently large class of stable derivators. For this we note that, in the proof of \autoref{thm:str-tria} verification of the axiom (Ex) relies on the essential surjectivity of the underlying diagram functor $\mathsf{dia}_{[n]}$, and similarly the axiom (wF) relies on the fullness of $\mathsf{dia}_{[n]}$. These two properties follow from the strongness assumption (since the categories $[n]$ are free). However, this assumption does not exclude many interesting derivators, since derivators associated to model categories are known to be strong (\cite[Prop.~2.15]{cisinski:derivable}, \cite[Thm.~9.8.5]{RB-cofibration}.
We consider the following generalization of \autoref{defn:triangle}.

\begin{defn}
Let $T$ be an additive category and $\Sigma\colon T\toiso T$ be an automorphism and $n,k\geq 1$. An $(n,k)$-\textbf{triangle} $(F,\phi)$ in $T$ consists of
\begin{enumerate}
\item a functor $F\colon\underline{\Lambda}_{n+k,k}\to T$ such that $F(f)\cong 0$ for $f\in\underline{\Lambda}_{n+k,k}$ non-injective.
\item a natural isomorphism $\phi\colon F\circ\mathsf{s}_2\toiso\Sigma^k\circ F$.
\end{enumerate}
A \textbf{morphism} of $(n,k)$-triangles $\psi\colon(F_1,\phi_1)\to(F_2,\phi_2)$ is a natural transformation $\psi\colon F_1\to F_2$ such that $(\Sigma^k\circ\psi)\circ\phi_1=\phi_2\circ(\psi\circ\mathsf{s}_2)$.
\end{defn}

An $(n,1)$-triangle is exactly an $n$-triangle. Similarly to the case $k=1$, we can use the structure of the derivators $\D_{n,k}$ to produce examples of $(n,k)$-triangles in $\D(\bbone)$.

\begin{eg}\label{eg:nk-triangle}
Let $n,k\geq 1$, $\D$ a stable derivator and $X\in\D_{n,k+1}(\bbone)\subset\D(\underline{\Lambda}_{n+k,k})$. By \autoref{cor:dnk-sigma} we obtain an isomorphism $\mathsf{s}_2^*(X)\cong\Sigma^k(X)$. We obtain an $(n,k)$-triangle $(F_X,\phi_X)$ in $\D(\bbone)$, where $F_X=\mathsf{dia}_{\underline{\Lambda}_{n+k,k}}(X)$ is the underlying diagram of $X$ and $\phi_X$ is the composition
\[
F_X\circ\mathsf{s}_2\toiso\mathsf{dia}_{\underline{\Lambda}_{n+k,k}}(\mathsf{s}_2^*(X))\toiso\mathsf{dia}_{\underline{\Lambda}_{n+k,k}}(\Sigma(X))\cong\Sigma\circ F_X.
\]
\end{eg}

The categories $[n]^k$ and $Sl_{n,k+1}$ are unfortunately not free for $k\geq2$ and $[n]\geq 1$. Moreover, we refer to \cite[Ex.~3.17]{bg:cubical} for an example showing that the underlying diagram functor $\mathsf{dia}_{\cube{2}}$ is not full in the case of the derivator of a field (even if we restrict to $\D_{1,3}$). Hence it is not reasonable to ask for the analogue of the weak functorliality axiom for $k\geq 2$. On the other hand we can show that the underlying diagram functor $\mathsf{dia}_{Sl_{n,k}}\colon sl\D_{n,k}\to\D(\bbone)^{Sl_{n,k}}_0$ is essentially surjective for a strong derivator $\D$ and $k=3$. Here $\D(\bbone)^{Sl_{n,k}}_0\subset\D(\bbone)^{Sl_{n,k}}$ denotes the full subcategory spanned by those functors $F\colon Sl_{n,k}\to\D(\bbone)$, such that $F(f)=0$ whenever $f$ is non-injective.

\begin{lem}\label{lem:full}
Let $F\colon C\to D$ be a full and essentially surjective functor between categories, and $A$ a finite, free category. Then
\[
F^A\colon C^A\to D^A
\]
is essentially surjective.
\end{lem}

\begin{proof}
Consider an object $G\colon A\to D$ in $D^A$. Since $F$ is essentially surjective, for all $a\in A$ there is an object $c_a\in C$ and an isomorphism $\phi_a\colon F(c_a)\to G(a)$. By assumption the is a finite set $B$ of morphisms in $A$ such that $A$ is generated freely by $B$. Since $F$ is full, for all $b\colon a\to a'$ in $B$ there is a morphism $c_b\in C$ such that $F(c_b)$ is the composition
\[
F(c_a)\xrightarrow{\phi_a}G(a)\xrightarrow{G(b)}G(a')\xrightarrow{\phi_{a'}^{-1}}F(c_{a'}).
\]
Since the category $A$ is free, the assignment $a\mapsto c_a, b\mapsto c_b$ defines a unique functor $\tilde{G}\colon A\to C$ and by construction, the collection of isomorphisms $\lbrace \phi_a\vert a\in A\rbrace$ defines a natural isomorphism $F^A(\tilde{G})\toiso G$.
\end{proof}

\begin{cor}\label{cor:ess-sur}
Let $\D$ be a strong derivator and $A,B$ finite, free categories. Then the underlying diagram functor
\[
\mathsf{dia}_{A\times B}\colon\D(A\times B)\to\D(\bbone)^{A\times B}
\]
is essentially surjective.
\end{cor}

\begin{proof}
The underlying diagram functor $\mathsf{dia}_{A\times B}$ is isomorphic to the composite
\[
\D(A\times B)=\D^A(B)\xrightarrow{\mathsf{dia}_B}\D^A(\bbone)^B=\D(A)^B\xrightarrow{\mathsf{dia}_A^B}(\D(\bbone)^A)^B\cong\D(\bbone)^{A\times B}.
\]
Since $\D$ is strong and $A,B$ are finte, free categories the functors $\mathsf{dia}_A$ and $\mathsf{dia}_B$ are full and essentially surjective. By \autoref{lem:full} the functor $\mathsf{dia}_A^B$ is essentially surjective. Hence $\mathsf{dia}_{A\times B}$ is as a composition of essentially surjective functors itself essentially surjective.
\end{proof}

\begin{cor}\label{cor:ess-sur-II}
Let $\D$ be a strong derivator. Then the retriction of the underlying diagram functor
\[
\mathsf{dia}_{Sl_{n,k}}\colon sl\D_{n,3}\to\D(\bbone)^{Sl_{n,3}}_0
\]
is essentially surjective.
\end{cor}

\begin{proof}
Let $X\in\D(\bbone)^{Sl_{n,3}}_0$. By the universal property of the zero-object, we can extend $X$ to an object $Y\in\D(\bbone)^{Sl_{n,3}^{\square}}$ with $Y\vert_{Sl_{n,3}}=X$ and $Y\vert_{Sl_{n,3}^{\square}\setminus Sl_{n,3}}=0$. Since $Sl_{n,3}^{\square}\cong[n]\times[n]$ there is by \autoref{cor:ess-sur} an object $Z\in sl\hat{\D}_{n,3}$ with $\mathsf{dia}_{Sl_{n,3}^{\square}}(Z)=Y$. Since underlying diagram functors are compatible with inverse images we obtain
\[
\mathsf{dia}_{Sl_{n,3}}((\mathsf{j}\vert_{Sl_{n,3}})^*(Z))\cong X.
\]
\end{proof}

\begin{rmk}
Using \autoref{cor:ess-sur-II} it should be straight forward to verify the analogues of the axioms (Ex) and (Rot) for the classes of $(n,2)$-triangles in the underlying category of a strong, stable derivator $\D$ defined by \autoref{eg:nk-triangle}. However, it is clear that the resulting formalism will be less useful that a strong triangulation due to the failure of the axiom (wF) in important examples.
\end{rmk}

In the following we point out that the situation becomes even worse if we pass to $(n,k)$-triangles for $k\geq 3$. More precisely, we show, that the underlying diagram functor $\mathsf{dia}_{Sl_{n,k}}\colon sl\D_{n,k}\to\D(\bbone)^{Sl_{n,k}}_0$ will in general not be essentially surjective for $k\geq 4$. The obstructions for this arise from non-trivial 3-fold Toda brackets. We recall the following alternative definition of 4-fold Toda brackets from \cite{CF-Toda}.

\begin{con}\label{con:alt-Toda}
Let $\mathcal{T}$ be a triangulated category.
\begin{enumerate}
\item Let $x_3\xrightarrow{u_3}x_2\xrightarrow{u_2}x_1\xrightarrow{u_1}x_0$ be a 3-simplex in $\mathcal{T}$. Then the 3-fold Toda bracket of the above sequence is the collection of composites $\beta\circ\Sigma\alpha\colon\Sigma x_3\to x_0$, where $\alpha$ and $\beta$ are maps making the following diagram, where the middle row is a distinguished triangle, commutative
\begin{equation}\label{eq:3-Toda}
\xymatrix{
x_3\ar[d]_{\alpha}\ar[r]^{u_3}&x_2\ar[d]_{\id}\\
Fu_2\ar[r]&x_2\ar[r]^{u_2}&x_1\ar[r]\ar[d]_{\id}&Cu_2\ar[d]_{\beta}\\
&&x_1\ar[r]^{u_1}&x_0.
}
\end{equation}
We observe that the 3-fold Toda bracket of $x_3\xrightarrow{u_3}x_2\xrightarrow{u_2}x_1\xrightarrow{u_1}x_0$ is non-empty iff $u_2\circ u_3=0$ and $u_1\circ u_2=0$.
\item Let $x_4\xrightarrow{u_4}x_3\xrightarrow{u_3}x_2\xrightarrow{u_2}x_1\xrightarrow{u_1}x_0$ be a 4-simplex in $\mathcal{T}$. Then the 4-fold Toda bracket of the above sequence is the union of the 3-fold Toda bracket associated to the sequences $\Sigma x_4\xrightarrow{\Sigma\alpha}Cu_3\xrightarrow{\beta}x_1\xrightarrow{u_1}x_0$ for all choices of morphisms $\alpha$ and $\beta$ making \eqref{eq:3-Toda} associated to $x_4\xrightarrow{u_4}x_3\xrightarrow{u_3}x_2\xrightarrow{u_2}x_1$ commutative. In particular, it follows from (i), that the 4-fold Toda bracket of $x_4\xrightarrow{u_4}x_3\xrightarrow{u_3}x_2\xrightarrow{u_2}x_1\xrightarrow{u_1}x_0$ is empty if 0 is not contained in the 3-fold Toda bracket of $x_4\xrightarrow{u_4}x_3\xrightarrow{u_3}x_2\xrightarrow{u_2}x_1$.
\end{enumerate}
\end{con}

\begin{rmk}
We refer to \cite[Ex.~5.3,Ex.~5.6,Ex.~5.7]{CF-Toda} for a detailed argument for the equivalence of \autoref{con:alt-Toda} and \autoref{defn:Toda}. Moreover, \autoref{con:alt-Toda} is generalized in \emph{loc.~cit.} to a definition for arbitrary higher Toda brackets and a general equivalence result to \autoref{defn:Toda} is established.
\end{rmk}

\begin{eg}\label{eg:not-es}
Let $\D$ be a stable derivator and $x_4\xrightarrow{u_4}x_3\xrightarrow{u_3}x_2\xrightarrow{u_2}x_1$ a sequence such that $u_3\circ u_4=0$, $u_2\circ u_3=0$ and such that 0 is not contained in its 3-fold Toda bracket. From the assumptions on the compositions and the universal property of the zero-object, we deduce the existence of a diagram $X\in\D(\bbone)^{\cube{3}}$ of the form
\[
\xymatrix{
&0\ar[rr]\ar[dd]&&0\ar[dd]\\
x_4\ar[dd]\ar[rr]^(.65){u_4}\ar[ur]&&x_3\ar[dd]^(.35){u_3}\ar[ur]\\
&0\ar[rr]&&x_1\\
0\ar[ur]\ar[rr]&&x_2.\ar[ur]_{u_2}
}
\]
We assume the existence of an object $Y\in\D(\cube{3})=\D^{\cube{3}}(\bbone)$ with $\mathsf{dia}_{\cube{3}}(Y)=X$. The structure of $X$ immediately implies that $Y\in\D^{\cube{3}}_{\rightarrow}$. Hence, $(\id_{\cube{2}}\times\mathrm{d}_2)_*(Y)\in\D^{\mathsf{T}_4}(\bbone)$ and by \autoref{thm:toda} the underlying diagram of $\mathsf{Toda}_4((\id_{\cube{2}}\times\mathrm{d}_2)_*(Y))$ lies in the 4-fold Toda bracket of $x_4\xrightarrow{u_4}x_3\xrightarrow{u_3}x_2\xrightarrow{u_2}x_1\xrightarrow{0}0$. On the other hand the assumption on $x_4\xrightarrow{u_4}x_3\xrightarrow{u_3}x_2\xrightarrow{u_2}x_1$ ensures by \autoref{con:alt-Toda} that the 4-fold Toda bracket above is empty, which leads to a contradiction.

A more explicit example for a sequence of morphisms $x_4\xrightarrow{u_4}x_3\xrightarrow{u_3}x_2\xrightarrow{u_2}x_1$ satisfying the assumption above is given by the sequence $S^2\xrightarrow{\eta}S^1\xrightarrow{2}S^1\xrightarrow{\eta}S^0$ in the underlying category of the derivator $\mathscr{H}$ associated to the homotopy theory of spectra. In fact, the 3-fold Toda bracket of this sequence consists exactly of the elements $-2\nu,2\nu\colon S^3\to S^0$ \cite[V.(5.4)]{toda:brackets}.
\end{eg}

It is immediate from \autoref{eg:not-es}, that the classes of $(n,k)$-triangles defined by \autoref{eg:nk-triangle} for a stable derivator $\D$ and $k\geq 3$ will not satisfy the analogs of the existence axiom for important examples of stable derivators.

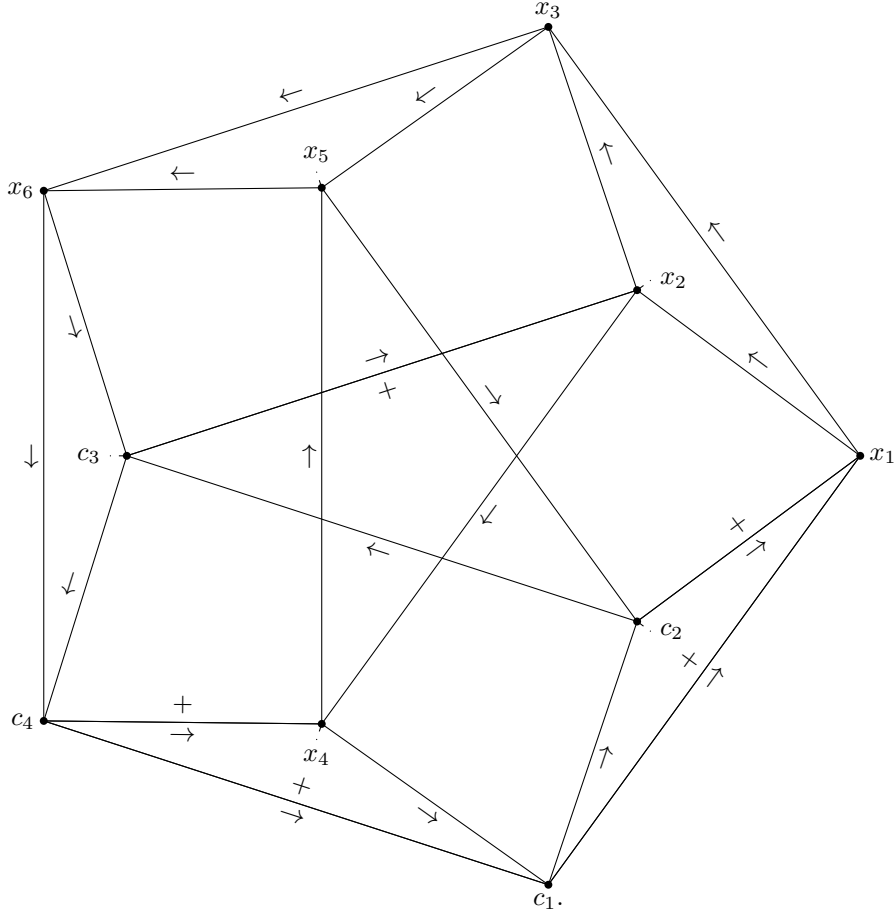
\begin{figure}
\begin{displaymath}
\begin{tikzpicture}[scale=0.75]
\coordinate[label=right:$x_1$] (A) at (0:8cm);
\coordinate[label=90:$x_3$] (B) at (72:8cm);
\coordinate[label=left:$x_6$] (C) at (144:8cm);
\coordinate[label=left:$c_4$] (D) at (216:8cm);
\coordinate[label=270:$c_1.$] (E) at (288:8cm);
\coordinate (F) at (36:5cm);
\coordinate (G) at (108:5cm);
\coordinate (H) at (180:5cm);
\coordinate (I) at (252:5cm);
\coordinate (J) at (324:5cm);
\coordinate[label=right:$x_2$] (K) at (36:5.3cm);
\coordinate[label=90:$x_5$] (L) at (108:5.3cm);
\coordinate[label=left:$c_3$] (M) at (180:5.3cm);
\coordinate[label=270:$x_4$] (N) at (252:5.3cm);
\coordinate[label=right:$c_2$] (O) at (324:5.3cm);
\draw (A) -- (B) node[midway,above,sloped]{$\leftarrow$};
\draw (B) -- (C) node[midway,above,sloped]{$\leftarrow$};
\draw (C) -- (D) node[midway,below,sloped]{$\rightarrow$};
\draw (D) -- (E) node[midway,below,sloped]{$\rightarrow$};
\draw (D) -- (E) node[midway,above,sloped]{$+$};
\draw (E) -- (A) node[midway,below,sloped]{$\rightarrow$};
\draw (E) -- (A) node[midway,above,sloped]{$+$};
\draw (F) -- (H) node[midway,above,sloped]{$\rightarrow$};
\draw (F) -- (H) node[midway,below,sloped]{$+$};
\draw (H) -- (J) node[midway,below,sloped]{$\leftarrow$};
\draw (J) -- (G) node[midway,above,sloped]{$\rightarrow$};
\draw (G) -- (I) node[midway,below,sloped]{$\leftarrow$};
\draw (I) -- (F) node[midway,below,sloped]{$\leftarrow$};
\draw (A) -- (J) node[midway,below,sloped]{$\rightarrow$};
\draw (A) -- (J) node[midway,above,sloped]{$+$};
\draw (B) -- (F) node[midway,above,sloped]{$\leftarrow$};
\draw (C) -- (G) node[midway,above,sloped]{$\leftarrow$};
\draw (D) -- (H) node[midway,above,sloped]{$\leftarrow$};
\draw (E) -- (I) node[midway,below,sloped]{$\rightarrow$};
\draw (F) -- (A) node[midway,above,sloped]{$\leftarrow$};
\draw (G) -- (B) node[midway,above,sloped]{$\leftarrow$};
\draw (H) -- (C) node[midway,below,sloped]{$\rightarrow$};
\draw (I) -- (D) node[midway,below,sloped]{$\rightarrow$};
\draw (I) -- (D) node[midway,above,sloped]{$+$};
\draw (J) -- (E) node[midway,below,sloped]{$\rightarrow$};
\draw[dashed] (F) -- (K);
\draw[dashed] (G) -- (L);
\draw[dashed] (H) -- (M);
\draw[dashed] (I) -- (N);
\draw[dashed] (J) -- (O);
\fill (A) circle (2pt);
\fill (B) circle (2pt);
\fill (C) circle (2pt);
\fill (D) circle (2pt);
\fill (E) circle (2pt);
\fill (F) circle (2pt);
\fill (G) circle (2pt);
\fill (H) circle (2pt);
\fill (I) circle (2pt);
\fill (J) circle (2pt);

\end{tikzpicture}
\end{displaymath}
\caption{The octahedral diagram for a 2-cofiber sequences. Here a morphism $x\xrightarrow{+}y$ denotes a morphism $x\to\Sigma^2y$. Moreover, we note that the rotation operation defined by $\mathsf{s}_3^*$ corresponds to the rotation around the center in the direction of the outer arrows.}
\label{fig:octaII}
\end{figure}

However, if we decide to work with coherent diagrams instead of just diagrams in some homotopy category, the results of \S\ref{sec:dnk} and \S\ref{sec:vertical} can be regarded as analogues (which are now consequences of having stable derivator rather than being axioms) of the axioms of a strong triangulation for \textbf{coherent} $(n,k)$-triangles in a stable derivator $\D$ (i.e. objects in $\D_{n,k+1}$).
\begin{itemize}
\item[(Ex')] \autoref{thm:slices}, \autoref{cor:dnk-sigma} (first part) and \autoref{cor:higher-Sdot-I} (second part),
\item[(wF')] \autoref{thm:slices},
\item[(Rot')] \autoref{cor:dnk-sym} and \autoref{cor:fracCY}.
\end{itemize}
We outline the analogy to \autoref{defn:str-tria} in the case of $(2,2)$-triangles. Consider an object $X\in sl\D_{2,3}$ with underlying diagram
\[
\xymatrix{
x_1\ar[r]\ar[d] & x_2 \ar[r]\ar[d] & x_3 \ar[d]\\
0 \ar[r] & x_4 \ar[r]\ar[d] & x_5 \ar[d]\\
& 0\ar[r] & x_6.
}
\]
Using \autoref{thm:slices} we obtain an object $\tilde{X}\in\D_{2,3}$. Moreover, by applying \autoref{eg:explicit} to the face morphisms $\mathrm{d}_1$, $\mathrm{d}_2$, $\mathrm{d}_3$ and $\mathrm{d}_4\colon\Lambda_3\to\Lambda_4$ we obtain (1,2)-triangles (i.e. 2-cofiber sequence) with bases
\[
\xymatrix{
x_4\ar[r]\ar[d]&x_5\ar[d]&x_2\ar[r]\ar[d]&x_3\ar[d]&x_1\ar[r]\ar[d]&x_3\ar[d]&x_1\ar[r]\ar[d]&x_2\ar[d]\\
0\ar[r]&x_6,&0\ar[r]&x_6,&0\ar[r]&x_5,&0\ar[r]&x_4,
}
\]
respectively. Let $c_4, c_3, c_2$ and $c_1$ denote the respective 2-cones, then \autoref{eg:explicit} applied to $\mathrm{d}_0\colon\Lambda_3\to\Lambda_4$ leads to a 2-cofiber sequence with base 
\[
\xymatrix{
c_1\ar[r]\ar[d]&c_2\ar[d]\\
0\ar[r]&c_3
}
\]
and 2-cone $c_4$. Hence, a (2,2)-triangle can be regarded as the analogue of an octahedral diagram for 2-cofiber sequences. Moreover, by applying the analogue of the procedure used in \autoref{rmk:str-tria} we can rewrite $\tilde{X}$ as displayed in \autoref{fig:octaII}.

Furthermore, \autoref{thm:main} implies that the datum of an $(n,k)$-triangle is equivalent to the datum of $(k,n)$-triangle.

In the following we indicate some relations between the calculus of $(n,k)$-triangles for some $k\geq 2$ fixed and $k+2$-angulated categories introduced by Geiss--Keller--Oppermann \cite{GKO-angulated}. First, we recall the basic definitions.

\begin{defn}
Let $n\geq 3$, $T$ be an additive category and $\tilde{\Sigma}\colon T\toiso T$ be an automorphism. A diagram in $T$ of the form
\[
x_1\xrightarrow{u_1} x_2\xrightarrow{u_2}\cdots\xrightarrow{u_{n-1}} x_n\xrightarrow{u_n}\tilde{\Sigma}x_1
\]
is called an $n-\tilde{\Sigma}$-sequence. A morphism of $n-\tilde{\Sigma}$-sequences is a commutative diagram
\[
\xymatrix{
x_1\ar[r]^{u_1}\ar[d]_{\phi_1}&x_2\ar[r]^{u_2}\ar[d]_{\phi_2}&\cdots\ar[r]^{u_{n-1}}&x_n\ar[r]^{u_n}\ar[d]_{\phi_n}&\tilde{\Sigma}x_1\ar[d]_{\tilde{\Sigma}\phi_1}\\
y_1\ar[r]_{v_1}&y_2\ar[r]_{v_2}&\cdots\ar[r]_{v_{n-1}}&y_n\ar[r]_{v_n}&\tilde{\Sigma}y_1,
}
\]
where both rows are $n-\tilde{\Sigma}$-sequences. A morphism of $n-\tilde{\Sigma}$-sequences is called a weak isomorphism if there is $1\leq i\leq n$ such that $\phi_i$ and $\phi_{i+1}$ (with $\phi_{n+1}=\tilde{\Sigma}\phi_1$) are isomorphisms.
\end{defn}

\begin{defn}
Let $n\geq 3$, $T$ be an additive category and $\tilde{\Sigma}\colon T\toiso T$ be an automorphism. An $n$-angulation on $(T,\tilde{\Sigma})$ consists of a class $T_n$ of $n-\tilde{\Sigma}$-sequences, whose elements are called $n$-angles, which is closed under weak isomorphisms, such that
\begin{enumerate}
\item for every object $x\in T$ there is a trivial $n$-angle of the form
\[
x\xrightarrow{\id}x\to 0\to\cdots\to 0\to\tilde{\Sigma}x,
\]
and every morphism $x\to y\in T$ is the first morphism of an $n$-angle,
\item if $x_1\xrightarrow{u_1} x_2\xrightarrow{u_2}\cdots\xrightarrow{u_{n-1}} x_n\xrightarrow{u_n}\tilde{\Sigma}x_1$ is an $n$-angle, then also its rotation
\[
x_2\xrightarrow{u_2} x_3\xrightarrow{u_3}\cdots\xrightarrow{u_{n-1}} x_n\xrightarrow{u_n}\tilde{\Sigma}x_1\xrightarrow{(-1)^n\tilde{\Sigma}u_1}\tilde{\Sigma}x_2
\]
is an $n$-angle,
\item every commutative diagram
\[
\xymatrix{
x_1\ar[r]^{u_1}\ar[d]_{\phi_1}&x_2\ar[r]^{u_2}\ar[d]_{\phi_2}&x_3\ar[r]^{u_3}\ar@{-->}[d]_{\phi_3}&\cdots\ar[r]^{u_{n-1}}&x_n\ar[r]^{u_n}\ar@{-->}[d]_{\phi_n}&\tilde{\Sigma}x_1\ar[d]_{\tilde{\Sigma}\phi_1}\\
y_1\ar[r]_{v_1}&y_2\ar[r]_{v_2}&y_3\ar[r]_{v_3}&\cdots\ar[r]_{v_{n-1}}&y_n\ar[r]_{v_n}&\tilde{\Sigma}y_1,
}
\]
where the rows are $n$-angles can be completed to a morphism of $n$-angles,
\item the morphism of $n$-angles above can be chosen, such that
\[
\resizebox{.9\hsize}{!}{$
x_2\oplus y_1\xrightarrow{\begin{vmatrix}-u_2&0\\\phi_2&v_1\end{vmatrix}}x_3\oplus y_2\xrightarrow{\begin{vmatrix}-u_3&0\\\phi_3&v_2\end{vmatrix}}\cdots\xrightarrow{\begin{vmatrix}-u_n&0\\\phi_n&v_{n-1}\end{vmatrix}}\tilde{\Sigma}x_1\oplus y_n\xrightarrow{\begin{vmatrix}-\tilde{\Sigma}u_1&0\\\tilde{\Sigma}\phi_n&v_n\end{vmatrix}}\tilde{\Sigma}x_2\oplus\tilde{\Sigma}y_1$}
\]
is an $n$-angle.
\end{enumerate}
\end{defn}

We refer to \cite{BT-axioms} for a different but equivalent axiomatization of $n$-angulated categories.

Given a stable derivator $\D$, we consider its underlying category $\D(\bbone)$ together with the automorphism $\tilde{\Sigma}=\Sigma^{n-2}$. Then the obvious candidates for $n$-angles in $\D(\bbone)$ are the underlying diagrams of $(1,n-2)$-triangles, i.e. $(n-1)$-cofiber sequences (c.f. \autoref{rmk:highercof}). In fact, using the methods from \S\ref{sec:Toda}, in particular \autoref{prop:filtered}, it is not hard to show that the resulting $n-\tilde{\Sigma}$-sequences canonically extend to a diagram as described in \cite[Thm.~1]{GKO-angulated} and \cite[Def.~5.15(iii)]{OT-cluster}. Let $T_n$ denote the class of the $n-\tilde{\Sigma}$-sequences defined by the $(n-1)$-cofiber sequences, as described above. We may ask whether $(\D(\bbone),\Sigma^{n-2},T_n)$ defines an $n$-angulation.

\noindent We will show that the answer is in general "no", by comparing the axioms of an $n$-angulation with the properties of coherent $(1,n-2)$-triangles.
\begin{itemize}
\item On the one hand the existence axiom in the context of $n$-angulated categories (axiom (i)) is weaker as the existence for $(1,n-2)$-triangles. For a morphism $x\to y$ in $\D(\bbone)$ we can find $X\in\D([1])\cong sl\D_{1,2}(\bbone)$ with $\mathsf{dia}_{[1]}(X)=(x\to y)$. Then $(\mathsf{d}^h[1])^{n-3}(X)$ and $(\mathsf{d}^h[-1])^{n-3}(X)$ both give rise to $n$-angles extending $x\to y$. The resulting $n$-angles are weakly isomorphic but in general not isomorphic.
\item On the other hand the weak functoriality axiom for $n$-angulated categories (axiom (iii)) does not follow from the corresponding property of $(1,n-2)$-triangles. In general there are too many different extensions of morphisms to $n$-angles. For instance, given a non-invertible morphism $x\to y$ in $\D(\bbone)$, with $X\in\D([1])$ as above. Since $\mathsf{d}^h[1](X)$ and $\mathsf{d}^h[-1](X)$ both define extensions of $(x\to y)$ to 4-angles, the identity on $(x\to y)$ should extend to a morphism of $4$-angles $\mathsf{d}^h[1](X)\to\mathsf{d}^h[-1](X)$. We assume that this morphism can be covered be a morphism in $\D_{1,3}(\bbone)$. By applying \autoref{thm:main} we obtain a morphism $\mathsf{d}^v[1](X)\to\mathsf{d}^v[-1](X)$, which has by \autoref{eg:explicit} an underlying diagram of the form
\[
\xymatrix{
x\ar[r]\ar@{=}[d]&y\ar@{=}[r]\ar[d]&y\ar@{=}[d]\\
x\ar@{=}[r]&x\ar[r]&y.
}
\]
This would imply that $(x\to y)$ is invertible and therefore contradict our assumption.
\end{itemize}
This is of course not surprising. They key examples of $n$-angulated categories arise as $(n-2)$-cluster tilting subcategories of triangulated categories (\cite[Thm.~1]{GKO-angulated}). Roughly speaking these are subcategories of triangulated categories, which are closed under bicartesian $(n-1)$-cubes but not necessarily under bicartesian squares. For a precise definition we refer to \cite[Def.~1.1]{iyama-higherA}. By restricting to an $(n-2)$-cluster tilting subcategory we make the set of possible extensions of a morphism to an $n$-angle smaller in a way such that the weak functoriality for $n$-angles holds.

\noindent Important examples of $(n-2)$-cluster tilting subcategories were constructed by Iyama \cite{iyama-higherA} as subcategories of the derived categories of $(n-2)$-Auslander-algebras associated to Dynkin quivers of type A. The following related questions and open problems seem to be interesting for future research.
\begin{itemize}
\item Are there analogues of higher cluster tilting subcategories for general stable homotopy theories?
\item Are there explicit descriptions of these, at least in case of Dynkin quivers of type A?
\item Do (coherent) $(k,n-2)$-triangles in $n$-cluster tilting subcategories satisfy special properties for $k\geq 2$?
\item Do these properties descent to an axiomatization of higher $n$-angulated categories?
\end{itemize}

\appendix
\section{Conjugation for 2-categories}
\label{sec:appendix}

In this short appendix we describe in which way a family of autoequivalences acts on a specific 2-category of pseudofunctors.
Of course there many sources in the literature (e.g. \cite{benabou:intro}) for an introduction to the theory of bicategories, but to fix the notation we start by recalling the elementary definitions. Since we restrict to the case of pseudofunctors between 2-categories from the beginning, we refer to \cite{renaudin} for a similar but more detailed exposition.

\begin{defn}
Let $\mathscr{C}$ and $\mathscr{C}'$ be 2-categories. A peudofunctor $\mathscr{F}:\mathscr{C}\rightarrow\mathscr{C}'$ consists of:
\begin{enumerate}
\item for $c \in Ob(\mathscr{C})$ an object $\mathscr{F}(c) \in Ob(\mathscr{C}')$,
\item for $c_1,c_2 \in Ob(\mathscr{C})$ a functor 
$$\mathscr{F}_{c_1,c_2}: \mathscr{C}(c_1,c_2) \rightarrow \mathscr{C}'(c_1,c_2).$$
\item for $c_1,c_2,c_3 \in Ob(\mathscr{C})$ a natural isomorphism
\[
\resizebox{.95\hsize}{!}{$
\mathrm{m}^{\mathscr{F}}_{c_1,c_2,c_3}: (-\circ -)\circ (\mathscr{F}_{c_1,c_2}\times\mathscr{F}_{c_2,c_3}) \Rightarrow \mathscr{F}_{c_1,c_3}\circ (-\circ -):\mathscr{C}(c_1,c_2)\times\mathscr{C}(c_2,c_3)\rightarrow\mathscr{C}'(\mathscr{F}(c_1),\mathscr{F}(c_3))
$}
\]
\item for $c \in Ob(\mathscr{C})$ a natural isomorphism
$$\mathrm{u}^{\mathscr{F}}_{c}: id_{\mathscr{F}(c)} \Rightarrow \mathscr{F}_{c,c}\circ id_c: \bbone \rightarrow \mathscr{C}'(\mathscr{F}(c),\mathscr{F}(c))$$
\end{enumerate}
such that the equalities of 2-isomorphisms:
\begin{enumerate}
\item (associativity) for all $(f,g,h) \in \mathscr{C}(c_1,c_2)\times\mathscr{C}(c_2,c_3)\times\mathscr{C}(c_3,c_4)$:
$$\mathrm{m}^{\mathscr{F}}_{g\circ f,h}\circ(\mathscr{F}(h)\mathrm{m}^{\mathscr{F}}_{g,f})=\mathrm{m}^{\mathscr{F}}_{f,h\circ g}\circ(\mathrm{m}^\mathscr{F}_{g,h}\mathscr{F}(f)):\mathscr{F}(h)\circ\mathscr{F}(g)\circ\mathscr{F}(f)\rightarrow\mathscr{F}(h\circ g\circ f)$$
\item (left unitality) for all $f \in \mathscr{C}(c_1,c_2)$:
$$id_{\mathscr{F}(f)}=\mathrm{m}^{\mathscr{F}}_{id_{c_1},f}\circ(\mathrm{u}_{c_1}^{\mathscr{F}}\mathscr{F}(f)):\mathscr{F}(f)\rightarrow\mathscr{F}(f)$$
\item (right unitality) for all $f \in \mathscr{C}(c_1,c_2)$:
$$id_{\mathscr{F}(f)}=\mathrm{m}^{\mathscr{F}}_{f,id_{c_2}}\circ(\mathscr{F}(f)\mathrm{u}_{c_2}^{\mathscr{F}}):\mathscr{F}(f)\rightarrow\mathscr{F}(f)$$
\end{enumerate}
hold true, with the simplified notation for $(f,g)\in \mathscr{C}(c_1,c_2)\times\mathscr{C}(c_2,c_3)$:
$$\mathrm{m}^{\mathscr{F}}_{f,g}:=(\mathrm{m}^{\mathscr{F}}_{c_1,c_2,c_3})_{(f,g)}:\mathscr{F}(g)\circ\mathscr{F}(f)\xrightarrow{\simeq}\mathscr{F}(g\circ f).$$
\end{defn}

\begin{defn}
Let $\mathscr{F},\mathscr{F}':\mathscr{C}\rightarrow\mathscr{C}'$ be pseudofunctors. A pseudonatural transformation $\alpha:\mathscr{F}\rightarrow\mathscr{F}'$ consists of:
\begin{enumerate}
\item for $c \in Ob{\mathscr{C}}$ a 1-morphism $\alpha_c:\mathscr{F}(c)\rightarrow\mathscr{F}'(c)$,
\item for $c_1,c_2 \in Ob(\mathscr{C})$ a natural isomorphism
$$\lambda^{\alpha}_{c_1,c_2}:(\alpha_{c_2})_*\circ\mathscr{F}_{c_1,c_2} \Rightarrow (\alpha_{c_1})^*\circ\mathscr{F}'_{c_1,c_2}:\mathscr{C}(c_1,c_2)\rightarrow\mathscr{C}'(\mathscr{F}(c_1),\mathscr{F}'(c_2))$$
\end{enumerate}
such that the equalities of 2-isomorphisms:
\begin{enumerate}
\item (associativity) for all $(f,g) \in\mathscr{C}(c_1,c_2)\times\mathscr{C}(c_2,c_3)$:
$$\lambda^{\alpha}_{g\circ f}\circ(id_{\alpha_{c_3}}\mathrm{m}^{\mathscr{F}}_{f,g})=(\mathrm{m}^{\mathscr{F}'}_{f,g}id_{\alpha_{c_1}})\circ(id_{\mathscr{F}'(g)}\lambda^{\alpha}_{f})\circ(\lambda^{\alpha}_{g}id_{\mathscr{F}(f)}):$$
$$\alpha_{c_3}\circ\mathscr{F}(g)\circ\mathscr{F}(f)\rightarrow\mathscr{F}'(g\circ f)\circ\alpha_{c_1}$$
\item (unitality) for all $c \in Ob(\mathscr{C})$:
$$\mathrm{u}^{\mathscr{F}'}_cid_{\alpha_c}=\lambda^{\alpha}_{id_c}\circ(id_{\alpha_c}\mathrm{u}^{\mathscr{F}}_c):\alpha_c \rightarrow \mathscr{F}(id_c)\circ\alpha_c$$
\end{enumerate}
hold true, with the simplified notation for $f \in \mathscr{C}(c_1,c_2)$:
$$\lambda^{\alpha}_f:=(\lambda^{\alpha}_{c_1,c_2})_f:a_{c_2}\circ\mathscr{F}(f)\xrightarrow{\cong}\mathscr{F}'(f)\circ\alpha_{c_1}$$
\end{defn}

\begin{defn}
Let $\mathscr{F},\mathscr{F}'\colon\mathscr{C}\to\mathscr{C}'$  be pseudofunctors and $\alpha,\alpha'\colon\mathscr{F}\to\mathscr{F}'$ be pseudonatural transformations. A modification $\theta\colon\alpha\to\alpha'$ consists of 
\begin{enumerate}
\item for $c\in Ob(\mathscr{C})$ a 2-morphism $\theta_c\colon\alpha_c\to\alpha'_c$,
\item such that for all $f\in\mathscr{C}(c_1,c_2)$ there is an equality of 2-morphisms
\[
(\id_{\mathscr{F}'(f)}\theta_{c_1})\circ\lambda^{\alpha}_f=\lambda^{\alpha'}_f\circ(\theta_{c_2}\id_{\mathscr{F}(f)}\colon\alpha_{c_2}\circ\mathscr{F}(f)\to\mathscr{F}'(f))\circ\alpha'_{c_1}.
\]
\end{enumerate}
\end{defn}

\begin{prop}\label{prop:conjugation}
Let $\mathscr{F}\colon\mathscr{C}\rightarrow\mathscr{C}'$ be a pseudofunctor between 2-categories.
Given for all $c\in Ob(\mathscr{C})$ an object $c_{\mathsf{S}}\in\mathscr{C}'$ and an equivalence
\[
(\mathsf{S}_c\colon\mathscr{F}(c)\to c_{\mathsf{S}},\mathsf{S}_c^{\vee}\colon c_{\mathsf{S}}\to\mathscr{F}(c),\eta_c\colon\id_{ c_{\mathsf{S}}}\toiso\mathsf{S}_c\circ\mathsf{S}_c^{\vee},\epsilon_c\colon\mathsf{S}_c^{\vee}\circ\mathsf{S}_c\toiso\id_{\mathscr{F}(c)}),
\]
then there is
\begin{enumerate}
\item a pseudofunctor $\mathscr{F}[\mathsf{S}]\colon\mathscr{C}\rightarrow\mathscr{C}'$ defined by
\begin{enumerate}
\item $\mathscr{F}[\mathsf{S}](c)=c_{\mathsf{S}}$
\item $\mathscr{F}[\mathsf{S}]_{c_1,c_2}:=\mathsf{S}_{c_2}\circ\mathscr{F}_{c_1,c_2}\circ\mathsf{S}_{c_1}^{\vee}$
\item $\mathrm{m}^{\mathscr{F}[\mathsf{S}]}_{c_1,c_2,c_3}:=\mathrm{m}^{\mathscr{F}}_{c_1,c_2,c_3}\circ\epsilon_{c_2}$
\item $\mathrm{u}^{\mathscr{F}[\mathsf{S}]}_{c}=\mathrm{u}^{\mathscr{F}}_{c}\circ\eta_c$
\end{enumerate}
\item and a pseudonatural equivalence $\alpha[\mathsf{S}]\colon\mathscr{F}\rightarrow\mathscr{F}[\mathsf{S}]$ defined by
\begin{enumerate}
\item $\alpha[\mathsf{S}]_c=\mathsf{S}_c$
\item $\lambda^{\alpha}_{c_1,c_2}=(\epsilon_{c_1})^{-1}$.
\end{enumerate}
\end{enumerate}
\end{prop}

\begin{proof}
Let $c,c_1,c_2,c_3\in\mathscr{C}$. Since functors and natural isomorphisms are closed under composition we conclude that $\mathscr{F}[\mathsf{S}]_{c_1,c_2}$ is a functor and $\mathrm{m}^{\mathscr{F}[\mathsf{S}]}_{c_1,c_2,c_3}$, $\mathrm{u}^{\mathscr{F}[\mathsf{S}]}_{c}$ and $\lambda^{\alpha}_{c_1,c_2}$ are natural isomorphisms. Hence it is sufficent to check the associativity and unitality conditions for $\mathscr{F}[\mathsf{S}]$ and $\alpha[\mathsf{S}]$.
\begin{enumerate}
\item For the associativity of $\mathscr{F}[\mathsf{S}]$, we have to show that for $f\in\mathscr{C}(c_1,c_2)$, $g\in\mathscr{C}(c_2,c_3)$ and $h\in\mathscr{C}(c_3,c_4)$the pastings
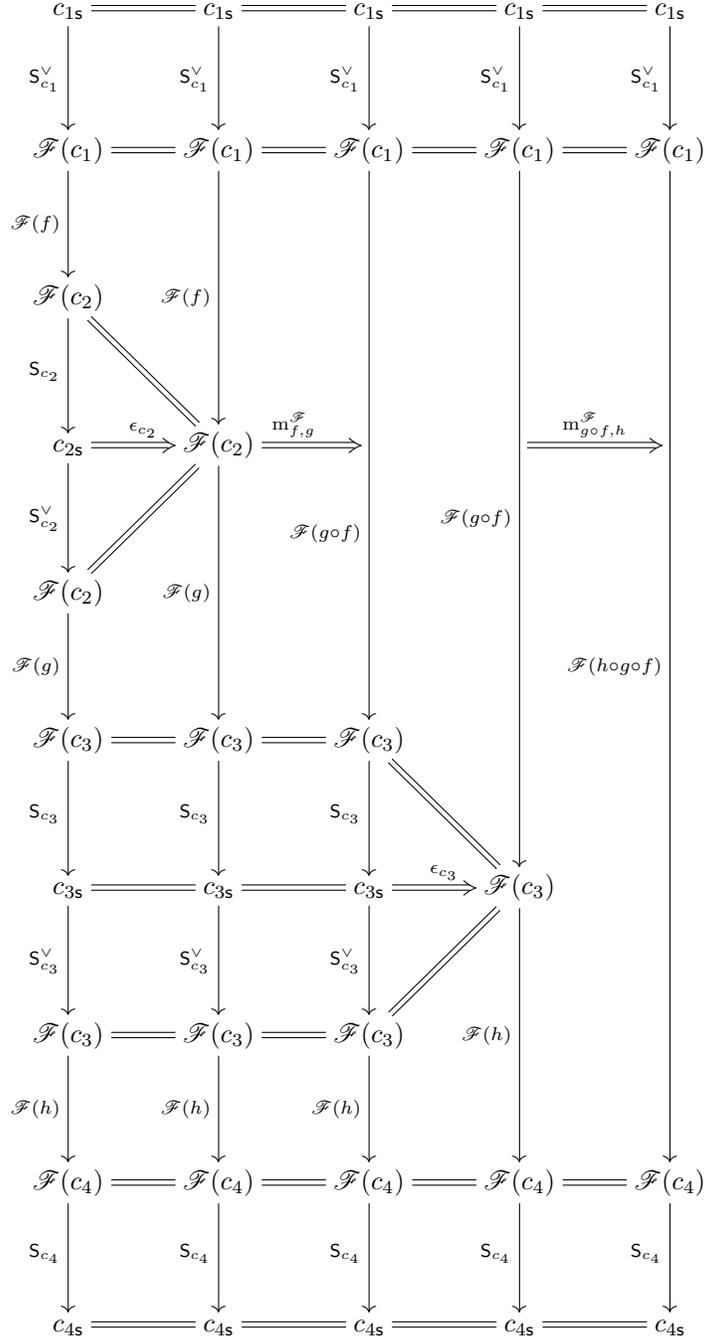
\[
\xymatrixrowsep{1.4cm}
\xymatrix{
c_{1\mathsf{s}}\ar[d]_{\mathsf{S}^{\vee}_{c_1}}\ar@{=}[r]&c_{1\mathsf{s}}\ar[d]_{\mathsf{S}^{\vee}_{c_1}}\ar@{=}[r]&c_{1\mathsf{s}}\ar[d]_{\mathsf{S}^{\vee}_{c_1}}\ar@{=}[r]&c_{1\mathsf{s}}\ar[d]_{\mathsf{S}^{\vee}_{c_1}}\ar@{=}[r]&c_{1\mathsf{s}}\ar[d]_{\mathsf{S}^{\vee}_{c_1}}\\
\mathscr{F}(c_1)\ar@{=}[r]\ar[d]_{\mathscr{F}(f)}&\mathscr{F}(c_1)\ar@{=}[r]\ar[dd]_{\mathscr{F}(f)}&\mathscr{F}(c_1)\ar@{=}[r]\ar[dddd]_(.65){\mathscr{F}(g\circ f)}&\mathscr{F}(c_1)\ar@{=}[r]\ar[ddddd]_{\mathscr{F}(g\circ f)}&\mathscr{F}(c_1)\ar[ddddddd]_{\mathscr{F}(h\circ g\circ f)}\\
\mathscr{F}(c_2)\ar[d]_{\mathsf{S}_{c_2}}\ar@{=}[rd]\\
c_{2\mathsf{s}}\ar[d]_{\mathsf{S}^{\vee}_{c_2}}\ar@{=>}[r]^{\epsilon_{c_2}}&\mathscr{F}(c_2)\ar[dd]_{\mathscr{F}(g)}\ar@{=>}[r]^{\mathrm{m}^{\mathscr{F}}_{f,g}}&&\ar@{=>}[r]^{\mathrm{m}^{\mathscr{F}}_{g\circ f,h}}&\\
\mathscr{F}(c_2)\ar[d]_{\mathscr{F}(g)}\ar@{=}[ru]\\
\mathscr{F}(c_3)\ar[d]_{\mathsf{S}_{c_3}}\ar@{=}[r]&\mathscr{F}(c_3)\ar[d]_{\mathsf{S}_{c_3}}\ar@{=}[r]&\mathscr{F}(c_3)\ar[d]_{\mathsf{S}_{c_3}}\ar@{=}[rd]\\
c_{3\mathsf{s}}\ar[d]_{\mathsf{S}^{\vee}_{c_3}}\ar@{=}[r]&c_{3\mathsf{s}}\ar[d]_{\mathsf{S}^{\vee}_{c_3}}\ar@{=}[r]&c_{3\mathsf{s}}\ar[d]_{\mathsf{S}^{\vee}_{c_3}}\ar@{=>}[r]^{\epsilon_{c_3}}&\mathscr{F}(c_3)\ar[dd]_{\mathscr{F}(h)}\\
\mathscr{F}(c_3)\ar[d]_{\mathscr{F}(h)}\ar@{=}[r]&\mathscr{F}(c_3)\ar[d]_{\mathscr{F}(h)}\ar@{=}[r]&\mathscr{F}(c_3)\ar[d]_{\mathscr{F}(h)}\ar@{=}[ru]\\
\mathscr{F}(c_4)\ar@{=}[r]\ar[d]_{\mathsf{S}_{c_4}}&\mathscr{F}(c_4)\ar@{=}[r]\ar[d]_{\mathsf{S}_{c_4}}&\mathscr{F}(c_4)\ar@{=}[r]\ar[d]_{\mathsf{S}_{c_4}}&\mathscr{F}(c_4)\ar@{=}[r]\ar[d]_{\mathsf{S}_{c_4}}&\mathscr{F}(c_4)\ar[d]_{\mathsf{S}_{c_4}}\\
c_{4\mathsf{s}}\ar@{=}[r]&c_{4\mathsf{s}}\ar@{=}[r]&c_{4\mathsf{s}}\ar@{=}[r]&c_{4\mathsf{s}}\ar@{=}[r]&c_{4\mathsf{s}}
}
\]
and
\[
\xymatrixrowsep{1.4cm}
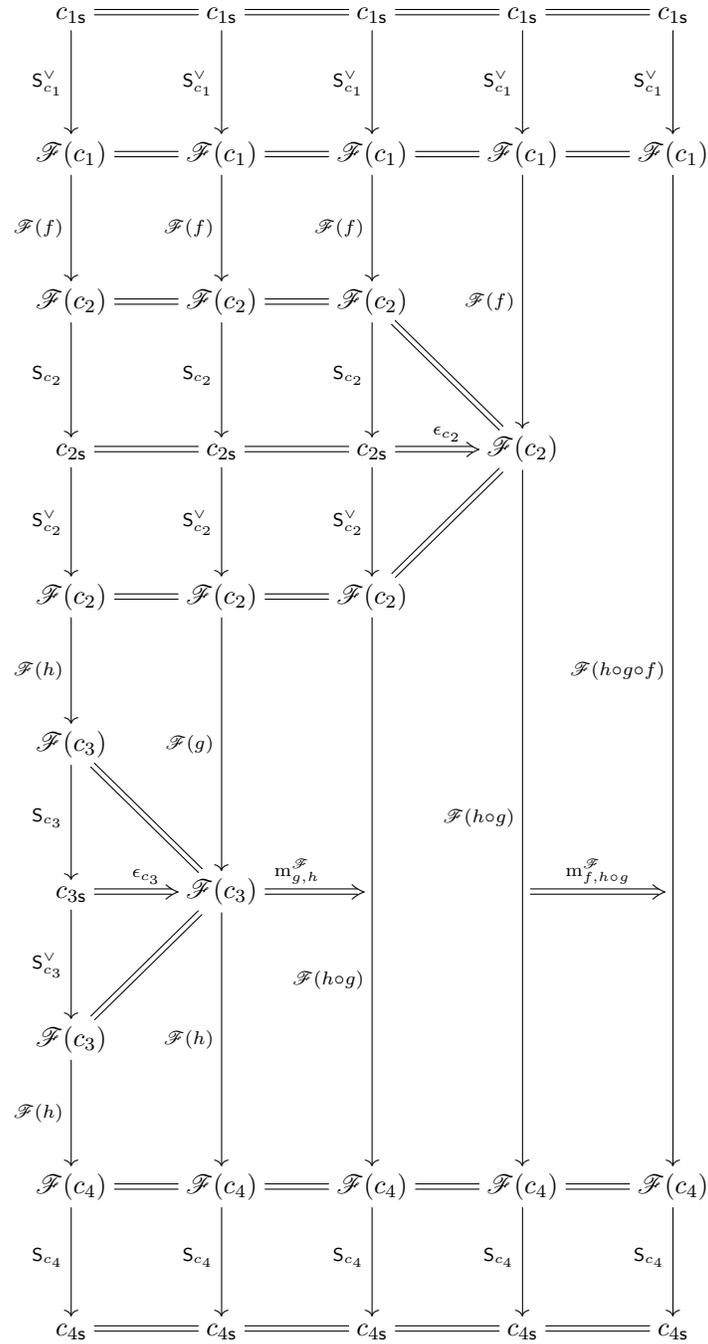
\xymatrix{
c_{1\mathsf{s}}\ar[d]_{\mathsf{S}^{\vee}_{c_1}}\ar@{=}[r]&c_{1\mathsf{s}}\ar[d]_{\mathsf{S}^{\vee}_{c_1}}\ar@{=}[r]&c_{1\mathsf{s}}\ar[d]_{\mathsf{S}^{\vee}_{c_1}}\ar@{=}[r]&c_{1\mathsf{s}}\ar[d]_{\mathsf{S}^{\vee}_{c_1}}\ar@{=}[r]&c_{1\mathsf{s}}\ar[d]_{\mathsf{S}^{\vee}_{c_1}}\\
\mathscr{F}(c_1)\ar@{=}[r]\ar[d]_{\mathscr{F}(f)}&\mathscr{F}(c_1)\ar@{=}[r]\ar[d]_{\mathscr{F}(f)}&\mathscr{F}(c_1)\ar@{=}[r]\ar[d]_{\mathscr{F}(f)}&\mathscr{F}(c_1)\ar@{=}[r]\ar[dd]_{\mathscr{F}(f)}&\mathscr{F}(c_1)\ar[ddddddd]_{\mathscr{F}(h\circ g\circ f)}\\
\mathscr{F}(c_2)\ar[d]_{\mathsf{S}_{c_2}}\ar@{=}[r]&\mathscr{F}(c_2)\ar[d]_{\mathsf{S}_{c_2}}\ar@{=}[r]&\mathscr{F}(c_2)\ar[d]_{\mathsf{S}_{c_2}}\ar@{=}[rd]\\
c_{2\mathsf{s}}\ar[d]_{\mathsf{S}^{\vee}_{c_2}}\ar@{=}[r]&c_{2\mathsf{s}}\ar[d]_{\mathsf{S}^{\vee}_{c_2}}\ar@{=}[r]&c_{2\mathsf{s}}\ar[d]_{\mathsf{S}^{\vee}_{c_2}}\ar@{=>}[r]^{\epsilon_{c_2}}&\mathscr{F}(c_2)\ar[ddddd]_{\mathscr{F}(h\circ g)}\\
\mathscr{F}(c_2)\ar[d]_{\mathscr{F}(h)}\ar@{=}[r]&\mathscr{F}(c_2)\ar[dd]_{\mathscr{F}(g)}\ar@{=}[r]&\mathscr{F}(c_2)\ar[dddd]_(.65){\mathscr{F}(h\circ g)}\ar@{=}[ru]\\
\mathscr{F}(c_3)\ar[d]_{\mathsf{S}_{c_3}}\ar@{=}[rd]\\
c_{3\mathsf{s}}\ar[d]_{\mathsf{S}^{\vee}_{c_3}}\ar@{=>}[r]^{\epsilon_{c_3}}&\mathscr{F}(c_3)\ar[dd]_{\mathscr{F}(h)}\ar@{=>}[r]^{\mathrm{m}^{\mathscr{F}}_{g,h}}&&\ar@{=>}[r]^{\mathrm{m}^{\mathscr{F}}_{f,h\circ g}}&\\
\mathscr{F}(c_3)\ar[d]_{\mathscr{F}(h)}\ar@{=}[ru]\\
\mathscr{F}(c_4)\ar@{=}[r]\ar[d]_{\mathsf{S}_{c_4}}&\mathscr{F}(c_4)\ar@{=}[r]\ar[d]_{\mathsf{S}_{c_4}}&\mathscr{F}(c_4)\ar@{=}[r]\ar[d]_{\mathsf{S}_{c_4}}&\mathscr{F}(c_4)\ar@{=}[r]\ar[d]_{\mathsf{S}_{c_4}}&\mathscr{F}(c_4)\ar[d]_{\mathsf{S}_{c_4}}\\
c_{4\mathsf{s}}\ar@{=}[r]&c_{4\mathsf{s}}\ar@{=}[r]&c_{4\mathsf{s}}\ar@{=}[r]&c_{4\mathsf{s}}\ar@{=}[r]&c_{4\mathsf{s}}
}
\]
agree. By contracting some identity cell this amounts to proving the equality of the pastings
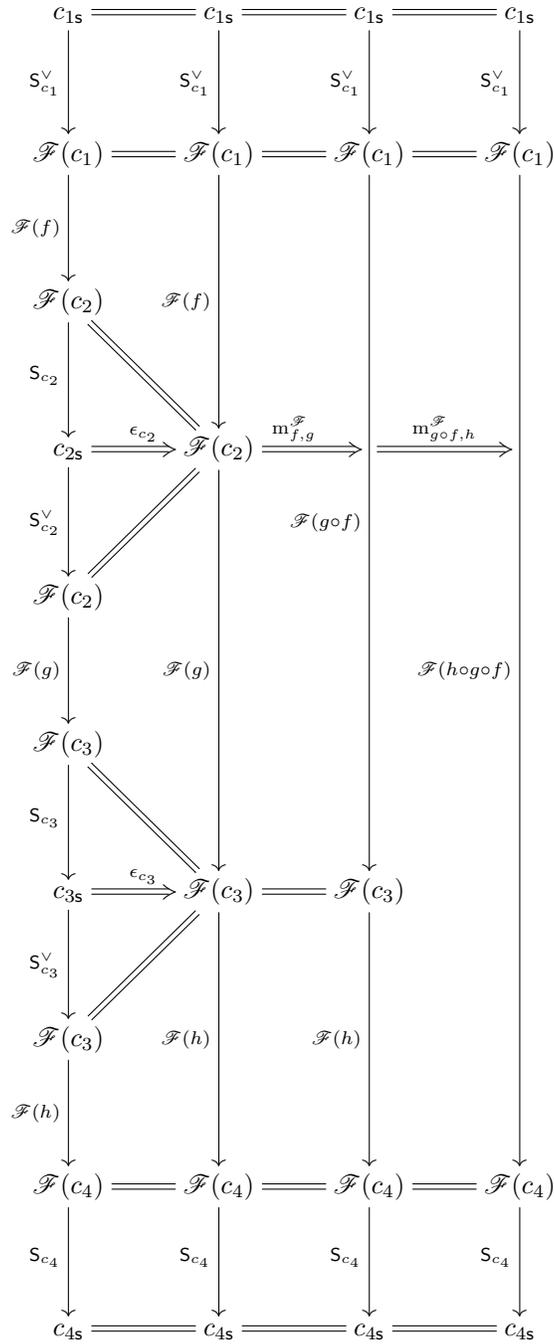
\[
\xymatrixrowsep{1.4cm}
\xymatrix{
c_{1\mathsf{s}}\ar[d]_{\mathsf{S}^{\vee}_{c_1}}\ar@{=}[r]&c_{1\mathsf{s}}\ar[d]_{\mathsf{S}^{\vee}_{c_1}}\ar@{=}[r]&c_{1\mathsf{s}}\ar[d]_{\mathsf{S}^{\vee}_{c_1}}\ar@{=}[r]&c_{1\mathsf{s}}\ar[d]_{\mathsf{S}^{\vee}_{c_1}}\\
\mathscr{F}(c_1)\ar@{=}[r]\ar[d]_{\mathscr{F}(f)}&\mathscr{F}(c_1)\ar@{=}[r]\ar[dd]_{\mathscr{F}(f)}&\mathscr{F}(c_1)\ar@{=}[r]\ar[ddddd]_{\mathscr{F}(g\circ f)}&\mathscr{F}(c_1)\ar[ddddddd]_{\mathscr{F}(h\circ g\circ f)}\\
\mathscr{F}(c_2)\ar[d]_{\mathsf{S}_{c_2}}\ar@{=}[rd]\\
c_{2\mathsf{s}}\ar[d]_{\mathsf{S}^{\vee}_{c_2}}\ar@{=>}[r]^{\epsilon_{c_2}}&\mathscr{F}(c_2)\ar[ddd]_{\mathscr{F}(g)}\ar@{=>}[r]^{\mathrm{m}^{\mathscr{F}}_{f,g}}&\ar@{=>}[r]^{\mathrm{m}^{\mathscr{F}}_{g\circ f,h}}&\\
\mathscr{F}(c_2)\ar[d]_{\mathscr{F}(g)}\ar@{=}[ru]\\
\mathscr{F}(c_3)\ar[d]_{\mathsf{S}_{c_3}}\ar@{=}[rd]\\
c_{3\mathsf{s}}\ar[d]_{\mathsf{S}^{\vee}_{c_3}}\ar@{=>}[r]^{\epsilon_{c_3}}&\mathscr{F}(c_3)\ar[dd]_{\mathscr{F}(h)}\ar@{=}[r]&\mathscr{F}(c_3)\ar[dd]_{\mathscr{F}(h)}\\
\mathscr{F}(c_3)\ar[d]_{\mathscr{F}(h)}\ar@{=}[ru]\\
\mathscr{F}(c_4)\ar@{=}[r]\ar[d]_{\mathsf{S}_{c_4}}&\mathscr{F}(c_4)\ar@{=}[r]\ar[d]_{\mathsf{S}_{c_4}}&\mathscr{F}(c_4)\ar@{=}[r]\ar[d]_{\mathsf{S}_{c_4}}&\mathscr{F}(c_4)\ar[d]_{\mathsf{S}_{c_4}}\\
c_{4\mathsf{s}}\ar@{=}[r]&c_{4\mathsf{s}}\ar@{=}[r]&c_{4\mathsf{s}}\ar@{=}[r]&c_{4\mathsf{s}}
}
\]
and
\[
\xymatrixrowsep{1.4cm}
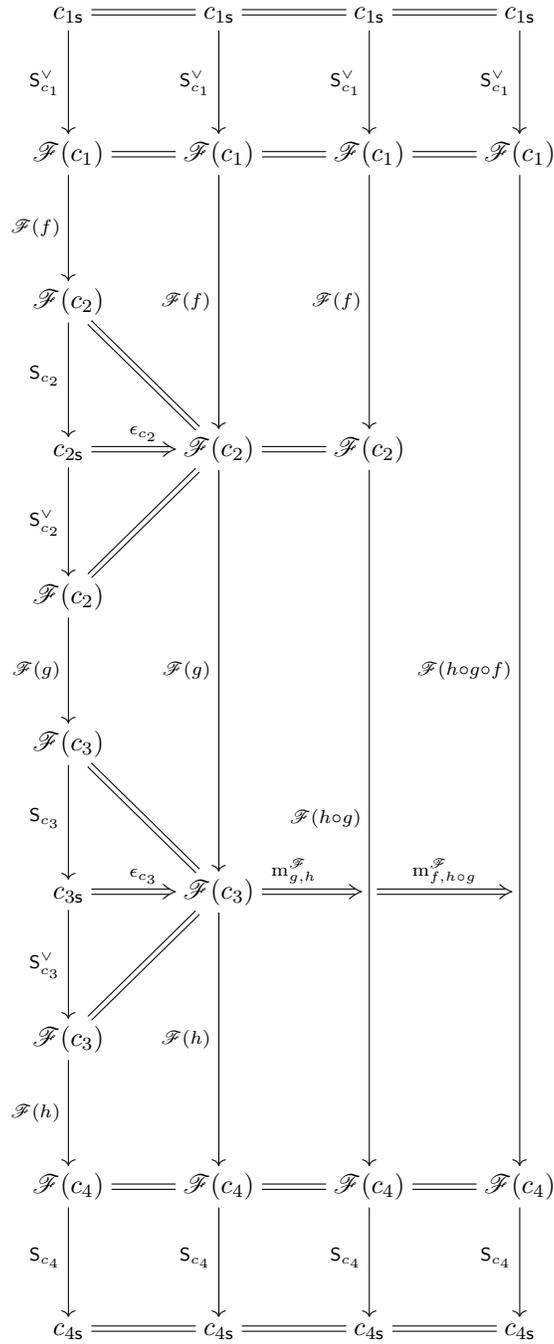
\xymatrix{
c_{1\mathsf{s}}\ar[d]_{\mathsf{S}^{\vee}_{c_1}}\ar@{=}[r]&c_{1\mathsf{s}}\ar[d]_{\mathsf{S}^{\vee}_{c_1}}\ar@{=}[r]&c_{1\mathsf{s}}\ar[d]_{\mathsf{S}^{\vee}_{c_1}}\ar@{=}[r]&c_{1\mathsf{s}}\ar[d]_{\mathsf{S}^{\vee}_{c_1}}\\
\mathscr{F}(c_1)\ar@{=}[r]\ar[d]_{\mathscr{F}(f)}&\mathscr{F}(c_1)\ar@{=}[r]\ar[dd]_{\mathscr{F}(f)}&\mathscr{F}(c_1)\ar@{=}[r]\ar[dd]_{\mathscr{F}(f)}&\mathscr{F}(c_1)\ar[ddddddd]_{\mathscr{F}(h\circ g\circ f)}\\
\mathscr{F}(c_2)\ar[d]_{\mathsf{S}_{c_2}}\ar@{=}[rd]\\
c_{2\mathsf{s}}\ar[d]_{\mathsf{S}^{\vee}_{c_2}}\ar@{=>}[r]^{\epsilon_{c_2}}&\mathscr{F}(c_2)\ar@{=}[r]\ar[ddd]_{\mathscr{F}(g)}&\mathscr{F}(c_2)\ar[ddddd]_{\mathscr{F}(h\circ g)}\\
\mathscr{F}(c_2)\ar[d]_{\mathscr{F}(g)}\ar@{=}[ru]\\
\mathscr{F}(c_3)\ar[d]_{\mathsf{S}_{c_3}}\ar@{=}[rd]\\
c_{3\mathsf{s}}\ar[d]_{\mathsf{S}^{\vee}_{c_3}}\ar@{=>}[r]^{\epsilon_{c_3}}&\mathscr{F}(c_3)\ar[dd]_{\mathscr{F}(h)}\ar@{=>}[r]^{\mathrm{m}^{\mathscr{F}}_{g,h}}&\ar@{=>}[r]^{\mathrm{m}^{\mathscr{F}}_{f,h\circ g}}&\\
\mathscr{F}(c_3)\ar[d]_{\mathscr{F}(h)}\ar@{=}[ru]\\
\mathscr{F}(c_4)\ar@{=}[r]\ar[d]_{\mathsf{S}_{c_4}}&\mathscr{F}(c_4)\ar@{=}[r]\ar[d]_{\mathsf{S}_{c_4}}&\mathscr{F}(c_4)\ar@{=}[r]\ar[d]_{\mathsf{S}_{c_4}}&\mathscr{F}(c_4)\ar[d]_{\mathsf{S}_{c_4}}\\
c_{4\mathsf{s}}\ar@{=}[r]&c_{4\mathsf{s}}\ar@{=}[r]&c_{4\mathsf{s}}\ar@{=}[r]&c_{4\mathsf{s}}
}
\]
But this follows, by comparing the subpastings of the inner cells in the second and third columns, since $\mathscr{F}$ is a pseudofunctor.
\item For the left unitality of $\mathscr{F}[\mathsf{S}]$ we have to show that for $f\in\mathscr{C}(c_1,c_2)$ the following pasting is equal to the identity.
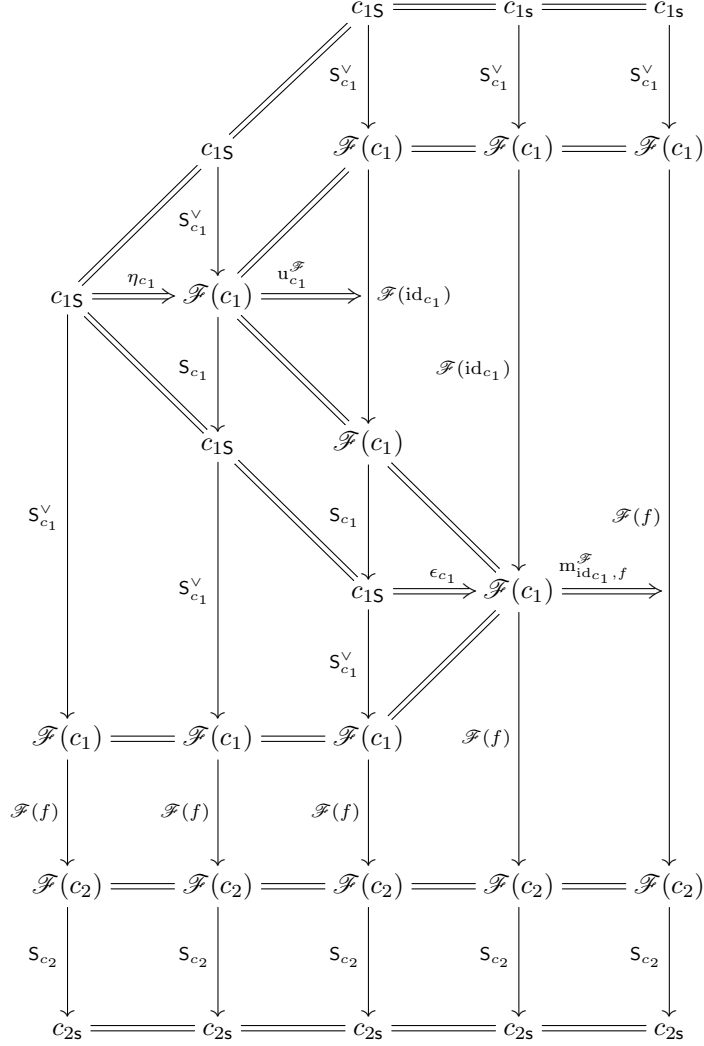
\[
\xymatrixrowsep{1.4cm}
\xymatrix{
&&c_{1\mathsf{S}}\ar[d]_{\mathsf{S}^{\vee}_{c_1}}\ar@{=}[r]&c_{1\mathsf{s}}\ar[d]_{\mathsf{S}^{\vee}_{c_1}}\ar@{=}[r]&c_{1\mathsf{s}}\ar[d]_{\mathsf{S}^{\vee}_{c_1}}\\
&c_{1\mathsf{S}}\ar[d]_{\mathsf{S}^{\vee}_{c_1}}\ar@{=}[ur]&\mathscr{F}(c_1)\ar[dd]^{\mathscr{F}(\id_{c_1})}\ar@{=}[r]&\mathscr{F}(c_1)\ar[ddd]_{\mathscr{F}(\id_{c_1})}\ar@{=}[r]&\mathscr{F}(c_1)\ar[ddddd]_{\mathscr{F}(f)}&\\
c_{1\mathsf{S}}\ar[ddd]_{\mathsf{S}^{\vee}_{c_1}}\ar@{=}[ur]\ar@{=}[rd]\ar@{=>}[r]^{\eta_{c_1}}&\mathscr{F}(c_1)\ar[d]_{\mathsf{S}_{c_1}}\ar@{=}[ur]\ar@{=}[rd]\ar@{=>}[r]^{\mathrm{u}^{\mathscr{F}}_{c_1}}&&&\\
&c_{1\mathsf{S}}\ar[dd]_{\mathsf{S}^{\vee}_{c_1}}\ar@{=}[rd]&\mathscr{F}(c_1)\ar[d]_{\mathsf{S}_{c_1}}\ar@{=}[rd]\\
&&c_{1\mathsf{S}}\ar[d]_{\mathsf{S}^{\vee}_{c_1}}\ar@{=>}[r]^{\epsilon_{c_1}}&\mathscr{F}(c_1)\ar@{=>}[r]^{\mathrm{m}^{\mathscr{F}}_{\id_{c_1},f}}\ar[dd]_{\mathscr{F}(f)}&\\
\mathscr{F}(c_1)\ar[d]_{\mathscr{F}(f)}\ar@{=}[r]&\mathscr{F}(c_1)\ar[d]_{\mathscr{F}(f)}\ar@{=}[r]&\mathscr{F}(c_1)\ar[d]_{\mathscr{F}(f)}\ar@{=}[ur]&&\\
\mathscr{F}(c_2)\ar[d]_{\mathsf{S}_{c_2}}\ar@{=}[r]&\mathscr{F}(c_2)\ar[d]_{\mathsf{S}_{c_2}}\ar@{=}[r]&\mathscr{F}(c_2)\ar[d]_{\mathsf{S}_{c_2}}\ar@{=}[r]&\mathscr{F}(c_2)\ar[d]_{\mathsf{S}_{c_2}}\ar@{=}[r]&\mathscr{F}(c_2)\ar[d]_{\mathsf{S}_{c_2}}\\
c_{2\mathsf{s}}\ar@{=}[r]&c_{2\mathsf{s}}\ar@{=}[r]&c_{2\mathsf{s}}\ar@{=}[r]&c_{2\mathsf{s}}\ar@{=}[r]&c_{2\mathsf{s}}
}
\]
For this, we note that the pasting of the second cells in the second, third and fourth column is the identity, since $\mathscr{F}$ is a pseudofunctor. The two remaining non-trivial cells compose to the identity because of one of the triangle identities of the eqiuvalence $(\mathsf{S}_{c_1},\mathsf{S}^{\vee}_{c_1})$. The proof of the right unitality is completely dual to this case.
\item For the associativity condition of $\alpha[\mathsf{S}]$ we have to show that for $f\in\mathscr{C}(c_1,c_2)$ and $g\in\mathscr{C}(c_2,c_3)$ the pasting
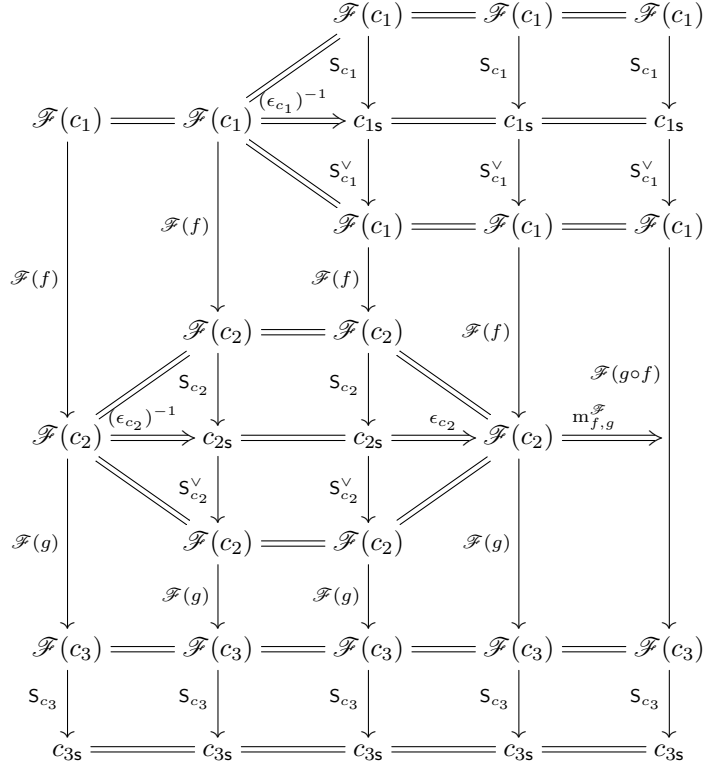
\[
\xymatrix{
&&\mathscr{F}(c_1)\ar[d]_{\mathsf{S}_{c_1}}\ar@{=}[r]&\mathscr{F}(c_1)\ar[d]_{\mathsf{S}_{c_1}}\ar@{=}[r]&\mathscr{F}(c_1)\ar[d]_{\mathsf{S}_{c_1}}\\
\mathscr{F}(c_1)\ar[ddd]_{\mathscr{F}(f)}\ar@{=}[r]&\mathscr{F}(c_1)\ar[dd]_{\mathscr{F}(f)}\ar@{=>}[r]^{(\epsilon_{c_1})^{-1}}\ar@{=}[ur]\ar@{=}[rd]&c_{1\mathsf{s}}\ar[d]_{\mathsf{S}^{\vee}_{c_1}}\ar@{=}[r]&c_{1\mathsf{s}}\ar@{=}[r]\ar[d]_{\mathsf{S}^{\vee}_{c_1}}&c_{1\mathsf{s}}\ar[d]_{\mathsf{S}^{\vee}_{c_1}}\\
&&\mathscr{F}(c_1)\ar[d]_{\mathscr{F}(f)}\ar@{=}[r]&\mathscr{F}(c_1)\ar[dd]_{\mathscr{F}(f)}\ar@{=}[r]&\mathscr{F}(c_1)\ar[dddd]_(.35){\mathscr{F}(g\circ f)}\\
&\mathscr{F}(c_2)\ar[d]_{\mathsf{S}_{c_2}}\ar@{=}[r]&\mathscr{F}(c_2)\ar@{=}[rd]\ar[d]_{\mathsf{S}_{c_2}}&&&\\
\mathscr{F}(c_2)\ar[dd]_{\mathscr{F}(g)}\ar@{=>}[r]^{(\epsilon_{c_2})^{-1}}\ar@{=}[ur]\ar@{=}[rd]&c_{2\mathsf{s}}\ar[d]_{\mathsf{S}^{\vee}_{c_2}}\ar@{=}[r]&c_{2\mathsf{s}}\ar[d]_{\mathsf{S}^{\vee}_{c_2}}\ar@{=>}[r]^{\epsilon_{c_2}}&\mathscr{F}(c_2)\ar@{=>}[r]^{\mathrm{m}^{\mathscr{F}}_{f,g}}\ar[dd]_{\mathscr{F}(g)}&\\
&\mathscr{F}(c_2)\ar[d]_{\mathscr{F}(g)}\ar@{=}[r]&\mathscr{F}(c_2)\ar@{=}[ur]\ar[d]_{\mathscr{F}(g)}&&&\\
\mathscr{F}(c_3)\ar[d]_{\mathsf{S}_{c_3}}\ar@{=}[r]&\mathscr{F}(c_3)\ar[d]_{\mathsf{S}_{c_3}}\ar@{=}[r]&\mathscr{F}(c_3)\ar[d]_{\mathsf{S}_{c_3}}\ar@{=}[r]&\mathscr{F}(c_3)\ar[d]_{\mathsf{S}_{c_3}}\ar@{=}[r]&\mathscr{F}(c_3)\ar[d]_{\mathsf{S}_{c_3}}\\
c_{3\mathsf{s}}\ar@{=}[r]&c_{3\mathsf{s}}\ar@{=}[r]&c_{3\mathsf{s}}\ar@{=}[r]&c_{3\mathsf{s}}\ar@{=}[r]&c_{3\mathsf{s}}
}
\]
agrees with
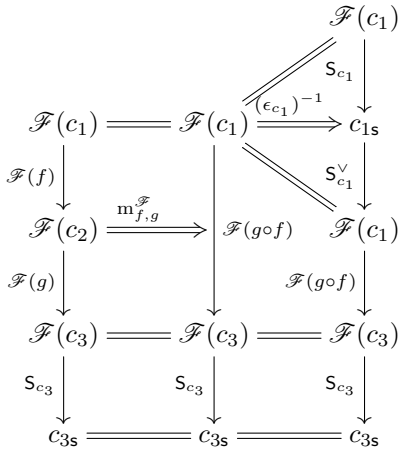
\[
\xymatrix{
&&\mathscr{F}(c_1)\ar[d]_{\mathsf{S}_{c_1}}\\
\mathscr{F}(c_1)\ar[d]_{\mathscr{F}(f)}\ar@{=}[r]&\mathscr{F}(c_1)\ar[dd]^{\mathscr{F}(g\circ f)}\ar@{=>}[r]^{(\epsilon_{c_1})^{-1}}\ar@{=}[ur]\ar@{=}[rd]&c_{1\mathsf{s}}\ar[d]_{\mathsf{S}^{\vee}_{c_1}}\\
\mathscr{F}(c_2)\ar@{=>}[r]^{\mathrm{m}^{\mathscr{F}}_{f,g}}\ar[d]_{\mathscr{F}(g)}&&\mathscr{F}(c_1)\ar[d]_{\mathscr{F}(g\circ f)}\\
\mathscr{F}(c_3)\ar[d]_{\mathsf{S}_{c_3}}\ar@{=}[r]&\mathscr{F}(c_3)\ar[d]_{\mathsf{S}_{c_3}}\ar@{=}[r]&\mathscr{F}(c_3)\ar[d]_{\mathsf{S}_{c_3}}\\
c_{3\mathsf{s}}\ar@{=}[r]&c_{3\mathsf{s}}\ar@{=}[r]&c_{3\mathsf{s}}
}
\]
But this is immediate since the hexagon in the center of the first diagram composes to the identity.
\item For the unitality of $\alpha[\mathsf{S}]$ we have to show that for $c\in\mathscr{C}$ the pastings
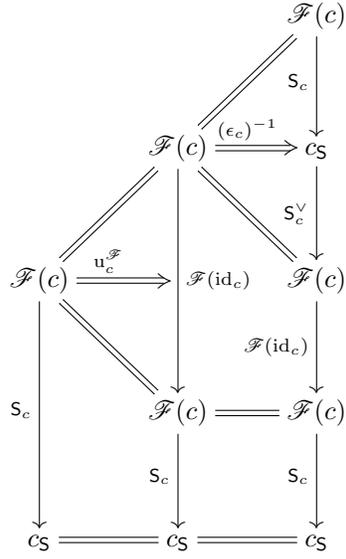
\[
\xymatrixrowsep{1.2cm}
\xymatrix{
&&\mathscr{F}(c)\ar[d]_{\mathsf{S}_{c}}\\
&\mathscr{F}(c)\ar[dd]^{\mathscr{F}(\id_c)}\ar@{=>}[r]^{(\epsilon_{c})^{-1}}\ar@{=}[ur]\ar@{=}[rd]&c_{\mathsf{S}}\ar[d]_{\mathsf{S}^{\vee}_{c}}\\
\mathscr{F}(c)\ar[dd]_{\mathsf{S}_{c}}\ar@{=}[ur]\ar@{=}[rd]\ar@{=>}[r]^{\mathrm{u}^{\mathscr{F}}_c}&&\mathscr{F}(c)\ar[d]_{\mathscr{F}(\id_c)}\\
&\mathscr{F}(c)\ar[d]_{\mathsf{S}_{c}}\ar@{=}[r]&\mathscr{F}(c)\ar[d]_{\mathsf{S}_{c}}\\
c_{\mathsf{S}}\ar@{=}[r]&c_{\mathsf{S}}\ar@{=}[r]&c_{\mathsf{S}}
}
\]
and
\[
\xymatrixrowsep{1.2cm}
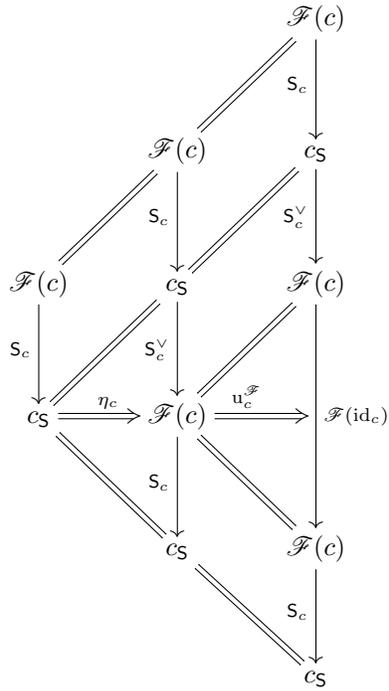
\xymatrix{
&&\mathscr{F}(c)\ar[d]_{\mathsf{S}_{c}}\\
&\mathscr{F}(c)\ar@{=}[ur]\ar[d]_{\mathsf{S}_{c}}&c_{\mathsf{S}}\ar[d]_{\mathsf{S}^{\vee}_{c}}\\
\mathscr{F}(c)\ar@{=}[ur]\ar[d]_{\mathsf{S}_{c}}&c_{\mathsf{S}}\ar[d]_{\mathsf{S}^{\vee}_{c}}\ar@{=}[ur]&\mathscr{F}(c)\ar[dd]^{\mathscr{F}(\id_c)}\\
c_{\mathsf{S}}\ar@{=>}[r]^{\eta_c}\ar@{=}[ur]\ar@{=}[dr]&\mathscr{F}(c)\ar@{=>}[r]^{\mathrm{u}^{\mathscr{F}}_c}\ar[d]_{\mathsf{S}_{c}}\ar@{=}[ur]\ar@{=}[dr]&\\
&c_{\mathsf{S}}\ar@{=}[dr]&\mathscr{F}(c)\ar[d]_{\mathsf{S}_{c}}\\
&&c_{\mathsf{S}}
}
\]
agree. We apply one of the triangle identities of the equivalence $(\mathsf{S}_{c},\mathsf{S}^{\vee}_{c})$ to the first column of the second pasting to obtain
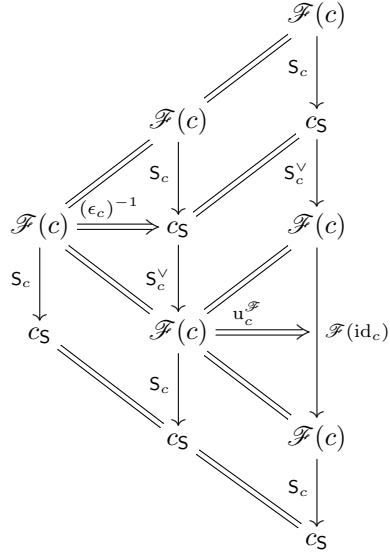
\[
\xymatrix{
&&\mathscr{F}(c)\ar[d]_{\mathsf{S}_{c}}\\
&\mathscr{F}(c)\ar@{=}[ur]\ar[d]_{\mathsf{S}_{c}}&c_{\mathsf{S}}\ar[d]_{\mathsf{S}^{\vee}_{c}}\\
\mathscr{F}(c)\ar@{=>}[r]^{(\epsilon_{c})^{-1}}\ar@{=}[ur]\ar@{=}[dr]\ar[d]_{\mathsf{S}_{c}}&c_{\mathsf{S}}\ar[d]_{\mathsf{S}^{\vee}_{c}}\ar@{=}[ur]&\mathscr{F}(c)\ar[dd]^{\mathscr{F}(\id_c)}\\
c_{\mathsf{S}}\ar@{=}[dr]&\mathscr{F}(c)\ar@{=>}[r]^{\mathrm{u}^{\mathscr{F}}_c}\ar[d]_{\mathsf{S}_{c}}\ar@{=}[ur]\ar@{=}[dr]&\\
&c_{\mathsf{S}}\ar@{=}[dr]&\mathscr{F}(c)\ar[d]_{\mathsf{S}_{c}}\\
&&c_{\mathsf{S}}
}
\]
This pasting is seen to be equal to the first diagram by contracting identity cells.
\end{enumerate}

\end{proof}

\bibliographystyle{alpha}
\bibliography{introduction}

\end{document}